\documentclass[11pt]{article}
\usepackage[margin=1in]{geometry}
\usepackage{amsmath,amssymb,amsthm,mathtools}
\usepackage{booktabs}
\usepackage{graphicx}
\usepackage{array}
\usepackage{tikz}
\usetikzlibrary{arrows.meta}
\usepackage[round]{natbib}
\usepackage[hidelinks,hypertexnames=false]{hyperref}
\usepackage[capitalise]{cleveref}

\usepackage{aliascnt}
\newtheorem{theorem}{Theorem}
\newaliascnt{proposition}{theorem}
\newtheorem{proposition}[proposition]{Proposition}
\aliascntresetthe{proposition}
\newaliascnt{lemma}{theorem}
\newtheorem{lemma}[lemma]{Lemma}
\aliascntresetthe{lemma}
\newaliascnt{corollary}{theorem}
\newtheorem{corollary}[corollary]{Corollary}
\aliascntresetthe{corollary}
\theoremstyle{definition}
\newaliascnt{definition}{theorem}
\newtheorem{definition}[definition]{Definition}
\aliascntresetthe{definition}

\theoremstyle{remark}
\newaliascnt{remark}{theorem}
\newtheorem{remark}[remark]{Remark}
\aliascntresetthe{remark}
\crefname{theorem}{Theorem}{Theorems}
\crefname{proposition}{Proposition}{Propositions}
\crefname{lemma}{Lemma}{Lemmas}
\crefname{corollary}{Corollary}{Corollaries}
\crefname{definition}{Definition}{Definitions}
\crefname{assumption}{Assumption}{Assumptions}
\crefname{condition}{Condition}{Conditions}
\crefname{remark}{Remark}{Remarks}
\crefname{appendix}{Appendix}{Appendices}
\crefname{subappendix}{Appendix}{Appendices}

\newcommand{\R}{\mathbb{R}}
\newcommand{\Prob}{\mathbb{P}}
\newcommand{\E}{\mathbb{E}}
\newcommand{\Var}{\operatorname{Var}}
\newcommand{\Cov}{\operatorname{Cov}}
\newcommand{\cG}{\mathcal{G}}
\newcommand{\cK}{\mathcal{K}}
\newcommand{\cM}{\mathcal{M}}
\newcommand{\cN}{\mathcal{N}}
\newcommand{\cX}{\mathcal{X}}
\newcommand{\MP}{\mathcal{MP}}
\newcommand{\Xeq}{\mathcal{X}_{\mathrm{eq}}}
\newcommand{\Xcanc}{\mathcal{X}_{\mathrm{canc}}}
\newcommand{\pa}{\operatorname{pa}}
\newcommand{\adj}{\operatorname{adj}}

\newcommand{\Leb}{\operatorname{Leb}}
\newcommand{\dtv}{d_{\mathrm{TV}}}
\newcommand{\KL}{\mathrm{KL}}
\newcommand{\ind}{\mathbf{1}}
\newcommand{\eps}{\varepsilon}

\newcommand{\orl}{\mathrm{or}}
\newcommand{\thadj}{\theta^{\mathrm{adj}}}
\newcommand{\abar}{\bar a}
\newcommand{\Phibar}{\bar\Phi}
\newcommand{\sigl}{\underline{\sigma}}
\newcommand{\sigu}{\bar{\sigma}}
\newcommand{\indep}{\perp\!\!\!\perp}
\newcommand{\DS}{\textup{\textsc{ds}}}
\newcommand{\DSA}{\ensuremath{\DS_A}}
\newcommand{\DSOp}{\ensuremath{\DS_O}}
\newcommand{\DSeff}{\ensuremath{\DS_{\mathrm{eff}}}}
\newcommand{\DSci}{\ensuremath{\DS^{\mathrm{CI}}}}
\newcommand{\DSciO}{\ensuremath{\DS^{\mathrm{CI}}_O}}
\newcommand{\DSlrO}{\ensuremath{\DS^{\mathrm{LR}}_O}}
\newcommand{\DSlreff}{\ensuremath{\DS^{\mathrm{LR}}_{\mathrm{eff}}}}
\newcommand{\DSdir}{\ensuremath{\DS^{P_0}_{\mathrm{eff}}}}
\newcommand{\DSsw}{\ensuremath{\DS^{\mathrm{sw}}_{\mathrm{eff}}}}
\newcommand{\DSV}{\textup{\textsc{dsv}}}
\newcommand{\DSVp}{\textup{\textsc{dsv}\textsuperscript{+}}}

\title{How long must an honest confidence interval be after causal discovery?\\ Benign and malignant weak edges}
\author{Ruiqi Zhao\\[2pt]
{\normalsize School of Data Science, The Chinese University of Hong Kong, Shenzhen}\\
{\normalsize\texttt{124090916@link.cuhk.edu.cn}}}
\date{November 2026}

\begin{document}
\maketitle

\begin{abstract}
After causal discovery, a confidence interval for a total effect valid uniformly over the model class must pay for structural uncertainty. We measure this price by the worst-case ratio of expected length to that of the efficient interval knowing the graph, for linear Gaussian models with a known superset of the adjacencies and $k$-triangle faithfulness. The global price is of order $\sqrt n$ once one edge can be malignant, so we study it locally, near a generic distribution with one edge absent. Our main result classifies this weak edge by Meek closures. It is benign if deleting it leaves the effect identifiable in the resulting MPDAG and Markov equivalent explanations agree; the local price is then bounded, and otherwise, for small $k$, of order $\sqrt n$. At benign edges the price is of order $\mathcal R+K$, with $\mathcal R$ a ratio of efficient lengths and $K$ a cross-model constant. It equals one, so that honesty is free, exactly when $K=0$, generically decided by a d-separation. For two Markov classes it is determined within a factor $3.40$, and its upper bound is approached without knowing the base point or the weak edge. A further graphical condition gives the rate uniformly in a second edge's margin.
\end{abstract}

\section{Introduction}\label{sec:intro}

A common way to estimate a causal effect from observational data is to learn a causal graph first, for example with the PC algorithm or greedy equivalence search, and then to estimate the effect by adjusting for the parents of the treatment in the learned graph. The resulting confidence interval treats the learned graph as known. It is well known that this is not valid: when an edge is weak, the learned graph is wrong with non-vanishing probability, and the interval undercovers \citep{strieder2021confidence}. Under the Markov and faithfulness assumptions, \citet{robins2003uniform} showed that tests and estimators of a causal effect can be pointwise but not uniformly consistent, and that no consistent confidence intervals exist. This holds when latent confounders are possible, and also when the causal order is unknown even without confounders; see \citet[Section~2.1]{zhang2003strong} for the linear Gaussian, causally sufficient case. \citet{kalisch2007estimating} showed that the PC algorithm is uniformly consistent under strong faithfulness, and \citet{uhler2013geometry} showed that the margins needed by strong faithfulness exclude a large part of the parameter space.

Several recent papers construct valid intervals or confidence sets after causal discovery: by resampling the discovery step \citep{chang2026post}, by differentially private selection \citep{gradu2025valid}, by selective inference on the selection event of FCI \citep{guo2026price}, by inverting tests over graphs \citep{strieder2023confidence,strieder2024dual,strieder2024identifying}, by weighting the functionals of several candidate models \citep{bhattacharya2026robust}, and by confidence sets of causal orderings, of causal models or of the equivalence class itself \citep{wang2026confidence,edgar2026inferential,battey2026post,kang2026interpretable}. On the identification side, \citet{teh2025towards} give a sufficient condition under which the effect is the same across a candidate set of Markov equivalence classes. We review these papers at the end of this section. To our knowledge, none of them establishes a minimax lower bound on the expected length of confidence intervals for a total effect, uniformly over a class that explicitly contains near-unfaithful distributions, together with a graphical criterion for which weak edges force such intervals to be long.

This paper asks: \emph{over a class of linear Gaussian models that contains near-unfaithful distributions, how short can a confidence interval for a total effect be if it is valid uniformly over the class?} We measure length relative to the efficient interval of an analyst who knows the true graph, and call the worst-case ratio the \emph{price} of honesty. The class allows every edge to be weak. For small $k$, its global price is therefore at least of order $\sqrt n$ as soon as one edge can be malignant in the sense below (\cref{rem:global}), and the informative answer is local. We study the local minimax price, as $n\to\infty$, near a generic distribution at which one edge of the true graph is absent; all our results on general graphs are of this kind. Locally, then, the adjacencies are known except for one edge. This is the generic situation at the boundary of a model: among distributions at which an edge of the skeleton is absent, those at which two are absent form a set of lower dimension. Near that set we obtain the rate uniformly in the margin of the second edge (\cref{thm:uniform}). The local price is the cost of honesty once the data have localized the problem to a neighbourhood. The neighbourhood, the threshold on $k$ and the sample size from which our bounds apply depend on the base point. Throughout, statements about the price on general graphs hold near a generic base point and for $k$ below a threshold that depends on it; we do not repeat these qualifiers.

\paragraph{Benign and malignant weak edges.} The simplest instance is already informative. Let $Z$ be either a confounder or a mediator of the effect of $X$ on $Y$, and let $W\to X$ be an instrument-like edge. If the edge between $X$ and $Z$ is weak, the two candidate graphs give effects that differ by an amount proportional to the weak coefficient, so the effect is continuous at the crossing point and no margin is needed (\cref{cor:benign-toy}). If instead $W\to X$ is weak, the two graphs become Markov equivalent while their effects differ by a fixed amount, and honest intervals cannot shrink (\cref{thm:rate-lower}). The general criterion says exactly this: deleting the benign edge keeps the effect identifiable; deleting the malignant one does not.

\begin{table}[t]
\centering
\caption{The three kinds of weak edge. $H$ is the true graph without the weak edge and $\MP(H)$ its MPDAG; the second criterion uses the MPDAGs of the Markov classes of its explanations (\cref{thm:K1p-graph}). Local prices are near a generic base point, for small $k$. The values for benign edges are approached by procedures that do not know the base point.}
\label{tab:results}
\small
\begin{tabular}{>{\raggedright\arraybackslash}p{0.14\textwidth}>{\raggedright\arraybackslash}p{0.28\textwidth}>{\raggedright\arraybackslash}p{0.3\textwidth}>{\raggedright\arraybackslash}p{0.16\textwidth}}
\toprule
Weak edge & Graphical criterion & Local price of intervals & Results\\
\midrule
Benign, free & no bad path in $\MP(H)$ or in the MPDAG of any Markov class, and no explanation meets the d-separation of \cref{prop:Kpos} & exactly $1$ & \cref{thm:free}\\
Benign, with a price & as above, but some explanation meets it & generically $>1$, and $\asymp\mathcal R+K$; within $3.40$ for two Markov classes & \cref{thm:K,cor:RK,thm:free}\\
Malignant & bad path in $\MP(H)$ & $\gtrsim\sqrt n$ at the base point & \cref{prop:malignant}\\
Non-identified class & no bad path in $\MP(H)$, one in the MPDAG of some Markov class & $\gtrsim\sqrt n$ on every fixed neighbourhood & \cref{prop:K1p-nec}\\
\bottomrule
\end{tabular}
\end{table}

\paragraph{Main result.} The skeleton is a known superset of the true adjacencies (a missing edge is an edge with zero coefficient), and the candidate graphs are its orientations that satisfy background knowledge: that the outcome is a sink, together with optional required directed edges such as a time order. Let $e$ be the edge of the true graph $G$ that is absent at the base point, and $H=G-e$. Our main result (\cref{thm:dichotomy}, \cref{tab:results}) is that every such weak edge falls into one of three classes, and that Meek closures decide which.
\begin{itemize}
\item \emph{Benign}: the effect is identifiable in the MPDAG of $H$ \citep{perkovic2020identifying}, and also in an MPDAG built from $H$ and each Markov class of explanations of the base point (\cref{thm:K1p-graph}). The local price is then bounded and of order $\mathcal R+K$ (\cref{thm:K,cor:RK}). Here $\mathcal R\ge1$ is the largest ratio of the efficient standard deviations of two explanations, and $K$ is a cross-model constant: the first-order effect discrepancy between two explanations, divided by the square root of the Fisher information that separates them. Honesty is free, with local minimax price exactly one, if and only if $K=0$, which a d-separation in the explanations decides (\cref{thm:free,prop:Kpos}). On five nodes this is the case for $8186$ of the $8426$ benign pairs $(G,e)$. For some benign edges, Markov equivalent candidates that do not explain the base point are excluded only by triangle faithfulness, so the bounded price rests on identification supplied by the class rather than by the data (\cref{rem:ktf-identifies}).
\item \emph{Malignant}: the effect is not identifiable in the MPDAG of $H$. The price at the base point is at least of order $\sqrt n$ (\cref{prop:malignant}).
\item \emph{Locally non-identified Markov class}: the effect is identifiable in the MPDAG of $H$ but not in the MPDAG of some Markov class. The price is at least of order $\sqrt n$ on every fixed neighbourhood (\cref{prop:K1p-nec}).
\end{itemize}
With required edges the classification needs two conditions, which hold for a time order (\cref{lem:tier}).

\paragraph{Other results.}
\begin{enumerate}
\item \textbf{The constant and adaptation} (\cref{cor:RK,def:ds,thm:coverage,cor:compact}). For two Markov classes the local minimax price is determined within the factor $3.40$ at $\alpha=0.05$, for every graph, dimension and eigenvalue range. \emph{Double screening} (\cref{def:ds}) knows neither the base point nor the weak edge. It estimates the structure on a vanishing fraction of the sample and directs one test per candidate at the alternative explanations; given the direction, the test is exact. It comes arbitrarily close to the upper bound, uniformly over compact sets of generic base points, and a variant comes arbitrarily close to one where honesty is free (\cref{thm:free}). With adjustment intervals it is honest in finite samples.
\item \textbf{The structure of $K$} (\cref{prop:canonical,thm:K-lambda}). Each pair of explanations gives a canonical limit experiment. As the base point approaches a malignant stratum at margin $\lambda$ of a second edge, $K$ is either of order $1/\lambda$ or bounded, and three graphical conditions decide which.
\item \textbf{Uniformity in the margin} (\cref{thm:uniform}). Under one more graphical condition, near a distribution at which both edges are absent, the minimax expected length is of order $n^{-1/2}+\min\{A,1/(\lambda\sqrt n)\}$ or $n^{-1/2}$, according to these conditions, uniformly in $\lambda$ and in the bound $A$ on the weak coefficient. Double screening attains it. There is no plateau that does not depend on $\lambda$ and $A$, unlike in the four-variable model below, where it comes from a known range of the noise variances.
\item \textbf{A four-variable model} (\cref{thm:rate-lower,thm:rate-upper}). When the mediator-or-confounder status of one variable is uncertain, we obtain the minimax expected length uniformly over the margin on the weak orientation edge and over bounds on two other coefficients. Besides the parametric term and a term inversely proportional to the margin times $\sqrt n$, it has a plateau, of order $n^{-1/4}$ for unit noise variances, that comes from the known range of the noise variances. Selecting a graph and reporting its interval undercovers at the benign crossing for every pre-test level.
\item \textbf{Sets versus intervals} (\cref{prop:union}). The union of the Wald intervals of all candidate graphs is an honest confidence set with price bounded in $n$. The $\sqrt n$ lower bounds are therefore a property of intervals: an interval that covers two points contains the segment between them.
\end{enumerate}
\Cref{tab:evidence}, at the beginning of the appendix, lists the evidence behind each result.

\paragraph{Inference after causal discovery.} \citet{chang2026post} resample the discovery step of the PC algorithm and screen the resampled graphs. Their intervals have asymptotically valid coverage of the true causal effect for each fixed distribution, and they note that the inflation of the width is difficult to characterize theoretically. \citet{gradu2025valid} use differentially private selection and target a projection parameter that depends on the selected graph. \citet{guo2026price} invert a truncated Gaussian pivot on the selection event of FCI. With known noise variance, this gives exact selective coverage of the selected adjustment functional when the adjustment set contains the outcome's Markov blanket, a condition they note generically fails. It also gives asymptotic coverage of the structural effect at a fixed distribution with a positive faithfulness gap (their Theorem~5). They compute the width ratio $1+\Theta(n^{-1/4})$ of one sample-splitting allocation relative to a known-graph oracle (their Proposition~4). \citet{strieder2023confidence} invert likelihood ratio tests in models with equal error variances and obtain finite-sample validity without regularity conditions, at the price of conservativeness and of enumerating graphs. \citet{strieder2024dual} simplify the test inversion by a dual likelihood, and \citet{strieder2024identifying} study identification and confidence regions under partial homoscedasticity. \citet{bhattacharya2026robust} avoid selection by weighting the functionals identified by several candidate models with data-driven measures of their validity; they target the pooled functional and bound its distance to the true effect. \citet{wang2026confidence} construct confidence sets for causal orderings. \citet{edgar2026inferential} construct confidence sets of causal models and study which models enter them when the parameters are at the borderline of detectability, with inclusion probabilities that are stable, slowly decaying or rapidly decaying; see also \citet{battey2026post}. \citet{kang2026interpretable} construct confidence regions for the CPDAG itself from perturbed runs of the PC algorithm. \citet{malinsky2024cautious} proposes a cautious constraint-based selection that limits wrongly deleted edges and favours denser graphs when the graph is used to estimate an effect.

\paragraph{Related work.} The minimax theory of confidence interval length goes back to \citet{pratt1961length}; lower bounds through a modulus of continuity are due to \citet{donoho1991geometrizingII} and \citet{donoho1994statistical}, and \citet{low1997nonparametric}, \citet{cai2004adaptation} and \citet{cai2017confidence} characterized the length of honest, adaptive and high-dimensional intervals. Our lower bounds are of this two-point or Pratt type (\cref{sec:setup}), except a dual bound in the four-variable model, which uses a numerically computed least favourable distribution as in \citet{elliott2015nearly}. In this language, the price is governed by the local modulus of the effect over the class; the new part is the computation of the modulus for causal graphical models, its graphical characterization and its attainment without knowledge of the base point. Intervals after model selection are studied by \citet{leeb2005model}, \citet{potscher2009confidence} and \citet{bachoc2020uniformly}, uniformly valid intervals with invalid instruments by \citet{guo2023causal}, inference where identification degenerates by \citet{andrews2012estimation} and \citet{bei2026inference}, impossibility results for non-differentiable functionals by \citet{hirano2012impossibility}, and partially identified parameters by \citet{imbens2004confidence} and \citet{stoye2009more}; at a malignant base point our identified set is finite. Locally, each pair of explanations gives a canonical experiment in which the target is one of two linear functionals on two lines through a common subspace (\cref{prop:canonical}). This is close to the local misspecification framework of \citet{hjort2003frequentist}, to intervals that use uncertain prior information \citep{kabaila2009confidence} and to adaptation to an unknown bias \citep{armstrong2025adapting}, and our weak-coefficient bound $A$ in \cref{thm:uniform} plays a role similar to the bound on the control coefficients in the bias-aware intervals of \citet{armstrong2020bias}. There are three differences. The union of the two models is not convex, so the affine-procedure theory of \citet{donoho1994statistical} does not apply; the target is held to the oracle of whichever model is true; and which pairs occur, and their constants, are determined by the graph. For bias-aware intervals with a known bound on the bias, \citet{freeman2026local} shows that the induced tests lose no local power, a bounded explicit amount, or all of it, according to the order of the bound relative to the standard error. Our benign and malignant cases have a similar structure, but the bias comes from structural uncertainty, and its size and its detectability are determined by the graph. \citet{zhang2003strong} and \citet{spirtes2014uniformly} replaced faithfulness by strong and $k$-triangle faithfulness to obtain uniformly consistent causal discovery, and we use the detectability result of \citet{zhang2008detection} to show that explanations requiring cancellation violate triangle faithfulness. \citet{evans2020model} showed that models with the same adjacencies are tangent at the point of total independence, in contrast with \cref{thm:E}. The graphical criteria build on \citet{meek1995causal}, \citet{perkovic2017interpreting}, \citet{perkovic2020identifying} and the IDA idea of \citet{maathuis2009estimating}; efficient adjustment is studied by \citet{henckel2022graphical} and \citet{guo2022efficient}.

\paragraph{Organization.} \Cref{sec:setup} sets up the problem and proves the result on sets. \Cref{sec:toy} studies the four-variable model. \Cref{sec:general} gives the general classification. \Cref{sec:sim} reports simulations and \cref{sec:discussion} discusses limitations. Proofs are in the appendices.

\section{Setup}\label{sec:setup}

\subsection{Models, candidate graphs and targets}

Let $V=(V_1,\dots,V_p)$ be a random vector. For a directed acyclic graph (DAG) $D$ on $\{1,\dots,p\}$, the linear Gaussian model $\cM_D$ is the set of distributions $N(\mu,\Sigma)$ with $\Sigma$ positive definite and
\[
V_j=\mu_j'+\sum_{i\in\pa_D(j)}\lambda_{ij}V_i+\eps_j,\qquad \eps_j\sim N(0,\omega_j)\ \text{independent},
\]
with arbitrary coefficients $\lambda_{ij}\in\R$ (zero is allowed) and noise variances $\omega_j>0$. Equivalently, $\cM_D$ is the set of positive definite Gaussian distributions that satisfy the local Markov property of $D$. Two DAGs are Markov equivalent, $D\simeq D'$, if they have the same skeleton and the same unshielded colliders; Markov equivalent DAGs define the same model.

We fix a treatment $X$ and an outcome $Y$, and a skeleton $S$. Background knowledge $\cK$ is a set of required orientations; in \cref{sec:general} $\cK$ consists of the requirement that $Y$ is a sink and of optional required directed edges. The candidate class $\cG=\cG(S,\cK)$ is the set of acyclic orientations of $S$ that satisfy $\cK$. Because coefficients may vanish, $S$ is a known \emph{superset} of the adjacencies of the true causal DAG: an edge of $S$ with zero coefficient is an absent edge. Our general results are local near distributions at which exactly one edge of $S$ is absent or weak (\cref{sec:general}); several simultaneously weak or absent edges are not covered (\cref{sec:discussion}).

For $D\in\cG$ and $P\in\cM_D$, the total effect $\theta_D(P)$ of $X$ on $Y$ is the sum over directed paths from $X$ to $Y$ in $D$ of the products of the coefficients. When $Y\notin\pa_D(X)$, the parents of $X$ form a valid adjustment set, so
\begin{equation}\label{eq:adj}
\theta_D(P)=\thadj_D(P):=\bigl[\Sigma_{AA}^{-1}\Sigma_{AY}\bigr]_X,\qquad A=\{X\}\cup\pa_D(X),
\end{equation}
the coefficient of $X$ in the population regression of $Y$ on $A$ \citep{maathuis2009estimating}. The right-hand side of \eqref{eq:adj} is defined for every positive definite $\Sigma$, also when $P\notin\cM_D$; we use it as the estimand of the candidate $D$ in all procedures.

The data are $n$ independent copies of $V$ from some $P\in\cM_{G^*}$, where the true graph $G^*\in\cG$ is unknown. The target is $\theta_{G^*}(P)$. A parameter class $\mathcal P$ is a set of pairs $(G^*,P)$ with $G^*\in\cG$ and $P\in\cM_{G^*}$. The same distribution $P$ may belong to several models, in which case the pairs $(G^*,P)$ and $(G',P)$ may have different targets; this is the source of the difficulty.

\subsection{Honesty, length and price}

A confidence set is a measurable map $C_n$ from the data to Borel subsets of $\R$, possibly using external randomization. It is \emph{honest} at level $1-\alpha$ over $\mathcal P$ if
\[
\inf_{(G^*,P)\in\mathcal P}\Prob_P\bigl(\theta_{G^*}(P)\in C_n\bigr)\ \ge\ 1-\alpha
\]
for every $n$ (finite-sample honesty), or asymptotically honest if the $\liminf_n$ of the left side is at least $1-\alpha$. Its expected length is $\E_P\Leb(C_n)$.

The benchmark is an analyst who knows the true graph. Let $s_{\mathrm{eff}}^2(G^*,P)$ be the efficient (Cram\'er--Rao) variance bound for $\theta_{G^*}$ in the regular parametric model $\cM_{G^*}$ at $P$, and let the \emph{oracle length} be
\[
\orl_n(G^*,P)=2z_{1-\alpha/2}\,s_{\mathrm{eff}}(G^*,P)\,n^{-1/2},
\]
the length of the efficient Wald interval. The \emph{price} of $C_n$ over $\mathcal P$ is
\[
\Pi_n(C_n;\mathcal P)=\sup_{(G^*,P)\in\mathcal P}\frac{\E_P\Leb(C_n)}{\orl_n(G^*,P)},
\]
and the minimax price is the infimum over honest $C_n$. Because the oracle is attached to the pair $(G^*,P)$, a distribution that lies in two models is held to both oracle lengths.\footnote{\citet{guo2026price} use ``price'' for the width ratio of one sample-splitting procedure relative to a known-graph oracle at a fixed distribution. Here the price is a worst case over a class and is minimized over all honest procedures.}

The parent-adjustment interval of the true graph, the $t$-interval for \eqref{eq:adj} with $D=G^*$, is what a practitioner would compute; it is not efficient in general \citep{henckel2022graphical,guo2022efficient}. The two benchmarks differ by a bounded factor, so none of the rates below depends on the choice.

\begin{lemma}[Oracle comparison]\label{lem:oracle}
Let $Y$ be a sink of $G^*$, $X\to Y\in G^*$, and let the eigenvalues of $\Sigma_P$ lie in $[\sigl^2,\sigu^2]$. Let $s_{\mathrm{pa}}^2=\sigma^2_{Y\mid A}/\sigma^2_{X\mid\pa_{G^*}(X)}$, $A=\{X\}\cup\pa_{G^*}(X)$, be the asymptotic variance of the parent-adjustment estimator. Then $\sigl/\sigu\le s_{\mathrm{eff}}\le s_{\mathrm{pa}}\le\sigu/\sigl$; in particular $1\le s_{\mathrm{pa}}/s_{\mathrm{eff}}\le(\sigu/\sigl)^2$.
\end{lemma}

The proof is in \cref{app:oracle}. In the four-variable model of \cref{sec:toy}, where the eigenvalue bound is not available, the ratio is explicit (\cref{lem:oracle-toy}).

\subsection{Sets versus intervals}

The first observation is that for confidence \emph{sets} structural uncertainty is cheap.

\begin{proposition}[Union of candidate intervals]\label{prop:union}
Let $\mathcal P$ be the class of pairs $(G^*,P)$ with $G^*\in\cG$, $P\in\cM_{G^*}$ and eigenvalues of $\Sigma$ in $[\sigl^2,\sigu^2]$. For each distinct adjustment set $A\in\{\{X\}\cup\pa_D(X):D\in\cG\}$ let $W_A$ be the level-$(1-\alpha)$ $t$-interval for the $X$-coefficient of the regression of $Y$ on $A$, and let $C_n^\cup=\bigcup_A W_A$. Suppose $Y$ is a sink in every $D\in\cG$ and $S$ contains the edge $X-Y$. Then for every $n\ge p+3$ and every $(G^*,P)\in\mathcal P$,
\begin{enumerate}
\item[(i)] $\Prob_P(\theta_{G^*}(P)\in C_n^\cup)\ge 1-\alpha$;
\item[(ii)] $\E_P\Leb(C_n^\cup)\le N_{\mathrm{pa}}\,(\sigu/\sigl)^2\,\kappa_n\,\orl_n(G^*,P)$,
\end{enumerate}
where $N_{\mathrm{pa}}$ is the number of distinct adjustment sets and $\kappa_n\to1$ is explicit (\cref{app:union}).
\end{proposition}

The proof (\cref{app:union}) uses only that $W_{A^*}$ is an exact $t$-interval for the true graph and that the regression of $Y$ on any fixed set is exactly normal given the design under joint Gaussianity, whether or not the candidate is correct (\cref{lem:M}). A disconnected confidence region of this kind is reported, for example, by the procedures of \citet{strieder2023confidence}.

For an interval the situation is different. If two pairs $(G_1,P)$ and $(G_2,Q)$ in $\mathcal P$ have targets $\theta_1\ne\theta_2$, an honest interval must cover both targets with high probability under nearby data distributions, and an interval that covers both points contains the segment between them. The following two-point bound is standard \citep{pratt1961length,donoho1994statistical,cai2004adaptation}.
\begin{lemma}[Two-point bound]\label{lem:two-point}
If $C_n$ is an interval that is honest at level $1-\alpha$ over a class containing $(G_1,P)$ and $(G_2,Q)$, then
\[
\E_P\Leb(C_n)\ \ge\ |\theta_{G_1}(P)-\theta_{G_2}(Q)|\,\bigl(1-2\alpha-\dtv(P^n,Q^n)\bigr)_+ .
\]
\end{lemma}
\begin{proof}
Let $A_P=\{\theta_{G_1}(P)\in C_n\}$ and $A_Q=\{\theta_{G_2}(Q)\in C_n\}$. Honesty gives $P^n(A_P)\ge1-\alpha$ and $P^n(A_Q)\ge Q^n(A_Q)-\dtv(P^n,Q^n)\ge 1-\alpha-\dtv(P^n,Q^n)$, so $P^n(A_P\cap A_Q)\ge 1-2\alpha-\dtv(P^n,Q^n)$. On $A_P\cap A_Q$ the interval contains both targets.
\end{proof}
All our lower bounds on lengths of intervals, except a dual bound in the four-variable model (\cref{thm:lfp}), are of this form or refine it by Pratt's integral; the dual bound holds for arbitrary measurable confidence sets and hence for intervals. From now on, unless stated otherwise, a confidence set is an interval and the price is that of intervals.

\paragraph{Interpretation.} \Cref{lem:two-point} reduces lower bounds to a between-model modulus of continuity: the largest difference of targets between pairs in the class at Hellinger distance of order $n^{-1/2}$ \citep{donoho1991geometrizingII,donoho1994statistical,low1997nonparametric,cai2004adaptation}. In this language, a malignant edge (\cref{prop:malignant}) is one at which the modulus does not vanish as the distance tends to zero: the class contains two pairs with the same distribution and different targets. At a benign edge the local minimax length is of the parametric order (\cref{thm:K}), the behaviour of a locally linear modulus; the key step is that the error of a wrong explanation is locally Lipschitz in the statistics that test it (\cref{lem:K4E}, via \cref{thm:E}). The four-variable model with known noise variances has a square-root modulus, which produces the $n^{-1/4}$ plateau of \cref{sec:toy-rate}. What is new is the computation of these moduli for causal graphical models: we read them off graphs, we include the identification that triangle faithfulness adds (\cref{rem:ktf-identifies}), we attain them with honest procedures that do not know the base point, and at benign edges we determine the slope within an explicit factor (\cref{cor:RK}). At a malignant base point the identified set is a finite set of candidate effects; \cref{prop:union} reports this set, as in the partial identification literature \citep{imbens2004confidence,stoye2009more}, while an interval must cover its convex hull. We study intervals because a single interval is what is reported for a scalar effect, and because the comparison with sets isolates the part of the price that is due to convexity.

\section{A four-variable model: confounder or mediator}\label{sec:toy}

\subsection{The model}\label{sec:toy-model}

Let $V=(W,X,Z,Y)$ with skeleton $W-X$, $X-Z$, $Z-Y$, $X-Y$, and background knowledge $\cK=\{W\to X,\ Y\text{ is a sink}\}$. The candidate class has two graphs:
\begin{align*}
C\ (\text{$Z$ confounds}):&\quad X=wW+aZ+\eps_X,\quad Z=\eps_Z,\quad Y=\beta X+bZ+\eps_Y,\quad \theta_C=\beta;\\
M\ (\text{$Z$ mediates}):&\quad X=wW+\eps_X,\quad Z=aX+\eps_Z,\quad Y=\beta X+bZ+\eps_Y,\quad \theta_M=\beta+ab,
\end{align*}
with $W=\eps_W$ and independent centred Gaussian noises with variances $\sigma_W^2,\sigma_X^2,\sigma_Z^2,\sigma_Y^2$. $C$ has the unshielded collider $W\to X\leftarrow Z$ and $M$ does not, so the two graphs are not Markov equivalent. For fixed constants $\bar A>0$, $\bar w>0$ and $\kappa\ge1$, and for $\lambda\in(0,\bar w]$, $A\in(0,\bar A]$ and $B>0$, let $\cM(\lambda,A,B,\kappa)$ be the class of pairs $(g,P)$, $g\in\{C,M\}$, with $|w|\in[\lambda,\bar w]$, $|a|\le A$, $|b|\le B$, $\beta\in\R$ and all noise variances in $[1,\kappa]$.

Two kinds of crossing occur. At $a=0$ the two graphs give the same distribution and the same effect; the effect is continuous across this crossing. At $w=0$ the two graphs become Markov equivalent (the collider disappears) while their effects still differ by $ab$. In the terminology of \cref{sec:general}, the edge $X-Z$ is a benign weak edge and the edge $W\to X$ is a malignant one: deleting $X-Z$ leaves an MPDAG in which the effect is identifiable, while deleting $W\to X$ leaves the undirected edge $X-Z$ on the path $X-Z\to Y$. The model is the instance of \cref{sec:general} with $R=\{W\to X\}$, a required directed edge. The classification of $X-Z$ depends on it: without it, $W\to X$ is still malignant, but $X-Z$ is not benign, because the chains $W\to X\to Z$ and $Z\to X\to W$, each with $Z\to Y$ and $X\to Y$, are Markov equivalent explanations of the base point with effects $\beta+ab$ and $\beta$, so \eqref{eq:K1p} fails and $X-Z$ falls in the third case of \cref{thm:dichotomy}. The orientation of $W-X$ can also come from a second collider, as in the graph with $W\to X\leftarrow U$ added, in which $X-Z$ is benign even with $R=\emptyset$.

The quantity that governs the malignant side is the largest $a$ for which a mediator model can be matched exactly by a confounder model with noise variances in $[1,\kappa]$:
\begin{equation}\label{eq:abar}
\abar(\kappa)^2=\begin{cases}1-1/\kappa,&1\le\kappa\le2,\\ \kappa/4,&\kappa\ge2.\end{cases}
\end{equation}
\Cref{app:abar} derives \eqref{eq:abar} as the value of an elementary optimization over the noise variances; the maximum is attained at $(\sigma_X^2,\sigma_Z^2)=(2,\kappa/2)$ for $\kappa\ge2$ and $(\kappa,1)$ for $\kappa<2$. With unit variances ($\kappa=1$) no exact matching is possible and $\abar=0$.

The oracle of \cref{sec:setup} is explicit in this model. Because $\beta$ is unbounded, the eigenvalue bound of \cref{lem:oracle} is not available, but the ratio to the parent-adjustment interval does not depend on $\beta$:
\begin{lemma}[Oracle in the four-variable model]\label{lem:oracle-toy}
Under both $C$ and $M$, $s_{\mathrm{pa}}^2/s_{\mathrm{eff}}^2=1+w^2\sigma_W^2/\sigma_X^2\in[1,1+\bar w^2\kappa]$; moreover $s_{\mathrm{eff}}^2=\sigma_Y^2/(w^2\sigma_W^2+\sigma_X^2)$ under $C$ and $s_{\mathrm{eff}}^2=(\sigma_Y^2+b^2\sigma_Z^2)/(w^2\sigma_W^2+\sigma_X^2)$ under $M$.
\end{lemma}

\subsection{The rate of the minimax length}\label{sec:toy-rate}

Let
\[
R_n=n^{-1/2}+\min\Bigl\{\frac{B}{\lambda\sqrt n},\ B\min\bigl(A,\sqrt{\abar(\kappa)^2+n^{-1/2}}\bigr)\Bigr\}.
\]

\begin{theorem}[Lower bound]\label{thm:rate-lower}
Fix $\bar w$, $\bar A$, $\kappa\ge1$ and $\alpha<1/2$. There are $c>0$ and $n_0$ such that for all $n\ge n_0$, $\lambda\in(0,\bar w]$, $A\in(0,\bar A]$ and $B>0$, every interval that is honest at level $1-\alpha$ over $\cM(\lambda,A,B,\kappa)$ satisfies
\[
\sup_{(g,P)\in\cM(\lambda,A,B,\kappa)}\E_P\Leb(C_n)\ \ge\ c\,R_n .
\]
\end{theorem}

The proof (\cref{app:rate-lower}) uses \cref{lem:two-point} with three pairs: two parameters of one graph that differ in $\beta$ by $n^{-1/2}$; a mediator model with $w=\lambda$ against the confounder model with the same coefficients; and, when $\kappa>1$, a mediator model with least favourable noise variances against a confounder model that matches all covariances except $\Cov(W,Z)$ and $\Cov(W,Y)$. In the last two pairs the conditional law of $Y$ given $(W,X,Z)$ is the same under both members, so the total variation distance does not depend on $b$, while the effects differ by $|ab|$.

The matching upper bound is attained by a procedure that screens each candidate in two ways.

\begin{definition}[Procedures {\normalfont\DSV} and {\normalfont\DSVp}]\label{def:dsv}
Let $\gamma_1+\gamma_2+\alpha'=\alpha$. \DSV: (1) \emph{Conditional independence screening}: reject $C$ if the $t$-test of $W\indep Z$ rejects at level $\gamma_1/2$, and reject $M$ if the $t$-test of $W\indep Z\mid X$ rejects at level $\gamma_1/2$. (2) \emph{Variance screening}: for each candidate $g$ and each node $j$, let $\hat\sigma^2_{j,g}$ be the residual variance of the regression of $j$ on $\pa_g(j)$ with intercept and $m_{j,g}=n-|\pa_g(j)|-1$ degrees of freedom; reject $g$ if some $\hat\sigma^2_{j,g}$ lies outside $[\chi^2_{m,\eta/2}/m,\ \kappa\chi^2_{m,1-\eta/2}/m]$ with $\eta=\gamma_2/8$. (3) Output the convex hull of the level-$(1-\alpha')$ $t$-intervals \eqref{eq:adj} of the candidates that are not rejected, and the empty set if both are rejected.
\DSVp: if $A\le1/(\bar w\sqrt n)$, output $\hat\theta_C\pm(t_{n-4,1-\alpha/2}\,\widehat{\mathrm{se}}_C+AB)$, where $\hat\theta_C$ is the coefficient of $X$ in the regression of $Y$ on $(1,X,W,Z)$; otherwise output \DSV.
\end{definition}

The fallback in \DSVp\ uses that the regression of $Y$ on $(X,W,Z)$ has $X$-coefficient $\beta$ under both graphs, while the targets differ from $\beta$ by at most $|ab|\le AB$. Like the oracle and the variance screen, \DSVp\ uses the parameters of the class as inputs.

\begin{theorem}[Upper bound]\label{thm:rate-upper}
\DSV\ and \DSVp\ are honest at level $1-\alpha$ over $\cM(\lambda,A,B,\kappa)$ for every $n\ge5$. Fix $(\bar w,\bar A,\kappa,\alpha,\gamma_1,\gamma_2)$. There are $C<\infty$ and $n_0$ such that for all $n\ge n_0$, $\lambda\in(0,\bar w]$, $A\in(0,\bar A]$ and $B>0$,
\[
\sup_{(g,P)\in\cM(\lambda,A,B,\kappa)}\E_P\Leb(C_n^{\DSVp})\ \le\ C\,R_n ,
\]
and the same bound holds for \DSV\ whenever $A\ge1/(\bar w\sqrt n)$.
\end{theorem}

Coverage is exact in finite samples because, under joint Gaussianity, the true model's residual variances are exactly scaled $\chi^2$ and its conditional independence $t$-statistics exactly $t$-distributed. The length bound (\cref{app:rate-upper}) controls $2\Delta\,\Prob(\text{wrong model kept})$, with $\Delta\le\kappa B|a|$: the conditional independence screen bounds it by $B/(\lambda\sqrt n)$, and the variance screen by $B\sqrt{\abar^2+n^{-1/2}}$, because if $a^2$ exceeds $\abar(\kappa)^2$ the wrong model's fitted noise variances leave $[1,\kappa]$ (\cref{lem:Vprime}).

Together, \cref{thm:rate-lower,thm:rate-upper} give $L_n^*\asymp R_n$ for each fixed $\kappa$, uniformly over $\lambda\in(0,\bar w]$, $A\in(0,\bar A]$ and $B>0$. Three regimes appear.
\begin{itemize}
\item If $\lambda$ is fixed, $L_n^*\asymp n^{-1/2}(1+B/\lambda)$: the length is of parametric order, but its constant grows like $B/\lambda$ and does not improve with $n$.
\item If $\lambda\sqrt n$ is bounded, the length does not shrink: it is of order $B\min(A,\abar(\kappa))$ when $\kappa>1$ (for $\kappa=1$ see the next item). The known range of the noise variances limits the malignant region through $\abar(\kappa)$.
\item With unit noise variances ($\kappa=1$) the plateau is $Bn^{-1/4}$; this is the doubly weak corner $w\asymp n^{-1/2}$, $a\asymp n^{-1/4}$.
\end{itemize}
The plateau is due to the known range of the noise variances: on general graphs, where only the eigenvalues of $\Sigma$ are restricted, \cref{thm:uniform} gives the same dependence on $(n,\lambda,A)$ with $\min(A,\abar)$ replaced by $A$. The constants $c,n_0$ of \cref{thm:rate-lower} depend on $\kappa$ and are not uniform as $\kappa\downarrow1$, where the threshold $w_0(\kappa)$ of the third pair in its proof tends to zero. The weak adjacency edge alone is harmless:
\begin{corollary}[Benign weak adjacency]\label{cor:benign-toy}
Fix $w_0\in(0,\bar w]$. Over $\cM(w_0,A,B,\kappa)$, with no margin on $a$, \DSVp\ is honest for every $n$ and its price is at most $C(1+B)$ for $n\ge n_0$, with $C,n_0$ depending only on $(w_0,\bar w,\bar A,\kappa,\alpha,\gamma_1,\gamma_2)$.
\end{corollary}
\begin{proof}
Take $\lambda=w_0$ in \cref{thm:rate-upper}: $R_n\le n^{-1/2}(1+B/w_0)$. By \cref{lem:oracle-toy}, $s_{\mathrm{eff}}^2\ge1/(\kappa(1+\bar w^2\kappa))$, so the oracle length is at least of order $n^{-1/2}$.
\end{proof}
The price is bounded in $n$, not in the effect scale: the factor $1+B$ is necessary, as the second pair in the proof of \cref{thm:rate-lower}, read with the confounder member as the reference distribution, shows.

\subsection{The price constant at the benign crossing}\label{sec:toy-constant}

With unit noise variances, the local limit experiment at $a=0$ is the canonical experiment of \cref{prop:canonical} with parameters $(1,b,0,1/\sqrt{1+w^2})$ (\cref{prop:limit}). At $w=b=1$ and $\alpha=0.05$ a computer-assisted dual bound gives its minimax price $\Pi^*(1,1)\ge1.2743$, against $1.2019$ from Pratt's two-point bound, and an interval procedure has price about $1.422$ numerically, with sets almost no shorter; so at the benign crossing the restriction to intervals costs little, unlike the malignant case (\cref{prop:union}). \Cref{app:toy-constant} gives the bound at fourteen values of $(w,b)$.

\subsection{Selecting a graph is not honest}\label{sec:toy-naive}

Selecting the graph whose constraint test has the larger $p$-value, or keeping one graph unless a pre-test of its constraint rejects at level $\alpha_t$, and reporting its $95\%$ interval undercovers in the local limit experiment, for $b\ne0$ at every $\alpha_t$; at $w=b=1$ the worst-case limiting coverage is at most $0.89$ and $0.73$ (\cref{prop:naive,prop:naive-more} and \cref{tab:naive-full} in \cref{app:naive}; see also \citealp{leeb2005model,kabaila2006large,potscher2009confidence}).

\section{General graphs}\label{sec:general}

\subsection{Setting}\label{sec:gen-setting}

Throughout this section the background knowledge $\cK$ consists of the requirement that $Y$ is a sink and of a possibly empty set $R$ of required directed edges of $S$ (a known orientation, a time order), the skeleton $S$ contains the edge $X-Y$, and $\cG=\cG(S,\cK)$ is the set of \emph{all} acyclic orientations of $S$ that satisfy $\cK$, assumed non-empty. All graphs in $\cG$ have the same skeleton, so for $D\in\cG$ and an edge $e$ of $S$ we write $D-e$ for $D$ with the edge between the endpoints of $e$ removed, whatever its orientation in $D$. For a DAG $D$ on $S$ or on $S-e$, $[D]_Y$ denotes the set of DAGs on the same skeleton that are Markov equivalent to $D$ and satisfy $\cK$ (with $e$ removed from $R$ if $D$ lives on $S-e$); for $R=\emptyset$, these are the DAGs equivalent to $D$ in which $Y$ is a sink. The weak edge $e$ is required to be \emph{admissible}: every $D\in[H]_Y$ becomes a member of $\cG$ when $e$ is added, with the orientation required by $R$ if $e\in R$. This is automatic if $e\notin R$. We call $\cK$ \emph{tier knowledge} if there is an ordered partition $T_1,\dots,T_m$ of $V\setminus\{Y\}$, $m\ge1$, such that $R$ is the set of edges of $S$ between different tiers, each oriented towards the later tier; $m=1$ gives $R=\emptyset$.

\paragraph{Parameter class.} Following \citet{spirtes2014uniformly}, a set $W\subseteq V\setminus\{u,v\}$ is \emph{legal} for a triangle $\{u,v,c\}$ of a DAG if $c\notin W$ when $c$ is a non-collider on $u-c-v$ and $c\in W$ when $c$ is a collider. For $k>0$ let $\cM_k$ be the class of pairs $(G^*,P)$ with $G^*\in\cG$, $P\in\cM_{G^*}$, all eigenvalues of $\Sigma$ in $[\sigl^2,\sigu^2]$, and
\[
|\rho_P(u,v\mid W)|\ \ge\ k\,|\lambda^{\mathrm{std}}_{G^*}(u\to v;P)|
\]
for every edge $u\to v$ of a triangle $\{u,v,c\}$ of $G^*$ and every legal $W$. Here $\lambda^{\mathrm{std}}_{G^*}(u\to v;P)=\beta\,\sigma_u/\sigma_v$, $\beta$ is the coefficient of $u$ in the regression of $v$ on $\pa_{G^*}(v)$, and $\sigma_u,\sigma_v$ are marginal standard deviations. This is the $k$-triangle-faithfulness condition, which forbids exact and near cancellation along triangles but allows arbitrarily weak edges and arbitrary cancellations along longer paths.

\paragraph{Weak edge and base point.} Fix $G\in\cG$ and an edge $e$ of $G$, and let $H=G-e$. Let $\Theta_H=\R^{|H|}\times(0,\infty)^p$ be the structural parameters (coefficients and noise variances) of $H$ and $P(\phi)$ the corresponding Gaussian law. A base point is $P_0\in\cM_H$ that is faithful to $H$; it is the limit of models in $\cM_G$ as the coefficient of $e$ tends to zero. ``Almost every base point'' means $P_0=P(\phi_0)$ for Lebesgue-almost every $\phi_0\in\Theta_H$. The candidates that explain $P_0$ are
\begin{gather*}
\cX=\{G'\in\cG:P_0\in\cM_{G'}\},\qquad \Xcanc=\cX\setminus\Xeq,\\
\Xeq=\{G'\in\cG:G'-e\text{ is Markov equivalent to }H\}.
\end{gather*}
Explanations in $\Xcanc$ need a cancellation; those in $\Xeq$ differ from $G$ only in the orientation of $e$ and in the orientations that $e$ leaves undetermined. The local price of intervals $C_n$ near $P_0$ is
\[
\Pi^{\mathrm{loc}}_n(C_n;r)=\sup\Bigl\{\frac{\E_P\Leb(C_n)}{\orl_n(G^*,P)}:(G^*,P)\in\cM_k,\ \|\Sigma_P-\Sigma_{P_0}\|\le r\Bigr\}.
\]

\begin{lemma}[Explanation sets are graphical]\label{lem:G}
If $P_0$ is faithful to $H$, then $G'\in\cX$ if and only if every d-separation of $G'$ is a d-separation of $H$. Moreover $\Xeq\subseteq\cX$. Hence $\cX$, $\Xeq$ and $\Xcanc$ depend only on $(G,e)$, not on the faithful base point.
\end{lemma}

\subsection{Two identification conditions and their graphical criteria}\label{sec:gen-criterion}

Let $\MP(H)$ be the maximally oriented partially directed acyclic graph (MPDAG) that represents $[H]_Y$ \citep{meek1995causal,perkovic2017interpreting}. In an MPDAG, a \emph{bad path} is a proper possibly causal path from $X$ to $Y$ that starts with an undirected edge; that is, a path $\langle X=V_0,\dots,V_k=Y\rangle$ on which $X$ occurs only once, the first edge is undirected, and there is no edge $V_i\leftarrow V_j$ with $i<j$ \citep{perkovic2020identifying}. The first condition concerns the explanations of the base point:
\begin{equation}\tag{K2}\label{eq:K2}
\thadj_{G'}(P)=\theta_H(P)\quad\text{for all }G'\in\Xeq\text{ and all }P\in\cM_H .
\end{equation}
It says that all explanations of the limit point that differ only in the orientation of $e$ estimate the same effect, identically on $\cM_H$. By analyticity, \eqref{eq:K2} holds on all of $\cM_H$ as soon as it holds on a non-empty relatively open subset. The following lemma is a translation of \citet[Theorem~3.6 and Proposition~3.2]{perkovic2020identifying}.

\begin{lemma}[Graphical criterion for \eqref{eq:K2}]\label{prop:graph}
Let $e$ be admissible.
(a) $\{G'-e:G'\in\Xeq\}=[H]_Y$.
(b) \eqref{eq:K2} holds if and only if $\theta_D=\theta_H$ on $\cM_H$ for every $D\in[H]_Y$.
(c) \eqref{eq:K2} holds if and only if the total effect of $X$ on $Y$ is identifiable in $\MP(H)$, that is, if and only if $\MP(H)$ has no bad path.
\end{lemma}

The second condition concerns explanations that are Markov equivalent to each other. Partition $\Xeq$ into the classes $S_c=\{G'\in\Xeq:G'\simeq G_1\}$ of Markov equivalent members.
\begin{equation}\tag{K1$'$}\label{eq:K1p}
\text{for every class }S_c\text{ and all }G_1,G_2\in S_c:\quad \thadj_{G_2}=\theta_{G_1}\ \text{on }\cM_{G_1}.
\end{equation}
Identifiability of the effect in $\MP(G')$ for every $G'\in\Xeq$ implies \eqref{eq:K1p} (\cref{cor:K1}) but is not necessary, because $S_c$ is in general a proper subset of $[G_1]_Y$ (\cref{app:certify}). With $e=\{u,v\}$ and $\mathrm{CC}_D=\{x:u\to x\leftarrow v\text{ in }D\}$ the set of common children of $u$ and $v$ in $D$, $S_c=\{D\in[G_1]_Y:\mathrm{CC}_D=\mathrm{CC}_{G_1}\}$ (\cref{lem:C}). Since ``$x\notin\mathrm{CC}_D$'' is a disjunction of two orientations, $S_c$ need not be representable by an MPDAG, and Perkovi\'c's criterion does not apply to it directly. The next theorem shows that, under \eqref{eq:K2}, it does apply after the conjunctive part of the constraint and the orientations of $\MP(H)$ are added as background knowledge.

\begin{theorem}[Graphical criterion for \eqref{eq:K1p}]\label{thm:K1p-graph}
Assume that $e$ is admissible and that \eqref{eq:K2} holds. For a class $S_c$ with representative $G_1$, let $\MP''(S_c)$ be the MPDAG of $[G_1]_Y$ with the background knowledge
\[
\{u\to x,\ v\to x:\ x\in\mathrm{CC}_{G_1}\}\ \cup\ \{\text{directed edges of }\MP(H)\}.
\]
$\MP''(S_c)$ does not depend on the representative. \eqref{eq:K1p} holds if no $\MP''(S_c)$ has a bad path; and if $\MP(H)$ is a chain graph, the converse holds, so that \eqref{eq:K1p} holds if and only if no $\MP''(S_c)$ has a bad path.
\end{theorem}

The proof (\cref{app:K1p}) shows that the members of the MPDAG that are not in $S_c$ create an extra common child $x$ of $u$ and $v$; replacing $e$ by the detour $u-x-v$, or by a shortcut through $x$, turns any bad path into a bad path of $\MP(H)$, which \eqref{eq:K2} excludes. Both conditions are thus decided by one Meek closure per class, given the classes of $\Xeq$. The chain hypothesis can fail with general required edges (\cref{app:certify}), but not for the most common kind of them.

\begin{lemma}[Tier knowledge]\label{lem:tier}
If $\cK$ is tier knowledge, every edge of $S$ is admissible, and for every $G\in\cG$ and every edge $e$ of $G$, $\MP(H)$ is a chain graph with chordal chain components.
\end{lemma}

With $Y$ in a last tier of its own, $\MP(H)$ is a tiered MPDAG in the sense of \citet{bang2023wiser}, and the chain property is their Theorem~1 and Corollary~1. Our proof (\cref{app:K1p}) identifies $[H]_Y$ with an interventional Markov equivalence class of \citet{hauser2012characterization}, which also covers some required edges that do not come from a time order.

\subsection{The classification}\label{sec:gen-dichotomy}

The two criteria classify every admissible weak edge. Case 1 uses three quantities defined in \cref{sec:gen-benign}, where it is proved: the cross-model constant $K$, which compares the first-order effect discrepancy between two explanations with the information that separates them; the ratio $\mathcal R$ of the efficient lengths of two explanations; and the threshold $k^*$ on triangle faithfulness below which the lower bounds apply. Cases 2 and 3 are proved below, and \cref{app:dichotomy} combines the three.

\begin{theorem}[Classification of weak edges]\label{thm:dichotomy}
Consider linear Gaussian models with a known skeleton $S$, background knowledge $\cK$ as above and the candidate class $\cG$ of all orientations of $S$ satisfying $\cK$; assume that the weak edge $e$ is admissible and that $\MP(H)$ is a chain graph, which both hold for tier knowledge, in particular for $R=\emptyset$ (\cref{lem:tier}). Let $G\in\cG$, $e$ an edge of $G$ and $H=G-e$, and consider intervals and $\alpha<1/2$. There is a Lebesgue-null set $\cN'\subset\Theta_H$ such that for every $\phi_0\notin\cN'$, with $P_0=P(\phi_0)$ and the eigenvalues of $\Sigma_{P_0}$ in the interior of $[\sigl^2,\sigu^2]$, exactly one of the following holds.
\begin{enumerate}
\item \emph{Benign.} $\MP(H)$ has no bad path, and no $\MP''(S_c)$ has a bad path. Then for $k$ below a threshold $k^*(P_0)>0$ the local minimax price is bounded in $n$ and of order $\mathcal R(P_0)+K(P_0)$, within the explicit factors of \cref{cor:RK}, and equal to one if and only if $K(P_0)=0$ (\cref{thm:free}), and double screening attains this order for every $k>0$ without knowing $P_0$.
\item \emph{Malignant.} $\MP(H)$ has a bad path. For every $k<k_0(P_0)$, the threshold of \cref{prop:malignant} for a pair $G_1,G_2\in\Xeq$ with $\theta_{G_1}(P_0)\ne\theta_{G_2}(P_0)$, every interval that is honest over $\cM_k$ has price at least $c(P_0)\sqrt n$ at $P_0$.
\item \emph{Locally non-identified Markov class.} $\MP(H)$ has no bad path, but some $\MP''(S_c)$ has one. For every $k<k_1(P_0)$, the threshold of \cref{prop:K1p-nec}, and every fixed small $r>0$, every interval that is honest over $\cM_k$ has local price at least $c(P_0,k,r)\sqrt n$, where $c(P_0,k,r)\to0$ as $r\to0$.
\end{enumerate}
The thresholds $k_0(P_0),k_1(P_0),k^*(P_0)>0$ depend on the base point. Without the chain hypothesis, cases 1 and 2 hold as stated, and a weak edge for which $\MP(H)$ has no bad path but some $\MP''(S_c)$ has one is in case 1 or case 3 according to whether \eqref{eq:K1p} holds.
\end{theorem}

\paragraph{The non-identified side.}
The lower bounds are, mathematically, consequences of non-identification. The class contains two pairs with the same distribution and different targets, so \cref{lem:two-point} applies with zero total variation. The content is that this happens at generic points of the class, which $k$-triangle faithfulness does not exclude, and that whether it happens can be read off the graph.

\begin{proposition}[Malignant weak edges]\label{prop:malignant}
Let $P_0$ be faithful to $H$, and suppose there are $G_1,G_2\in\Xeq$ with $\Delta_0=|\theta_{G_1}(P_0)-\theta_{G_2}(P_0)|>0$. Let $k_0(P_0)$ be the smallest ratio $|\rho_{P_0}(u,v\mid W)|/|\lambda^{\mathrm{std}}_{G_i}(u\to v;P_0)|$ over triangle edges $u\to v\ne e$ of $G_1,G_2$ and legal $W$, with $k_0(P_0)=+\infty$ if there are no such edges. Then $k_0(P_0)>0$, and for $k<k_0(P_0)$:
\begin{enumerate}
\item[(i)] if the eigenvalues of $\Sigma_{P_0}$ lie in $[\sigl^2,\sigu^2]$, every interval that is honest over $\cM_k$ has $\E_{P_0}\Leb(C_n)\ge(1-2\alpha)\Delta_0$, so its price at $P_0$ is at least $c(P_0)\sqrt n$ with $c(P_0)=(1-2\alpha)\Delta_0\sigl/(2z_{1-\alpha/2}\sigu)$;
\item[(ii)] if the eigenvalues lie in the interior of $[\sigl^2,\sigu^2]$ and the class additionally requires $|\lambda^{\mathrm{std}}_e|\ge\lambda$, there are $t_1,L,C>0$ depending on $P_0$ such that for $\lambda\le t_1$, $\sup\E\Leb(C_n)\ge(\Delta_0-2L\lambda)(1-2\alpha-2C\lambda\sqrt n)_+$.
\end{enumerate}
\end{proposition}

The proof is \cref{lem:two-point} applied to $(G_1,P_0)$ and $(G_2,P_0)$; the point to check is that $k$-triangle faithfulness does not exclude them, which for triangles that contain $e$ follows from faithfulness of $P_0$ to $H$ (\cref{lem:M1}). Part (ii) gives nothing for fixed $\lambda$ and $\lambda\sqrt n\to\infty$; in that regime the relevant result is \cref{thm:K-lambda}.

When \eqref{eq:K2} holds but \eqref{eq:K1p} fails, the effects of two Markov equivalent members of $\Xeq$ agree on $\cM_H$ but not off it, and the price is large on every neighbourhood of $P_0$, although not necessarily at $P_0$ itself.

\begin{proposition}[Locally non-identified Markov class]\label{prop:K1p-nec}
Let $P_0$ be faithful to $H$ with the eigenvalues of $\Sigma_{P_0}$ in the interior of $[\sigl^2,\sigu^2]$, and let $G_1\simeq G_2$ in $\Xeq$ with $q=\thadj_{G_2}-\theta_{G_1}$ not identically zero on $\cM_{G_1}$. Let $k_1(P_0)>0$ be defined as $k_0(P_0)$ in \cref{prop:malignant}, with $G_1,G_2$. For every $k<k_1(P_0)$ there is $r_1>0$ such that for every $r\in(0,r_1]$, $c(r)=\sup\{|q(P)|:P\in\cM_{G_1},\|\Sigma_P-\Sigma_{P_0}\|\le r\}>0$ and every interval that is honest over $\cM_k$ satisfies, for $n\ge p+3$,
\[
\Pi^{\mathrm{loc}}_n(C_n;r)\ \ge\ \frac{(1-2\alpha)\,c(r)\,\sigl}{2z_{1-\alpha/2}\,\sigu}\,\sqrt n .
\]
\end{proposition}

The constant depends on $(P_0,k,r)$ and $c(r)\to0$ as $r\to0$; the statement is for a fixed neighbourhood. The proof (\cref{app:K1p-nec}) uses the pairs $(G_1,P)$ and $(G_2,P)$ for $P\in\cM_{G_1}=\cM_{G_2}$ near $P_0$, at which $q\ne0$.

\begin{remark}
The four-variable model of \cref{sec:toy} is the instance with $R=\{W\to X\}$, which is tier knowledge with tiers $\{W\}$ and $\{X,Z\}$; $X-Z$ is benign and $W\to X$ is malignant, in agreement with an exact rational computation (\cref{app:certify}); with $R=\emptyset$, $X-Z$ falls in case 3 (\cref{sec:toy-model}). On five nodes, of the $26988$ pairs $(G,e)$, $14302$ are malignant, $8426$ benign and $4260$ in case 3 (\cref{app:certify}). For a generic $P_0$, an edge is malignant exactly when the IDA estimates \citep{maathuis2009estimating} at $P_0$ computed over $[H]_Y$ are not all equal. By \cref{prop:union}, a confidence set that may be disconnected can avoid the $\sqrt n$ price in cases 2 and 3 by reporting the union of the candidate intervals.
\end{remark}

\begin{remark}[Local and global price]\label{rem:global}
The classification is local, and for a reason: the class $\cM_k$ allows any edge to be arbitrarily weak, so its global price is the largest of the local prices. If some edge of $S$ is malignant for some orientation, as $W\to X$ in \cref{sec:toy}, then the global price over $\cM_k$ is at least of order $\sqrt n$ for every $k$ below the threshold $k_0(P_0)$ of \cref{prop:malignant} at some base point $P_0$. A global bound would require every edge, and every set of simultaneously weak edges, to be benign; this set-wise condition is not implied by the single-edge criterion and is not studied here (\cref{sec:discussion}). The local price says, near a given point of the model, whether the structural uncertainty that remains after discovery costs a constant factor or a factor of order $\sqrt n$.
\end{remark}

\subsection{The benign side: the local minimax price}\label{sec:gen-benign}

The local length bound rests on a non-degeneracy property that holds for all graphs. Parametrize $\cM_{G^*}$ by structural coefficients and noise variances, and let $h_{G'}$ be the vector of all local Markov partial correlations of $G'$.

\begin{theorem}[First-order distinguishability]\label{thm:E}
Let $G^*,G'$ be DAGs with the same skeleton $S$ and $e\in S$, such that $G^*-e$ and $G'-e$ are Markov equivalent and $G^*$ and $G'$ are not. Then $\partial h_{G'}/\partial\lambda_e$ is not identically zero on $\{\lambda_e=0\}\subset\Theta_{G^*}$.
\end{theorem}

The proof (\cref{app:E}) compares the unshielded colliders that $e$ creates in $G^*$ and $G'$ and, in each of four cases, exhibits a local Markov constraint of $G'$ with non-zero derivative at a point of $\{\lambda_e=0\}$ with only one other non-zero coefficient. So the two models are never tangent along $\cM_H$ at generic points, whereas models with the same adjacencies are tangent at the point of total independence \citep[Theorem~4.2]{evans2020model}; this is why the benign side needs no second-order analysis.

\begin{lemma}[Cancellations violate triangle faithfulness]\label{lem:HC}
For every $G'\in\Xcanc$ there are a triangle edge $u-v$ of $G'$ with non-zero coefficient and a legal $W$ with $\rho_{P_0}(u,v\mid W)=0$.
\end{lemma}
The proof (\cref{app:HC}) reduces the lemma to \citet[Theorem~2]{zhang2008detection}: for a causally sufficient set of variables, as here, an undetectable failure of faithfulness, one in which the distribution is faithful to some other DAG, is a failure of triangle faithfulness.

\paragraph{The cross-model constant.} For a DAG $G'$ and a positive definite $P$, let $D_{G'}(P)=\min_{Q\in\cM_{G'}}\KL(P\|Q)$ and let $\Pi_{G'}P$ be the minimizer. For an ordered pair $(G^*,G')$ of non-equivalent members of $\Xeq$, let $P_\varepsilon\in\cM_{G^*}$ be obtained from the parameter of $P_0$ under $G^*$ by setting the coefficient of $e$ to $\varepsilon$, and define
\[
I_{G^*,G'}=\frac{d^2}{d\varepsilon^2}D_{G'}(P_\varepsilon)\Big|_{\varepsilon=0},\qquad
\delta_{G^*,G'}=\frac{d}{d\varepsilon}\Bigl[\theta_{G'}(\Pi_{G'}P_\varepsilon)-\theta_{G^*}(P_\varepsilon)\Bigr]_{\varepsilon=0},
\]
the Fisher information that separates the two explanations along the weak edge, and the first-order discrepancy between the target and the effect of the closest distribution of the other explanation. The \emph{cross-model constant} is
\[
K(P_0)=\max_{(G^*,G')}\ \frac{|\delta_{G^*,G'}|}{s_{\mathrm{eff}}(G^*,P_0)\,\sqrt{I_{G^*,G'}}},
\]
with $K(P_0)=0$ if $\Xeq$ is a single Markov equivalence class. Let also $\mathcal R(P_0)\ge1$ be the largest ratio $s_{\mathrm{eff}}(D,P_0)/s_{\mathrm{eff}}(D',P_0)$ of the efficient standard deviations of two explanations $D,D'\in\Xeq$; by \cref{lem:oracle} it is at most $(\sigu/\sigl)^2$. By \cref{thm:E}, $I_{G^*,G'}>0$ at base points outside the null set $\cN$ of \cref{thm:K}, so $K(P_0)<\infty$. The ratio does not depend on the choice of the curve (\cref{app:K}), and all quantities are explicit (\cref{lem:Dproj}): $\Pi_{G'}P$ is the product of the node-wise regressions of $P$ along $G'$, $D_{G'}(P)=\frac12\sum_j\log(\sigma^2_{j\mid\pa_{G'}(j)}/\sigma^2_{j\mid\mathrm{pred}(j)})$ for any topological order of $G'$, and $I_{G^*,G'}$ is a sum of quadratic forms in first derivatives of partial covariances.

\paragraph{Two test statistics.} For a DAG $D$ with a topological order, let $\mathrm{pred}(j)$ be the nodes before $j$, $T_D=n\sum_j\log(\hat\sigma^2_{j\mid\pa_D(j)}/\hat\sigma^2_{j\mid\mathrm{pred}(j)})$ the likelihood ratio statistic of $\cM_D$ against the saturated model, with residual variances of regressions with intercept, and $\hat\tau_D\in\R^q$ its whitened coordinates, the standardized coefficients of $\mathrm{pred}(j)\setminus\pa_D(j)$ in the regressions of $j$ on $\mathrm{pred}(j)$, where $q$ is the number of non-adjacent pairs of $S$; their population version $\mu_D$ vanishes exactly on $\cM_D$ (\cref{app:coverage}). For Markov classes $S^*\ne S'$ of $\Xeq$ and a member $D_{S'}$ of $S'$, the direction in which data from $S^*$ move the coordinates of $S'$ is $u_{S^*,S'}=\frac{d}{d\varepsilon}\mu_{D_{S'}}(P_\varepsilon)|_{\varepsilon=0}$, with $P_\varepsilon$ the curve of a member of $S^*$ in the definition of $I$.

\paragraph{Procedure.} The procedure estimates these directions on a small part of the sample and tests each candidate in them on the rest. For $\gamma\in(0,1)$ and $d\ge1$ let $\varpi_d(\gamma)=\sup_{x\ge0}x\,\Prob\bigl(\chi^2_d(x^2)\le\chi^2_{d,1-\gamma}\bigr)$, where $\chi^2_d(x^2)$ is non-central with non-centrality $x^2$, and $\bar\varpi_d=\max_{d'\le d}\varpi_{d'}$.

\begin{definition}[Double screening]\label{def:ds}
Let $\gamma_1,\gamma_2,\alpha'>0$ and $\gamma_3\ge0$ with $\gamma_1+\gamma_2+\gamma_3+\alpha'=\alpha$, let $\varepsilon_n\in(0,1)$ and $\tau_n>0$, and split the data into the first $n_1=\lceil\varepsilon_nn\rceil$ observations and the other $n_2=n-n_1$.
\begin{enumerate}
\item[(0)] \emph{Structure estimate}, on the first part, with sample covariance $\hat\Sigma^{(1)}$: find $(\hat D,\hat e)$ that minimizes $D_{D-e'}(\hat\Sigma^{(1)})$ over $D\in\cG$ and edges $e'$ of $S$, with ties broken by a fixed rule. If $n_1\ge p+1$ and the minimum is at most $\tau_n$, let $\hat H=\hat D-\hat e$ and $\hat\cX=\{D\in\cG:D-\hat e\simeq\hat H\}$; for each Markov class $S'$ of $\hat\cX$, fix a member $D_{S'}$ and a topological order by a rule that depends only on $\hat\cX$, and let $\hat B_{S'}$ be an orthonormal basis of the span of the vectors $\hat u_{S^*,S'}$ over the other classes $S^*$ of $\hat\cX$, defined as $u_{S^*,S'}$ with $(P_0,e)$ replaced by $(\Pi_{\hat H}\hat\Sigma^{(1)},\hat e)$, and empty if all of them vanish. Otherwise step (0) returns no estimate, and step (1b) is not applied.
\item[(1)] \emph{Model screen.} (a) Reject $D$ if $T_D$, computed on all data, exceeds the $1-\gamma_1$ quantile of its null law. (b) If $\gamma_3>0$ and $D$ is Markov equivalent to some $D_{S'}$ for which $\hat B_{S'}$ is defined and non-empty, reject $D$ if $|\hat B_{S'}^\top\hat\tau^{(2)}_{D_{S'}}|^2$, computed on the second part, exceeds the $1-\gamma_3$ quantile of its null law given $\hat B_{S'}$.
\item[(2)] \emph{$k$-TF screen.} Let $\mathcal T_D$ be the set of $k$-TF elements $(u\to v,W)$ of $D$, a triangle edge $u\to v$ of $D$ and a legal $W$, and $N_T=\sum_{D\in\cG}|\mathcal T_D|$. For each element of $\mathcal T_D$ form exact level-$(1-\eta)$ intervals, $\eta=\gamma_2/(4N_T)$: the $t$-interval $I_\beta$ for the coefficient of $u$ in the regression of $v$ on $\pa_D(v)$, $\chi^2$-intervals $[\sigma^2_{u,L},\sigma^2_{u,U}]$ and $[\sigma^2_{v,L},\sigma^2_{v,U}]$ for the marginal variances, and the exact interval $I_\rho$ for the partial correlation obtained by inverting the distribution of the sample partial correlation. Reject $D$ if $k\,\mathrm{dist}(0,I_\beta)\sqrt{\sigma^2_{u,L}/\sigma^2_{v,U}}>\max_{r\in I_\rho}|r|$ for some element.
\item[(3)] \emph{Interval.} Output the convex hull, over the candidates that are not rejected, of their intervals, and the empty set if all candidates are rejected. For \DSA\ the interval of $D$ is the level-$(1-\alpha')$ $t$-interval for the coefficient of $X$ in the regression of $Y$ on $X$ and $A(D)$, where $A$ assigns to each $D$ a valid adjustment set $A(D)\not\ni Y$; for \DSeff\ it is the Wald interval $\hat\theta_D\pm z_{1-\alpha'/2}\hat s_Dn^{-1/2}$ of the efficient estimator, $\hat\theta_D=\theta_D(\Pi_D\hat\Sigma)$ and $\hat s_D=s_{\mathrm{eff}}(D,\Pi_D\hat\Sigma)$.
\end{enumerate}
\end{definition}

The procedure uses the candidate class, $k$ and the levels, but not the base point, the weak edge or $H$. We write \DSOp\ for $A(D)=O(X,Y,D)$, the optimal adjustment set of \citet{henckel2022graphical}, and use $\varepsilon_n=\tau_n=n^{-1/3}$ in \cref{sec:sim}. The model screen plays the role of the first step of the Bonferroni two-step methods of \citet{romano2014practical} and \citet{mccloskey2017bonferroni}: a pre-test restricts the candidates, and its level is spent from that of the final interval. Like the honest intervals of nonparametric regression, which need a smoothness constant \citep{armstrong2018optimal}, the procedure takes $k$ as an input.

\begin{theorem}[Null laws and honesty]\label{thm:coverage}
(i) For $P\in\cM_D$ and $n\ge p+1$, $T_D$ has the law of $\sum_{j=1}^pn\log(1+U_j/V_j)$ and $\hat\tau_D$ that of $(N_j/\sqrt{V_j/(n-j)})_j$, with independent $U_j\sim\chi^2_{q_j}$, $q_j=|\mathrm{pred}(j)\setminus\pa_D(j)|$, $N_j\sim N(0,I_{q_j})$ and $V_j\sim\chi^2_{n-j}$ ($U_j=|N_j|^2$), where $j$ is the position in the order; these laws do not depend on $P$, and the quantiles can be computed by quadrature or Monte Carlo.
(ii) For every $k>0$, every rule $A$ and every $n$ with $n_2\ge p+1$ and $n\ge p+4$, \DSA\ is honest at level $1-\alpha$ over $\cM_k$, and \DSeff\ is asymptotically honest at level $1-\alpha$ over $\cM_k$.
\end{theorem}

On the event that all tests and exact intervals of the true graph succeed, the true graph passes the $k$-TF screen because $k\,\mathrm{dist}(0,I_\beta)\sqrt{\sigma^2_{u,L}/\sigma^2_{v,U}}\le k|\lambda^{\mathrm{std}}|\le|\rho|\le\max_{r\in I_\rho}|r|$, so its interval is part of the output. Step (0) does not affect honesty: given the first part of the sample, the directed test of the true graph has exact level $\gamma_3$ by (i), because the two parts are independent (\cref{app:coverage}).

\paragraph{Two further constants.} Let $m$ be the number of Markov classes of $\Xeq$. Members of one class $S'$ have the same model, hence the same $I_{G^*,G'}$, and by \eqref{eq:K1p} the same $\delta_{G^*,G'}$; writing $I_{G^*,S'}$ and $\delta_{G^*,S'}$, let
\[
K_+(P_0)=\max_{G^*\in\Xeq}\sum_{S'\not\ni G^*}\frac{|\delta_{G^*,S'}|}{s_{\mathrm{eff}}(G^*,P_0)\sqrt{I_{G^*,S'}}},
\]
so that $K\le K_+\le(m-1)K$, and $K_+=K$ if $m=2$; $K_+$ enters only the upper bounds. Finally, let
\[
k^*(P_0)=\min\frac{|\rho_{P_0}(u,v\mid W)|}{|\lambda^{\mathrm{std}}_D(u\to v;P_0)|},
\]
the minimum over $D\in\Xeq$, triangle edges $u\to v\ne e$ of $D$ and legal $W$, where a ratio with zero denominator is $+\infty$; it is positive because $u$ and $v$ are adjacent in $D-e\simeq H$. The lower bounds below hold for $k<k^*(P_0)$.

\begin{theorem}[Benign weak edges: the local price]\label{thm:K}
Assume \eqref{eq:K2} and \eqref{eq:K1p}. There is a relatively closed real-analytic subset $\cN\subset\Theta_H$ of Lebesgue measure zero such that for every $\phi_0\in\Theta_H\setminus\cN$, with $P_0=P(\phi_0)$, the following hold.
\begin{enumerate}
\item[(i)] If the eigenvalues of $\Sigma_{P_0}$ lie in the interior of $[\sigl^2,\sigu^2]$, then for every $k\in(0,k^*(P_0))$, every $r>0$ and every sequence of intervals that are honest, or asymptotically honest, over $\cM_k$,
\[
\liminf_n\Pi^{\mathrm{loc}}_n(C_n;r)\ \ge\ \max\Bigl\{\Psi_\alpha\bigl(K(P_0)\bigr),\ c_P\,\mathcal R(P_0)\Bigr\},\qquad \Psi_\alpha(K)\ \ge\ \max\Bigl\{c_P,\ \tfrac12c_P+c_KK\Bigr\}.
\]
\item[(ii)] Let $m\ge2$, fix the levels with $\gamma_3>0$, and let $\varepsilon_n\to0$, $\tau_n\to0$ and $n\varepsilon_n\tau_n/\log n\to\infty$. For every $k>0$, \DSeff\ satisfies
\[
\lim_{r\downarrow0}\limsup_{n\to\infty}\Pi^{\mathrm{loc}}_n(C^{\mathrm{eff}}_n;r)\ \le\ \mathcal R(P_0)\,\frac{z_{1-\alpha'/2}+e_m}{z_{1-\alpha/2}}\ +\ \frac{\bar\varpi_{m-1}(\gamma_3)}{2z_{1-\alpha/2}}\,K_+(P_0),
\]
with $e_2=\sqrt{2/\pi}$ and $e_m=\sqrt{2\log m}$ for $m\ge3$. For \DSA\ the same holds with $\mathcal R^A\ge\mathcal R$, $\sqrt{2\log|\Xeq|}$ and $K^A_+$, defined in \cref{rem:Oset}, in place of $\mathcal R$, $e_m$ and $K_+$ (\cref{thm:explicit}(d)).
\end{enumerate}
\end{theorem}

Here $\Psi_\alpha$ is the explicit non-decreasing function of \cref{app:K}, $c_P=(a\Phi(a)+\varphi(a))/z_{1-\alpha/2}$ with $a=z_{1-\alpha}$ is Pratt's constant for one parametric model, and $c_K=\sup_d d\,(\Phi(z_{1-\alpha}-d)-\alpha)/(2z_{1-\alpha/2})$; at $\alpha=0.05$, $c_P=0.8498$, $c_K=0.1952$ and $\Psi_\alpha(K)=c_P$ for $K\le1$. The lower bound is a Pratt-type integral bounded by Neyman--Pearson from alternatives inside $\cM_{G^*}$ and from $\Pi_{G'}P_\varepsilon$ at $\varepsilon\asymp(nI)^{-1/2}$; the term $c_P\mathcal R$ comes from the convention of \cref{sec:setup} that a distribution in two models is held to both oracle lengths, and by \eqref{eq:R-identity} below it exceeds $c_P$ only if $K>0$. The two slopes have the same form, a distance times a type II error maximized over the distance: in the lower bound the error of a one-sided level-$\alpha$ test less $\alpha$, in the upper bound that of the directed test. The information that defines $K$ is exactly the non-centrality of both tests: $|u_{S^*,S'}|^2=I_{G^*,S'}$, and along the weak edge at coefficient $c/\sqrt n$, $T_{S'}$ is asymptotically $\chi^2_q(c^2I_{G^*,S'})$ and the directed statistic $\chi^2_{d}(c^2I_{G^*,S'})$ with $d\le m-1$ (\cref{app:K}). In the language of \cref{sec:setup}, $\max\{\mathcal R,K\}$ is, up to constants that do not depend on the graph, the slope at zero of the local modulus of continuity of the target over the class, in units of $s_{\mathrm{eff}}$.

\begin{corollary}[The local minimax price]\label{cor:RK}
Under the assumptions of \cref{thm:K}, let the eigenvalues of $\Sigma_{P_0}$ lie in the interior of $[\sigl^2,\sigu^2]$ and $k<k^*(P_0)$, and let $\Pi^*_k(P_0)$ be the infimum of $\lim_{r\downarrow0}\limsup_n\Pi^{\mathrm{loc}}_n(C_n;r)$ over the sequences of intervals that are asymptotically honest over $\cM_k$. If $m\ge2$, then for every $\gamma_3\in(0,\alpha)$
\[
\max\bigl\{1,\ c_P\mathcal R,\ \tfrac12c_P+c_KK\bigr\}\ \le\ \Pi^*_k(P_0)\ \le\ a_m\mathcal R+b_mK_+\ \le\ a_m\mathcal R+(m-1)\,b_mK,
\]
with $a_m=(z_{1-(\alpha-\gamma_3)/2}+e_m)/z_{1-\alpha/2}$ and $b_m=\bar\varpi_{m-1}(\gamma_3)/(2z_{1-\alpha/2})$, and the ratio of the two sides is at most $a_m/c_P+(m-1)b_m/c_K$ for all values of $(\mathcal R,K)$. For every $\eta>0$, \DSeff\ with $\gamma_1+\gamma_2$ small enough, which does not depend on $P_0$, has local price at most $(a_m+\eta)\mathcal R+b_mK_+$ at every base point to which the corollary applies with $m\ge2$. If $m=1$, then $\mathcal R=1$, $K=0$ and $\Pi^*_k(P_0)=1$; the lower bound $1$ and the case $K=0$ are \cref{thm:free}.
\end{corollary}

At $\alpha=0.05$ and $\gamma_3=0.03$, for two Markov classes, the local minimax price lies between $\max\{1,\,0.849\,\mathcal R,\,0.42+0.195K\}$ and $1.595\,\mathcal R+0.296K$ (lower constants rounded down, upper ones up), within the factor $3.40$ for every graph, every dimension and every eigenvalue range.\footnote{This $\gamma_3$ minimizes the factor for $m=2$. For $m=3$ the factor is $5.76$ at the same $\gamma_3$ and $5.70$ at its own optimum $\gamma_3\approx0.039$; since $m$ is not known in advance, a procedure fixes one $\gamma_3$ (\cref{app:certify}).} The upper bound, which \cref{thm:explicit}(c) gives for a procedure that knows the base point, is attained to within $\eta\mathcal R$ by one that knows neither the base point nor the weak edge; since the two bounds differ by a factor up to $3.40$, this does not show that adaptation is free, only that it is not needed to reach the bound. It rests on the directed test of step (1b), estimated on a vanishing fraction of the sample and exact given its direction: the likelihood ratio test alone has $q$ degrees of freedom, while near $P_0$ the deviation from a wrong explanation has at most $m-1$ dimensions, and without step (1b) the factor grows with $q$: at $\gamma_1=\gamma_2=0.01$ it is $3.66$ for $q=1$ and $4.95$ for $q=15$ (\cref{thm:explicit}). Since $\mathcal R\le(\sigu/\sigl)^2$, the order is $1+K$ for a fixed eigenvalue range, but $\mathcal R$ can be large while $K$ is bounded, by \eqref{eq:R-identity}.

\begin{remark}[Sample size and a fixed neighbourhood]\label{rem:samplesize}
The bounds of \cref{thm:K,cor:RK} are asymptotic: the sample size from which they hold depends on $P_0$, through the separation of the other edges, the curvature of $D_H$ and $I_{G^*,S'}$, and on $k$, which enters only there and in the exponential rate, of order $nk^2$, at which cancellation explanations are rejected (\cref{lem:B}). With a conditional independence screen in place of step (1) and parent adjustment, a bound $C(1+K(P_0))$ holds on a fixed neighbourhood of $P_0$ for $n\ge n_1$, with $C$ independent of $P_0$ and $k$ (\cref{prop:CI}). For the upper bounds, the dependence on $P_0$ is locally uniform (\cref{cor:compact}).
\end{remark}

Because double screening does not depend on the base point, this dependence is locally uniform. Write $\Pi^{\mathrm{loc}}_n(C_n;r,\phi)$ for the local price around $P(\phi)$, and let $U(\phi)$ be the right side of \cref{thm:K}(ii) at $P(\phi)$.

\begin{corollary}[Uniformity in the base point]\label{cor:compact}
Under the assumptions of \cref{thm:K}(ii), let $\Phi\subset\Theta_H\setminus\cN$ be compact, with the eigenvalues of $\Sigma_{P(\phi)}$ in the interior of $[\sigl^2,\sigu^2]$ for $\phi\in\Phi$. Then $\min_\Phi k^*>0$, so the lower bounds of \cref{thm:K}(i) and \cref{cor:RK} hold at every point of $\Phi$ for one $k$, and for every $k>0$ and $\eta>0$ there are $r>0$ and $n_0$ with
\[
\Pi^{\mathrm{loc}}_n(C^{\mathrm{eff}}_n;r,\phi)\ \le\ U(\phi)+\eta\qquad\text{for all }\phi\in\Phi\text{ and }n\ge n_0 .
\]
The same holds for \DSA\ with its bound, and for the screen of \cref{prop:CI} with the bound $C(1+\max_\Phi K)$. For the lower bounds only the threshold on $k$ is uniform; they are statements for each base point as $n\to\infty$.
\end{corollary}
\begin{proof}
On $\Theta_H\setminus\cN$, $I_{G^*,G'}>0$ and $s_{\mathrm{eff}}$, $\delta$ and $I$ are continuous, so $\mathcal R$, $K$, $K_+$ and $U$ are continuous. The partial correlations in the definition of $k^*$ are continuous and, by faithfulness to $H$, which holds off $\cN$ (\cref{app:benign}), non-zero on $\Phi$, and $|\lambda^{\mathrm{std}}|\le(\sigu/\sigl)^2$, so $k^*$ is bounded below on $\Phi$ by a positive constant. The local price depends on $\phi$ only through the ball $\{\|\Sigma_P-\Sigma_{P(\phi)}\|\le r\}$, and the procedure does not depend on $\phi$. Fix $\phi\in\Phi$. By \cref{thm:K}(ii) there are $r_\phi$ and $n_\phi$ with $\Pi^{\mathrm{loc}}_n(C_n;r_\phi,\phi)\le U(\phi)+\eta/2$ for $n\ge n_\phi$. If $\|\Sigma_{P(\phi')}-\Sigma_{P(\phi)}\|<r_\phi/2$, the ball of radius $r_\phi/2$ around $\Sigma_{P(\phi')}$ lies in the ball of radius $r_\phi$ around $\Sigma_{P(\phi)}$, so $\Pi^{\mathrm{loc}}_n(C_n;r_\phi/2,\phi')\le\Pi^{\mathrm{loc}}_n(C_n;r_\phi,\phi)$, which is at most $U(\phi')+\eta$ if also $|U(\phi')-U(\phi)|<\eta/2$. These two conditions define an open neighbourhood of $\phi$; finitely many of them cover $\Phi$, and $r$, $n_0$ are the smallest $r_\phi/2$ and the largest $n_\phi$ over them. The other procedures are treated in the same way.
\end{proof}

The argument needs a procedure that does not depend on the base point; it does not apply to the screen of \cref{thm:explicit}(c), whose directed tests use $P_0$, but it applies to those of \cref{thm:explicit}(a) and (b). Near the boundary of $\Theta_H\setminus\cN$ the constants degenerate; where the boundary is a second vanishing edge that satisfies the condition \eqref{eq:M} below, \cref{thm:uniform} gives the rate uniformly.

\begin{remark}[Identification by triangle faithfulness]\label{rem:ktf-identifies}
A candidate that is Markov equivalent to $G^*$ but not in $\Xeq$ fits every distribution in $\cM_{G^*}$; near $P_0$ it is rejected only because it would need a near cancellation along a triangle (\cref{lem:HC}). Under \eqref{eq:K1p} without the stronger identifiability condition, the bounded price therefore rests on identification supplied by the class rather than by the data, on a neighbourhood that shrinks as $k\to0$. A misspecified $k$ can invalidate coverage; in the simulations of \cref{app:sim} this happened almost only through empty output.
\end{remark}

\subsection{The structure of the cross-model constant}\label{sec:gen-constant}

\Cref{thm:K,cor:RK} say that the local price at a benign base point is of order $\mathcal R(P_0)+K(P_0)$, or $1+K(P_0)$ for a fixed eigenvalue range. This subsection describes the local limit experiment behind $K$, shows that the price is exactly one when $K=0$ and only then, decides from the graphs when this happens, and decides when $K$ blows up.

\paragraph{The canonical limit experiment.} The constant $K$ is one number extracted from the local limit experiment at a benign base point, and that experiment has a simple form. Let $L_0=T_{P_0}\cM_H$, let $L_1,L_2$ be the tangent spaces of $\cM_{G^*}$ and $\cM_{G'}$ at $P_0$ in the Fisher metric, and let $a_i$ be the gradient of the target in $L_i$, so that $s_{\mathrm{eff},i}=|a_i|$. By \cref{thm:E}, $L_i=L_0\oplus\R n_i$ with unit vectors $n_i\perp L_0$ and, after a choice of the sign of $n_2$, $\rho=\langle n_1,n_2\rangle\in[0,1)$, and by \eqref{eq:K2}, $a_1$ and $a_2$ have the same projection onto $L_0$, of length $a_0$; let $\alpha_i=\langle a_i,n_i\rangle$. Let $K_{1\to2}$ and $K_{2\to1}$ be the ratio in the definition of $K$ for the ordered pairs $(G^*,G')$ and $(G',G^*)$. Both are explicit functions of $(a_0,\alpha_1,\alpha_2,\rho)$ (\cref{app:K}), and they determine the ratio of the efficient lengths:
\begin{equation}\label{eq:R-identity}
\frac{s_{\mathrm{eff},1}^2}{s_{\mathrm{eff},2}^2}=\frac{1-K_{2\to1}^2}{1-K_{1\to2}^2}\qquad\text{unless }K_{1\to2}=K_{2\to1}=1 ,
\end{equation}
and $K_{1\to2}>1$ if and only if $K_{2\to1}>1$.

\begin{proposition}[Canonical experiment]\label{prop:canonical}
In the setting of \cref{thm:K}, let $(G^*,G')$ be an ordered non-equivalent pair in $\Xeq$ with the parameters $(a_0,\alpha_1,\alpha_2,\rho)$ above. Let $\Pi^*_{\mathrm{can}}(a_0,\alpha_1,\alpha_2,\rho)$ be the minimax price over randomized confidence sets of the following experiment: observe $Z_0\sim N(\eta,1)$ and $Z=(Z_1,Z_2)\sim N(\tau n_g,I_2)$, independent, with $n_1=(1,0)$ and $n_2=(\rho,\sqrt{1-\rho^2})$, where the model $g\in\{1,2\}$, $\eta\in\R$ and $\tau\in\R$ are unknown, the target is $a_0\eta+\alpha_g\tau$, and the oracle length is $2z_{1-\alpha/2}\sqrt{a_0^2+\alpha_g^2}$. Then:
\begin{enumerate}
\item[(i)] (limit experiment) in the local limit experiment of the class $\cM_k$ at $P_0$ restricted to the pair $(G^*,G')$, the subexperiment in which the coordinates of $L_0$ orthogonal to the target are fixed has minimax price $\Pi^*_{\mathrm{can}}$, and so has the limit experiment of the pair on $L_1+L_2$;
\item[(ii)] (finite samples) for every $k<k^*(P_0)$, every $r>0$ and every sequence of intervals that is asymptotically honest over $\cM_k$, $\liminf_n\Pi^{\mathrm{loc}}_n(C_n;r)\ge\Pi^*_{\mathrm{can}}(a_0,\alpha_1,\alpha_2,\rho)$.
\end{enumerate}
The four-variable limit experiment of \cref{prop:limit} is the case $(a_0,\alpha_1,\alpha_2,\rho)=(1,b,0,\sqrt v)$ after division by $a_0$, with $v=1/(1+w^2)$, $G^*=M$ and $G'=C$.
\end{proposition}

Since randomized confidence sets include intervals, $\Pi^*_{\mathrm{can}}$ of each pair lower-bounds the price constant at $P_0$. It depends on three parameters, of which $K_{1\to2}$, $K_{2\to1}$ and $\mathcal R$ are functions, and not on $K$ alone; on a grid of canonical parameters the lower bound of \cref{thm:K} is within a factor $1.56$ of its discretized value (\cref{app:certify}).

\paragraph{When honesty is free.} In the canonical form, $K_{1\to2}=K_{2\to1}=0$ if and only if $\alpha_1=\alpha_2=0$, because $|\rho|<1$; so $K(P_0)>0$ exactly when the target has a non-zero component normal to $\cM_H$ in the tangent space of some explanation. This component has a closed form, and its non-vanishing is decided by a d-separation.

\begin{proposition}[When the cross-model constant is positive]\label{prop:Kpos}
Let $e$ be admissible, assume \eqref{eq:K2}, and let $\Xeq$ have at least two Markov classes. For $G\in\Xeq$ write $e=u\to v$ in $G$, $D=G-e$ and $W=\pa_G(v)\setminus\{u\}$, and let $D^+$ be $D$ with an added vertex $\xi$ and edge $\xi\to X$. For $\phi_0$ outside a Lebesgue-null subset of $\Theta_H$, $K(P_0)>0$ if and only if some $G\in\Xeq$ satisfies
\begin{enumerate}
\item[(a)] $v=Y$ or $v$ is an ancestor of $Y$ in $D$, and
\item[(b)] $\xi$ and $u$ are d-connected given $W$ in $D^+$.
\end{enumerate}
If no member satisfies (a) and (b), then $K(P_0)=0$ at every base point outside $\cN$. Under \eqref{eq:K1p}, either all or none of the members of a Markov class satisfy (a) and (b), so one member per class suffices.
\end{proposition}

For a member $G$, the normal component is, up to a non-zero factor, $\tau_{v\to Y}\Cov(\varepsilon_X,r_u)$, where $\tau_{v\to Y}$ is the total effect of $v$ on $Y$, $\varepsilon_X$ the noise of $X$ in $D$, and $r_u$ the residual of the regression of $u$ on $W$ (\cref{app:K}). The information that the weak edge adds to $D$ must carry exogenous variation of $X$, and $v$ must affect $Y$. The criterion is not that $e$ lies on a directed path from $X$ to $Y$; the two conditions coincide on five nodes but differ in both directions on six (\cref{app:K,app:certify}).

When $K(P_0)=0$, all explanations have, at $P_0$, the efficient influence function of the model $\cM_H$ without the weak edge: the efficient estimator that ignores the weak edge has no first-order bias under any explanation, and $\mathcal R(P_0)=1$. The price is then one, and otherwise it exceeds one. We call $(\Xeq,e)$ \emph{free} if $\Xeq$ is a single Markov class or no member of $\Xeq$ satisfies (a) and (b) of \cref{prop:Kpos}.

\begin{theorem}[When honesty is free]\label{thm:free}
In the setting of \cref{cor:RK}:
\begin{enumerate}
\item[(i)] $\Pi^*_k(P_0)\ge1$, with equality if and only if $K(P_0)=0$. If $(\Xeq,e)$ is free, $K(P_0)=0$ at every base point outside $\cN$; otherwise $K(P_0)>0$ for $\phi_0$ outside a Lebesgue-null subset of $\Theta_H$.
\item[(ii)] If $K(P_0)>0$, then $\Pi^*_k(P_0)\ge1+\delta_{\mathrm{can}}(a_0,\alpha_1,\alpha_2,\rho)$ for the canonical parameters of each ordered non-equivalent pair with $K_{1\to2}>0$, where $\delta_{\mathrm{can}}$ is given by an explicit formula and is positive unless $\alpha_1=\alpha_2=0$ (\cref{lem:strict}).
\item[(iii)] If $K(P_0)=0$, then $\mathcal R(P_0)=1$, and the limiting local price of \DSlreff, double screening without the directed test (\cref{thm:explicit}), and for $m\ge2$ that of \DSeff\ under the conditions of \cref{thm:K}(ii), is at most $z_{1-\alpha'/2}/z_{1-\alpha/2}$, where $\alpha'$ is the level of the final interval.
\end{enumerate}
For every $\eta>0$ there are levels such that a variant \DSsw\ of \DSeff, which computes steps (1)--(3) on the second part of the sample and omits the directed test when the explanations estimated in step~(0) are free, has limiting local price at most $1+\eta$ at every base point outside $\cN$ at which $(\Xeq,e)$ is free, and at most $(a_m+\eta)\mathcal R(P_0)+b_mK_+(P_0)$ at every other base point to which \cref{cor:RK} applies. It is asymptotically honest over $\cM_k$, and its version with adjustment intervals is honest in finite samples (\cref{app:K}).
\end{theorem}

The lower bound in (i) is the local asymptotic minimax bound of one regular model, transferred to finite samples as in \cref{prop:canonical}(ii) (in one Gaussian model the usual interval is even admissible, \citealp{joshi1969admissibility}); the bound $\Psi_\alpha(K)$ of \cref{thm:K}(i), which equals $c_P<1$ for $K\le1$, holds at single points of the neighbourhood. That a vanishing first-order discrepancy costs nothing is the analogue of the case $\omega=0$ in the local misspecification framework of \citet{hjort2003frequentist}, in which the estimators of all submodels have the same limit distribution. The new part is (ii): when $K>0$, the efficient intervals of two explanations differ on a set of data that is typical under both. The proof (\cref{app:K}) restricts a Bayes argument with diffuse priors to this set, where the two Bayes intervals differ by an amount proportional to $|\alpha_1n_1-\alpha_2n_2|$. The function $\delta_{\mathrm{can}}$ shows that the price is not one, not by how much: on the $27$ canonical parameters of \cref{app:certify} numerical lower bounds for it lie between $1.6\times10^{-8}$ and $4.3\times10^{-4}$, while the discretized program exceeds one by between $0.001$ and $0.87$, and along one ray it decays like $K^4$ as $K\to0$. Screening is still needed when $K=0$: the discrepancy of a wrong explanation is then of second order in the weak coefficient, which is not negligible for coefficients larger than $n^{-1/4}$, and the screen removes it at a cost that vanishes with its level. In all benign cases of \cref{tab:simB}, $K(P_0)=0$, and the median ratio $1.11$ there agrees with the factor of (iii) for the levels used. On five nodes with only $Y$ a sink, $(\Xeq,e)$ is free for $8186$ of the $8426$ benign pairs $(G,e)$, $5954$ of them with a single Markov class, and on six nodes for $1358485$ of $1448757$ (\cref{app:certify}). In the four-variable model with unit variances, $K=\max\{b/w,\ b\sqrt{1+w^2}/(w\sqrt{1+b^2})\}$ (\cref{app:K}), which is $b/w$ when $b\ge w$, so the benign edge $X-Z$ is free only if $b=0$.

\paragraph{Approaching a malignant stratum.} The constant explains the term $B/(\lambda\sqrt n)$ of \cref{thm:rate-lower}, and it shows that, along this approach, $1/\lambda$ is the only possible blow-up of $K$. Let $f$ be an edge of $H$, let $\phi_\lambda$ be $\phi_0$ with the coefficient of $f$ in $\Theta_H$ replaced by $\lambda$, let $P_\lambda=P(\phi_\lambda)$, and let $\phi_0'\in\Theta_{H-f}$ denote the other coordinates of $\phi_0$. At $\lambda=0$ the base point lies on the lower stratum $\cM_{H-f}$, and whether a pair of explanations can be separated there to first order, and whether their effects then differ to first order, is decided by graphical conditions. For an ordered non-equivalent pair $(G^*,G')$ in $\Xeq$, consider:
\begin{enumerate}
\item[(F)] the v-structures of $G^*$ and $G'$ differ in exactly one, and its two arms are $e$ and $f$; we write it $x_e\to z\leftarrow x_f$, with $e=x_e-z$ and $f=z-x_f$, and call $G_\wedge$ the one of $G^*,G'$ that contains it and $G_\vee$ the other;
\item[(C)] $\mathrm{CC}_{G_\wedge}(z,x_f)=\mathrm{CC}_H(z,x_f)$, where $\mathrm{CC}_D(s,t)$ denotes the set of common children of $s$ and $t$ in $D$;
\item[(D)] $X$ is an endpoint of $e$, and, with $w$ the other endpoint and $U=\pa_{G_\wedge}(z)\setminus\{x_e,x_f\}$: $e\notin R$, some member of $[H]_Y$ has $\pa(x_e)=U$, and $Y$ and $w$ are d-connected given $U\cup\{X\}$ in $H-f$.
\end{enumerate}
Condition (F) is equivalent to the statement that every witness configuration of the proof of \cref{thm:E} (\cref{app:E}) uses the edge $f$, and (C) is equivalent to $(G_\wedge-e)-f\simeq H-f$: deleting $f$ from the two equivalent graphs $G_\wedge-e$ and $H$ leaves equivalent graphs (\cref{app:K}). All three conditions can be read off the graphs; for $R=\emptyset$ the parent-set requirement of (D) is a local condition on the CPDAG of $H$, and $U$ does not depend on the choice of $G_\wedge$ within its class (\cref{app:K}).

\begin{theorem}[Price near a malignant stratum]\label{thm:K-lambda}
In the setting of \cref{thm:K}, assume that the eigenvalues of $\Sigma_{P_\lambda}$ at $\lambda=0$ lie in the interior of $[\sigl^2,\sigu^2]$. For $\phi_0'$ outside a Lebesgue-null subset of $\Theta_{H-f}$ there are $\lambda_0>0$ and $0<c\le C$ such that for $0<|\lambda|<\lambda_0$, $\phi_\lambda\notin\cN$ and:
\begin{enumerate}
\item[(i)] $K(P_\lambda)\le C/|\lambda|$;
\item[(ii)] if some ordered non-equivalent pair $(G^*,G')$ in $\Xeq$ has $I_{G^*,G'}(P_{\lambda=0})=0$ and $\delta_{G^*,G'}(P_{\lambda=0})\ne0$, then $K(P_\lambda)\ge c/|\lambda|$; otherwise $K(P_\lambda)\le C$;
\item[(iii)] for each pair, $I_{G^*,G'}(P_{\lambda=0})=0$ if and only if (F) and (C) hold;
\item[(iv)] for each pair that satisfies (F) and (C), $\delta_{G^*,G'}(P_{\lambda=0})\ne0$ if (D) holds and $\delta_{G^*,G'}(P_{\lambda=0})=0$ otherwise.
\end{enumerate}
In particular, $K(P_\lambda)\ge c/|\lambda|$ if some ordered non-equivalent pair satisfies (F), (C) and (D), and $K(P_\lambda)\le C$ otherwise. Hence, for each such $\lambda$ and every $k\in(0,k^*(P_\lambda))$, the local minimax price at $P_\lambda$ is, as $n\to\infty$, of order $1/|\lambda|$ in the first case of (ii) and of order one in the second, with constants that depend on $\phi_0'$ but not on $\lambda$.
\end{theorem}

So whether $K$ blows up like $1/\lambda$ along this approach is decided by the graphs alone. The proof (\cref{app:K}) works at $P_{\lambda=0}\in\cM_{H-f}$: (F) and (C) hold exactly when members of the two classes can be chosen to coincide after deleting $f$, a statement about Markov classes proved by covered edge reversals \citep{chickering1995transformational} (\cref{lem:onev}), and then the $e$-curve stays in $\cM_{G'}$, so $I$ vanishes; if (F) or (C) fails, a sparse base point gives $I>0$; and $\delta$ is then computed from members of the two classes that differ in $e$ alone, which exist exactly when $e\notin R$ and $x_e$ can have the parent set $U$ (\cref{lem:Econf}). This requirement cannot be dropped: on six nodes it is the only failing part of (D) in some configurations with bounded $K(P_\lambda)$ (\cref{app:certify}). The theorem is a statement for each fixed $\lambda$: the threshold $k^*(P_\lambda)$ and the sample size from which the bounds apply may depend on $\lambda$, and the simulations of \cref{app:sim} show that this sample size can grow quickly as the margin shrinks. In the four-variable model, with $e=X-Z$ and $f=W\to X$, $C$ and $M$ satisfy (F), (C) and (D), $\delta\to\pm b$, and (ii) gives the factor $b/w$ of the explicit computation; \cref{app:K} gives a five-node instance without required edges (\cref{fig:instances}).

On all skeletons with five nodes, the criteria of \cref{thm:K-lambda,prop:Kpos} agree with exact rational computation in every case (\cref{app:certify}).

\subsection{Uniformity in the margin}\label{sec:gen-uniform}

\Cref{thm:K-lambda} holds for each fixed margin: its threshold $k^*(P_\lambda)$ and the sample size from which it applies depend on $\lambda$. In the four-variable model the rate is uniform in the margin (\cref{thm:rate-lower,thm:rate-upper}). The same holds on general graphs, near a base point at which both $e$ and $f$ are absent, under one graphical condition on $f$.

Let $f$ be an edge of $H$, let $H_2=H-f$, and let $P_{00}=P(\phi_{00})$ with $\phi_{00}\in\Theta_{H_2}$, faithful to $H_2$ and with the eigenvalues of $\Sigma_{P_{00}}$ in the interior of $[\sigl^2,\sigu^2]$. For a candidate $D$, an edge $g$ of $D$ and a positive definite $P$, let $\beta^D_g(P)$ be the coefficient of $g$ in the regression parametrization of $P$ under $D$. For $r,\lambda,A>0$ let
\[
\cM_k(\lambda,A;r)=\bigl\{(G^*,P)\in\cM_k:\ \|\Sigma_P-\Sigma_{P_{00}}\|\le r,\ |\beta^{G^*}_e(P)|\le A,\ |\beta^{G^*}_f(P)|\ge\lambda\bigr\},
\]
so that $e$ is weak, with coefficient at most $A$, and $f$ has margin $\lambda$; in the four-variable model $e=X-Z$ and $f=W\to X$, with the bounds $A$ on $|a|$ and $\lambda$ on $|w|$. Let $L_n^*(\lambda,A;r)$ be the infimum of $\sup_{\cM_k(\lambda,A;r)}\E_P\Leb(C_n)$ over intervals that are honest at level $1-\alpha$ over $\cM_k(\lambda,A;r)$. The condition on $f$ is
\begin{equation}\tag{M}\label{eq:M}
D\in[H_2]_Y\ \text{and}\ D+f\ \text{acyclic and allowed by }\cK\quad\Longrightarrow\quad D+f\in[H]_Y,
\end{equation}
where $D+f$ is $D$ with $f$ added in either orientation.
It says that if $f$ were the weak edge at a base point on $\cM_{H_2}$, its explanations would form a single Markov class; it then follows that $f$ is benign there (\cref{app:uniform}). It holds in the four-variable model, where $[H_2]_Y$ has a single member and $f=W\to X$ is required.

\begin{theorem}[Uniformity in the margin]\label{thm:uniform}
Let $e$ be admissible, assume \eqref{eq:K2} and \eqref{eq:K1p}, let $f$ be an edge of $H$ that satisfies \eqref{eq:M}, and assume that $H$ contains a directed path from $X$ to $Y$, which holds unless $e$ is the edge $X-Y$. Let $\iota=1$ if some ordered non-equivalent pair in $\Xeq$ satisfies (F), (C) and (D), and $\iota=0$ otherwise. For $\phi_{00}$ outside a Lebesgue-null subset of $\Theta_{H_2}$ there is $k_0>0$ such that for every $k\in(0,k_0)$ there are $r,\lambda_0,A_0>0$, $n_0$ and $0<c\le C$ with
\[
c\,R_n\ \le\ L_n^*(\lambda,A;r)\ \le\ \sup_{\cM_k(\lambda,A;r)}\E_P\Leb\bigl(C_n^{\DSA}\bigr)\ \le\ C\,R_n,\qquad R_n=n^{-1/2}+\iota\min\Bigl\{A,\frac1{\lambda\sqrt n}\Bigr\},
\]
for all $n\ge n_0$, $\lambda\in(0,\lambda_0]$ and $A\in(0,A_0]$, where \DSA\ is double screening (\cref{def:ds}) with any adjustment rule and any levels, which is honest in finite samples over all of $\cM_k$ (\cref{thm:coverage}). The constants depend on $\phi_{00}$, $k$, $\alpha$ and the levels, but not on $n$, $\lambda$ or $A$.
\end{theorem}

Since the oracle length is of order $n^{-1/2}$ on the class (\cref{lem:oracle}), the price is of order $1+\iota\min\{A\sqrt n,1/\lambda\}$, and three regimes appear. For fixed $\lambda$ it is of order $1/\lambda$, as in \cref{thm:K-lambda}, now with a threshold on $k$ and a sample size from which the bounds apply that do not depend on $\lambda$. If $\lambda\sqrt n\lesssim1/A$ it is of order $A\sqrt n$: at $\lambda=0$ the pair of explanations coincides on all of $\cM_{G_\wedge-f}$, so that $e$ is locally non-identified relative to $H_2$, as in the third case of \cref{thm:dichotomy}, and $A$ plays the role of the radius. If $\iota=0$ the price is bounded uniformly in $\lambda$ and $A$. There is no further term, in particular no plateau of order $n^{-1/4}$ that does not depend on $\lambda$ and $A$. The plateau $B\min(A,\sqrt{\abar(\kappa)^2+n^{-1/2}})$ of the four-variable model, and its $n^{-1/4}$ for unit variances, come from the known range $\kappa$ of the noise variances; $\cM_k$ restricts only the eigenvalues of $\Sigma$, and locally it behaves like $\kappa=\infty$, where $\min(A,\abar(\kappa))=A$. Double screening attains the rate without knowing $\lambda$, $A$, $P_{00}$ or which edges are weak, and, as in \cref{cor:compact}, its upper bound holds uniformly over compact sets of base points $\phi_{00}$ outside the exceptional set, which is relatively closed: $k_0$, a minimum of finitely many continuous functions that are positive off this set, is bounded below on such a compact set, and the covering argument of \cref{cor:compact} applies to the class $\cM_k(\lambda,A;r)$, which depends on $\phi_{00}$ only through the centre of the ball. The path condition is used only for the term $n^{-1/2}$ of the lower bound: if it fails, then $e=X\to Y$, by \eqref{eq:K2} no member of $[H]_Y$ has a directed path from $X$ to $Y$, the effect equals the coefficient of $e$ on the class, $\iota=0$ because $Y$ is a sink, so that neither $e$ nor $f$ can be reversed and no pair satisfies (F), and the intersection of \DSA\ with $[-A,A]$ is honest over the class with length $O(\min\{A,n^{-1/2}\})$.

Three observations make the expansions behind \cref{thm:K-lambda} uniform (\cref{app:uniform}). First, by \cref{lem:HC} applied at $P_{00}$, the true graphs near $P_{00}$ lie in $\cX_2=\{G^*\in\cG:G^*-e-f\simeq H_2\}$, which \eqref{eq:M} places in $\Xeq$; the graphs that explain $P_\lambda$ for small $\lambda$ but are not in $\cX_2$ are cancellation explanations of $P_{00}$ and are rejected by the $k$-TF screen. Second, if two members of $\cX_2$ have v-structures that differ only in $x_e\to z\leftarrow x_f$, deleting $f$ leaves Markov equivalent graphs. Third, when $e$ and $f$ lie in a common triangle, no pair satisfies (F), and the lower bound uses only laws at which the coefficient of $e$ vanishes, where the triangle conditions involving $e$ are void. So, in the coordinates $(\varepsilon,\lambda')$ of the coefficients of $e$ and $f$ in the true graph, the projection distance factorizes, for the pairs that satisfy (F), as $D_{G'}=\varepsilon^2\lambda'^2\psi$ with $\psi$ bounded away from zero, while the effect discrepancy is $O(|\varepsilon|)$, and $O(|\varepsilon\lambda'|)$ when (D) fails; the rest is the two-point bound and an exponential bound on keeping a separated candidate.

Condition \eqref{eq:M} cannot simply be dropped, and it is restrictive. On six nodes with only $Y$ a sink, $496$ configurations satisfy (F), (C) and (D) with $e$ benign in $H$. In $368$ of them $f$ is not benign in $H_2$, that is, not benign as the weak edge of $H$ at base points on $\cM_{H_2}$, so that the doubly weak base point is itself one at which the weak edge $f$ is not benign. Of the other $128$, $64$ satisfy \eqref{eq:M}; in the remaining $64$ some $H'\neq H$ with $H'-f\simeq H_2$ explains points near $P_{00}$, and $e$ is not benign in $H'$. On five nodes all four such configurations have $f$ not benign in $H_2$, so without required edges $\iota=1$ under \eqref{eq:M} first occurs on six nodes; required edges, such as $W\to X$ in the four-variable model, make \eqref{eq:M} easier to satisfy. On a grid of $(n^{-1/2},\lambda,A)$ spanning three to four decades in each coordinate, the modulus of continuity of the effect over the class is between $0.29$ and $2.62$ times $R_n$ in the $64$ configurations with \eqref{eq:M} and in the four-variable model, and more than $2.2\times10^3$ times $R_n$ in each of the $64$ without it (\cref{app:certify}). The other positions of $e$ and $f$, which have $\iota=0$, are checked on five nodes in \cref{app:certify}.

\section{Simulations}\label{sec:sim}

We compare, on the same data, \emph{Naive}, which selects one graph and reports its $95\%$ $t$-interval; double screening (\cref{def:ds}) with $\gamma_1=0.001$, $\gamma_2=\gamma_3=0.01$, $\alpha'=0.029$ and $\varepsilon_n=\tau_n=n^{-1/3}$, in the finite-sample honest version \DSOp\ and the efficient version \DSeff; \emph{Union}, the union of all candidate $95\%$ intervals reported as a set (\cref{prop:union}); and \emph{Hull}, the convex hull of the same intervals without screening. In the four-variable model the honest procedure is \DSV\ of \cref{def:dsv}. Lengths are ratios to the efficient interval $2z_{1-\alpha/2}s_{\mathrm{eff}}n^{-1/2}$ of the true graph, so that their supremum is the price of \cref{sec:setup}. Coverages are minima over many cells, each with Monte Carlo error, so a minimum somewhat below $0.95$ is expected. The partial correlation interval of the $k$-TF screen uses the Fisher $z$-transform. The simulations check the implementations and illustrate the theory; they do not prove it. \Cref{app:sim} gives the designs and the results for the procedures of \cref{app:K}, including the conditional independence screen, and for a misspecified $k$.

\paragraph{Four-variable model.} With the least favourable noise variances for $\kappa=4$, $w\in\{0.05,0.2,1\}$, $a\in\{0,0.1,0.3,1\}$, $n\in\{500,2000,8000\}$ and $20000$ replications per cell (\cref{tab:simA} in \cref{app:sim}), Naive undercovers down to $0.489$, while the minimum coverage of \DSV\ is $0.962$. The length ratio of \DSV\ grows with $n$ for $w=0.05$, the malignant regime, and is bounded for $w=1$, decreasing from $2.74$ to $2.40$; the union set has bounded length throughout, as \cref{prop:union} predicts.

\paragraph{Five-node graphs.} From the sampling protocol of \cref{app:certify} we draw $20$ benign and $20$ malignant weak edges according to \cref{prop:graph}, with coefficients $\pm U(0.5,1.5)$, noise variances $U(0.5,1.5)$, a true model that is $k$-triangle faithful with ratio at least $2k$, $k=0.05$, and weak-edge coefficient $c/\sqrt n$, $c\in\{0,1,3\}$; Naive selects the candidate with the highest likelihood (\cref{app:sim}). \Cref{tab:simB} shows the results. Naive fails badly on malignant edges (coverage $0.25$ to $0.59$) and undercovers on benign ones at moderate $n$. On benign edges the median length ratio of double screening is about $1.11$ for all $n$, the ratio $z_{1-\alpha'/2}/z_{1-\alpha/2}$ of the quantiles used: in the median case only candidates equivalent to the true graph survive and no price is paid for the structural uncertainty. In all $20$ benign cases $K(P_0)=0$ (\cref{app:sim}), and this ratio agrees with the limiting factor of \cref{thm:free}(iii) for the levels used; the larger maxima at $n\le2000$ are finite-sample effects. The largest benign ratio comes from one case in which the weak edge is $X\to Y$ itself and the candidate class has $18$ graphs. Screening is what keeps an interval short: Hull is also honest, but on benign edges at $c=0$ its median ratio grows from $2.35$ to $4.61$ and its maximum from $7.88$ to $26.64$ as $n$ grows from $500$ to $8000$, because it keeps candidates whose effects differ by a fixed amount however clearly the data reject them; a set avoids this without screening (\cref{prop:union}), as the Union columns show. On malignant edges the ratios of double screening grow from about $4$ to $10$ at $c=0$, and the maximum roughly doubles with each fourfold increase of $n$, consistent with the predicted $\sqrt n$ scaling. The structure estimate of step (0) picked the weak edge and the class of $H$ in $0.84$, $0.95$ and $0.97$ of the replications on average at $n=500$, $2000$ and $8000$ (rounded down); a wrong pick only changes which candidates receive the directed test, so it does not affect honesty. The minimum coverage over the $353$ retained cells is $0.944$ for both versions.

\begin{table}[t]
\centering
\caption{Five-node graphs: minimum coverage over the retained cases, and median (maximum) length ratio relative to the efficient interval; for Hull and Union only the median is shown. \DSOp\ and \DSeff\ are double screening with optimal adjustment sets and with efficient intervals; Hull is the convex hull of all candidate $95\%$ intervals, without screening; Union is the union of the same intervals, reported as a set. Benign (ben.) and malignant (mal.) edges: of the $20$ cases per edge type, those whose true model violates $k$-triangle faithfulness at the given $c$ are dropped; $17$ to $20$ cases are retained per row. $2000$ replications per cell (coverage standard error about $0.005$). Coverages are rounded down and ratios up. Script \texttt{s12b\_table.py}.}
\label{tab:simB}
\footnotesize
\setlength{\tabcolsep}{3pt}
\begin{tabular}{lrrccccccccc}
\toprule
& & & \multicolumn{5}{c}{minimum coverage} & \multicolumn{4}{c}{length ratio}\\
\cmidrule(lr){4-8}\cmidrule(lr){9-12}
edge & $n$ & $c$ & Naive & \DSOp & \DSeff & Hull & Union & \DSOp & \DSeff & Hull & Union\\
\midrule
ben. & 500 & 0 & 0.798 & 0.949 & 0.948 & 0.947 & 0.943 & 1.12 (7.72) & 1.11 (7.60) & 2.35 & 2.35\\
ben. & 500 & 3 & 0.803 & 0.955 & 0.954 & 0.951 & 0.942 & 1.12 (7.73) & 1.11 (7.62) & 2.45 & 2.45\\
ben. & 2000 & 0 & 0.870 & 0.949 & 0.948 & 0.950 & 0.943 & 1.11 (10.68) & 1.11 (10.61) & 2.85 & 2.47\\
ben. & 2000 & 3 & 0.861 & 0.953 & 0.953 & 0.945 & 0.945 & 1.11 (10.80) & 1.11 (10.73) & 2.60 & 2.60\\
ben. & 8000 & 0 & 0.943 & 0.953 & 0.952 & 0.960 & 0.943 & 1.11 (2.64) & 1.11 (2.63) & 4.61 & 2.47\\
ben. & 8000 & 3 & 0.938 & 0.944 & 0.944 & 0.955 & 0.938 & 1.11 (2.67) & 1.11 (2.67) & 4.59 & 2.69\\
mal. & 500 & 0 & 0.258 & 0.972 & 0.971 & 0.971 & 0.948 & 3.95 (7.17) & 3.73 (7.15) & 4.35 & 3.20\\
mal. & 500 & 3 & 0.532 & 0.954 & 0.952 & 0.971 & 0.939 & 3.42 (7.06) & 3.36 (7.07) & 4.29 & 3.31\\
mal. & 2000 & 0 & 0.251 & 0.974 & 0.973 & 0.970 & 0.945 & 5.94 (13.32) & 5.70 (13.31) & 7.70 & 3.60\\
mal. & 2000 & 3 & 0.593 & 0.954 & 0.955 & 0.969 & 0.950 & 4.83 (9.76) & 4.77 (9.76) & 7.51 & 3.80\\
mal. & 8000 & 0 & 0.258 & 0.973 & 0.973 & 0.970 & 0.946 & 10.02 (25.63) & 9.79 (25.62) & 14.40 & 3.70\\
mal. & 8000 & 3 & 0.524 & 0.960 & 0.959 & 0.969 & 0.945 & 8.39 (25.03) & 8.28 (25.04) & 14.25 & 4.09\\
\bottomrule
\end{tabular}
\end{table}

\paragraph{The price on all five-node classes with $K>0$.} To check the order $\mathcal R+K$ and the attainment without knowledge of the base point, we ran the procedures on all $96$ distinct pairs $(\Xeq,e)$ on five nodes with $K>0$, at one random base point each, and on the ten configurations of \cref{tab:simK}, with each Markov class of $\Xeq$ as the true graph, weak-edge coefficient $c/\sqrt n$ with $c\le32$, $n=32000$ and $1000$ replications per cell; all have two Markov classes. At every one of the $212$ pairs of a base point and a true graph, the price of \DSeff\ was below the bound of \cref{thm:K}(ii), at between $0.55$ and $0.84$ of it, and that of \DSOp\ below its bound, at between $0.50$ and $0.79$ of it; the price of \DSeff\ lies between $0.46$ and $0.95$ times $\mathcal R+K$. Estimating the directions cost almost nothing: the price of \DSeff\ was at most $1.072$ times that of the same procedure with the directions computed from the base point. The minimum coverage over the $2287$ cells in the class was $0.936$ for \DSeff\ and $0.938$ for \DSOp, with means $0.967$ and $0.968$, and a rerun of the nine lowest cells gave between $0.955$ and $0.968$ (\cref{app:sim}).

\paragraph{Uniformity in the margin.} On a six-node configuration that satisfies (F), (C), (D) and \eqref{eq:M} (\cref{app:sim}), with a base point on $\cM_{H_2}$ and the coefficients $\lambda$ of $f$ and $c/\sqrt n$ of $e$ in the true graph, we ran \DSOp\ for $n\in\{2000,8000\}$, $\lambda\in\{0.4,0.2,0.1,0.05,0.025,0.0125,0\}$ and $|c|\le16$, with $2000$ replications per cell. The largest length ratio over $|c|\le A\sqrt n$ behaves qualitatively as $\min\{A\sqrt n,1/\lambda\}$ predicts. For $\lambda=0.4$ it stays between $2.70$ and $2.78$ for every $A\sqrt n\ge4$; for $\lambda\le0.05$ it hardly depends on $\lambda$, down to $\lambda=0$, and grows with $A\sqrt n$, from about $2.0$ at $A\sqrt n=1$ to about $6.0$ at $n=2000$ and $7.0$ at $n=8000$ for $A\sqrt n=16$; and at fixed $(\lambda,A\sqrt n)$ the values at the two sample sizes differ by a factor of at most $1.18$. Hull, which does not screen, has ratios between $6.61$ and $13.57$ at $n=2000$ and between $12.23$ and $23.92$ at $n=8000$ over all cells in the class. The minimum coverage of \DSOp\ over the $291$ cells in the class is $0.958$.

\section{Discussion}\label{sec:discussion}

Locally near generic base points, the price of honest inference after causal discovery depends on which edge is weak, and the dependence can be read off the graph (\cref{thm:dichotomy}). Suppose that deleting the weak edge leaves the effect identifiable in the resulting MPDAG, and that the Markov equivalent explanations of the base point agree on the effect. Then honest intervals cost a bounded factor, of order $\mathcal R+K$ for an explicit cross-model constant $K$ and the ratio $\mathcal R$ of efficient lengths. In the limit they cost nothing at all exactly when $K=0$; a d-separation decides this, and it is the common case in our enumerations of all graphs with only $Y$ a sink (\cref{thm:free}). For two Markov classes the local minimax price is determined within a factor that does not depend on the graph. Its upper bound is approached by a procedure that knows neither the base point nor the weak edge, uniformly over compact sets of generic base points (\cref{thm:K,cor:RK,cor:compact}). Near a distribution at which a second edge is also absent, and under a graphical condition, the rate is uniform in the margin of that edge (\cref{thm:uniform}); it has no plateau of order $n^{-1/4}$ that does not depend on the margin and on the size of the weak coefficient. If either condition on the weak edge fails, then for every fixed $k$ below a threshold that depends on the base point, honest intervals cannot shrink at the parametric rate. This holds at the base point in the malignant case, and on every fixed neighbourhood of it in the third case. On the malignant side it is a consequence of non-identification at the base point, which $k$-triangle faithfulness does not exclude. The constant $k$, which plays a central role for uniform consistency \citep{spirtes2014uniformly}, does not enter the constants of the price on the benign side, only the neighbourhood and the sample size from which they apply.

\paragraph{Confidence sets of models.} \citet{edgar2026inferential} study likelihood-ratio confidence sets of causal models near a point at which two models intersect: in their Proposition~5.1, a model that is $c$-equivalent to the true one in the sense of \citet{evans2020model} is retained with probability $1-\alpha+o(1)$ within $O(n^{-1/(2c)})$ of the point and with probability $o(1)$ farther away, and Markov equivalent models are always retained. Our results concern the length of an interval for a functional, and the two pictures are complementary: a retained explanation costs length only through the difference of its effect from the target, so the three classes of \cref{thm:dichotomy} are the combinations of separability and effect difference. Markov equivalent explanations, always retained, are harmless exactly when \eqref{eq:K1p} holds; explanations that coincide with the truth at the base point but have a different effect cost a factor of order $\sqrt n$; and explanations in $\Xeq$ that are not Markov equivalent are separated at first order (\cref{thm:E}), so they cost a bounded factor governed by $K$, the ratio of the first-order effect difference to the first-order separation. In the four-variable model with unit noise variances the two explanations differ at second order when the margin vanishes, and the plateau $Bn^{-1/4}$ is attained at $a\asymp n^{-1/4}$, the exponent of their Proposition~5.1 with $c=2$. Without a known variance range they coincide when the margin vanishes (\cref{lem:Vf}), and this intermediate regime disappears.

\paragraph{Limitations.} The results are for linear Gaussian models, a skeleton that is a known superset of the true adjacencies, and locally a single weak edge, or two under \eqref{eq:M} (\cref{thm:uniform}); more weak edges, or two without \eqref{eq:M}, are not covered, and the global price over a class that allows any edge to be weak is of order $\sqrt n$ as soon as one edge is malignant (\cref{rem:global}). With required edges other than a time order, the classification assumes admissibility and, for one converse, that $\MP(H)$ is a chain graph; the criterion held without the chain assumption in every case we examined (\cref{app:certify}). The procedures need $k$, and in the four-variable model the variance range $\kappa$, as inputs, much as honest intervals in nonparametric regression need a known smoothness constant \citep{armstrong2018optimal}; a $k$ larger than the truth allows can make double screening screen out the true graph as $n$ grows, and in our simulations it then returned the empty set rather than a wrong interval, but we have no guarantee of this. The bounds on general graphs are local: the threshold $k^*(P_0)$ is bounded below, and the neighbourhood and the sample size from which the upper bounds apply are uniform, over compact sets of generic base points, but they degenerate near non-generic ones, and \cref{thm:uniform}, although uniform in the margin, is local around $P_{00}$ and requires \eqref{eq:M}, which without required edges holds in half of the six-node configurations with $\iota=1$ and $f$ benign in $H_2$. The factor $3.40$ is the ratio of two asymptotic bounds: it bounds how far the local minimax price is from either, but says nothing about the price at a given sample size. Sharper values ($1.2743$ and $1.422$) are available only in the local limit experiment of the four-variable model, the first computer-assisted and the second numerical. Finally, double screening enumerates the candidate graphs and their triangle constraints, which limits it to small graphs; an implementation that works on the MPDAG directly, and a $k$-triangle-faithfulness screen that tests fewer conditioning sets, remain to be developed.

\paragraph{Open questions.} What is the price when more than two edges are weak at the same time, or two without \eqref{eq:M}, where the weak edge can lose its benign status in another explanation of the doubly weak base point? \Cref{thm:uniform} gives the rate uniformly in the margin but not the constant: is the factor of \cref{cor:RK} uniform as well? With a known range of the noise variances, do intermediate rates such as the $n^{-1/4}$ of the four-variable model occur on general graphs? Can the criteria of \cref{prop:Kpos,thm:K-lambda} be read off $\MP(H)$ directly, as the identification conditions are? Is $K_+$ or $K$ the right constant for several Markov classes, and can a lower bound that combines $\mathcal R$ and $K$, or a sharper upper bound, reduce the factor $3.40$? One-sided screens are not available, because the sign of the weak coefficient is unknown. Is adaptation to an unknown $k$ possible, and can the classification be extended to latent variables or to non-Gaussian models \citep{wang2022uniformly}?

\section*{Use of AI tools}
AI assistants took part in discussing the choice of research direction; the author decided the main line of the work and carried it forward. AI assistants, mainly Claude (Anthropic), with DeepSeek and Codex (OpenAI) for review, carried out much of the mathematical derivation, wrote and checked the code, drafted and revised the text, and reviewed the proofs and analyses. The author has personally checked the proofs and results and takes full responsibility for the content.

\section*{Code availability}
Code and stored outputs that reproduce all numerical results are available at \url{https://github.com/Zrq23/honest-length-causal-discovery}; each table states the script that produced it.

\clearpage
\bibliographystyle{plainnat}
\bibliography{refs}

\begin{thebibliography}{68}
\providecommand{\natexlab}[1]{#1}
\providecommand{\url}[1]{\texttt{#1}}
\expandafter\ifx\csname urlstyle\endcsname\relax
  \providecommand{\doi}[1]{doi: #1}\else
  \providecommand{\doi}{doi: \begingroup \urlstyle{rm}\Url}\fi

\bibitem[Anderson(2003)]{anderson2003introduction}
T.~W. Anderson.
\newblock \emph{An Introduction to Multivariate Statistical Analysis}.
\newblock Wiley, Hoboken, NJ, third edition, 2003.

\bibitem[Andersson et~al.(1997)Andersson, Madigan, and
  Perlman]{andersson1997characterization}
Steen~A. Andersson, David Madigan, and Michael~D. Perlman.
\newblock A characterization of {M}arkov equivalence classes for acyclic
  digraphs.
\newblock \emph{The Annals of Statistics}, 25\penalty0 (2):\penalty0 505--541,
  1997.
\newblock \doi{10.1214/aos/1031833662}.

\bibitem[Andrews and Cheng(2012)]{andrews2012estimation}
Donald W.~K. Andrews and Xu~Cheng.
\newblock Estimation and inference with weak, semi-strong, and strong
  identification.
\newblock \emph{Econometrica}, 80\penalty0 (5):\penalty0 2153--2211, 2012.
\newblock \doi{10.3982/ECTA9456}.

\bibitem[Armstrong and Koles{\'a}r(2018)]{armstrong2018optimal}
Timothy~B. Armstrong and Michal Koles{\'a}r.
\newblock Optimal inference in a class of regression models.
\newblock \emph{Econometrica}, 86\penalty0 (2):\penalty0 655--683, 2018.
\newblock \doi{10.3982/ECTA14434}.

\bibitem[Armstrong et~al.(2020)Armstrong, Koles{\'a}r, and
  Kwon]{armstrong2020bias}
Timothy~B. Armstrong, Michal Koles{\'a}r, and Soonwoo Kwon.
\newblock Bias-aware inference in regularized regression models.
\newblock arXiv:2012.14823, 2020.

\bibitem[Armstrong et~al.(2025)Armstrong, Kline, and
  Sun]{armstrong2025adapting}
Timothy~B. Armstrong, Patrick Kline, and Liyang Sun.
\newblock Adapting to misspecification.
\newblock \emph{Econometrica}, 93\penalty0 (6):\penalty0 1981--2005, 2025.
\newblock \doi{10.3982/ECTA21991}.

\bibitem[Bachoc et~al.(2020)Bachoc, Preinerstorfer, and
  Steinberger]{bachoc2020uniformly}
Fran{\c{c}}ois Bachoc, David Preinerstorfer, and Lukas Steinberger.
\newblock Uniformly valid confidence intervals post-model-selection.
\newblock \emph{The Annals of Statistics}, 48\penalty0 (1):\penalty0 440--463,
  2020.
\newblock \doi{10.1214/19-AOS1815}.

\bibitem[Bang and Didelez(2023)]{bang2023wiser}
Christine~W. Bang and Vanessa Didelez.
\newblock Do we become wiser with time? {O}n causal equivalence with tiered
  background knowledge.
\newblock In \emph{Proceedings of the Thirty-Ninth Conference on Uncertainty in
  Artificial Intelligence}, volume 216 of \emph{Proceedings of Machine Learning
  Research}, pages 119--129, 2023.

\bibitem[Battey et~al.(2026)Battey, {Garc{\'i}a Rasines}, and
  Tang]{battey2026post}
H.~S. Battey, D.~{Garc{\'i}a Rasines}, and Y.~Tang.
\newblock Post-reduction inference for confidence sets of models.
\newblock \emph{Biometrika}, 113\penalty0 (4):\penalty0 asag045, 2026.
\newblock \doi{10.1093/biomet/asag045}.

\bibitem[Bei and Navjeevan(2026)]{bei2026inference}
Xinyue Bei and Manu Navjeevan.
\newblock Inference under first-order degeneracy.
\newblock arXiv:2602.07377, 2026.

\bibitem[Bhattacharya et~al.(2026)Bhattacharya, Ocelli, and
  Westling]{bhattacharya2026robust}
Rohit Bhattacharya, Ina Ocelli, and Ted Westling.
\newblock Robust weighted triangulation of causal effects under model
  uncertainty.
\newblock In \emph{Proceedings of the 42nd Conference on Uncertainty in
  Artificial Intelligence}, volume 337 of \emph{Proceedings of Machine Learning
  Research}, pages 567--584, 2026.

\bibitem[Cai and Guo(2017)]{cai2017confidence}
T.~Tony Cai and Zijian Guo.
\newblock Confidence intervals for high-dimensional linear regression: Minimax
  rates and adaptivity.
\newblock \emph{The Annals of Statistics}, 45\penalty0 (2):\penalty0 615--646,
  2017.
\newblock \doi{10.1214/16-AOS1461}.

\bibitem[Cai and Low(2004)]{cai2004adaptation}
T.~Tony Cai and Mark~G. Low.
\newblock An adaptation theory for nonparametric confidence intervals.
\newblock \emph{The Annals of Statistics}, 32\penalty0 (5):\penalty0
  1805--1840, 2004.
\newblock \doi{10.1214/009053604000000049}.

\bibitem[Chang et~al.(2026)Chang, Guo, and Malinsky]{chang2026post}
Ting-Hsuan Chang, Zijian Guo, and Daniel Malinsky.
\newblock Post-selection inference for causal effects after causal discovery.
\newblock \emph{Biometrika}, 113\penalty0 (1):\penalty0 asaf073, 2026.
\newblock \doi{10.1093/biomet/asaf073}.

\bibitem[Chickering(1995)]{chickering1995transformational}
David~Maxwell Chickering.
\newblock A transformational characterization of equivalent {B}ayesian network
  structures.
\newblock In \emph{Proceedings of the Eleventh Conference on Uncertainty in
  Artificial Intelligence}, pages 87--98, 1995.

\bibitem[Chickering(2002)]{chickering2002optimal}
David~Maxwell Chickering.
\newblock Optimal structure identification with greedy search.
\newblock \emph{Journal of Machine Learning Research}, 3:\penalty0 507--554,
  2002.

\bibitem[Donoho(1994)]{donoho1994statistical}
David~L. Donoho.
\newblock Statistical estimation and optimal recovery.
\newblock \emph{The Annals of Statistics}, 22\penalty0 (1):\penalty0 238--270,
  1994.
\newblock \doi{10.1214/aos/1176325367}.

\bibitem[Donoho and Liu(1991)]{donoho1991geometrizingII}
David~L. Donoho and Richard~C. Liu.
\newblock Geometrizing rates of convergence, {II}.
\newblock \emph{The Annals of Statistics}, 19\penalty0 (2):\penalty0 633--667,
  1991.
\newblock \doi{10.1214/aos/1176348114}.

\bibitem[Edgar and Battey(2026)]{edgar2026inferential}
Charlotte Edgar and Heather Battey.
\newblock On inferential equivalence classes of causal models.
\newblock arXiv:2608.08327, 2026.

\bibitem[Elliott et~al.(2015)Elliott, M{\"u}ller, and
  Watson]{elliott2015nearly}
Graham Elliott, Ulrich~K. M{\"u}ller, and Mark~W. Watson.
\newblock Nearly optimal tests when a nuisance parameter is present under the
  null hypothesis.
\newblock \emph{Econometrica}, 83\penalty0 (2):\penalty0 771--811, 2015.

\bibitem[Evans(2020)]{evans2020model}
Robin~J. Evans.
\newblock Model selection and local geometry.
\newblock \emph{The Annals of Statistics}, 48\penalty0 (6):\penalty0
  3513--3544, 2020.
\newblock \doi{10.1214/19-AOS1940}.

\bibitem[Freeman(2026)]{freeman2026local}
Hugo Freeman.
\newblock Local asymptotic power of honest confidence intervals.
\newblock arXiv:2607.10558v3, 2026.

\bibitem[Gradu et~al.(2025)Gradu, Zrnic, Wang, and Jordan]{gradu2025valid}
Paula Gradu, Tijana Zrnic, Yixin Wang, and Michael~I. Jordan.
\newblock Valid inference after causal discovery.
\newblock \emph{Journal of the American Statistical Association}, 120\penalty0
  (550):\penalty0 1127--1138, 2025.
\newblock \doi{10.1080/01621459.2024.2402089}.

\bibitem[Guo et~al.(2026)Guo, Wu, and Yiu]{guo2026price}
Dongxin Guo, Jikun Wu, and Siu~Ming Yiu.
\newblock The price of valid inference after causal discovery.
\newblock In \emph{Proceedings of the 42nd Conference on Uncertainty in
  Artificial Intelligence}, volume 337 of \emph{Proceedings of Machine Learning
  Research}, pages 1804--1822, 2026.

\bibitem[Guo and Perkovi{\'c}(2022)]{guo2022efficient}
F.~Richard Guo and Emilija Perkovi{\'c}.
\newblock Efficient least squares for estimating total effects under linearity
  and causal sufficiency.
\newblock \emph{Journal of Machine Learning Research}, 23\penalty0
  (104):\penalty0 1--41, 2022.

\bibitem[Guo(2023)]{guo2023causal}
Zijian Guo.
\newblock Causal inference with invalid instruments: post-selection problems
  and a solution using searching and sampling.
\newblock \emph{Journal of the Royal Statistical Society Series B}, 85\penalty0
  (3):\penalty0 959--985, 2023.
\newblock \doi{10.1093/jrsssb/qkad049}.

\bibitem[Hauser and B{\"u}hlmann(2012)]{hauser2012characterization}
Alain Hauser and Peter B{\"u}hlmann.
\newblock Characterization and greedy learning of interventional {M}arkov
  equivalence classes of directed acyclic graphs.
\newblock \emph{Journal of Machine Learning Research}, 13:\penalty0 2409--2464,
  2012.

\bibitem[Henckel et~al.(2022)Henckel, Perkovi{\'c}, and
  Maathuis]{henckel2022graphical}
Leonard Henckel, Emilija Perkovi{\'c}, and Marloes~H. Maathuis.
\newblock Graphical criteria for efficient total effect estimation via
  adjustment in causal linear models.
\newblock \emph{Journal of the Royal Statistical Society Series B: Statistical
  Methodology}, 84\penalty0 (2):\penalty0 579--599, 2022.
\newblock \doi{10.1111/rssb.12451}.

\bibitem[Hirano and Porter(2012)]{hirano2012impossibility}
Keisuke Hirano and Jack~R. Porter.
\newblock Impossibility results for nondifferentiable functionals.
\newblock \emph{Econometrica}, 80\penalty0 (4):\penalty0 1769--1790, 2012.

\bibitem[Hjort and Claeskens(2003)]{hjort2003frequentist}
Nils~Lid Hjort and Gerda Claeskens.
\newblock Frequentist model average estimators.
\newblock \emph{Journal of the American Statistical Association}, 98\penalty0
  (464):\penalty0 879--899, 2003.
\newblock \doi{10.1198/016214503000000828}.

\bibitem[Imbens and Manski(2004)]{imbens2004confidence}
Guido~W. Imbens and Charles~F. Manski.
\newblock Confidence intervals for partially identified parameters.
\newblock \emph{Econometrica}, 72\penalty0 (6):\penalty0 1845--1857, 2004.
\newblock \doi{10.1111/j.1468-0262.2004.00555.x}.

\bibitem[Joshi(1969)]{joshi1969admissibility}
V.~M. Joshi.
\newblock Admissibility of the usual confidence sets for the mean of a
  univariate or bivariate normal population.
\newblock \emph{The Annals of Mathematical Statistics}, 40\penalty0
  (3):\penalty0 1042--1067, 1969.
\newblock \doi{10.1214/aoms/1177697608}.

\bibitem[Kabaila and Giri(2009)]{kabaila2009confidence}
Paul Kabaila and Khageswor Giri.
\newblock Confidence intervals in regression utilizing prior information.
\newblock \emph{Journal of Statistical Planning and Inference}, 139\penalty0
  (10):\penalty0 3419--3429, 2009.
\newblock \doi{10.1016/j.jspi.2009.03.018}.

\bibitem[Kabaila and Leeb(2006)]{kabaila2006large}
Paul Kabaila and Hannes Leeb.
\newblock On the large-sample minimal coverage probability of confidence
  intervals after model selection.
\newblock \emph{Journal of the American Statistical Association}, 101\penalty0
  (474):\penalty0 619--629, 2006.
\newblock \doi{10.1198/016214505000001140}.

\bibitem[Kalisch and B{\"u}hlmann(2007)]{kalisch2007estimating}
Markus Kalisch and Peter B{\"u}hlmann.
\newblock Estimating high-dimensional directed acyclic graphs with the
  {PC}-algorithm.
\newblock \emph{Journal of Machine Learning Research}, 8:\penalty0 613--636,
  2007.

\bibitem[Kang et~al.(2026)Kang, Mazumder, Lembo, and
  Taeb]{kang2026interpretable}
Jian Kang, Arnav Mazumder, Giovanni~Carter Lembo, and Armeen Taeb.
\newblock Interpretable confidence regions for causal discovery.
\newblock arXiv:2611.17771, 2026.

\bibitem[Kiefer(1957)]{kiefer1957invariance}
J.~Kiefer.
\newblock Invariance, minimax sequential estimation, and continuous time
  processes.
\newblock \emph{The Annals of Mathematical Statistics}, 28\penalty0
  (3):\penalty0 573--601, 1957.

\bibitem[Laurent and Massart(2000)]{laurent2000adaptive}
B{\'e}atrice Laurent and Pascal Massart.
\newblock Adaptive estimation of a quadratic functional by model selection.
\newblock \emph{The Annals of Statistics}, 28\penalty0 (5):\penalty0
  1302--1338, 2000.

\bibitem[Le~Cam(1986)]{lecam1986asymptotic}
Lucien Le~Cam.
\newblock \emph{Asymptotic Methods in Statistical Decision Theory}.
\newblock Springer, New York, 1986.

\bibitem[Leeb and P{\"o}tscher(2005)]{leeb2005model}
Hannes Leeb and Benedikt~M. P{\"o}tscher.
\newblock Model selection and inference: Facts and fiction.
\newblock \emph{Econometric Theory}, 21\penalty0 (1):\penalty0 21--59, 2005.
\newblock \doi{10.1017/S0266466605050036}.

\bibitem[Low(1997)]{low1997nonparametric}
Mark~G. Low.
\newblock On nonparametric confidence intervals.
\newblock \emph{The Annals of Statistics}, 25\penalty0 (6):\penalty0
  2547--2554, 1997.
\newblock \doi{10.1214/aos/1030741084}.

\bibitem[Maathuis et~al.(2009)Maathuis, Kalisch, and
  B{\"u}hlmann]{maathuis2009estimating}
Marloes~H. Maathuis, Markus Kalisch, and Peter B{\"u}hlmann.
\newblock Estimating high-dimensional intervention effects from observational
  data.
\newblock \emph{The Annals of Statistics}, 37\penalty0 (6A):\penalty0
  3133--3164, 2009.
\newblock \doi{10.1214/09-AOS685}.

\bibitem[Malinsky(2024)]{malinsky2024cautious}
Daniel Malinsky.
\newblock A cautious approach to constraint-based causal model selection.
\newblock arXiv:2404.18232, 2024.

\bibitem[McCloskey(2017)]{mccloskey2017bonferroni}
Adam McCloskey.
\newblock Bonferroni-based size-correction for nonstandard testing problems.
\newblock \emph{Journal of Econometrics}, 200\penalty0 (1):\penalty0 17--35,
  2017.

\bibitem[Meek(1995)]{meek1995causal}
Christopher Meek.
\newblock Causal inference and causal explanation with background knowledge.
\newblock In \emph{Proceedings of the Eleventh Conference on Uncertainty in
  Artificial Intelligence}, pages 403--410, 1995.

\bibitem[Mityagin(2020)]{mityagin2020zero}
Boris~S. Mityagin.
\newblock The zero set of a real analytic function.
\newblock \emph{Mathematical Notes}, 107\penalty0 (3):\penalty0 529--530, 2020.
\newblock \doi{10.1134/S0001434620030189}.

\bibitem[Perkovi{\'c}(2020)]{perkovic2020identifying}
Emilija Perkovi{\'c}.
\newblock Identifying causal effects in maximally oriented partially directed
  acyclic graphs.
\newblock In \emph{Proceedings of the 36th Conference on Uncertainty in
  Artificial Intelligence}, volume 124 of \emph{Proceedings of Machine Learning
  Research}, pages 530--539, 2020.
\newblock Revised proof of Lemma~A.3 in arXiv:1910.02997v3.

\bibitem[Perkovi{\'c} et~al.(2017)Perkovi{\'c}, Kalisch, and
  Maathuis]{perkovic2017interpreting}
Emilija Perkovi{\'c}, Markus Kalisch, and Marloes~H. Maathuis.
\newblock Interpreting and using {CPDAGs} with background knowledge.
\newblock In \emph{Proceedings of the Thirty-Third Conference on Uncertainty in
  Artificial Intelligence}, 2017.

\bibitem[P{\"o}tscher(2009)]{potscher2009confidence}
Benedikt~M. P{\"o}tscher.
\newblock Confidence sets based on sparse estimators are necessarily large.
\newblock \emph{Sankhy{\=a}: The Indian Journal of Statistics, Series A},
  71\penalty0 (1):\penalty0 1--18, 2009.
\newblock arXiv:0711.1036.

\bibitem[Pratt(1961)]{pratt1961length}
John~W. Pratt.
\newblock Length of confidence intervals.
\newblock \emph{Journal of the American Statistical Association}, 56\penalty0
  (295):\penalty0 549--567, 1961.

\bibitem[Robins et~al.(2003)Robins, Scheines, Spirtes, and
  Wasserman]{robins2003uniform}
James~M. Robins, Richard Scheines, Peter Spirtes, and Larry Wasserman.
\newblock Uniform consistency in causal inference.
\newblock \emph{Biometrika}, 90\penalty0 (3):\penalty0 491--515, 2003.
\newblock \doi{10.1093/biomet/90.3.491}.

\bibitem[Romano et~al.(2014)Romano, Shaikh, and Wolf]{romano2014practical}
Joseph~P. Romano, Azeem~M. Shaikh, and Michael Wolf.
\newblock A practical two-step method for testing moment inequalities.
\newblock \emph{Econometrica}, 82\penalty0 (5):\penalty0 1979--2002, 2014.

\bibitem[Rotnitzky and Smucler(2020)]{rotnitzky2020efficient}
Andrea Rotnitzky and Ezequiel Smucler.
\newblock Efficient adjustment sets for population average causal treatment
  effect estimation in graphical models.
\newblock \emph{Journal of Machine Learning Research}, 21\penalty0
  (188):\penalty0 1--86, 2020.

\bibitem[Spirtes and Zhang(2014)]{spirtes2014uniformly}
Peter Spirtes and Jiji Zhang.
\newblock A uniformly consistent estimator of causal effects under the
  $k$-triangle-faithfulness assumption.
\newblock \emph{Statistical Science}, 29\penalty0 (4):\penalty0 662--678, 2014.
\newblock \doi{10.1214/13-STS429}.

\bibitem[Spirtes et~al.(2000)Spirtes, Glymour, and
  Scheines]{spirtes2000causation}
Peter Spirtes, Clark Glymour, and Richard Scheines.
\newblock \emph{Causation, Prediction, and Search}.
\newblock MIT Press, Cambridge, MA, second edition, 2000.

\bibitem[Stoye(2009)]{stoye2009more}
J{\"o}rg Stoye.
\newblock More on confidence intervals for partially identified parameters.
\newblock \emph{Econometrica}, 77\penalty0 (4):\penalty0 1299--1315, 2009.
\newblock \doi{10.3982/ECTA7347}.

\bibitem[Strieder and Drton(2023)]{strieder2023confidence}
David Strieder and Mathias Drton.
\newblock Confidence in causal inference under structure uncertainty in linear
  causal models with equal variances.
\newblock \emph{Journal of Causal Inference}, 11\penalty0 (1):\penalty0
  20230030, 2023.
\newblock \doi{10.1515/jci-2023-0030}.

\bibitem[Strieder and Drton(2024{\natexlab{a}})]{strieder2024dual}
David Strieder and Mathias Drton.
\newblock Dual likelihood for causal inference under structure uncertainty.
\newblock In \emph{Proceedings of the Third Conference on Causal Learning and
  Reasoning}, volume 236 of \emph{Proceedings of Machine Learning Research},
  pages 1--17, 2024{\natexlab{a}}.

\bibitem[Strieder and Drton(2024{\natexlab{b}})]{strieder2024identifying}
David Strieder and Mathias Drton.
\newblock Identifying total causal effects in linear models under partial
  homoscedasticity.
\newblock In \emph{Proceedings of the 12th International Conference on
  Probabilistic Graphical Models}, volume 246 of \emph{Proceedings of Machine
  Learning Research}, pages 213--230, 2024{\natexlab{b}}.

\bibitem[Strieder et~al.(2021)Strieder, Freidling, Haffner, and
  Drton]{strieder2021confidence}
David Strieder, Tobias Freidling, Stefan Haffner, and Mathias Drton.
\newblock Confidence in causal discovery with linear causal models.
\newblock In \emph{Proceedings of the Thirty-Seventh Conference on Uncertainty
  in Artificial Intelligence}, volume 161 of \emph{Proceedings of Machine
  Learning Research}, pages 1217--1226, 2021.

\bibitem[Teh et~al.(2025)Teh, Sadeghi, and Soo]{teh2025towards}
Kai~Z. Teh, Kayvan Sadeghi, and Terry Soo.
\newblock Towards robust causal effect identification beyond {M}arkov
  equivalence.
\newblock ICML 2025 Workshop on Scaling Up Intervention Models;
  arXiv:2506.15561, 2025.

\bibitem[Uhler et~al.(2013)Uhler, Raskutti, B{\"u}hlmann, and
  Yu]{uhler2013geometry}
Caroline Uhler, Garvesh Raskutti, Peter B{\"u}hlmann, and Bin Yu.
\newblock Geometry of the faithfulness assumption in causal inference.
\newblock \emph{The Annals of Statistics}, 41\penalty0 (2):\penalty0 436--463,
  2013.
\newblock \doi{10.1214/12-AOS1080}.

\bibitem[van~der Vaart(1998)]{vaart1998asymptotic}
A.~W. van~der Vaart.
\newblock \emph{Asymptotic Statistics}.
\newblock Cambridge University Press, Cambridge, 1998.

\bibitem[Wang and Spirtes(2022)]{wang2022uniformly}
Shuyan Wang and Peter Spirtes.
\newblock A uniformly consistent estimator of non-{G}aussian causal effects
  under the $k$-triangle-faithfulness assumption.
\newblock In \emph{Proceedings of the First Conference on Causal Learning and
  Reasoning}, volume 177 of \emph{Proceedings of Machine Learning Research},
  pages 861--876, 2022.

\bibitem[Wang et~al.(2026)Wang, Kolar, and Drton]{wang2026confidence}
Y.~Samuel Wang, Mladen Kolar, and Mathias Drton.
\newblock Confidence sets for causal orderings.
\newblock \emph{Journal of the American Statistical Association}, 121\penalty0
  (553):\penalty0 690--703, 2026.
\newblock \doi{10.1080/01621459.2025.2542552}.

\bibitem[Witte et~al.(2020)Witte, Henckel, Maathuis, and
  Didelez]{witte2020efficient}
Janine Witte, Leonard Henckel, Marloes~H. Maathuis, and Vanessa Didelez.
\newblock On efficient adjustment in causal graphs.
\newblock \emph{Journal of Machine Learning Research}, 21\penalty0
  (246):\penalty0 1--45, 2020.

\bibitem[Zhang and Spirtes(2003)]{zhang2003strong}
Jiji Zhang and Peter Spirtes.
\newblock Strong faithfulness and uniform consistency in causal inference.
\newblock In \emph{Proceedings of the Nineteenth Conference on Uncertainty in
  Artificial Intelligence}, pages 632--639, 2003.

\bibitem[Zhang and Spirtes(2008)]{zhang2008detection}
Jiji Zhang and Peter Spirtes.
\newblock Detection of unfaithfulness and robust causal inference.
\newblock \emph{Minds and Machines}, 18\penalty0 (2):\penalty0 239--271, 2008.
\newblock \doi{10.1007/s11023-008-9096-4}.

\end{thebibliography}

\appendix
\crefalias{section}{appendix}
\crefalias{subsection}{subappendix}
\clearpage
\begin{table}[h]
\centering
\caption{Evidence behind the main results. ``Proof'' means a complete proof in this paper under the stated conditions; ``computer-assisted'' means a calculation with an explicit error allowance under a floating-point error model, not a formally verified computation; ``numerical'' means a computation without an error bound; ``Monte Carlo'' refers to simulation.}
\label{tab:evidence}
\small
\begin{tabular}{p{0.42\textwidth}p{0.52\textwidth}}
\toprule
Result & Evidence\\
\midrule
Sets: bounded price (\cref{prop:union}) & Proof, exact in finite samples\\
Four-variable rate (\cref{thm:rate-lower,thm:rate-upper}) & Proofs; exact finite-sample honesty\\
Dual lower bound in the four-variable model (\cref{thm:lfp}) & Proof; the value $1.2743$ computer-assisted\\
Interval procedure, price about $1.422$ & Tails in closed form, middle numerical\\
Select-then-estimate (\cref{prop:naive,prop:naive-more}) & Proofs; coverage values and the sign of $\xi_C$ numerical\\
Criteria for \eqref{eq:K2}, \eqref{eq:K1p} (\cref{prop:graph,thm:K1p-graph,lem:tier,lem:tier-cut}) & Proofs; with required edges other than a time order under admissibility and, for one converse, a chain-graph hypothesis\\
First-order distinguishability, non-identified sides (\cref{thm:E,prop:malignant,prop:K1p-nec}) & Proofs\\
Honesty, local minimax price, adaptive attainment (\cref{thm:coverage,thm:K,cor:RK}; \cref{prop:CI,thm:explicit}) & Proofs; constants by numerical quadrature; simulated on all five-node classes with $K>0$\\
Canonical experiment, $K>0$, near a malignant stratum (\cref{prop:canonical,prop:Kpos,thm:K-lambda}) & Proofs; canonical prices by a linear program; the graphical criteria also checked in exact arithmetic on five and six nodes\\
When honesty is free (\cref{thm:free,lem:strict}) & Proofs; $\delta_{\mathrm{can}}$ evaluated by quadrature, and free pairs counted on five and six nodes\\
Uniformity in the margin (\cref{thm:uniform}) & Proof; the modulus of continuity computed on a grid for all six-node configurations with (F), (C), (D) and \eqref{eq:M}; simulated on one of them\\
Simulations (\cref{sec:sim}) & Monte Carlo\\
\bottomrule
\end{tabular}
\end{table}

\section*{Notation}\label{app:notation}
\addcontentsline{toc}{section}{Notation}

\begin{center}
\footnotesize
\renewcommand{\arraystretch}{1.0}
\begin{tabular}{@{}>{\raggedright\arraybackslash}p{0.19\textwidth}>{\raggedright\arraybackslash}p{0.55\textwidth}>{\raggedright\arraybackslash}p{0.2\textwidth}@{}}
\toprule
Symbol & Meaning & Where\\
\midrule
$X$, $Y$, $\theta_D(P)$ & treatment, outcome; total effect of $X$ on $Y$ under the DAG $D$ & \cref{sec:setup}\\
$\cM_D$ & linear Gaussian models along $D$ (zero coefficients allowed) & \cref{sec:setup}\\
$D\simeq D'$ & $D$ and $D'$ are Markov equivalent & \cref{sec:setup}\\
$\thadj_D$ & coefficient of $X$ in the regression of $Y$ on $X$ and $\pa_D(X)$ & \eqref{eq:adj}\\
$s_{\mathrm{eff}}$, $\orl_n$ & efficient standard deviation of $\theta_{G^*}$ in $\cM_{G^*}$; efficient interval length $2z_{1-\alpha/2}s_{\mathrm{eff}}n^{-1/2}$ & \cref{sec:setup}\\
price, $\Pi_n$ & worst-case ratio of expected length to $\orl_n$ over a class & \cref{sec:setup}\\
$\sigl^2$, $\sigu^2$ & bounds on the eigenvalues of $\Sigma$ & \cref{sec:setup}\\
$S$, $\cK$, $R$, $\cG$ & skeleton; background knowledge ($Y$ a sink and the required edges $R$); candidate class of all orientations of $S$ satisfying $\cK$ & \cref{sec:gen-setting}\\
tier knowledge & $R$ the edges between ordered tiers of $V\setminus\{Y\}$, oriented forward & \cref{sec:gen-setting}\\
$G$, $e$, $H$ & a candidate, its weak edge, and $H=G-e$ & \cref{sec:gen-setting}\\
$[D]_Y$ & DAGs Markov equivalent to $D$ that satisfy $\cK$ & \cref{sec:gen-setting}\\
$P_0$, $\Theta_H$ & base point, a distribution in $\cM_H$ faithful to $H$; parameter space of $H$ & \cref{sec:gen-setting}\\
$\cX$, $\Xeq$, $\Xcanc$ & candidates that explain $P_0$; those with $G'-e$ equivalent to $H$; those in $\cX$ not in $\Xeq$ & \cref{sec:gen-setting}\\
$\cM_k$, $k$, $\lambda^{\mathrm{std}}$, $\rho$ & $k$-triangle-faithful class and its constant; standardized coefficient; partial correlation given a legal set & \cref{sec:gen-setting}\\
$\MP(H)$, bad path & MPDAG of $[H]_Y$; proper possibly causal path from $X$ to $Y$ starting with an undirected edge & \cref{sec:gen-criterion}\\
$S_c$, $\mathrm{CC}_D$, $\MP''(S_c)$ & Markov classes of $\Xeq$; common children of the endpoints of $e$; MPDAG of \cref{thm:K1p-graph} & \cref{sec:gen-criterion}\\
\eqref{eq:K2}, \eqref{eq:K1p} & identification of the effect across $\Xeq$, and within each class $S_c$ & \cref{sec:gen-criterion}\\
\DSA, \DSOp, \DSeff & double screening with adjustment rule $A$, with optimal adjustment sets, with efficient intervals & \cref{def:ds}\\
$\gamma_1$, $\gamma_2$, $\gamma_3$, $\alpha'$; $\varepsilon_n$, $\tau_n$ & levels of the likelihood ratio, $k$-TF and directed tests and of the candidate intervals; split fraction and threshold of the structure estimate & \cref{def:ds}\\
$I$, $\delta$ & Fisher information separating two explanations along $e$; first-order effect discrepancy & \cref{sec:gen-benign}\\
$K$, $K_+$, $m$ & cross-model constant; its sum over Markov classes; number of Markov classes of $\Xeq$ & \cref{sec:gen-benign}\\
$\mathcal R$, $k^*$ & largest ratio of efficient standard deviations in $\Xeq$ (calligraphic, not the set $R$); threshold on $k$ & \cref{sec:gen-benign}\\
$T_D$, $\hat\tau_D$, $\mu_D$, $u_{S^*,S'}$ & likelihood ratio statistic of $\cM_D$; its whitened coordinates and their population version; direction of the directed test & \cref{sec:gen-benign}\\
$q$, $\varpi_d(\gamma)$, $\bar\varpi_d$, $e_m$ & non-adjacent pairs of $S$; slope functions; intercept term & \cref{sec:gen-benign}; \cref{thm:K}\\
$K^A$, $K^A_+$, $\mathcal R^A$ & analogues of $K$, $K_+$, $\mathcal R$ for an adjustment rule $A$ & \cref{rem:Oset}\\
\DSci, \DSlrO, \DSlreff, \DSdir & variants used in the proofs and simulations: conditional independence screen; no directed test; directions known & \cref{prop:CI,thm:explicit}\\
$\Pi^*_k(P_0)$, $a_m$, $b_m$ & local minimax price over $\cM_k$; constants of its upper bound ($m\ge2$) & \cref{cor:RK}\\
free, $\delta_{\mathrm{can}}$, \DSsw & $K\equiv0$ by \cref{prop:Kpos}; bound on $\Pi^*_{\mathrm{can}}-1$; switching variant & \cref{thm:free}\\
$\cN$ & Lebesgue-null set of base points excluded in \cref{thm:K} & \cref{app:benign}\\
$\Psi_\alpha$, $c_P$, $c_K$ & lower bound function and its constants & \cref{thm:K}\\
$\lambda$, (F), (C), (D) & margin of an edge $f$ of $H$; conditions of \cref{thm:K-lambda} & \cref{sec:gen-constant}\\
$H_2$, $P_{00}$, $\cX_2$, \eqref{eq:M}, $\cM_k(\lambda,A;r)$ & $H-f$; base point on $\cM_{H_2}$; $\{G^*:G^*-e-f\simeq H_2\}$; condition on $f$; class with $|\beta_e|\le A$, $|\beta_f|\ge\lambda$ & \cref{sec:gen-uniform}\\
$w$, $b$, $\kappa$, $\Pi^*$ & four-variable model: coefficients, variance range, local price constant & \cref{sec:toy}\\
\DSV, \DSVp & honest procedures of the four-variable model & \cref{def:dsv}\\
Naive, Union, Hull, Oracle & select-then-estimate; union of the candidate intervals as a set; their convex hull; interval of the true graph & \cref{sec:sim}\\
\bottomrule
\end{tabular}
\end{center}

\section{Proofs for the four-variable model}\label{app:toy}

Throughout, $\sigma_W^2,\sigma_X^2,\sigma_Z^2,\sigma_Y^2$ are the noise variances, constants $c,C$ depend only on the fixed quantities $(\bar w,\bar A,\kappa,\alpha,\gamma_1,\gamma_2)$, and $\eta=\gamma_2/8$.

\subsection{The matching threshold}\label{app:abar}

At $w=0$, a mediator model with noise variances $\sigma_X^2,\sigma_Z^2$ has $\Var Z=\sigma_Z^2+a^2\sigma_X^2$ and $\Var(X\mid Z)=\sigma_X^2\sigma_Z^2/(\sigma_Z^2+a^2\sigma_X^2)$. A confounder model reproduces its covariance of $(X,Z)$ exactly if and only if its noise variances are $\Var Z$ for $Z$ and $\Var(X\mid Z)$ for $X$, with coefficient $a\sigma_X^2/\Var Z$. With variances restricted to $[1,\kappa]$ this requires $\sigma_Z^2+a^2\sigma_X^2\le\kappa$ and $\sigma_X^2\sigma_Z^2\ge\sigma_Z^2+a^2\sigma_X^2$, that is
\[
a^2\le\min\Bigl\{\frac{\kappa-\sigma_Z^2}{\sigma_X^2},\ \frac{\sigma_Z^2(\sigma_X^2-1)}{\sigma_X^2}\Bigr\}.
\]
The two bounds are equal at $\sigma_Z^2=\kappa/\sigma_X^2$, where $a^2\le\kappa(\sigma_X^2-1)/\sigma_X^4$. Over $\sigma_X^2\in[1,\kappa]$ this is maximized at $\sigma_X^2=2$ with value $\kappa/4$ if $\kappa\ge2$, and at $\sigma_X^2=\kappa$ with value $1-1/\kappa$ if $\kappa<2$. This is \eqref{eq:abar}; we write $(m^*,s^*)=(2,\kappa/2)$ or $(\kappa,1)$ for the least favourable values of $(\sigma_X^2,\sigma_Z^2)$.

\subsection{Proof of \texorpdfstring{\cref{thm:rate-lower}}{the lower bound}}\label{app:rate-lower}

\paragraph{Tools.} We use \cref{lem:two-point} with pairs $(g_1,P)$, $(g_2,Q)$ in the class. All members below have $\beta$, $b$, $\sigma_Y^2$ such that $Y=\beta X+bZ+\eps_Y$, so the joint law factors as the law of $(W,X,Z)$ times the conditional law of $Y$ given $(W,X,Z)$. If the conditional laws coincide, $\KL(P\|Q)=\KL(P_{WXZ}\|Q_{WXZ})$. For centred Gaussian laws of $(W,X,Z)$ with covariance matrices $\Sigma_P,\Sigma_Q$, writing $E=\Sigma_Q^{-1/2}(\Sigma_P-\Sigma_Q)\Sigma_Q^{-1/2}$,
\[
\KL(P_{WXZ}\|Q_{WXZ})=\tfrac12\bigl(\mathrm{tr}E-\log\det(I+E)\bigr)\le\tfrac12\|E\|_F^2\le\frac{\|\Sigma_P-\Sigma_Q\|_F^2}{2\lambda_{\min}(\Sigma_Q)^2}\qquad\text{if }\|E\|\le\tfrac12 .
\]
In every member, $(W,X,Z)$ is generated by a lower triangular structural system with coefficients bounded by $\bar w$ and $C_A:=\bar A(\bar w^2+\kappa)$ and noise variances in $[1,\kappa]$, so $\lambda_{\min}(\Sigma_Q)\ge c_\Sigma>0$ with $c_\Sigma$ depending only on $(\bar w,\bar A,\kappa)$. By Pinsker's inequality $\dtv(P^n,Q^n)\le\sqrt{n\KL(P\|Q)/2}$. In each pair we choose the free parameter so that $\dtv(P^n,Q^n)\le(1-2\alpha)/2$, and \cref{lem:two-point} gives $\sup\E\Leb(C_n)\ge\tfrac12(1-2\alpha)|\theta_{g_1}(P)-\theta_{g_2}(Q)|$. All constants $c,c_1,C$ below depend only on $(\bar w,\bar A,\kappa,\alpha)$.

\paragraph{Pair 1: the $n^{-1/2}$ term.} Both members are $M$ with unit noise variances, $w=\lambda$, $a=0$ and $b=0$, and $\beta=0$ for $P$, $\beta=d/\sqrt n$ for $Q$. Both lie in $\cM(\lambda,A,B,\kappa)$ for every $A>0$ and $B>0$. Here the conditional laws of $Y$ differ, but only through the mean $\beta X$: $(W,X,Z)$ has the same law under both and $n\KL(P\|Q)=n\,\E(X^2)(d/\sqrt n)^2/2=d^2(1+\lambda^2)/2\le d^2(1+\bar w^2)/2$, symmetric in $P,Q$. With $d=d(\bar w,\alpha)$ small, $\dtv\le(1-2\alpha)/2$, and the targets differ by $d/\sqrt n$. This gives $\mathrm{LB}_1=c\,n^{-1/2}$.

\paragraph{Pair 2: unit variances.} $P$ is $M$ and $Q$ is $C$, both with unit noise variances, $w=\lambda$, the same $a\ge0$, $b=B$ and $\beta=0$; both lie in the class whenever $a\le A$. The $Y$-equations coincide, so only $(W,X,Z)$ contributes to the divergence. The covariance matrices are
\[
\Sigma_P=\begin{pmatrix}1&\lambda&a\lambda\\ \lambda&\lambda^2+1&a(\lambda^2+1)\\ a\lambda&a(\lambda^2+1)&a^2(\lambda^2+1)+1\end{pmatrix},\qquad
\Sigma_Q=\begin{pmatrix}1&\lambda&0\\ \lambda&\lambda^2+a^2+1&a\\ 0&a&1\end{pmatrix},
\]
so all entries of $\Sigma_P-\Sigma_Q$ are listed by $\Delta\Var X=-a^2$, $\Delta\Cov(X,Z)=a\lambda^2$, $\Delta\Var Z=a^2(1+\lambda^2)$, $\Delta\Cov(W,Z)=a\lambda$, and $\|\Sigma_P-\Sigma_Q\|_F\le C(a^2+\lambda a)$. The targets are $\theta_M=aB$ and $\theta_C=0$. Take $a=\min\{A,c_1n^{-1/4},c_1/(\lambda\sqrt n)\}$. Then $a^2\le c_1^2n^{-1/2}$ and $\lambda a\le c_1n^{-1/2}$, so $\|E\|\le Cc_1n^{-1/2}\le\tfrac12$ for $n\ge n_0$, and $n\KL\le Cc_1^2$. With $c_1$ small this gives
\[
\mathrm{LB}_2=c\,B\min\{A,\,n^{-1/4},\,1/(\lambda\sqrt n)\}.
\]

\paragraph{Pair 3: exact matching, $\kappa>1$.} If $\kappa=1$ then $\abar=0$ and this pair is not used. Let $(m^*,s^*)$ be the least favourable variances of \cref{app:abar} and let $a\in[0,\abar/2]$. $P$ is $M$ with $\sigma_W^2=\sigma_Y^2=1$, $(\sigma_X^2,\sigma_Z^2)=(m^*,s^*)$, $w=\lambda$, $b=B$, $\beta=0$. Under $P$, $\Var X=\lambda^2+m^*$, $\Cov(X,Z)=a\Var X$, $\Var Z=a^2\Var X+s^*$, $\Cov(W,X)=\lambda$ and $\Cov(W,Z)=a\lambda$. $Q$ is $C$ with $\sigma_W^2=\sigma_Y^2=1$, $w=\lambda$, $b=B$, $\beta=0$ and
\[
\sigma_Z'^2=\Var_PZ,\qquad a'=a\frac{\Var_PX}{\Var_PZ},\qquad \sigma_X'^2=\Var_P(X\mid Z)-\lambda^2=\Var_PX-\frac{\Cov_P(X,Z)^2}{\Var_PZ}-\lambda^2 .
\]
Under $Q$, $\Var Z=\sigma_Z'^2$, $\Cov(X,Z)=a'\sigma_Z'^2$, $\Var X=\lambda^2+a'^2\sigma_Z'^2+\sigma_X'^2$, $\Cov(W,X)=\lambda$ and $\Cov(W,Z)=0$. Substituting the definitions, $\Var Z$, $\Cov(X,Z)$, $\Var X$ and $\Cov(W,X)$ agree with $P$, so $\Sigma_P-\Sigma_Q$ has the single non-zero pair of entries $\Cov(W,Z)=a\lambda$. The $Y$-equations coincide, so $\KL(P\|Q)=\KL(P_{WXZ}\|Q_{WXZ})$ and $\|\Sigma_P-\Sigma_Q\|_F=\sqrt2\,a\lambda$. The targets are $\theta_M=aB$ and $\theta_C=0$.

\emph{Feasibility.} $Q$ must have $\sigma_Z'^2,\sigma_X'^2\in[1,\kappa]$ and $|a'|\le A$. We have $\sigma_Z'^2\ge s^*\ge1$ and $\sigma_X'^2\le\Var_P(X\mid Z)-\lambda^2\le m^*\le\kappa$. The remaining two constraints, $\sigma_Z'^2\le\kappa$ and $\sigma_X'^2\ge1$, hold at $\lambda=0$ with uniform slack over $a\in[0,\abar/2]$: for $\kappa\ge2$, $(m^*,s^*)=(2,\kappa/2)$ gives $\sigma_Z'^2=2a^2+\kappa/2\le5\kappa/8$ and $\Var(X\mid Z)=\kappa/(\kappa/2+2a^2)\ge8/5$; for $1<\kappa<2$, $(m^*,s^*)=(\kappa,1)$ gives $\sigma_Z'^2=1+a^2\kappa\le1+(\kappa-1)/4<\kappa$ and $\Var(X\mid Z)=\kappa/(1+a^2\kappa)\ge4\kappa/(\kappa+3)>1$. Denote the slacks by $\delta_Z(\kappa)=\kappa-\max_a\sigma_Z'^2|_{\lambda=0}>0$ and $\delta_X(\kappa)=\min_a\sigma_X'^2|_{\lambda=0}-1>0$. The dependence on $\lambda$ is explicit: with $x=\lambda^2$, $\sigma_Z'^2=s^*+a^2(m^*+x)$ increases by $a^2x\le\abar^2x/4$, and
\[
\frac{d}{dx}\,\sigma_X'^2=\frac{d}{dx}\Bigl\{\frac{(m^*+x)s^*}{s^*+a^2(m^*+x)}-x\Bigr\}=\frac{(s^*)^2}{[s^*+a^2(m^*+x)]^2}-1\in[-1,0],
\]
so $\sigma_X'^2$ decreases by at most $x$. Hence, for $\kappa>1$ (so that $\abar(\kappa)>0$), both constraints hold for all $a\in[0,\abar/2]$ and $\lambda\le w_0(\kappa):=\min\{1,\sqrt{\delta_X(\kappa)},\sqrt{4\delta_Z(\kappa)/\abar(\kappa)^2}\}$, which depends only on $\kappa$. For $\kappa=1$, $w_0$ is not defined and Pair 3 is not used. Finally $|a'|\le a(\lambda^2+m^*)/s^*\le a(\bar w^2+\kappa)$, so $a\le A/(\bar w^2+\kappa)$ guarantees $|a'|\le A$ and also $a\le A$.

For $\lambda\le w_0$ take $a=\min\{\abar/2,\ A/(\bar w^2+\kappa),\ c_1/(\lambda\sqrt n)\}$. Then $a\lambda\le c_1n^{-1/2}$, so $\|E\|\le\tfrac12$ for $n\ge n_0$ and $n\KL\le Cc_1^2$, and with $c_1$ small
\[
\mathrm{LB}_3=c\,B\min\{A,\,\abar,\,1/(\lambda\sqrt n)\}\qquad(\lambda\le w_0),
\]
where $c$ absorbs the factors $\tfrac12$ and $(\bar w^2+\kappa)^{-1}$.

\paragraph{Combination.} Write $z=1/(\lambda\sqrt n)$ and $\varrho_n=\sqrt{\abar^2+n^{-1/2}}$; then $R_n=n^{-1/2}+B\min\{z,A,\varrho_n\}$. Set $\mathrm{LB}_3=0$ when $\kappa=1$ or $\lambda>w_0(\kappa)$, and when $\kappa>1$ let $n\ge n_0\ge w_0(\kappa)^{-4}$. The three cases below are exhaustive and mutually exclusive.
\begin{itemize}
\item $\kappa>1$ and $\lambda\le w_0$: $\max(\mathrm{LB}_2,\mathrm{LB}_3)\ge c\,B\min\{A,\max(n^{-1/4},\abar),z\}\ge\tfrac c2B\min\{A,\varrho_n,z\}$, using $\max(n^{-1/4},\abar)\ge\varrho_n/\sqrt2$.
\item $\kappa>1$ and $\lambda>w_0$: $z<1/(w_0\sqrt n)\le n^{-1/4}$, so $\mathrm{LB}_2=cB\min\{A,z\}\ge cB\min\{A,\varrho_n,z\}$.
\item $\kappa=1$: $\abar=0$, $\varrho_n=n^{-1/4}$ and $\mathrm{LB}_2=cB\min\{A,\varrho_n,z\}$; Pair 3 and $w_0$ are not used.
\end{itemize}
In all cases $\max(\mathrm{LB}_1,\mathrm{LB}_2,\mathrm{LB}_3)\ge\tfrac c2\max\{n^{-1/2},B\min\{z,A,\varrho_n\}\}\ge\tfrac c4R_n$, using $\max(x,y)\ge(x+y)/2$. The constants and $n_0$ depend only on $(\bar w,\bar A,\kappa,\alpha)$, and the bound holds for every $A\in(0,\bar A]$, $B>0$ and $\lambda\in(0,\bar w]$. \qed

\subsection{Proof of \texorpdfstring{\cref{thm:rate-upper}}{the upper bound}}\label{app:rate-upper}

\paragraph{Coverage.} Under the true model $g^*$, each residual variance satisfies $m\hat\sigma^2_{j,g^*}/\sigma_j^2\sim\chi^2_m$ by \cref{lem:M}, with $\sigma_j^2\in[1,\kappa]$, so each of the $4$ variance checks of $g^*$ rejects it with probability at most $\eta$ ($\eta=\gamma_2/8$ is conservative). The conditional independence $t$-statistic of $g^*$ is exactly $t$-distributed and the interval of $g^*$ is exact by (M3). A union bound gives coverage at least $1-\gamma_1-\gamma_2-\alpha'=1-\alpha$.

\paragraph{The fallback of {\normalfont\DSVp}.} Let $A\le1/(\bar w\sqrt n)$. Under both graphs the regression of $Y$ on $(X,W,Z)$ has $X$-coefficient $\beta$, residual variance $\sigma_Y^2$ and Gaussian errors, so by (M3) of \cref{lem:M} the interval $\hat\theta_C\pm t_{n-4,1-\alpha/2}\widehat{\mathrm{se}}_C$ covers $\beta$ with probability exactly $1-\alpha$; the target is $\beta$ under $C$ and $\beta+ab$ with $|ab|\le AB$ under $M$, so the widened interval covers it with probability at least $1-\alpha$, for every $n\ge5$. The residual variance of $X$ given $(W,Z)$ is $\sigma_X^2\ge1$ under $C$ and $\sigma_X^2\sigma_Z^2/(\sigma_Z^2+a^2\sigma_X^2)\ge1/(1+\bar A^2\kappa)$ under $M$, so by (M2) the expected length is at most $Cn^{-1/2}+2AB$. Since $A\le1/(\bar w\sqrt n)\le1/(\lambda\sqrt n)$ and $A\le n^{-1/4}$ for $n\ge\bar w^{-4}$, $\min\{B/(\lambda\sqrt n),B\min(A,\sqrt{\abar^2+n^{-1/2}})\}=BA$, so the length is at most $CR_n$. It remains to treat \DSV\ for $A>1/(\bar w\sqrt n)$.

\paragraph{Step 0: $B/\sqrt n\lesssim R_n$ when $A>1/(\bar w\sqrt n)$.} Since $\lambda\le\bar w$, $A>1/(\bar w\sqrt n)$ and $\sqrt{\abar^2+n^{-1/2}}\ge n^{-1/4}\ge1/(\bar w\sqrt n)$ for $n\ge\bar w^{-4}$, $R_n\ge B/(\bar w\sqrt n)$.

\paragraph{Decomposition.} Let $g'$ be the wrong model. Then $\Leb(C_n)\le\Leb(W_C)+\Leb(W_M)+2\ind\{g'\text{ kept}\}|\hat\theta_{g'}-\theta_{g^*}|$. For the regressions of $Y$ on $\{X,W,Z\}$ and on $\{X,W\}$, the residual variance of $X$ given the other regressors is at least $1/(1+\bar A^2\kappa)$ and the residual variance of $Y$ is at most $\kappa+B^2\kappa(1+\bar A^2\kappa)$. By (M1) and (M2), $\E\Leb(W_g)$ and $\E|\hat\theta_g-\thadj_g(P)|$ are at most $C(1+B)/\sqrt n\le CR_n$, by Step 0. Hence
\[
\E\Leb(C_n)\le CR_n+2\Delta\,\Prob(g'\text{ kept}),\qquad\Delta=|\thadj_{g'}(P)-\theta_{g^*}(P)|.
\]

\paragraph{The bias $\Delta$.} If the truth is $M$, the regression of $Y$ on $(X,W,Z)$ has $X$-coefficient $\beta$ while $\theta_M=\beta+ab$, so $\Delta=|ab|$. If the truth is $C$, the regression of $Y$ on $(X,W)$ has $X$-coefficient $\beta+b\Cov(Z,X\mid W)/\Var(X\mid W)=\beta+ba\sigma_Z^2/(\sigma_X^2+a^2\sigma_Z^2)$, so $\Delta\le\kappa|ab|$. In both cases $\Delta\le\kappa B|a|$, and it remains to show
\begin{equation}\label{eq:star}
B|a|\,\Prob(g'\text{ kept})\le C\min\Bigl\{\frac{B}{\lambda\sqrt n},\ B\min\bigl(A,\sqrt{\abar^2+n^{-1/2}}\bigr)\Bigr\}+CB/\sqrt n .
\end{equation}
Trivially $B|a|\Prob\le BA$.

\paragraph{The conditional independence screen.} If the truth is $M$, the constraint of $C$ has population value $\rho(W,Z)=aw\sigma_W/\sqrt{\sigma_Z^2+a^2(\sigma_X^2+w^2\sigma_W^2)}$; if the truth is $C$, the constraint of $M$ has $|\rho(W,Z\mid X)|=|wa|\sigma_W\sigma_Z/\sqrt{(\sigma_X^2+a^2\sigma_Z^2)(w^2\sigma_W^2+\sigma_X^2)}$. In both cases $|\rho|\ge c_1|wa|$ with $c_1$ depending on $(\kappa,\bar A,\bar w)$. By \cref{lem:T}, with $x=c_1\lambda|a|\sqrt{(m+1)/2}$,
\[
B|a|\Prob(g'\text{ passes})\le\frac{B}{c_1\lambda\sqrt{(m+1)/2}}\sup_{x\ge0}x\Phibar(x-2\bar c)+B|a|\eps_n,\qquad\eps_n=e^{-(m+1)/16}+e^{-m/2}.
\]
The first term is at most $CB/(\lambda\sqrt n)$. The second is at most $BA\eps_n\le B\bar A\eps_n$, which for $n\ge n_0$ is below each of $CB/(\lambda\sqrt n)$ (as $\eps_n$ decays exponentially and $\lambda\le\bar w$), $CB\min(A,n^{-1/4})$ (it is at most $BA$, and at most $Bn^{-1/4}$ once $\bar A\eps_n\le n^{-1/4}$) and $BA$.

\paragraph{The variance screen.} Let $\tau_j^2(P)$ be the population residual variances of the wrong model's node-wise regressions and $D(P)=\max_j\mathrm{dist}(\tau_j^2,[1,\kappa])$, $d=\min(D,1/2)$.

\begin{lemma}\label{lem:V}
There are $c,C$ depending on $\kappa$ and $m_0$ depending also on $\eta$ such that for $m\ge m_0$, where $m=n-3$ is the smallest of the degrees of freedom of the node-wise regressions (the bound below is monotone in the degrees of freedom),
\[
\Prob(g'\text{ passes the variance screen})\le2\exp\Bigl\{-cm\bigl(d/\kappa-C\sqrt{\log(1/\eta)/m}\bigr)_+^2\Bigr\}.
\]
\end{lemma}
\begin{proof}
By \cref{lem:M}, $m\hat\sigma^2_{j,g'}/\tau_j^2\sim\chi^2_m$ whether or not $g'$ is correct. The acceptance probability $\pi(\tau^2)$ of a node is monotone in $\tau^2$ outside the acceptance interval: with $g(x)=xf(x)$ for the density $f$ of $\chi^2_m/m$, which is unimodal with mode $1$, $\pi'(\tau^2)=-[g(U/\tau^2)-g(L/\tau^2)]/\tau^2$ has a constant sign on each side. So it suffices to treat $D=d\le1/2$. If $\tau_j^2\ge\kappa+d$, acceptance requires $\chi^2_m/m\le q_{\mathrm{hi}}/(1+d/\kappa)\le1-(d/(2\kappa)-C\sqrt{\log(1/\eta)/m})$, using $q_{\mathrm{hi}}\le1+2\sqrt{\log(2/\eta)/m}+2\log(2/\eta)/m$ \citep{laurent2000adaptive}; the lower tail of \citet[Lemma~1]{laurent2000adaptive} bounds this probability. The case $\tau_j^2\le1-d$ is symmetric.
\end{proof}

\begin{lemma}\label{lem:Vprime}
Let $C_\kappa=\max\{3,(1+\kappa)/2\}$. For every $P$ in the class, $D(P)\ge\min\{1/2,\ (a^2-\abar(\kappa)^2-C_ww^2)_+/C_\kappa\}$, with $C_w=0$ if the truth is $M$ and $C_w=2\kappa$ if the truth is $C$.
\end{lemma}
\begin{proof}
\emph{Truth $M$.} Write $u=\sigma_Z^2+a^2\sigma_X^2$. Under the wrong model $C$, the fitted variance of $Z$ is $\Var Z=u+a^2w^2\sigma_W^2\ge u$ and that of $X$ given $(W,Z)$ is $\sigma_X^2\sigma_Z^2/u\le\kappa$; those of $W$ and $Y$ are correct. So $D\ge\max\{(u-\kappa)_+,(1-\sigma_X^2\sigma_Z^2/u)_+\}$. Suppose $D\le\delta\le1/2$. Then $a^2\le\min\{(\kappa+\delta-\sigma_Z^2)/\sigma_X^2,\ \sigma_Z^2(\sigma_X^2-1+\delta)/((1-\delta)\sigma_X^2)\}$. The two bounds are equal at $\sigma_Z^2=(\kappa+\delta)(1-\delta)/\sigma_X^2$, giving $a^2\le(\kappa+\delta)(\sigma_X^2-1+\delta)/\sigma_X^4$. For $\kappa\ge2$ the maximum over $\sigma_X^2$ is $(\kappa+\delta)/(4(1-\delta))\le\kappa/4+\delta(1+\kappa)/2$. For $\kappa<2$ and $\kappa\le2(1-\delta)$ the maximum is at $\sigma_X^2=\kappa$ and exceeds $1-1/\kappa$ by $[\delta(2\kappa-1)+\delta^2]/\kappa^2\le3\delta$; for $\kappa\in(2-2\delta,2)$ it is $(\kappa+\delta)/(4(1-\delta))$, which exceeds $1-1/\kappa\ge(1-2\delta)/(2-2\delta)$ by at most $5\delta/(4(1-\delta))\le3\delta$. In all cases $a^2\le\abar^2+C_\kappa\delta$; taking the contrapositive proves the claim.

\emph{Truth $C$.} Write $\eps=w^2\sigma_W^2\le\kappa w^2$. Under the wrong model $M$, the fitted variance of $X$ given $W$ is $\sigma_X^2+a^2\sigma_Z^2$ and that of $Z$ given $X$ is $\sigma_Z^2(\sigma_X^2+\eps)/(\sigma_X^2+\eps+a^2\sigma_Z^2)$. If $D\le\delta\le1/2$, then $a^2\le\min\{(\kappa+\delta-\sigma_X^2)/\sigma_Z^2,\ (\sigma_X^2+\eps)(\sigma_Z^2-1+\delta)/((1-\delta)\sigma_Z^2)\}$. Using $\min\{x,y+z\}\le\min\{x,y\}+z$ for $z\ge0$, the right side is at most the bound of the first case with the roles of $\sigma_X^2$ and $\sigma_Z^2$ exchanged, which is $\le\abar^2+C_\kappa\delta$, plus $\eps(\sigma_Z^2-1+\delta)/((1-\delta)\sigma_Z^2)\le2\eps\le2\kappa w^2$.
\end{proof}

\paragraph{Proof of \eqref{eq:star} for the plateau.} Let $\zeta=a^2-\abar^2$. If $\zeta\le K_0n^{-1/2}$, then $|a|\le\sqrt{\abar^2+K_0n^{-1/2}}$ and \eqref{eq:star} holds. Otherwise distinguish two cases.
(a) $C_ww^2\le \zeta/2$. By \cref{lem:Vprime}, $d\ge\min(1/2,\zeta/(2C_\kappa))$. If $\zeta/(2C_\kappa)\ge1/2$, \cref{lem:V} gives an exponentially small acceptance probability and $B|a|\Prob\le BA\eps'_n$, absorbed as before. Otherwise, for $K_0$ large, \cref{lem:V} gives $\Prob\le2e^{-c''ns^2}$, and with $x=\sqrt ns$, $B|a|\Prob\le B\sqrt{\abar^2+n^{-1/2}}\,\sqrt{1+x}\,2e^{-c''x^2}\le CB\sqrt{\abar^2+n^{-1/2}}$.
(b) $C_ww^2>\zeta/2$, which can only happen when the truth is $C$, for every size of $\zeta$. Then $|w|>\sqrt{\zeta/(2C_w)}$ and $|a|\ge\sqrt \zeta$, so $\sqrt n|wa|\ge x/\sqrt{2C_w}$, and the conditional independence bound gives $B|a|\Prob\le B\sqrt{\abar^2+n^{-1/2}}\sqrt{1+x}\,\Phibar(c'x-2\bar c)+BA\eps_n\le CB\sqrt{\abar^2+n^{-1/2}}$. The variance screen is not used in this case; it may accept with probability close to one, since all fitted variances of the wrong model can be feasible here.
Combining the three bounds proves \eqref{eq:star} and the theorem. \qed

\subsection{The local limit experiment and the dual bound}\label{app:toy-constant}

In this subsection the noise variances are $1$ and we consider the local parametrization $\beta=\beta_0+\tilde\beta/\sqrt n$, $a=c_g t/\sqrt n$ with $c_M=1$, $c_C=1+w^2$; the oracle is the efficient interval of \cref{sec:setup}.

\begin{proposition}[Local limit experiment]\label{prop:limit}
Let $v=1/(1+w^2)$. Along these local parameters the experiment converges to observing independent $(S_\beta,P,E)$, each with variance $v$, with means $(\tilde\beta,t,0)$ under $M$ and $(\tilde\beta,t,-wt)$ under $C$. The local target is $\psi=\tilde\beta+bt$ under $M$ and $\psi=\tilde\beta$ under $C$, and the oracle lengths are $2z_{1-\alpha/2}s$ under $C$ and $2z_{1-\alpha/2}s\sqrt{1+b^2}$ under $M$, $s=\sqrt v$.
\end{proposition}

Here $S_\beta$ is the efficient score for $\beta$, $P$ the score for the regression of $Z$ on $X$, and $E=D-wP$ with $D=n^{-1/2}\sum W_iZ_i$. The log-likelihood ratios are exactly quadratic in the local parameters (\cref{app:limit}). This is the canonical experiment of \cref{prop:canonical} with parameters $(1,b,0,\sqrt v)$. In the limit experiment the minimax price is
\[
\Pi^*(w,b)=\inf_{\varphi\ \mathrm{honest}}\ \sup_{g,t,\tilde\beta}\ \frac{\E\,\Leb(\varphi)}{\orl_g},
\]
over all measurable, possibly randomized, confidence sets $\varphi$.

\begin{theorem}[Dual lower bound]\label{thm:lfp}
For a finite set $K$ of points $(g_j,t_j)$, weights $\pi\in\Delta(K)$ and multipliers $\mu\in[0,\infty)^K$, let $f_j$ be the density of $(P,E)$ at $(g_j,t_j)$, $\psi_j$ the local target at $\tilde\beta=0$, and
\[
R(\pi,\mu)=(1-\alpha)\sum_j\mu_j+\int\min\Bigl\{0,\ \sum_j\frac{\pi_j f_j(p,e)}{\orl_{g_j}}-\sum_j\mu_j f_j(p,e)\,\phi_v(u-\psi_j)\Bigr\}\,dp\,de\,du .
\]
(i) In the limit experiment, every honest confidence set has price at least $R(\pi,\mu)$. (ii) In finite samples, under conditions (S1)--(S4) of \cref{app:lfp}, every sequence of uniformly honest confidence sets satisfies $\liminf_n\sup\E\Leb(C_n)/L^{\mathrm{or}}_n\ge R(\pi,\mu)$.
\end{theorem}

Part (i) averages over $\tilde\beta\sim U[-M,M]$ and lets $M\to\infty$; part (ii) uses the exact quadratic expansion of the likelihood ratio. At $w=b=1$ and $\alpha=0.05$, a witness $(\hat\pi,\hat\mu)$ from the dual of a discretized linear program, evaluated by box bounds with an explicit allowance for floating-point rounding (\cref{app:certify-lfp}), gives $\Pi^*(1,1)\ge R(\hat\pi,\hat\mu)\ge1.27438$; this is a computer-assisted bound under a floating-point error model, not a formally verified computation. The Pratt bound at the crossing point gives $\max\{1,c_P\sqrt{1+b^2}\}=1.2019$, and the discretized linear program about $1.284$ (an approximation). The interval procedure of \cref{app:bias2i} has price about $1.422$, and a union-fallback set procedure about $1.421$, both numerical (\cref{tab:price-curve}).

\subsection{Proof of \texorpdfstring{\cref{prop:limit}}{the limit experiment}}\label{app:limit}

At the base point ($a=0$, unit variances, $w\ne0$), let $\eps_{Xi}=X_i-wW_i$, $\eps_{Yi}=Y_i-\beta_0X_i-bZ_i$, $S_{\beta n}=n^{-1/2}\sum X_i\eps_{Yi}$, $S_{Mn}=n^{-1/2}\sum X_iZ_i$ and $S_{Cn}=n^{-1/2}\sum Z_i\eps_{Xi}$. With $\beta=\beta_0+u/\sqrt n$ and $a=\eta/\sqrt n$, expanding the Gaussian log-likelihoods gives, without remainder,
\begin{align*}
\log\Lambda_{M,n}&=uS_{\beta n}+\eta S_{Mn}-\tfrac{u^2+\eta^2}{2n}\textstyle\sum X_i^2,\\
\log\Lambda_{C,n}&=uS_{\beta n}+\eta S_{Cn}-\tfrac{u^2}{2n}\textstyle\sum X_i^2-\tfrac{\eta^2}{2n}\sum Z_i^2 .
\end{align*}
Under the base point $(S_{\beta n},S_{Mn},S_{Cn})\Rightarrow N(0,I)$ with $I=\begin{psmallmatrix}q&0&0\\0&q&1\\0&1&1\end{psmallmatrix}$, $q=1+w^2$, and the remainders are uniformly $o_P(1)$ over bounded local parameters because $n^{-1}\sum X_i^2\to q$ and $n^{-1}\sum Z_i^2\to1$. In the coordinates $S_\beta/q$, $P=S_M/q$, $E=(S_M-qS_C)/(qw)$ the limit is the experiment of \cref{prop:limit} (the first coordinate is the statistic called $S_\beta$ there); also $E=D-wP$ with $D=n^{-1/2}\sum Z_iW_i$, and $P-wE=S_C$ exactly.

\subsection{Proof of \texorpdfstring{\cref{thm:lfp}}{the dual bound}}\label{app:lfp}

Write $A(p,e)=\sum_j\pi_jf_j/\orl_{g_j}$, $g_j=\mu_jf_j(p,e)\phi_v(u-\psi_j)$, $G=\sum_jg_j$ and $J=\int\min\{0,A-G\}$; $J$ is finite since $-G\le\min\{0,A-G\}\le0$ and $\int G=\sum\mu_j$.

\emph{(i).} Let $\varphi$ be honest with price $r<\infty$, and let $q_M$ be the density of $U[-M,M]$. Averaging the price and coverage constraints against $q_M(\tilde\beta)\pi_j$ and $q_M(\tilde\beta)\mu_j$, using the Pratt identity $\E\Leb=\int\Prob(\psi_0\in\varphi)\,d\psi_0$ and Tonelli separately for the two non-negative terms, gives $r\ge(1-\alpha)\sum\mu_j+J_M$, where after the change of variables $u=\psi_0-S_\beta$,
\[
J_M=\int\min\bigl\{0,\ A(p,e)a_M(S_\beta)-\textstyle\sum_jg_j(p,e,u)\,q_M(S_\beta+u-\psi_j)\bigr\},\qquad a_M=q_M*\phi_v .
\]
Since $\int a_M=1$, $J=\int\min\{0,Aa_M-Ga_M\}$, and the $1$-Lipschitz property of $x\mapsto\min\{0,x\}$ gives $|J_M-J|\le\sum_j\int g_j\|q_M(\cdot+u-\psi_j)-a_M\|_1$. Writing $a_M(S_\beta)=\E q_M(S_\beta-Z')$, $Z'\sim N(0,v)$, the $L^1$ norm is at most $\E|u-\psi_j+Z'|/M$, and integrating gives $|J_M-J|\le(2s/(M\sqrt\pi))\sum\mu_j$. Let $M\to\infty$.

\emph{(ii): conditions.} Fix a base point $\theta_0$ with $a=0$, $w\ne0$, fixed $b,\beta_0$ and unit noise variances, and the local points $P_{n,(\tilde\beta,j)}$ of branch $g_j$ with $\beta=\beta_0+\tilde\beta/\sqrt n$ and $a=c_{g_j}t_j/\sqrt n$.
\begin{itemize}
\item[(S1)] The local target $\sqrt n(\psi-\beta_0)$ equals $\tilde\beta+\psi_j$ for every $n$ ($b$ fixed along the path).
\item[(S2)] For every $M<\infty$ there is $N(M)$ such that $P_{n,(\tilde\beta,j)}$ belongs to the class for all $n\ge N(M)$, $|\tilde\beta|\le M$ and $j\in K$. (A class with a fixed positive margin on $a$ does not satisfy this.)
\item[(S3)] The procedures are uniformly honest: $\liminf_n\inf_{P}\Prob_P(\psi(P)\in C_n)\ge1-\alpha$ over the class.
\item[(S4)] The oracle $L^{\mathrm{or}}_n=2z_{1-\alpha/2}\sigma_{\mathrm{eff}}/\sqrt n$ is the deterministic efficient length of the labelled branch, so that $\sqrt nL^{\mathrm{or}}_n$ converges uniformly on $[-M,M]\times K$ to $\orl_{g_j}$ and is bounded away from zero.
\end{itemize}

\emph{(ii): proof.} Let $r_n=\sup\E_P\Leb(C_n)/L^{\mathrm{or}}_n(g,P)$ over the class; if $r_n=\infty$ there is nothing to prove. Write $\vartheta$ for the local $\beta$ coordinate, $x$ for the local target coordinate, $D_n=\sqrt n(C_n-\beta_0)$, $L_{n,j}(\vartheta)=dP_{n,(\vartheta,j)}/dP_{n,0}$ for the likelihood ratio against the base point, and $\ell_{n,j}(\vartheta)=\sqrt n\,L^{\mathrm{or}}_n$ at $(\vartheta,j)$. Fix $M$ and $T>M+\max_j|\psi_j|$, and let $q_M(\vartheta)=\ind\{|\vartheta|\le M\}/(2M)$. For a randomized set, replace indicators $\ind\{x\in D_n\}$ below by membership probabilities in $[0,1]$; nothing changes.

\emph{Step 1 (averaging).} By (S2), for $n\ge N(M)$ every $(\vartheta,j)$ with $|\vartheta|\le M$ is in the class, so $r_n\ge\E_{n,(\vartheta,j)}\Leb(D_n)/\ell_{n,j}(\vartheta)$. By (S1) the local target is $\vartheta+\psi_j$, and by (S3) $\Prob_{n,(\vartheta,j)}(\vartheta+\psi_j\in D_n)\ge1-\alpha-\eps_n$ with $\eps_n=((1-\alpha)-\inf_{|\vartheta|\le M,j}\mathrm{cov}_n)_+\to0$. Averaging the first inequality against $\pi_jq_M(\vartheta)\,dh$ (total mass $1$) and subtracting the non-negative coverage slack weighted by $\mu_jq_M(\vartheta)$,
\[
r_n\ge(1-\alpha-\eps_n)\sum_j\mu_j+\sum_j\pi_j\!\int\! q_M(\vartheta)\frac{\E_{n,(\vartheta,j)}\Leb(D_n)}{\ell_{n,j}(\vartheta)}dh-\sum_j\mu_j\!\int\! q_M(\vartheta)\Prob_{n,(\vartheta,j)}(\vartheta+\psi_j\in D_n)\,dh .
\]

\emph{Step 2 (change of measure and truncation).} Change measure to $P_{n,0}$, use Tonelli separately for the two non-negative terms, bound $\Leb(D_n)\ge\int_{-T}^T\ind\{x\in D_n\}dx$, and substitute $x=\vartheta+\psi_j$ in the coverage term; since $q_M(x-\psi_j)=0$ for $|x|>T$, the coverage term is also an integral over $[-T,T]$. With
\[
A_n=\sum_j\pi_j\int q_M(\vartheta)\frac{L_{n,j}(\vartheta)}{\ell_{n,j}(\vartheta)}dh,\qquad B_n(x)=\sum_j\mu_jq_M(x-\psi_j)L_{n,j}(x-\psi_j),
\]
the last two terms of Step 1 equal $\E_{n,0}\int_{-T}^T\ind\{x\in D_n\}(A_n-B_n(x))\,dx\ge\E_{n,0}\Xi_{n}$, where
\[
\Xi_{n}=\int_{-T}^T\min\{0,A_n-B_n(x)\}\,dx,\qquad\text{so}\qquad r_n\ge(1-\alpha-\eps_n)\textstyle\sum_j\mu_j+\E_{n,0}\Xi_n .
\]

\emph{Step 3 (domination).} $0\le-\Xi_n\le\int B_n=\Upsilon_n:=\sum_j\mu_j\int q_M(\vartheta)L_{n,j}(\vartheta)\,dh$, and $\E_{n,0}\Upsilon_n=\sum_j\mu_j$ because $\E_{n,0}L_{n,j}(\vartheta)=1$ (the Gaussian laws are mutually absolutely continuous).

\emph{Step 4 (convergence).} By the exact expansion of \cref{app:limit}, the processes $\vartheta\mapsto L_{n,j}(\vartheta)$, $j\in K$, converge in distribution in $C([-M,M])^K$ under $P_{n,0}$ to the likelihood ratio processes $L_j(\vartheta)$ of the limit experiment, and by (S4) $\ell_{n,j}\to\orl_{g_j}$ uniformly with a positive lower bound. The map from $(L_j)_j$ to $(\Xi,\Upsilon)$ is continuous on $C([-M,M])^K$: it integrates continuous functions over bounded sets (the jumps of $q_M$ at $\pm M$ do not affect this), and $\min\{0,\cdot\}$ is $1$-Lipschitz. Hence $(\Xi_n,\Upsilon_n)\Rightarrow(\Xi,\Upsilon)$, and $\E\Upsilon=\sum_j\mu_j$ since $\E L_j(\vartheta)=1$. If $Y_n\ge0$, $Y_n\Rightarrow Y$ and $\E Y_n\to\E Y<\infty$, then $(Y_n)$ is uniformly integrable (because $\E(Y_n-K)_+=\E Y_n-\E(Y_n\wedge K)\to\E(Y-K)_+$ and $y\ind\{y>2K\}\le2(y-K)_+$). So $(\Upsilon_n)$, and hence $(\Xi_n)$, is uniformly integrable, and $\E_{n,0}\Xi_n\to\E\Xi$.

\emph{Step 5 (identification).} In the limit experiment the observation is $(S_\beta,P,E)$ with density $f_{\vartheta,j}(b,p,e)=\phi_v(b-\vartheta)f_j(p,e)$ and $L_j(\vartheta)=f_{\vartheta,j}/f_0$. Since $f_0>0$ and $\min\{0,y\}f_0=\min\{0,yf_0\}$, and $\int q_M(\vartheta)\phi_v(b-\vartheta)dh=a_M(b)$,
\[
\E\Xi=\int_{-T}^T\!\!\int\min\Bigl\{0,\ A(p,e)a_M(b)-\sum_j\mu_jq_M(x-\psi_j)\phi_v(b-x+\psi_j)f_j(p,e)\Bigr\}\,db\,dp\,de\,dx .
\]
With $u=x-b$, $\phi_v(b-x+\psi_j)=\phi_v(u-\psi_j)$ and $q_M(x-\psi_j)=q_M(b+u-\psi_j)$, so $\E\Xi=J_{M,T}$, where $J_{M,T}$ denotes the integral defining $J_M$ in part (i) restricted to $|x|\le T$. Since $T>M+\max_j|\psi_j|$, for $|x|>T$ and every $j$ we have $|x-\psi_j|>M$, hence $q_M(x-\psi_j)=0$, and the integrand reduces to $\min\{0,A(p,e)a_M(b)\}=0$ because $A\ge0$. Therefore the integrand of $J_M$ vanishes outside $|x|\le T$ and $J_{M,T}=J_M$, so $\E\Xi=J_M$. Letting $n\to\infty$ gives $\liminf_nr_n\ge(1-\alpha)\sum_j\mu_j+J_M\ge R(\pi,\mu)-(2s/(M\sqrt\pi))\sum_j\mu_j$ by part (i), and letting $M\to\infty$ proves (ii). \qed

\Cref{tab:price-curve} gives the dual lower bound $R(\hat\pi,\hat\mu)$ and the other quantities of \cref{sec:toy-constant} at fourteen values of $(w,b)$.

\begin{table}[t]
\centering
\caption{Price constant in the limit experiment, $\alpha=0.05$. ``$K$'' is the cross-model constant $\max\{b/w,\ b\sqrt{1+w^2}/(w\sqrt{1+b^2})\}$ of \cref{sec:gen-constant}. ``Crossing'' is the Pratt bound $\max\{1,c_P\sqrt{1+b^2}\}$, $c_P\approx0.8498863$, rounded down. ``Dual'' is the lower-bound functional $R(\hat\pi,\hat\mu)$ of \cref{thm:lfp}, rounded down; at $(w,b)=(1,1)$ the value is a computer-assisted bound under the floating-point error model of \cref{app:certify-lfp}, at the other points it is a midpoint evaluation without error bound. ``LP'' is the discretized linear program, an approximation and not a bound. ``Set'' and ``interval'' are numerical values (closed-form tails, numerical middle segment, rounded up after adding a heuristic $10^{-3}$ buffer) for a union-fallback set procedure and for the interval procedure; they are not bounds. Scripts: \texttt{d8\_lfp\_bound.py}, \texttt{d9\_rho\_curve.py}, \texttt{d11\_bias2\_grid.py}, \texttt{d17\_bias2i\_tail.py}.}
\label{tab:price-curve}
\begin{tabular}{rrcccccc}
\toprule
$w$ & $b$ & $K$ & crossing & dual (lower) & LP & set & interval\\
\midrule
0.25 & 1 & 4.00 & 1.201 & 1.826 & 1.83 & -- & --\\
0.5 & 1 & 2.00 & 1.201 & 1.474 & 1.48 & 1.691 & --\\
0.75 & 1 & 1.33 & 1.201 & 1.324 & 1.33 & -- & --\\
1 & 1 & 1.00 & 1.201 & 1.274$^\dagger$ & 1.28 & 1.421 & 1.423\\
1.5 & 1 & 0.85 & 1.201 & 1.249 & 1.27 & -- & --\\
2 & 1 & 0.79 & 1.201 & 1.242 & 1.27 & 1.353 & --\\
1 & 0.25 & 0.34 & 1.000 & 1.018 & 1.02 & -- & --\\
1 & 0.5 & 0.63 & 1.000 & 1.068 & 1.08 & 1.163 & --\\
1 & 1.5 & 1.50 & 1.532 & 1.592 & 1.61 & -- & --\\
1 & 2 & 2.00 & 1.900 & 1.954 & 1.97 & 2.131 & --\\
1 & 3 & 3.00 & 2.687 & 2.732 & 2.76 & -- & --\\
0.5 & 0.5 & 1.00 & 1.000 & 1.179 & 1.18 & -- & --\\
0.5 & 2 & 4.00 & 1.900 & 2.191 & 2.20 & -- & --\\
2 & 2 & 1.00 & 1.900 & 1.909 & 1.96 & -- & --\\
\bottomrule
\end{tabular}

\smallskip\footnotesize $^\dagger$Computer-assisted bound under the floating-point error model of \cref{app:certify-lfp}: $R(\hat\pi,\hat\mu)\ge1.27438$. Without the buffer, the interval value at $(1,1)$ is $1.4218$.
\end{table}

\subsection{The interval procedure in the limit experiment}\label{app:bias2i}

This procedure is defined and analysed in the limit experiment of \cref{prop:limit} only; we do not claim finite-sample honesty for a plug-in version.

Let $X=E/s$, $r=k_C/k_M$, $K_C=\{|D|\le k_C\}$, $K_M=\{|X|\le k_M\}$, $W_C=[-cs,cs]$, $W_M=bP\pm c_M\mathrm{sd}_M$, $\mathrm{ext}=\varkappa|b|s$, $\omega_M=\Phi((k_M-|X|)/\tau_M)$, and $I^+=[-cs-\omega_M(-bP+\mathrm{ext})_+,\ cs+\omega_M(bP+\mathrm{ext})_+]$ (in the coordinate $u=\psi-S_\beta$). The procedure reports $I^+$ on $K_C$, $W_M$ on $K_C^c\cap K_M$, and on $K_C^c\cap K_M^c$ it reports $W_C$ if $|D|\le r|X|$ and $W_M$ otherwise. By equivariance and the symmetry $(P,E,u)\mapsto(-P,-E,-u)$ it suffices to consider $t\ge0$. Let $q(y;\mu)=\Phi(y-\mu)-\Phi(-y-\mu)$, $q_M(t)=q(k_M;wt/s)$, $q_C(t)=q(k_C;wt)$, $\Phi_c=2\Phi(c)-1$, and $\Omega(t)=\Phi((k_M-wt/s)/\sqrt{1+\tau_M^2})$. For any finite grid $\mathcal X\subset[0,\infty)$:
{\small\begin{align*}
\mathrm{cov}_C(t)&\ge\Phi_c\bigl(1-q_M(t)-\min_{x\in\mathcal X}[2\Phibar(x)+q(x/r;wt/s)]\bigr),\\
\mathrm{cov}_M(t)&\ge\Phi_{c_M}-q_C(t)-\min_{x\in\mathcal X}[q(x;wt)+2\Phibar(x/r)],\\
\E\Leb_C&\le2cs+\Omega(t)\bigl(|b|(t+s\sqrt{2/\pi})+2\mathrm{ext}\bigr)+2\Phi(-k_C)(2c_M\mathrm{sd}_M-2cs)_+,\\
\E\Leb_M&\le2c_M\mathrm{sd}_M+q_C(t)(2cs+2\mathrm{ext}-2c_M\mathrm{sd}_M)_++|b|\sqrt{t^2+v}\sqrt{q_C(t)}\\&\quad+2\Phi(-k_M)(2cs-2c_M\mathrm{sd}_M)_+ .
\end{align*}}
The coverage bounds hold because $W_C$ can be lost only on $K_M\cup\{|D|>r|X|\}$ and $W_M$ only on $K_C\cup\{|D|\le r|X|\}$, and $\{|D|>r|X|\}\subseteq\{|D|>x\}\cup\{|X|<x/r\}$. They are non-decreasing in $t$. The length bounds are non-increasing for $t\ge T_0$ provided $y(T_0)=(k_M-wT_0/s)/\sqrt{1+\tau_M^2}<0$, $wT_0>k_C$ and $T_0\ge\sqrt v$, and the sufficient derivative conditions $a(|b|T_0+c_0)|y(T_0)|>|b|$ and $w\,a'(1-e^{-2k_C(a'+k_C)})>2T_0/(T_0^2+v)$ hold at $T_0$, with $a=w/(s\sqrt{1+\tau_M^2})$, $c_0=|b|s\sqrt{2/\pi}+2\mathrm{ext}$ and $a'=wT_0-k_C$; the second uses $|q_C'|/q_C\ge wa'(1-e^{-2k_C(a'+k_C)})$.

At $w=b=1$ with $(c,\varkappa,c_M,k_C,k_M,\tau_M)=(1.990,0.535,2.554,2.298,2.258,0.518)$, these conditions hold at $T_0=4.75$ with margins $30.19$ and $5.95$; for $t\ge T_0$ the closed-form bounds give coverage at least $0.9524$ and price at most $1.4064$. On $[0,T_0]$, numerical integration on a grid of step $0.03s$ and $t$-step $0.05$, refined to step $0.015s$ and $t$-step $0.005$ near the extremes, gives coverage at least $0.9504$ and price at most $1.4218$; this part has no certified error bound (script \texttt{d17\_bias2i\_tail.py}).

\subsection{Proof of \texorpdfstring{\cref{lem:oracle-toy}}{the oracle lemma for the four-variable model}}\label{app:oracle-toy}

Under $C$ the target is $\beta$, the structural parameters of the four equations vary independently, and $\beta$ appears only in the equation of $Y$, whose regressors are $(X,Z)$; so $s_{\mathrm{eff}}^2=\sigma_Y^2/\Var(X\mid Z)=\sigma_Y^2/(w^2\sigma_W^2+\sigma_X^2)$. The parent-adjustment set is $\{X,W,Z\}$, and $s_{\mathrm{pa}}^2=\sigma_Y^2/\Var(X\mid W,Z)=\sigma_Y^2/\sigma_X^2$. Under $M$ the target is $\beta+ab$, with $(\beta,b)$ in the equation of $Y$ (regressors $(X,Z)$) and $a$ in the equation of $Z$ (regressor $X$). The efficient variance is $\sigma_Y^2(1,a)\Sigma_{XZ}^{-1}(1,a)^\top+b^2\sigma_Z^2/\Var X$, and $\Sigma_{XZ}e_1=\Var X\cdot(1,a)^\top$ gives $(1,a)\Sigma_{XZ}^{-1}(1,a)^\top=1/\Var X$; so $s_{\mathrm{eff}}^2=(\sigma_Y^2+b^2\sigma_Z^2)/\Var X$ with $\Var X=w^2\sigma_W^2+\sigma_X^2$. The parent-adjustment set is $\{X,W\}$, and $s_{\mathrm{pa}}^2=\Var(Y\mid X,W)/\Var(X\mid W)=(b^2\sigma_Z^2+\sigma_Y^2)/\sigma_X^2$. In both cases the ratio is $\Var X/\sigma_X^2=1+w^2\sigma_W^2/\sigma_X^2$. \qed

\subsection{Select-then-estimate}\label{app:naive}

We state the result of \cref{sec:toy-naive} precisely. In the local limit experiment of \cref{prop:limit}, let $D=E+wP$ and $V=P-wE$; then $D$ and $V$ are independent with unit variance, the constraint $t$-statistics of $C$ and $M$ converge to $D$ and $E/s$, and the two parent-adjustment estimators converge to $\tilde\beta+U$ and $\tilde\beta+U+bV$ with $U\sim N(0,1)$ independent of $(D,V)$ (paragraph \emph{Limits} below). We consider three selection rules: \textup{\textsc{cmp}} selects $C$ if $|D|\le|E|/s$ (the larger $p$-value wins); \textup{\textsc{pre}}$_C(\alpha_t)$ selects $M$ only if the test of $C$'s constraint rejects at level $\alpha_t$; \textup{\textsc{pre}}$_M(\alpha_t)$ selects $C$ only if the test of $M$'s constraint rejects at level $\alpha_t$. Let $\mathrm{cov}_g(t)$ be the limiting coverage of the selected $95\%$ Wald interval at $(g,t)$; it does not depend on $\tilde\beta$.

\begin{proposition}[Select-then-estimate]\label{prop:naive}
Fix $w\ne0$ and $b$ with $|b|\le B$.
\begin{enumerate}
\item[(i)] $\liminf_n\inf_{(g,P)\in\cM(w_0,A,B,\kappa)}\Prob_P(\theta_g(P)\in C_n^{\mathrm{sel}})\le\inf_{g,t}\mathrm{cov}_g(t)$ whenever $w_0\le|w|$, so the procedure is honest over the benign class of \cref{cor:benign-toy} only if $\inf_{g,t}\mathrm{cov}_g(t)\ge1-\alpha$.
\item[(ii)] For $b\ne0$ and every $\alpha_t\in(0,1)$, the pre-test rules undercover as $|t|\to\infty$: under $g=C$, \textup{\textsc{pre}}$_C(\alpha_t)$ has $\lim_{|t|\to\infty}\mathrm{cov}_C(t)=(1-\alpha)(1-\alpha_t)$, and under $g=M$, \textup{\textsc{pre}}$_M(\alpha_t)$ has $\lim_{|t|\to\infty}\mathrm{cov}_M(t)=1-\alpha-\pi_M(\alpha_t,w,b)$ with $\pi_M>0$ the probability that the test of $M$'s constraint rejects while the mediator interval would have covered. Both limits are below $1-\alpha$. Under $C$, the test of $C$'s constraint has exact level $\alpha_t$ for every $n$ and every $a$.
\end{enumerate}
\end{proposition}

\Cref{tab:naive-full} gives the limiting worst-case coverage for five rules. The proofs of \cref{prop:naive,prop:naive-more} follow the table.

\begin{proposition}[Further properties of selection]\label{prop:naive-more}
In the setting of \cref{prop:naive}:
\begin{enumerate}
\item[(i)] for every rule, $\limsup_{|b|\to\infty}\inf_t\mathrm{cov}_M(t)\le\Prob_0(\text{$M$ selected})$, which is $1/2$ for \textup{\textsc{cmp}} and $\alpha_t$ for \textup{\textsc{pre}}$_C(\alpha_t)$;
\item[(ii)] for \textup{\textsc{cmp}} and fixed $(g,t)$, $\mathrm{cov}_g(t)=1-\alpha+z\phi(z)\xi_g(t)b^2+O(b^4)$ with an explicit $\xi_g(t)$; numerically $\min_t\xi_C(t)<0$ for $w\in\{0.5,1,2\}$, so \textup{\textsc{cmp}} undercovers for all sufficiently small $|b|\ne0$.
\end{enumerate}
\end{proposition}

\begin{table}[t]
\centering
\caption{Limiting worst-case coverage $\inf_{g,t}\mathrm{cov}_g(t)$ of the selected $95\%$ interval at $w=1$; the infimum is taken over a grid of $t\in[0,40]$ together with the limit $|t|\to\infty$ of \cref{prop:naive}(ii). One-dimensional integrals, cross-checked against adaptive quadrature (largest difference below $10^{-13}$) and Monte Carlo; since the entries are upper bounds on the worst-case coverage, they are rounded up to four decimals (scripts \texttt{d23b\_naive\_limit\_exact.py}, \texttt{d23c\_table3\_roundup.py}).}
\label{tab:naive-full}
\begin{tabular}{lcccccc}
\toprule
rule & $b=0.25$ & $0.5$ & $1$ & $2$ & $5$ & $20$\\
\midrule
\textup{\textsc{cmp}} & 0.9452 & 0.9330 & 0.8872 & 0.7817 & 0.6094 & 0.5052\\
\textup{\textsc{pre}}$_C(0.05)$ & 0.9025 & 0.8909 & 0.7244 & 0.4076 & 0.1375 & 0.0555\\
\textup{\textsc{pre}}$_C(0.001)$ & 0.9091 & 0.7822 & 0.4290 & 0.1083 & 0.0116 & 0.0016\\
\textup{\textsc{pre}}$_M(0.05)$ & 0.9033 & 0.8953 & 0.8073 & 0.7185 & 0.6724 & 0.6624\\
\textup{\textsc{pre}}$_M(0.001)$ & 0.8887 & 0.7365 & 0.4742 & 0.3164 & 0.2616 & 0.2515\\
\bottomrule
\end{tabular}

\smallskip\footnotesize At $w=b=1$, the best pre-test level over $\alpha_t\in\{10^{-4},\dots,0.95\}$ gives $0.7799$ for \textup{\textsc{pre}}$_C$ and $0.8595$ for \textup{\textsc{pre}}$_M$ (both at $\alpha_t=0.1$).
\end{table}

\paragraph{Limits.} In the notation of \cref{app:limit}, $D=E+wP=n^{-1/2}\sum W_iZ_i$ and $V=P-wE=S_C$ hold exactly. In the limit $P$ and $E$ are independent with variance $v$, so $D$ and $V$ are independent with unit variance; their means are $\E D=\E E+wt$ and $\E V=c_gt$. The regression of $Y$ on $(X,W,Z)$ has $X$-coefficient $\beta$ under both graphs, and the residual of $X$ on $(W,Z)$ is $\eps_X$ at the base point, so $\sqrt n(\hat\theta_C-\beta_0)\to\tilde\beta+U$ with $U$ the limit of $n^{-1/2}\sum\eps_{Yi}\eps_{Xi}$. The $X$-coefficient of the regression of $Y$ on $(X,W)$ is $\beta+b\gamma$, where $\gamma$ is the $X$-coefficient of $Z$ on $(X,W)$; $\gamma=a$ under $M$ and $\gamma=a/(1+a^2)$ under $C$ at unit variances, so $\sqrt n\gamma=c_gt+O(n^{-1})$, and $\sqrt n(\hat\gamma-\gamma)\to V-c_gt$. Hence $\sqrt n(\hat\theta_M-\beta_0)\to\tilde\beta+U+bV$. The same $U$ appears in both limits, because the residual of $X$ is $\eps_X$ in both regressions and the residual of $Y$ is $\eps_Y$ up to terms that vanish; and $U$ is independent of $(D,V)$, because $\E[\eps_Y\eps_XWZ]=\E[\eps_Y\eps_XZ\eps_X]=0$. The standardized standard errors converge to $1$ and $\sqrt{1+b^2}$ under both graphs. The constraint $t$-statistics converge to $D$ and to $E/s=s(D-wV)$, and the means of all limits under the local alternatives follow from the exact quadratic expansion of \cref{app:limit} and Le Cam's third lemma.

\paragraph{(i).} The limiting coverage events have boundaries of probability zero, so for each fixed $(g,t,\tilde\beta)$ the coverage converges to $\mathrm{cov}_g(t)$ by the continuous mapping theorem. For $n\ge(c_gt/A)^2$ the local point lies in $\cM(w_0,A,B,\kappa)$, so $\liminf_n\inf\le\mathrm{cov}_g(t)$ for every $(g,t)$, and taking the infimum over $(g,t)$ (through a sequence $t_k$ if it is approached only as $|t|\to\infty$) proves (i).

\paragraph{\Cref{prop:naive-more}(i).} Take $g=M$ and $t=K/|b|$. If $C$ is selected the interval $\hat\theta_C\pm z$ covers $\psi_M$ with probability $\Phi(z+K)-\Phi(K-z)$, independently of $(D,V)$. If $M$ is selected, coverage is $\{|U/b+V_0|\le z\sqrt{1+b^2}/|b|\}$ with $V_0=V-t\sim N(0,1)$, whose indicator converges almost surely to $\ind\{|V_0|\le z\}$ as $|b|\to\infty$. As $t\to0$ the law of $(D,V_0)$ converges to its law at $t=0$; bounded convergence gives $\limsup_{|b|\to\infty}\mathrm{cov}_M(K/|b|)\le\Phi(z+K)-\Phi(K-z)+\Prob_0(M\text{ selected},|V|\le z)$, and $K\to\infty$ gives (ii). Under \textup{\textsc{cmp}}, $(D,E/s)$ at $t=0$ is bivariate normal with equal variances and correlation $s$, so it is exchangeable and $\Prob_0(|D|>|E|/s)=1/2$.

\paragraph{(ii).} Under $C$, $Z=\eps_Z$ is independent of $W$ for every $a$, so the $t$-test of $W\indep Z$ is exactly $t$-distributed for every $n$; in the limit $D\sim N(0,1)$ for every $t$, and \textup{\textsc{pre}}$_C(\alpha_t)$ selects $M$ with probability $\alpha_t$, an event that depends only on $D$. If $C$ is selected, coverage $\{|U|\le z\}$ has probability $1-\alpha$, independently of $D$. If $M$ is selected, the interval is centred at $\tilde\beta+U+bV$ with $V\sim N((1+w^2)t,1)$ and has half-width $z\sqrt{1+b^2}$, so its coverage of $\tilde\beta$ tends to $0$ as $|t|\to\infty$ for fixed $b\ne0$. Hence $\mathrm{cov}_C(t)\to(1-\alpha)(1-\alpha_t)$. Under $M$, $E\sim N(0,v)$ for every $t$ and \textup{\textsc{pre}}$_M(\alpha_t)$ selects $C$ when $|E/s|>c$, $c=z_{1-\alpha_t/2}$; the interval $\hat\theta_C\pm z$ then misses $\psi_M=\tilde\beta+bt$ with probability tending to one. If $M$ is selected, coverage is $\{|U+bV_0|\le z\sqrt{1+b^2}\}$ with $V_0=V-t\sim N(0,1)$, whose joint law with $E/s$ does not depend on $t$; but $V_0$ and $E/s$ are correlated, with correlation $-ws$, so this event is not independent of the selection. Hence
\begin{align*}
\lim_{|t|\to\infty}\mathrm{cov}_M(t)&=\Prob\bigl(|E/s|\le c,\ |U+bV_0|\le z\sqrt{1+b^2}\bigr)=1-\alpha-\pi_M,\\
\pi_M&=\Prob\bigl(|E/s|>c,\ |U+bV_0|\le z\sqrt{1+b^2}\bigr),
\end{align*}
using that $U+bV_0\sim N(0,1+b^2)$, so that the unconditional coverage is exactly $1-\alpha$. The event in $\pi_M$ has positive probability because $(E/s,U+bV_0)$ has a non-degenerate Gaussian law. The limit is a one-dimensional integral over $V_0$, given $V_0$ the law of $E/s$ being $N(-wsV_0,1-w^2s^2)$; at $w=b=1$ and $\alpha_t=0.05$ it is $0.9093$, larger than $(1-\alpha)(1-\alpha_t)=0.9025$.

\paragraph{\Cref{prop:naive-more}(ii).} Write $\pi_C(V)=\Prob(C\text{ selected}\mid V)$, which is explicit because given $V$ the selection event of \textup{\textsc{cmp}} is an interval in $D$. Using $\Phi(z+x)-\Phi(x-z)=1-\alpha-z\phi(z)x^2+O(x^4)$ and $\sqrt{1+b^2}=1+b^2/2+O(b^4)$,
\[
\xi_M(t)=\E\bigl[-t^2\pi_C(V)+(1-\pi_C(V))(1-(V-t)^2)\bigr],\qquad \xi_C(t)=\E\bigl[(1-\pi_C(V))(1-V^2)\bigr],
\]
with $V\sim N(c_gt,1)$; the $O(b^4)$ remainder is uniform for $t$ in compacts. The values $\min_t\xi_C(t)=-2.686,-0.704,-0.208$ at $w=0.5,1,2$ are computed by adaptive quadrature and agree with $(\mathrm{cov}_C(t)-0.95)/(z\phi(z)b^2)$ at $b=0.005$ to four digits (script \texttt{d23b\_naive\_limit\_exact.py}). The same script checks the limit in \cref{prop:naive}(ii) against direct integration at $t=40$ and $t=200$ (largest difference $1.2\cdot10^{-14}$). \qed

\section{Proofs for general graphs}\label{app:general}

\subsection{Exact facts about Gaussian regressions}\label{app:regression}

\begin{lemma}\label{lem:M}
Let $V_1,\dots,V_n$ be i.i.d.\ $N(\mu,\Sigma)$ with $\Sigma$ positive definite. Fix a response $j$, a set $A\not\ni j$ and $a\in A$, and let $\hat\beta$ be the least-squares coefficients of $j$ on $(1,A)$, with residual variance $\hat\sigma^2_{j\mid A}$ on $m=n-|A|-1$ degrees of freedom. Let $\beta$ be the population coefficients, $\sigma^2_{j\mid A}$ the population residual variance and $\sigma^2_{a\mid A\setminus a}$ the residual variance of $a$ on $A\setminus\{a\}$. Then
\[
\hat\beta_a-\beta_a=\frac{\sigma_{j\mid A}}{\sigma_{a\mid A\setminus a}}\frac{Z}{\sqrt Q},\qquad \hat\sigma^2_{j\mid A}=\sigma^2_{j\mid A}\frac{R}{m},
\]
with $Z\sim N(0,1)$, $Q\sim\chi^2_{m+1}$, $R\sim\chi^2_m$ independent. Consequently, if the eigenvalues of $\Sigma$ lie in $[\sigl^2,\sigu^2]$: (M3) holds for $m\ge1$, and (M1), (M2), (M4) hold for $m\ge2$:
\begin{enumerate}
\item[(M1)] $\E|\hat\beta_a-\beta_a|\le(\sigu/\sigl)\sqrt{2/\pi}\,(m-1)^{-1/2}$;
\item[(M2)] the $t$-interval $\hat\beta_a\pm t\,\widehat{\mathrm{se}}$ has expected length $2t\,(\sigma_{j\mid A}/\sigma_{a\mid A\setminus a})\,c_m\,\E Q^{-1/2}\le 2t(\sigu/\sigl)(m-1)^{-1/2}$, where $c_m=\E\sqrt{\chi^2_m/m}$;
\item[(M3)] $(\hat\beta_a-\beta_a)/\widehat{\mathrm{se}}\sim t_m$ exactly;
\item[(M4)] the expected length is at least $2t(\sigl/\sigu)\,c_m\,(m+1)^{-1/2}$.
\end{enumerate}
None of this requires a graphical model for $\Sigma$.
\end{lemma}

\begin{proof}
Given the design, the conditional law of $j$ is linear and homoscedastic Gaussian with residual variance $\sigma^2_{j\mid A}$, so $\hat\beta$ is conditionally normal with covariance $\sigma^2_{j\mid A}(\tilde D^\top\tilde D)^{-1}$ and the residual sum of squares is $\sigma^2_{j\mid A}\chi^2_m$, independent of $(\hat\beta,\tilde D)$ \citep[Chapter~8]{anderson2003introduction}. The quantity $1/[(\tilde D^\top\tilde D)^{-1}]_{aa}$ is the residual sum of squares of column $a$ on $(1,A\setminus\{a\})$, distributed as $\sigma^2_{a\mid A\setminus a}\chi^2_{m+1}$ and independent of the conditional normal. (M1) and (M2) use $\E Q^{-1/2}\le(\E Q^{-1})^{1/2}=(m-1)^{-1/2}$, $c_m\le1$, $\sigma_{j\mid A}\le\sigu$ and $\sigma_{a\mid A\setminus a}\ge\sigl$; (M4) uses $\E Q^{-1/2}\ge(\E Q)^{-1/2}$.
\end{proof}

\begin{proof}[Proof of \cref{lem:oracle}]\label{app:oracle}
The parent-adjustment estimator is regular and asymptotically normal in $\cM_{G^*}$ with variance $s_{\mathrm{pa}}^2$, so $s_{\mathrm{eff}}\le s_{\mathrm{pa}}$. Since $\sigma^2_{Y\mid A}\le\Var Y\le\sigu^2$ and $\sigma^2_{X\mid\pa(X)}=1/[(\Sigma_{SS})^{-1}]_{XX}\ge\lambda_{\min}(\Sigma_{SS})\ge\sigl^2$ for $S=\{X\}\cup\pa(X)$, $s_{\mathrm{pa}}\le\sigu/\sigl$. For the lower bound on $s_{\mathrm{eff}}$, take the one-dimensional submodel that changes only the coefficient of the edge $X\to Y$ by $\delta$. The effect is the sum of path products and the direct edge contributes its coefficient, so $d\theta/d\delta=1$; only the equation of $Y$ changes (its intercept is adjusted to keep $\E Y$ fixed), so the Fisher information is $\Var X/\omega_Y$, with $\Var X\le\sigu^2$ and $\omega_Y=\Var(Y\mid\pa(Y))=1/(\Sigma^{-1})_{YY}\ge\sigl^2$ because $Y$ is a sink. The efficient variance bound is at least $(d\theta/d\delta)^2$ divided by this information, so $s_{\mathrm{eff}}^2\ge\sigl^2/\sigu^2$.
\end{proof}

\begin{proof}[Proof of \cref{prop:union}]\label{app:union}
(i) $C_n^\cup\supseteq W_{A^*}$ for $A^*=\{X\}\cup\pa_{G^*}(X)$, and $W_{A^*}$ covers $\theta_{G^*}(P)$ with probability exactly $1-\alpha$ by (M3) and \eqref{eq:adj}. (ii) $\Leb(C_n^\cup)\le\sum_A\Leb(W_A)$. Each $m_A\ge m_{\min}=n-p$ because $Y\notin A$; by (M2), monotonicity of $t$-quantiles in the degrees of freedom and of $\E\chi_{m+1}^{-1}$ in $m$, $\E\Leb(W_A)\le2t_{n-p,1-\alpha/2}(\sigu/\sigl)(n-p-1)^{-1/2}$. By \cref{lem:oracle}, $\orl_n\ge2z_{1-\alpha/2}(\sigl/\sigu)n^{-1/2}$. Dividing gives (ii) with $\kappa_n=(t_{n-p,1-\alpha/2}/z_{1-\alpha/2})\sqrt{n/(n-p-1)}\to1$.
\end{proof}

\begin{lemma}[Acceptance probability of a conditional independence test]\label{lem:T}
Let $\Sigma$ be positive definite, $(j,i,S)$ a constraint with population partial correlation $\rho=\rho_\Sigma(j,i\mid S)$, $m=n-|S|-2\ge2$, and $c$ the critical value of the two-sided $t$-test of the coefficient of $i$ in the regression of $j$ on $(1,i,S)$. Then
\[
\Prob(\text{the test does not reject})\le\Phibar\bigl(|\rho|\sqrt{(m+1)/2}-2c\bigr)+e^{-(m+1)/16}+e^{-m/2}.
\]
\end{lemma}
\begin{proof}
By \cref{lem:M} the $t$-statistic is $(Z+\nu)/\sqrt{R/m}$ with $\nu=(\rho/\sqrt{1-\rho^2})\sqrt Q$, using $\beta_i=\rho\sigma_{j\mid S}/\sigma_{i\mid S}$ and $\sigma^2_{j\mid i,S}=\sigma^2_{j\mid S}(1-\rho^2)$. By \citet[Lemma~1]{laurent2000adaptive}, $\Prob(Q<(m+1)/2)\le e^{-(m+1)/16}$ and $\Prob(R/m>4)\le e^{-m/2}$. On the complement, $|\nu|\ge|\rho|\sqrt{(m+1)/2}$ and non-rejection implies $|Z+\nu|\le2c$, whose conditional probability given $Q$ is at most $\Phibar(|\nu|-2c)$.
\end{proof}

\begin{lemma}[Sample correlations]\label{lem:R}
Let $(U_i,E_i)$, $i\le\nu$, be i.i.d.\ $N(0,I_2)$ with sample second moments $a,d,b$, and let $r(\rho)=(\rho a+sb)/\sqrt{a(\rho^2a+2\rho sb+s^2d)}$ with $s=\sqrt{1-\rho^2}$; $r(\rho)$ has the law of a sample correlation from $\nu$ pairs with correlation $\rho$. (R1) If $|a-1|,|d-1|,|b|\le\eps\le1/14$, then $|r(\rho)-\rho|\le5\eps$ for all $|\rho|<1$. (R2) For $\eps\le1/14$ and all $|\rho|<1$, $\Prob(|r(\rho)-\rho|>5\eps)\le8e^{-\nu\eps^2/16}$.
\end{lemma}
\begin{proof}
(R1) Write $a=1+\alpha_0$, $d=1+\delta_0$ and $T=1+\rho^2\alpha_0+s^2\delta_0+2\rho sb$; then $|T-1|\le2\eps$ and $r=\rho\sqrt{(1+\alpha_0)/T}+sb/\sqrt{(1+\alpha_0)T}$. The ratio $x=(1+\alpha_0)/T$ satisfies $|x-1|\le3\eps/(1-2\eps)\le3.5\eps$, so the first term differs from $\rho$ by at most $3.5\eps$; the second is at most $\eps/\sqrt{0.79}<1.13\eps$. (R2) \citet[Lemma~1]{laurent2000adaptive} with $x=\nu\eps^2/16$ for the upper tail and $x=\nu\eps^2/4$ for the lower tail bounds $\Prob(|a-1|>\eps)$ and $\Prob(|d-1|>\eps)$ by $2e^{-\nu\eps^2/16}$ each. For $b$, $UE=[((U+E)/\sqrt2)^2-((U-E)/\sqrt2)^2]/2$, so $b=(a'-d')/2$ with $a',d'$ independent $\chi^2_\nu/\nu$, and $\{|b|>\eps\}\subseteq\{|a'-1|>\eps\}\cup\{|d'-1|>\eps\}$. A union bound gives the constant $8$.
\end{proof}

\paragraph{The exact partial correlation interval.} Let $r$ be the sample partial correlation of $(u,v)$ given $W$ (with intercept), $\nu_W=n-|W|-1$. Its law is that of the sample correlation of $\nu_W+1$ bivariate Gaussian pairs and depends only on $\rho=\rho(u,v\mid W)$ \citep[Chapter~4]{anderson2003introduction}; write $F_\nu(\cdot;\rho)$ for its distribution function, which is continuous and non-increasing in $\rho$ for each fixed argument, because the sample correlation is stochastically increasing in $\rho$. For level $1-\eta$ let
\[
\rho_U(r)=\sup\{\rho':F_\nu(r;\rho')\ge\eta/2\},\qquad \rho_L(r)=\inf\{\rho':F_\nu(r;\rho')\le1-\eta/2\},
\]
over $\rho'\in[-1,1]$ with the limiting conventions $F_\nu(r;-1)=1$ and $F_\nu(r;1)=0$, and $I_\rho=[\rho_L(r),\rho_U(r)]$. Since $F_\nu(r;\rho)$ is uniform on $(0,1)$ when $|\rho|<1$, $\Prob(\rho\in I_\rho)=1-\eta$ exactly.

\begin{corollary}[Deterministic endpoints]\label{cor:Rprime}
Let $\eps_\nu=21\sqrt{\log(16/\eta)/\nu}\le5/14$. For every $r\in(-1,1)$, $\rho_U(r)\le r+\eps_\nu$ and $\rho_L(r)\ge r-\eps_\nu$; hence $\max_{\rho'\in I_\rho}|\rho'|\le|r|+\eps_\nu$.
\end{corollary}
\begin{proof}
Let $\rho'=r+\eps_\nu<1$ (otherwise there is nothing to prove) and let $\varrho$ denote the sample partial correlation under $\rho'$. By \cref{lem:R}(R2) with $\eps=\eps_\nu/5$, $F_\nu(r;\rho')\le\Prob_{\rho'}(|\varrho-\rho'|>\eps_\nu)\le8(\eta/16)^{441/400}<\eta/2$. By monotonicity, $F_\nu(r;\rho'')<\eta/2$ for all $\rho''\ge\rho'$, so $\rho_U(r)\le\rho'$. Symmetrically, with $\rho'=r-\eps_\nu$, $1-F_\nu(r;\rho')<\eta/2$, so $F_\nu(r;\rho'')>1-\eta/2$ for all $\rho''\le\rho'$ and $\rho_L(r)\ge\rho'$.
\end{proof}

\begin{lemma}[The $k$-TF screen rejects violators]\label{lem:B}
Let the eigenvalues of $\Sigma$ lie in $[\sigl^2,\sigu^2]$, $K=\sigu/\sigl$, and let a $k$-TF element $(u\to v,W)$ of a candidate satisfy $k|\lambda^{\mathrm{std}}|-|\rho|\ge\mu>0$ at the population level. There is $n_B(\mu)$, of order $\mu^{-2}\max(k^2K^4,1)\log(N_T/\gamma_2)$ for small $\mu$, such that for $n\ge n_B(\mu)$ the candidate passes this element of the screen with probability at most $16\exp(-c_Bn\mu^2)$, $c_B=\min\{1/(6272k^2K^4),1/(5\cdot10^5)\}$.
\end{lemma}
\begin{proof}
Note $|\lambda^{\mathrm{std}}|\le K^2$, so $\mu\le kK^2$. Let $\eps_s=\mu/(12kK^2)$, $\eps_\beta=\mu/(4kK)$ and $\eps_r=\min(\mu/25,5/14)$. Take $n_B(\mu)$ so large that $n\ge2p+2$ (all relevant degrees of freedom are at least $n/2$), $4\sqrt2K t'/\sqrt{m_\beta+1}\le\eps_\beta/2$ with $t'$ the quantile of $I_\beta$, the ratio $\chi^2_{n-1,\eta/2}/\chi^2_{n-1,1-\eta/2}\ge1-\eps_s$, $\eps_\nu\le\mu/5$ and $\eps_s\le1/2$. Consider the events $G_1=\{Q\ge(m_\beta+1)/2,\ R/m_\beta\le4,\ |Z|\le\eps_\beta\sqrt{m_\beta+1}/(4\sqrt2K)\}$ in the representation of \cref{lem:M}, $G_2=\{|s_u^2/\sigma_u^2-1|\le\eps_s,\ |s_v^2/\sigma_v^2-1|\le\eps_s\}$ and $G_3=\{|r-\rho|\le\eps_r\}$. On $G_1$, $\mathrm{dist}(0,I_\beta)\ge|\beta|-\eps_\beta$; on $G_2$, $\sqrt{\sigma^2_{u,L}/\sigma^2_{v,U}}\ge(\sigma_u/\sigma_v)(1-3\eps_s)$; on $G_3$ and by \cref{cor:Rprime}, $\max_{I_\rho}|\rho'|\le|\rho|+\eps_r+\mu/5$. Hence on $G_1\cap G_2\cap G_3$,
\[
k\,\mathrm{dist}(0,I_\beta)\sqrt{\sigma^2_{u,L}/\sigma^2_{v,U}}\ge k|\lambda^{\mathrm{std}}|-\mu/2\ge|\rho|+\mu/2>\max_{I_\rho}|\rho'|,
\]
and the candidate is rejected. By \cref{lem:M}, \citet[Lemma~1]{laurent2000adaptive} and \cref{lem:R}, the complements have probabilities at most $e^{-(m_\beta+1)/16}+e^{-m_\beta/2}+2e^{-\eps_\beta^2(m_\beta+1)/(64K^2)}$, $4e^{-(n-1)\eps_s^2/16}$ and $8e^{-\nu\eps_r^2/400}$. Each exponent is at least $c_Bn\mu^2$ (for the truncated case of $\eps_r$, $\nu(5/14)^2/400\ge n/6272\ge n\mu^2/(6272k^2K^4)$), and the coefficients sum to $16$.
\end{proof}

\subsection{Coverage}\label{app:coverage}

\paragraph{Whitened coordinates.} Let $q_j=|\mathrm{pred}(j)\setminus\pa_D(j)|$, so that $\sum_jq_j=q$; then $T_D=\sum_jn\log(1+|\hat\tau_{D,j}|^2/(n-j))$, and replacing the member of $S^*$ that defines the curve $P_\varepsilon$ changes $u_{S^*,S'}$ by a non-zero factor, while replacing $D_{S'}$ or its topological order changes it by an isometry. For a DAG $D$ with a topological order, let $E(j)=\mathrm{pred}(j)\setminus\pa_D(j)$, let $\hat\gamma_j$ be the coefficients of $E(j)$ in the regression of $j$ on $(1,\mathrm{pred}(j))$, $M_j$ the cross-product matrix of the residuals of the columns of $E(j)$ after regression on $(1,\pa_D(j))$, and $\hat\sigma^2_j$ the residual sum of squares of the regression of $j$ on $(1,\mathrm{pred}(j))$ divided by $n-j$. Then $\hat\tau_D=(M_j^{1/2}\hat\gamma_j/\hat\sigma_j)_j\in\R^q$, and its population version is $\mu_D(P)=(\Sigma_{E(j)E(j)\cdot\pa_D(j)}^{1/2}\gamma_j/\sigma_{j\mid\mathrm{pred}(j)})_j$, with the population coefficients $\gamma_j$ and partial covariances; since $\sigma^2_{j\mid\pa_D(j)}=\sigma^2_{j\mid\mathrm{pred}(j)}+\gamma_j^\top\Sigma_{E(j)E(j)\cdot\pa_D(j)}\gamma_j$, $\mu_D$ vanishes exactly on $\cM_D$.

\begin{proof}[Proof of \cref{thm:coverage}]
(i) The maximum likelihood estimates of $\cM_D$ and of the saturated model factorize along $D$ and along the complete DAG with the same order, and consist of node-wise regressions with intercept, so $T_D$ is twice the log likelihood ratio, and $\hat\sigma^2_{j\mid\pa(j)}/\hat\sigma^2_{j\mid\mathrm{pred}(j)}=1+|\hat\tau_{D,j}|^2/(n-j)$. If $P\in\cM_D$, then conditionally on the data of $\mathrm{pred}(j)$ the data of $j$ are $a\mathbf 1+X_{\pa(j)}\beta+\varepsilon$ with $\varepsilon\sim N(0,\sigma^2I_n)$, so $M_j^{1/2}\hat\gamma_j=M_j^{-1/2}\tilde X_{E(j)}^\top(I-H_j)\varepsilon\sim N(0,\sigma^2I_{q_j})$, where $\tilde X_{E(j)}$ are the columns of $E(j)$ and $H_j$ is the projection onto the span of $(\mathbf 1,X_{\pa(j)})$, and the residual sum of squares of the regression on $(\mathbf 1,X_{\mathrm{pred}(j)})$ is $\sigma^2\chi^2_{n-j}$. The two are functions of the projections of $\varepsilon$ onto orthogonal subspaces, so they are conditionally independent, and their conditional law does not depend on the design. Hence $(N_j,V_j)$ is independent of the data of $\mathrm{pred}(j)$, of which $(N_i,V_i)$, $i<j$, are functions, and by induction all are independent; $U_j=|N_j|^2$.

(ii) Let $(G^*,P)\in\cM_k$ and let $E_n$ be the intersection of: $G^*$ passes step (1); all $4|\mathcal T_{G^*}|$ exact intervals of the elements of $\mathcal T_{G^*}$ cover their parameters; and the interval of $G^*$ covers $\theta_{G^*}(P)$. By (i), step (1a) rejects $G^*$ with probability at most $\gamma_1$ (exactly $\gamma_1$ if $q\ge1$). Step (1b) tests $G^*$ only if step (0) returns an estimate, which requires $n_1\ge p+1$, and $G^*\simeq D_{S'}$ for a class $S'$ of $\hat\cX$; given the first part of the sample, $\hat B_{S'}$ and $D_{S'}$ are fixed, the second part is an independent sample of size $n_2\ge p+1$ from $P\in\cM_{G^*}=\cM_{D_{S'}}$, and by (i) the conditional law of $\hat\tau^{(2)}_{D_{S'}}$, hence of the test statistic, is the null law, so the conditional rejection probability is $\gamma_3$; taking expectations, step (1b) rejects $G^*$ with probability at most $\gamma_3$. The exact intervals of step (2) fail with total probability at most $\gamma_2$. For \DSA\ the interval of $G^*$ is exact by (M3) and validity of $A(G^*)$; for \DSeff\ its non-coverage tends to $\alpha'$ uniformly over $\cM_k$ (proof of \cref{thm:explicit}(b)). A union bound gives $\Prob(E_n^c)\le\gamma_1+\gamma_2+\gamma_3+\alpha'=\alpha$ for \DSA, and $\alpha+o(1)$ uniformly for \DSeff. On $E_n$, for every element of $\mathcal T_{G^*}$, $k\,\mathrm{dist}(0,I_\beta)\sqrt{\sigma^2_{u,L}/\sigma^2_{v,U}}\le k|\beta|\sigma_u/\sigma_v=k|\lambda^{\mathrm{std}}|\le|\rho|\le\max_{I_\rho}|\rho'|$, so $G^*$ is not rejected and its interval, which contains the target, is part of the output. For the conditional independence screen of \cref{prop:CI}, step (1) instead rejects $G^*$ with probability at most $\gamma_1$ by a union bound, because under $P$ the local Markov constraints of $G^*$ hold and each $t$-statistic is exactly $t$-distributed.
\end{proof}

\subsection{Explanation sets and the graphical criterion}\label{app:graph}

\begin{proof}[Proof of \cref{lem:G}]
If every d-separation of $G'$ holds in $H$, then, $P_0$ being Markov to $H$, $P_0$ satisfies the global Markov property of $G'$, so $P_0\in\cM_{G'}$. Conversely $P_0\in\cM_{G'}$ implies that the d-separations of $G'$ hold in $P_0$, and by faithfulness they are d-separations of $H$. If $G'-e$ is equivalent to $H$, then $P_0\in\cM_H=\cM_{G'-e}\subseteq\cM_{G'}$.
\end{proof}

\begin{lemma}\label{lem:Z}
If $G'\in\cX$ and the coefficient of $e$ in the regression parametrization of $P_0$ under $G'$ is zero, then $G'\in\Xeq$.
\end{lemma}
\begin{proof}
If $e=t\to h$ in $G'$ has zero coefficient, then $h\indep t\mid\pa_{G'}(h)\setminus\{t\}$ and the factorization of $P_0$ along $G'$ is a factorization along $G'-e$, so $P_0\in\cM_{G'-e}$. By faithfulness the d-separations of $G'-e$ are d-separations of $H$. The two graphs have the same skeleton, so by \citet[Theorem~4]{chickering2002optimal} $H$ can be transformed into $G'-e$ by covered edge reversals only, which preserve Markov equivalence \citep[Lemma~2]{chickering2002optimal}.
\end{proof}

\begin{proof}[Proof of \cref{prop:graph}]
(a) If $G'\in\Xeq$, then $G'-e$ is equivalent to $H$ and $Y$ is a sink in it. Conversely, let $D\in[H]_Y$. If $e\notin R$, take a topological order of $D$ in which $Y$ is last, which is possible because $Y$ has no children, and orient $e$ from the earlier to the later endpoint; if $e\in R$, orient it as required, which gives an acyclic graph by admissibility. The resulting $G'=D+e$ is acyclic, has skeleton $S$, satisfies $\cK$, so $G'\in\cG$; and $G'-e=D$ is equivalent to $H$, so $G'\in\Xeq$. (This uses that $\cG$ contains all such orientations.)

(b) Let $G'\in\Xeq$, $D=G'-e$ and $P\in\cM_H=\cM_D$. If $e=t\to h$ in $G'$, then $t$ is a non-descendant of $h$ in $D$, so $h\indep t\mid\pa_D(h)$ and the coefficient of $e$ in the regression parametrization of $P$ under $G'$ is zero; the other coefficients coincide with those under $D$. Hence $\theta_{G'}(P)=\theta_D(P)$, and $\thadj_{G'}(P)=\theta_{G'}(P)$ because $\pa_{G'}(X)$ is a valid adjustment set and $Y\notin\pa_{G'}(X)$. Together with (a) this proves (b).

(c) The required orientations $\{v\to Y:v\in\adj(Y)\}\cup(R\setminus\{e\})$ are consistent background knowledge (they hold in $H$), so by \citet{meek1995causal} the MPDAG obtained from the CPDAG of $H$ represents exactly $[H]_Y$. If the effect is identifiable in $\MP(H)$, \citet[Theorem~3.6]{perkovic2020identifying} gives a formula for $f(y\mid do(x))$ in terms of the observational density that is the same for all members of $[H]_Y$. For $P\in\cM_H$ and $D\in[H]_Y$, the interventional mean $\E_D[Y\mid do(X=x)]$ is affine in $x$ with slope $\theta_D(P)$, so all $\theta_D(P)$ coincide, which is (b). If the effect is not identifiable, the proof of \citet[Proposition~3.2, Appendix~B]{perkovic2020identifying} constructs $D_1,D_2\in[H]_Y$ and a zero-mean linear Gaussian structural equation model along $D_1$ with positive noise variances and unit marginal variances, hence a distribution $f\in\cM_{D_1}=\cM_H$, such that $\E_{D_1}[Y\mid do(X=1)]\ne\E_{D_2}[Y\mid do(X=1)]$. For a positive definite Gaussian $f\in\cM_D$ the interventional distribution under $D$ is unique (truncated factorization), and both means vanish at $x=0$, so $\theta_{D_1}(f)\ne\theta_{D_2}(f)$ and (b) fails.
\end{proof}

\begin{corollary}\label{cor:K1}
For $G'\in\Xeq$, the following are equivalent: (i) at every point of $\cM_{G'}$ that is faithful to $G'$, all candidates in $\cG$ that explain it give the same value of \eqref{eq:adj}; (ii) the total effect of $X$ on $Y$ is identifiable in $\MP(G')$.
\end{corollary}
\begin{proof}
At a point $P$ faithful to $G'$, a candidate $G''$ explains $P$ if and only if the d-separations of $G''$ hold in $G'$ (as in \cref{lem:G}); since the skeletons agree, \citet[Theorem~4, Lemma~2]{chickering2002optimal} show that then $G''$ is equivalent to $G'$. So the explaining candidates are $[G']_Y$. If (i) holds, then since faithful points are dense in $\cM_{G'}$ and the adjustment coefficients are continuous on the positive definite cone, all members of $[G']_Y$ give the same value on all of $\cM_{G'}$; if (ii) failed, the linear Gaussian counterexample in the proof of \cref{prop:graph}(c) would give a point of $\cM_{G'}$ where two members differ. Conversely (ii) implies (i) by \citet[Theorem~3.6]{perkovic2020identifying}.
\end{proof}

\subsection{The graphical criterion for \texorpdfstring{\eqref{eq:K1p}}{(K1')}}\label{app:K1p}

Throughout, $e=\{u,v\}$, $S_c\subseteq\Xeq$ is a Markov equivalence class with representative $G_1$, $C^*=\mathrm{CC}_{G_1}$, and for a partially directed graph $\mathcal P$ we write $\mathrm{arr}(\mathcal P)$ for its directed edges. We use the four orientation rules of \citet{meek1995causal} (an arrow is $\to$, an undirected edge is $-$, and adjacency ignores direction):
R1: $a\to b-c$, $a,c$ non-adjacent $\Rightarrow b\to c$;
R2: $a\to c\to b$, $a-b$ $\Rightarrow a\to b$;
R3: $a-c\to b$, $a-d\to b$, $a-b$, $c,d$ non-adjacent $\Rightarrow a\to b$;
R4: $a-b$, $c\to d\to b$, $a$ adjacent to $c$ and $d$, $c,b$ non-adjacent $\Rightarrow a\to b$.
By \citet[Theorem~4]{meek1995causal}, closing a CPDAG together with consistent background knowledge under R1--R4 gives the MPDAG whose directed edges are exactly the orientations common to all DAGs in the class that satisfy the background knowledge. A \emph{chain graph} has no partially directed cycle; its chain components are the connected components of its undirected part. We use two facts: (F1) the two endpoints of a directed edge lie in different chain components, since otherwise the edge and an undirected path between its endpoints form a partially directed cycle; (F2) the relation $\mathrm{comp}(a)\prec\mathrm{comp}(b)$ generated by directed edges $a\to b$ is acyclic.

\begin{lemma}\label{lem:C}
For $G'\simeq G_1$, $G'-e\simeq H$ if and only if $\mathrm{CC}_{G'}=C^*$. Hence $S_c=\{D\in[G_1]_Y:\mathrm{CC}_D=C^*\}$.
\end{lemma}
\begin{proof}
For a DAG $D$ containing $e$, an unshielded collider $a\to x\leftarrow b$ of $D-e$ either has $a,b$ non-adjacent in $S$, and is then an unshielded collider of $D$ with no arm equal to $e$, or has $\{a,b\}=e$, that is, $x\in\mathrm{CC}_D$. Conversely, unshielded colliders of $D$ with no arm equal to $e$ remain unshielded colliders of $D-e$. The two parts are disjoint. If $G'\simeq G_1$ the first parts agree, so $G'-e$ and $G_1-e\simeq H$ have the same unshielded colliders, and hence are equivalent, if and only if $\mathrm{CC}_{G'}=C^*$.
\end{proof}

\begin{lemma}\label{lem:P}
$\mathrm{arr}(\MP(H))\subseteq\mathrm{arr}(\mathcal P)$ for $\mathcal P$ the partially directed graph of orientations common to $S_c$, and for $\mathcal P=\MP''(S_c)$. If \eqref{eq:K2} holds, every bad path of such a $\mathcal P$ contains $e$.
\end{lemma}
\begin{proof}
Members of $S_c$, and members of the extension set of $\MP''(S_c)$, contain the directed edges of $\MP(H)$: the former because their restriction to $S-e$ lies in $[H]_Y$, the latter by construction. Let $p$ be a bad path of $\mathcal P$ that avoids $e$. Then $p$ lies in $S-e$. It has no backward edge in $\MP(H)$, because $\MP(H)$ has fewer directed edges than $\mathcal P$, and for the same reason its first edge is undirected in $\MP(H)$. So $p$ is a bad path of $\MP(H)$, which \eqref{eq:K2} excludes by \cref{prop:graph}(c).
\end{proof}

\begin{lemma}\label{lem:MPH}
If $R=\emptyset$, $\MP(H)$ is the CPDAG of $H$ with every edge incident to $Y$ oriented into $Y$. It is a chain graph in which $Y$ is a singleton chain component. If $R\ne\emptyset$ and $\MP(H)$ is a chain graph, then $Y$ is a singleton chain component because it has no outgoing and no undirected edge.
\end{lemma}
\begin{proof}
The CPDAG is closed under R1--R4. After orienting the edges at $Y$, a rule can fire anew only if one of its premise arrows is a new arrow $w\to Y$; requirements on undirected edges and adjacencies cannot become easier to meet. In R1 the new arrow is $a\to b$ with $b=Y$, but $b-c$ would be an undirected edge at $Y$. In R2, R3 and R4 the new arrow ends at $b=Y$ and the premise needs $a-b$ undirected, or (R2 with $c=Y$, R4 with $d=Y$) it needs an arrow out of $Y$. Neither exists, so the closure is unchanged. The CPDAG is a chain graph \citep{andersson1997characterization}, and a partially directed cycle through $Y$ would need an edge out of $Y$ or an undirected edge at $Y$. (With required edges the closure is the MPDAG of the CPDAG with $R\setminus\{e\}$ and the edges into $Y$ \citep{meek1995causal}; it need not be a chain graph, though it is for tier knowledge by \cref{lem:tier}, and the hypothesis is used only in \cref{lem:detour}, that is, only for the converse in \cref{thm:K1p-graph}.)
\end{proof}

The proof of \cref{lem:tier} uses a more general statement, which also covers some required edges that do not come from a time order.

\begin{lemma}[Cut knowledge]\label{lem:tier-cut}
$\MP(H)$ is a chain graph with chordal chain components whenever there are sets $I_1,\dots,I_r\subseteq V$ such that $R$ is the set of edges of $S$ not at $Y$ that have exactly one endpoint in some $I_j$ (the edges at $Y$ are oriented by the sink requirement anyway); for instance, $R$ may consist of all edges not at $Y$ at some nodes. Admissibility can fail in this case.
\end{lemma}

\begin{proof}[Proof of \cref{lem:tier,lem:tier-cut}]
\Cref{lem:tier}, admissibility. Put $T_{m+1}=\{Y\}$ and let $t(a)$ be the index of the tier of $a$. In every $D\in[H]_Y$ an edge between different tiers points to the later tier: by $R\setminus\{e\}$ if neither endpoint is $Y$, and because $Y$ is a sink otherwise. So $t$ does not decrease along directed paths of $D$. If $e\in R$, then $e=u\to v$ with $t(u)<t(v)$, $D$ has no directed path from $v$ to $u$, and $D+e$ is acyclic; it has skeleton $S$ and satisfies $\cK$, so $D+e\in\cG$. If $e\notin R$, admissibility is automatic (\cref{sec:gen-setting}).

\Cref{lem:tier-cut}. We use the notions of \citet{hauser2012characterization}. For $I\subseteq V$ and a DAG $D$, the intervention graph $D^{(I)}$ is $D$ without the edges whose head lies in $I$ (their Definition~5). An edge with both endpoints in $I$ is absent from $D^{(I)}$ and an edge with no endpoint in $I$ is present, whatever their orientation; an edge with exactly one endpoint in $I$ is present if and only if it points out of $I$. Hence two DAGs $D,D'$ on $S-e$ have intervention graphs $D^{(I)},D'^{(I)}$ with the same skeleton if and only if they orient alike every edge of $S-e$ with exactly one endpoint in $I$. Let $I_{r+1}=V\setminus\{Y\}$ and $\mathcal I=\{\emptyset,I_1,\dots,I_{r+1}\}$, which is a conservative family of targets because it contains $\emptyset$ (their Definition~6). The edges of $S-e$ with exactly one endpoint in some $I_j$, $j\le r$, are those of $R\setminus\{e\}$ and possibly some edges at $Y$, and those with exactly one endpoint in $I_{r+1}$ are the edges at $Y$. By their Theorem~10, (i)$\Leftrightarrow$(iv), a DAG $D$ on $S-e$ is $\mathcal I$-Markov equivalent to $H$ if and only if it is Markov equivalent to $H$ and orients $R\setminus\{e\}$ and the edges at $Y$ as $H$ does, that is, if and only if $D\in[H]_Y$, because $H\in[H]_Y$. So $[H]_Y$ is the $\mathcal I$-Markov equivalence class of $H$. By \citet[Theorem~4]{meek1995causal}, an edge is directed in $\MP(H)$ if and only if all members of $[H]_Y$ orient it alike, so $\MP(H)$ is the union of the members of $[H]_Y$ in the graph-theoretic sense, which is the $\mathcal I$-essential graph of $H$ (their Definition~11). By their Proposition~15 it is a chain graph whose chain components induce chordal subgraphs; their chain graphs, graphs without a directed cycle in the sense of their Appendix~A.1, are ours.

For \cref{lem:tier}, chain structure, apply \cref{lem:tier-cut} with $r=m-1$ and $I_j=T_1\cup\dots\cup T_j$: an edge between $T_i$ and $T_l$, $i<l\le m$, has exactly one endpoint in $I_i$, and an edge within a tier has none or both endpoints in each $I_j$, so the edges not at $Y$ with exactly one endpoint in some $I_j$ are those of $R$.

Admissibility can fail in \cref{lem:tier-cut}. Take $S$ with edges $X-b$, $b-c$, $X-c$, $X-Y$, the graph $G=\{X\to b,b\to c,X\to c,X\to Y\}$, $R=\{X\to b,b\to c\}$, the set of edges at $b$, and $e=b\to c$. Then $D=\{X\to b,c\to X,X\to Y\}$ is Markov equivalent to $H$ and lies in $[H]_Y$, but $D+e$ contains the cycle $c\to X\to b\to c$.
\end{proof}

\begin{lemma}\label{lem:D1}
Let $\mathcal P$ be as in \cref{lem:P}, let $D$ be a consistent extension of $\mathcal P$ (so $D\supseteq\mathrm{arr}(\MP(H))$), and let $x\in\mathrm{CC}_D\setminus C^*$. Then $x\ne Y$, $u,v\ne Y$, and $u-x$ and $v-x$ are undirected in $\MP(H)$.
\end{lemma}
\begin{proof}
If $x=Y$ then $u\to Y\leftarrow v$ in every candidate, so $Y\in C^*$; and $Y$ has no children, so $u,v\ne Y$. Since $u\to x$ and $v\to x$ in $D\supseteq\mathrm{arr}(\MP(H))$, $\MP(H)$ contains neither $x\to u$ nor $x\to v$. Suppose $v\to x$ is in $\MP(H)$. If also $u\to x$, then $u\to x\leftarrow v$ is an unshielded collider of every member of $[H]_Y$ ($u,v$ are non-adjacent in $H$), so $x\in\mathrm{CC}_{G_1}=C^*$, a contradiction. If $u-x$ is undirected, R1 applied to $v\to x-u$ forces $x\to u$ in $\MP(H)$, a contradiction. The case $u\to x$ is symmetric.
\end{proof}

\begin{lemma}[Detour]\label{lem:detour}
Assume that $\MP(H)$ is a chain graph, as it is for tier knowledge and in particular for $R=\emptyset$ (\cref{lem:tier}). Let $\mathcal P$ and $D$ be as in \cref{lem:D1}, and suppose that some $x\in\mathrm{CC}_D\setminus C^*$ exists. If $\mathcal P$ has a bad path, then $\MP(H)$ has a bad path. The lemma does not use \eqref{eq:K2}.
\end{lemma}
\begin{proof}
Let $p$ be a bad path of $\mathcal P$. Two facts transfer from $\mathcal P$ to $\MP(H)$ because $\mathrm{arr}(\MP(H))\subseteq\mathrm{arr}(\mathcal P)$: a path in $S-e$ without backward edges in $\mathcal P$ (including chords) has none in $\MP(H)$, and an edge of $S-e$ that is undirected in $\mathcal P$ is undirected in $\MP(H)$. If $p$ avoids $e$ it is a bad path of $\MP(H)$, as in \cref{lem:P}. Otherwise let $u=V_i$, $v=V_{i+1}$ in the order of $p$ (exchange the names if necessary), and take $x$ as in \cref{lem:D1}.
(1) If $x=V_j$ with $j<i$, let $p'=\langle V_0,\dots,V_j,V_{i+1},\dots,V_k\rangle$; its new edge $x-v$ is undirected, its vertices are a subsequence of those of $p$, so it has no backward edges, and its first edge is undirected (it is $X-v$ if $j=0$).
(2) If $x=V_j$ with $j>i+1$, let $p'=\langle V_0,\dots,V_i,V_j,\dots,V_k\rangle$, and argue in the same way with the edge $u-x$.
(3) If $x$ is not on $p$, let $p'=\langle V_0,\dots,u,x,v,\dots,V_k\rangle$; its new edges are undirected and its first edge is undirected. A backward chord would be either $x\to y$ with $y=V_j$, $j\le i-1$, or $y\to x$ with $y=V_j$, $j\ge i+2$. In the first case the segment of $p$ from $y$ to $u$ uses only undirected and forward edges, so $\mathrm{comp}(y)\preceq\mathrm{comp}(u)$ by (F1) and \cref{lem:MPH}; $u-x$ gives $\mathrm{comp}(u)=\mathrm{comp}(x)$; and $x\to y$ gives $\mathrm{comp}(x)\prec\mathrm{comp}(y)$ by (F1). This is a cycle, contradicting (F2). The second case is symmetric; there $y\ne Y$ because $Y$ has no children.
In each case $p'$ is a bad path of $\MP(H)$.
\end{proof}

\begin{proof}[Proof of \cref{thm:K1p-graph}]
\emph{Independence of the representative.} Members of $S_c$ have the same Markov equivalence class and, by \cref{lem:C}, the same $\mathrm{CC}$; the directed edges of $\MP(H)$ depend only on $H$. The background knowledge is consistent because $G_1$ satisfies it.
\emph{($\Leftarrow$).} Let $\mathrm{Ext}''$ be the members of $[G_1]_Y$ that satisfy the background knowledge, so that $\MP''(S_c)$ represents $\mathrm{Ext}''$ \citep[Theorem~4]{meek1995causal}. Members of $S_c$ satisfy it (\cref{lem:C,lem:P}), so $S_c\subseteq\mathrm{Ext}''$. If $\MP''(S_c)$ has no bad path, all members of $\mathrm{Ext}''$, hence of $S_c$, have the same effect on $\cM_{G_1}$ \citep[Theorem~3.6]{perkovic2020identifying}.
\emph{($\Rightarrow$).} Suppose $\MP''(S_c)$ has a bad path. If $\mathrm{Ext}''\ne S_c$, some $D\in\mathrm{Ext}''\setminus S_c$ has $\mathrm{CC}_D\supsetneq C^*$ (it contains $C^*$ by the background knowledge, and differs from $C^*$ by \cref{lem:C}), and \cref{lem:detour} gives a bad path of $\MP(H)$, contradicting \eqref{eq:K2}. So $\mathrm{Ext}''=S_c$, $S_c$ is represented by the MPDAG $\MP''(S_c)$, and the linear Gaussian counterexample in the proof of \citet[Proposition~3.2]{perkovic2020identifying}, used as in \cref{prop:graph}(c), gives two members of $S_c$ with different effects at a point of $\cM_{G_1}$.
\emph{Common orientations.} If $\MP(H)$ is a chain graph, \eqref{eq:K1p} also holds if and only if, for every $S_c$, the partially directed graph of orientations common to all members of $S_c$ has no bad path, and this graph is the MPDAG of the set of members of $[G_1]_Y$ that contain it. The same argument with $\mathcal P$ the graph of orientations common to $S_c$ and $\mathrm{Ext}=\{D\simeq G_1:D\supseteq\mathrm{arr}(\mathcal P)\}$ applies. $\mathrm{Ext}$ is represented by $\mathcal P$: the orientations common to $\mathrm{Ext}$ contain $\mathrm{arr}(\mathcal P)$ by definition, and are contained in it because $S_c\subseteq\mathrm{Ext}$.
\end{proof}

\subsection{Cancellation explanations}\label{app:HC}

\begin{proof}[Proof of \cref{lem:HC}]
Let $G'\in\Xcanc$ and let $G''$ be $G'$ with all edges whose coefficient in the regression parametrization of $P_0$ is zero removed.
\emph{Markov.} Removing a zero-coefficient parent from a regression leaves the other coefficients and the residual variance unchanged, because the residual is orthogonal to all original parents; so $P_0\in\cM_{G''}$.
\emph{Minimality.} Every edge of $G''$ has non-zero coefficient, so removing any edge violates a local Markov condition; since Markov models grow with edges, $P_0$ is Markov to no proper subgraph of $G''$, which is the minimality condition of \citet[p.~246]{zhang2008detection}.
\emph{Failure of faithfulness.} By \cref{lem:Z}, the coefficient of $e$ under $G'$ is non-zero, so $e\in G''$. If $P_0$ were faithful to $G''$, the d-separations of $G''$ and $H$ would coincide, so $G''$ and $H$ would be Markov equivalent and have the same skeleton, which is false.
\emph{Undetectability.} $P_0$ is faithful to the DAG $H$, which is the definition of an undetectable failure in \citet[p.~252]{zhang2008detection}: the distribution is not faithful to the true DAG but is faithful to some other DAG; no relation between $H$ and $G''$ is required.
With $G''$ in the role of the true DAG, \citet[Theorem~2]{zhang2008detection} shows that triangle faithfulness, as defined in \citet[p.~257]{zhang2008detection} with the same quantifiers as our legal sets, fails for $(G'',P_0)$: there are a triangle $\{u,v,c\}$ of $G''$ and a legal $W$ with $\rho_{P_0}(u,v\mid W)=0$, where the edge $u-v$ has non-zero coefficient. (Their proof, in Appendix~C of that paper, uses only Markov, minimality and faithfulness to some DAG, as here.) The triangle is also a triangle of $G'$ with the same orientations, and the coefficient of $u-v$ under $G'$ equals that under $G''$.
\end{proof}

\subsection{The local error bound}\label{app:K4E}

\begin{lemma}\label{lem:K4E}
Let $G^*,G'\in\Xeq$ be non-equivalent and assume \eqref{eq:K2}. Parametrize $\cM_{G^*}$ by $\phi=(\lambda,\omega)$, let $q(\phi)=\thadj_{G'}(P_\phi)-\theta_{G^*}(P_\phi)$ and $h(\phi)=h_{G'}(P_\phi)$, and let $\phi_0$ satisfy $\lambda_e(\phi_0)=0$ and $\partial_{\lambda_e}h(\phi_0)\ne0$. Then there are a neighbourhood $B$ of $\phi_0$ and $C_4<\infty$ with $|q|\le C_4\|h\|_\infty$ on $B$; as $B$ shrinks, $C_4\to|\partial_eq(\phi_0)|/\|\partial_eh(\phi_0)\|_\infty$.
\end{lemma}
\begin{proof}
On $\{\lambda_e=0\}$, $P_\phi\in\cM_{G^*-e}=\cM_H=\cM_{G'-e}\subseteq\cM_{G'}$, so $h=0$; and by \eqref{eq:K2} and \cref{prop:graph}(b), $\thadj_{G'}(P_\phi)=\theta_H(P_\phi)=\theta_{G^*}(P_\phi)$, so $q=0$. Both are analytic, so by Hadamard's lemma $h=\lambda_e\tilde h$ and $q=\lambda_e\tilde q$ with analytic $\tilde h,\tilde q$ and $\tilde h(\phi_0)=\partial_eh(\phi_0)\ne0$. On a small neighbourhood $\inf\|\tilde h\|_\infty>0$, and $|q|/\|h\|_\infty=|\tilde q|/\|\tilde h\|_\infty$ off the hyperplane.
\end{proof}
If instead $\partial_eh$ vanished identically on the hyperplane, a second application of Hadamard's lemma would give $h=\lambda_e^2\hat h$, and if moreover $\partial_eq\not\equiv0$ the error bound would fail at generic points; this is the second-order near-equivalence of \citet{evans2020model}. \Cref{thm:E} shows that it does not occur: $\partial_eh$ is not identically zero on the hyperplane, so its zero set there is a null set, which is excluded through $\cN$.

\subsection{The benign side: a bound with a local constant}\label{app:benign}

\begin{lemma}[Bounded price with a local constant]\label{lem:benign-local}
Assume \eqref{eq:K2} and \eqref{eq:K1p}, and let $\cN$ be the set defined in the next paragraph. For every $\phi_0\in\Theta_H\setminus\cN$ and every fixed $k>0$ there are $r_0>0$ and $n_1$, depending on $P_0$ and $k$, such that the procedure \DSci\ of \cref{prop:CI} satisfies $\Pi^{\mathrm{loc}}_n(C^{\mathrm{CI}}_n;r_0)\le C'$ for $n\ge n_1$, where $C'$ depends on $p$, $|\cG|$, $\sigu/\sigl$, $\alpha$, $\gamma_1$, $\gamma_2$, $N_C$ and the local constant $C_4(P_0)$ of \cref{lem:K4E}, but not on $k$.
\end{lemma}
The set $\cN$ is the exceptional set of \cref{thm:K}; the rest of this subsection proves \cref{lem:benign-local}, and \cref{app:K} refines its case (d) to obtain \cref{prop:CI}.

\paragraph{The exceptional set.} For a DAG $D$ let $F_D(\lambda,\omega)=(I-L_D)^{-1}\mathrm{diag}(\omega)(I-L_D)^{-\top}$; it is polynomial, and its inverse on $\cM_D$, given by the node-wise regressions, is rational without poles. For $D$ equivalent to $H$ (in particular $D=G^*-e$ for $G^*\in\Xeq$), $T_{H\to D}=F_D^{-1}\circ F_H$ is an analytic diffeomorphism $\Theta_H\to\Theta_D$, and $\Theta_D$ is identified with $\Theta_{G^*}\cap\{\lambda_e=0\}$. For each ordered non-equivalent pair $(G^*,G')$ in $\Xeq$ let $a_{G^*,G'}(\phi)=\|\partial_eh_{G'}(T_{H\to G^*-e}(\phi),0)\|_2^2$, which is analytic and, by \cref{thm:E} and surjectivity of $T_{H\to G^*-e}$, not identically zero. (It can vanish at non-generic points, for example at the point of total independence, so $\cN$ keeps these factors.) For each triple $(i,j,S)$ that is not d-separated in $H$, clearing denominators in the Gaussian partial covariance gives a non-zero polynomial $b_{ijS}$. Let $\cN=\{\phi:\prod b_{ijS}(\phi)\prod a_{G^*,G'}(\phi)=0\}$. As $\Theta_H$ is connected and the product is analytic and not identically zero, $\cN$ is a relatively closed analytic subset of Lebesgue measure zero \citep{mityagin2020zero}. For $\phi_0\notin\cN$, $P_0$ is faithful to $H$ and \cref{lem:K4E} applies to every pair.

\paragraph{Step 0: true graphs near $P_0$ lie in $\Xeq$.} Suppose $(G^*,P_m)\in\cM_k$ with $P_m\to P_0$. Since $\cM_{G^*}$ is relatively closed in the positive definite cone and the $k$-TF inequalities are closed conditions, $P_0\in\cM_{G^*}$ and $(G^*,P_0)$ satisfies $k$-TF. By \cref{lem:HC}, every $G'\in\Xcanc$ violates $k$-TF at $P_0$ for every $k>0$, so $G^*\in\Xeq$. As $\cG$ is finite, there is $r_0$ such that every $(G^*,P)\in\cM_k$ with $\|\Sigma_P-\Sigma_{P_0}\|\le r_0$ has $G^*\in\Xeq$. If $\Xcanc=\emptyset$ this step is void.

\paragraph{Step 1: length decomposition.} If the output is non-empty, $\Leb(C_n)\le\sum_{G'}\Leb(W_{G'})+2\sum_{G'}|\hat\theta_{G'}-\theta_{G^*}(P)|\ind\{G'\text{ kept}\}$, and
\[
\E\bigl[|\hat\theta_{G'}-\theta_{G^*}|\ind\{G'\text{ kept}\}\bigr]\le\E|\hat\theta_{G'}-\thadj_{G'}(P)|+\Delta_{G'}(P)\,\Prob(G'\text{ kept}),\qquad\Delta_{G'}=|\thadj_{G'}-\theta_{G^*}|.
\]
By (M1), (M2) the first and the interval terms are $O(n^{-1/2})$ uniformly over all $P$ with bounded eigenvalues, whether or not $P\in\cM_{G'}$; $\Delta_{G'}\le2\sigu/\sigl$ because a population regression coefficient is at most the ratio of residual standard deviations. By \cref{lem:oracle} the oracle length is at least $2z_{1-\alpha/2}(\sigl/\sigu)n^{-1/2}$.

\paragraph{Step 2: four cases.} (a) $G'\notin\cX$: some local Markov partial correlation of $G'$ is at least $d_0/2>0$ near $P_0$, and by \cref{lem:T} the probability of keeping $G'$ is exponentially small. (b) $G'\in\Xcanc$: by \cref{lem:HC} and continuity, some $k$-TF element of $G'$ has violation margin at least $m_0/2$ near $P_0$, and \cref{lem:B} makes the probability of keeping $G'$ at most $16e^{-c_Bnm_0^2/4}$ for $n\ge n_B(m_0/2)$. (c) $G'\in\Xeq$ equivalent to $G^*$: by \eqref{eq:K1p}, $\Delta_{G'}\equiv0$ on $\cM_{G^*}$. A candidate equivalent to $G^*$ but not in $\Xeq$ is in $\Xcanc$, because $P_0\in\cM_H\subseteq\cM_{G^*}=\cM_{G'}$ puts it in $\cX$; it is handled by (b). (d) $G'\in\Xeq$ not equivalent to $G^*$: by \cref{lem:K4E}, $\Delta_{G'}\le C_4\delta$ with $\delta=\|h_{G'}\|_\infty$; applying \cref{lem:T} to a constraint attaining $\delta$ and using $\sup_{x\ge0}x\Phibar(x-2c)\le2c+1$,
\[
\Delta_{G'}\Prob(G'\text{ kept})\le C_4\Bigl[(2\bar c+1)\sqrt{2/(m+1)}+\delta_0\bigl(e^{-(m+1)/16}+e^{-m/2}\bigr)\Bigr]=O(n^{-1/2}),
\]
with $\bar c$ the largest critical value. Summing over the finitely many candidates proves the lemma; the exponential terms are absorbed for $n\ge n_1$, and only they and $r_0$ depend on $k$.

\subsection{Proof of \texorpdfstring{\cref{thm:K}}{the price-constant theorem} and \texorpdfstring{\cref{thm:K-lambda}}{the theorem near a malignant stratum}}\label{app:K}

Throughout, $g_\Sigma(\Delta_1,\Delta_2)=\frac12\mathrm{tr}(\Sigma^{-1}\Delta_1\Sigma^{-1}\Delta_2)$ is the Fisher metric of centred Gaussian laws and $\|\Delta\|_g^2=g_\Sigma(\Delta,\Delta)$; means play no role and are omitted.

\begin{lemma}[Projection onto a DAG model]\label{lem:Dproj}
Let $P=N(0,\Sigma)$ be positive definite and $G'$ a DAG with a topological order and predecessor sets $\mathrm{pred}(j)$. The minimizer $\Pi_{G'}P$ of $\KL(P\|Q)$ over $Q\in\cM_{G'}$ is $\prod_jP(x_j\mid x_{\pa_{G'}(j)})$, and
\[
D_{G'}(P)=\tfrac12\sum_j\log\frac{\sigma^2_{j\mid\pa_{G'}(j)}}{\sigma^2_{j\mid\mathrm{pred}(j)}}=-\tfrac12\sum_j\log(1-R_j^2),
\]
where $R_j^2$ is the squared multiple partial correlation of $j$ with $\mathrm{rest}(j)=\mathrm{pred}(j)\setminus\pa_{G'}(j)$ given $\pa_{G'}(j)$. If $P_\varepsilon$ is a smooth curve with $P_0\in\cM_{G'}$ and $c_j(\varepsilon)=\Cov_{P_\varepsilon}(j,\mathrm{rest}(j)\mid\pa_{G'}(j))$, then $\frac{d^2}{d\varepsilon^2}D_{G'}(P_\varepsilon)|_0=\sum_j\dot c_j^\top M_j^{-1}\dot c_j/\sigma^2_{j\mid\pa_{G'}(j)}$ with $M_j=\Cov_{P_0}(\mathrm{rest}(j)\mid\pa_{G'}(j))$; this equals $\|v^\perp\|_g^2$, where $v^\perp$ is the $g$-orthogonal component of $v=\dot\Sigma(0)$ normal to $\cM_{G'}$ at $P_0$, and $\frac{d}{d\varepsilon}\Pi_{G'}P_\varepsilon|_0$ is the $g$-orthogonal projection $v^\parallel$ of $v$ onto the tangent space.
\end{lemma}
\begin{proof}
For $Q\in\cM_{G'}$ and $Q^*=\prod_jP(x_j\mid x_{\pa(j)})\in\cM_{G'}$, the chain rule gives $\KL(P\|Q)=\KL(P\|Q^*)+\sum_j\E_P\KL(P_{j\mid\pa(j)}\|Q_{j\mid\pa(j)})$, so $Q^*$ is the minimizer. $\KL(P\|Q^*)=\sum_jH(x_j\mid x_{\pa(j)})-H(P)$, and $\log\det\Sigma=\sum_j\log\sigma^2_{j\mid\mathrm{pred}(j)}$ for every order; the ratio of the two residual variances is $1-R_j^2$. Along the curve, $c_j(0)=0$ because $\mathrm{rest}(j)$ consists of non-descendants of $j$, and $R_j^2=\varepsilon^2\dot c_j^\top M_j^{-1}\dot c_j/\sigma^2_{j\mid\pa(j)}+O(\varepsilon^3)$; with $-\frac12\log(1-x)=\frac x2+O(x^2)$ this gives the second derivative.

For the geometric statements, the structural parametrization $\Theta_{G'}\to\cM_{G'}$ is injective with a rational inverse without poles (the node-wise regressions), so $\cM_{G'}$ is an embedded submanifold of the positive definite cone, and $P\mapsto\Pi_{G'}P$ is a rational map, hence smooth. Write $F(\Sigma,\Sigma')=\KL(N(0,\Sigma)\|N(0,\Sigma'))$. On the diagonal, $\partial_{\Sigma'}F=0$, $\partial^2_{\Sigma'\Sigma'}F=g$ and $\partial^2_{\Sigma'\Sigma}F=-g$; since $F$ vanishes on the diagonal, also $\partial^2_{\Sigma\Sigma}F=g$, so $F(\Sigma+\Delta_1,\Sigma+\Delta_2)=\frac12\|\Delta_1-\Delta_2\|_g^2+O(\|\Delta\|^3)$. For every $\varepsilon$, $\Pi_{G'}P_\varepsilon$ minimizes $F(\Sigma(\varepsilon),\cdot)$ over $\cM_{G'}$, so $\partial_{\Sigma'}F(\Sigma(\varepsilon),\Sigma_{\Pi_{G'}P_\varepsilon})[w]=0$ for every tangent vector $w$ at $\Pi_{G'}P_\varepsilon$. Differentiating at $\varepsilon=0$, where $\Pi_{G'}P_0=P_0$ and the term from the motion of the tangent space vanishes because $\partial_{\Sigma'}F=0$ on the diagonal, gives $g(\Pi'-v,w)=0$ for all $w\in T_{P_0}\cM_{G'}$. As $\Pi'\in T_{P_0}\cM_{G'}$, $\Pi'=v^\parallel$. Finally $D_{G'}(P_\varepsilon)=F(\Sigma(\varepsilon),\Sigma_{\Pi_{G'}P_\varepsilon})=\frac12\varepsilon^2\|v-v^\parallel\|_g^2+O(\varepsilon^3)$, so the second derivative is $\|v^\perp\|_g^2$.
\end{proof}

\begin{lemma}[Projection distance and local Markov tests]\label{lem:Dh}
Let the eigenvalues of $\Sigma_P$ lie in $[\sigl^2,\sigu^2]$ and let $\rho_{\max}$ be the largest absolute local Markov partial correlation $\rho_P(j,i\mid\pa_{G'}(j))$ of $G'$. Then $\frac12\rho_{\max}^2\le D_{G'}(P)$, and if all $R_j^2\le1/2$ for some topological order of $G'$, then $D_{G'}(P)\le p^2(\sigu/\sigl)^2\rho_{\max}^2$.
\end{lemma}
\begin{proof}
For a constraint $(j,i,\pa(j))$, $i$ is a non-descendant of $j$, so there is a topological order in which $i$ precedes $j$; in it, $\rho^2\le R_j^2$, and $D\ge-\frac12\log(1-R_j^2)\ge\frac12R_j^2$ because every term of $D$ is non-negative and $D$ does not depend on the order. For the upper bound use the order of the assumption, write $c=\Cov(j,\mathrm{rest}(j)\mid\pa(j))$ and $M=\Cov(\mathrm{rest}(j)\mid\pa(j))$, so that $c_i=\rho_i(\sigma^2_{j\mid\pa(j)}M_{ii})^{1/2}$ with $\rho_i$ a local Markov partial correlation. Since $M$ is a Schur complement of a principal submatrix of $\Sigma$, its eigenvalues lie in $[\sigl^2,\sigu^2]$, and $R_j^2=c^\top M^{-1}c/\sigma^2_{j\mid\pa(j)}\le|\mathrm{rest}(j)|(\sigu/\sigl)^2\rho_{\max}^2$. If $R_j^2\le1/2$ then $-\frac12\log(1-R_j^2)\le R_j^2$; sum over $j$.
\end{proof}

\begin{lemma}[Two discrepancies]\label{lem:Kadj}
Let $\delta^{\mathrm{adj}}_{G^*,G'}=\frac{d}{d\varepsilon}[\thadj_{G'}(P_\varepsilon)-\theta_{G^*}(P_\varepsilon)]_{\varepsilon=0}$. Then $|\delta^{\mathrm{adj}}_{G^*,G'}-\delta_{G^*,G'}|\le C_\theta\sqrt{I_{G^*,G'}}$ with $C_\theta=2(\sigu/\sigl)^2(1+\sigu/\sigl)$.
\end{lemma}
\begin{proof}
On $\cM_{G'}$, $\thadj_{G'}=\theta_{G'}$, so the difference of the two discrepancies is the derivative of $\thadj_{G'}(P_\varepsilon)-\thadj_{G'}(\Pi_{G'}P_\varepsilon)$, which by \cref{lem:Dproj} is $d\thadj_{G'}[v^\perp]$. Writing $\thadj_{G'}=e_X^\top\Sigma_{AA}^{-1}\Sigma_{AY}$ with $\beta_A=\Sigma_{AA}^{-1}\Sigma_{AY}$, $d\thadj_{G'}[\Delta]=e_X^\top\Sigma_{AA}^{-1}(\Delta_{AY}-\Delta_{AA}\beta_A)$, so $|d\thadj_{G'}[\Delta]|\le\sigl^{-2}(1+\sigu/\sigl)\|\Delta\|_{\mathrm{Frob}}$, using $|\beta_A|\le\sigu/\sigl$ and that $\Delta_{AY}$ and $\Delta_{AA}$ are disjoint blocks of $\Delta$. Finally $\|\Delta\|_{\mathrm{Frob}}\le\sqrt2\sigu^2\|\Delta\|_g$ and $\|v^\perp\|_g=\sqrt{I}$, so $|\delta^{\mathrm{adj}}-\delta|\le\sqrt2(\sigu/\sigl)^2(1+\sigu/\sigl)\sqrt I\le C_\theta\sqrt I$.
\end{proof}

\begin{proof}[Proof of \cref{thm:K}(i)]
Write $a_\alpha=z_{1-\alpha}$, $z=z_{1-\alpha/2}$, $J=a_\alpha\Phi(a_\alpha)+\varphi(a_\alpha)=\int_0^\infty\Phi(a_\alpha-u)\,du$, and $m_d=1-\alpha-\Phi(d-a_\alpha)$ for $d>0$. Define
\[
\Psi_\alpha(K)=\sup_{0<d<2a_\alpha}\frac1{2z}\Bigl[2J+\int_0^{dK}\bigl(m_d-\Phi(a_\alpha-u)\bigr)_+\,du\Bigr],
\]
so that $\Psi_\alpha(0)=J/z=c_P$ and $\Psi_\alpha$ is non-decreasing. Since $\int_0^{dK}(m_d-\Phi(a_\alpha-u))_+du\ge dKm_d-J$, we have $\Psi_\alpha(K)\ge\frac12c_P+c_KK$ with $c_K=\sup_ddm_d/(2z)$. (At $\alpha=0.05$: $c_P=0.8498863$, $c_K=0.19526$; script \texttt{d26\_lower\_bound\_psi.py}.)

\emph{Local asymptotic normality.} The model of centred Gaussian laws with unrestricted covariance matrix is a regular exponential family of fixed dimension, hence locally asymptotically normal at $P_0$ with Fisher metric $g$. The three laws that occur below have the same mean (the projection $\Pi_{G'}$ keeps the mean, because the regressions have intercepts), the targets depend on $\Sigma$ only, and the Fisher metric is block diagonal in the mean and the covariance, so the means play no role and the bounds below hold for all tests, including those that do not know the mean. For local parameters $h_n,h_n'$ that converge, with $\|h_n-h_n'\|_g\to d$, the most powerful level-$\alpha'$ test of $P_{h_n'}^n$ against $P_{h_n}^n$ has power converging to $\Phi(d-z_{1-\alpha'})$, and this limit is continuous in $\alpha'$. Asymptotic honesty is uniform over $\cM_k$ (\cref{sec:setup}), so the coverage of every pair in the class is at least $1-\alpha_n$ with $\limsup\alpha_n\le\alpha$; with $\alpha_n'=\max(\alpha_n,\alpha)\to\alpha$, honesty holds at level $\alpha_n'$, which is used in place of $\alpha$ below.

\emph{Pratt's integral.} Let $(G^*,P)\in\cM_k$ be near $P_0$ and let $C_n$ be an honest interval. Since $\E_P\Leb(C_n)=\int\Prob_P(\vartheta\in C_n)\,d\vartheta$, it suffices to bound $\Prob_P(\vartheta\in C_n)$ from below for each $\vartheta$.

\emph{(1) Alternatives inside $\cM_{G^*}$.} Let $P^{(t)}$, $|t|\le Ms_{\mathrm{eff}}n^{-1/2}$ for every fixed $M$, be a one-dimensional least favourable submodel of $\cM_{G^*}$ through $P$ with $\theta_{G^*}(P^{(t)})=\theta_{G^*}(P)+t$ exactly and Fisher information $s_{\mathrm{eff}}^{-2}$ at $t=0$; it exists because $d\theta_{G^*}\ne0$, and by the argument below for class membership $(G^*,P^{(t)})\in\cM_k$ for $n$ large. Consider the test of the null hypothesis $P^{(t)n}$ against the alternative $P^n$ that rejects when $\theta_{G^*}(P)+t\notin C_n$; its size under $P^{(t)n}$ is the non-coverage of the true target, at most $\alpha_n'$, and its power is $\Prob_P(\theta_{G^*}(P)+t\notin C_n)$. Hence $\Prob_P(\theta_{G^*}(P)+t\in C_n)\ge1-\beta_n(t)$, with $\beta_n(t)$ the maximal power at level $\alpha_n'$ of this testing problem. For $t=us_{\mathrm{eff}}n^{-1/2}$, $\beta_n\to\Phi(|u|-a_\alpha)$, so the bound tends to $\Phi(a_\alpha-|u|)$ for every $u$ (a one-sided Neyman--Pearson test for each sign).

\emph{(2) The alternative $\Pi_{G'}P$.} Take $P=P_a$ with $a=d/\sqrt{nI}$ and $Q=\Pi_{G'}P_a\in\cM_{G'}$. The local parameters of $P_a$ and $Q$ differ by $\sqrt n\,a\,v^\perp+o(1)$ by \cref{lem:Dproj}, and $\|v^\perp\|_g=\sqrt I$, so their $g$-distance converges to $d$. Put $\theta_1=\theta_{G^*}(P_a)$, $\theta_2=\theta_{G'}(Q)$ and take $\theta_2>\theta_1$ (the other case is symmetric); then $\theta_2-\theta_1=|\delta|a+O(a^2)$ because the discrepancy vanishes at $a=0$ by \eqref{eq:K2}. The event $\{\theta_2\notin C_n\}$ is a test of $Q^n$ against $P_a^n$ of size at most $\alpha_n'$, so $\Prob_{P_a}(\theta_2\in C_n)\ge1-\beta_n'$ with $\beta_n'\to\Phi(d-a_\alpha)$. For $\vartheta\in[\theta_1,\theta_2]$, an interval that covers both $\theta_1$ and $\theta_2$ covers $\vartheta$ (this is where $C_n$ being an interval is used; for sets that may be disconnected the bound fails, \cref{prop:union}), so
\[
\Prob_{P_a}(\vartheta\in C_n)\ \ge\ 1-\alpha_n-\beta_n'\ \to\ m_d .
\]

\emph{Oracle mismatch.} For $D,D'\in\Xeq$ the pairs $(D,P_0)$ and $(D',P_0)$ belong to $\cM_k$ for $k<k^*(P_0)$, by the first paragraph on class membership. Step (1) applied to $P=P_0$ in the model $\cM_D$ gives $\liminf_n\sqrt n\,\E_{P_0}\Leb(C_n)\ge2J\,s_{\mathrm{eff}}(D,P_0)=c_P\cdot2z\,s_{\mathrm{eff}}(D,P_0)$, and the price at $(D',P_0)$ divides by $\orl_n(D',P_0)=2z\,s_{\mathrm{eff}}(D',P_0)n^{-1/2}$. Maximizing over $D,D'$ gives $\liminf_n\Pi^{\mathrm{loc}}_n\ge c_P\mathcal R(P_0)$.

\emph{Integration.} Put $\vartheta=\theta_1+us_{\mathrm{eff}}n^{-1/2}$. By (1) and (2), $\Prob_{P_a}(\vartheta\in C_n)\ge\Phi(a_\alpha-|u|)-o(1)$ for every $u$, and for $0\le u\le dK$ also $\ge m_d-o(1)$, where $dK$ is the length of $[\theta_1,\theta_2]$ in units of $s_{\mathrm{eff}}n^{-1/2}$, up to a factor $1+o(1)$ by continuity of $s_{\mathrm{eff}}$. The integrand is non-negative and its pointwise limit inferior is at least $g(u)=\Phi(a_\alpha-|u|)$ for $u<0$ and $u>dK$ and $\max\{\Phi(a_\alpha-u),m_d\}$ for $0\le u\le dK$. By Fatou's lemma,
\[
\liminf_n\frac{\E_{P_a}\Leb(C_n)}{2z\,s_{\mathrm{eff}}n^{-1/2}}\ \ge\ \frac1{2z}\int g(u)\,du=\frac1{2z}\Bigl[2J+\int_0^{dK}\bigl(m_d-\Phi(a_\alpha-u)\bigr)_+du\Bigr].
\]
Since $\orl_n(G^*,P_a)=2z\,s_{\mathrm{eff}}(G^*,P_a)n^{-1/2}$ with $s_{\mathrm{eff}}$ continuous, this is a lower bound for the local price, and taking the supremum over $d\in(0,2a_\alpha)$ gives $\Psi_\alpha(K)$, with $K=|\delta|/(s_{\mathrm{eff}}\sqrt I)$ for the pair attaining $K(P_0)$. For $K=0$ the bound is $c_P$ (the usual Pratt bound), and the maximum of the two lower bounds in the theorem follows from $\Psi_\alpha(K)\ge c_P$ and $\Psi_\alpha(K)\ge\frac12c_P+c_KK$. Moreover $m_d=\Phi(a_\alpha-d)-\alpha<\Phi(a_\alpha-d)\le\Phi(a_\alpha-u)$ for $u\le d$, so the integral over $[0,dK]$ vanishes when $K\le1$ and $\Psi_\alpha(K)=c_P$ exactly for $K\le1$; numerically this holds up to $K\approx1.07$, and $\frac12c_P+c_KK$ exceeds $c_P$ only for $K>2.18$. Intuitively, an alternative inside $\cM_{G^*}$ at $g$-distance $d$ already gives $\Phi(a_\alpha-d)$, which exceeds what the alternative in $\cM_{G'}$ gives on the segment of length $dK\le d$.

\emph{Class membership.} Let $k<k^*(P_0)$ and $D\in\Xeq$. Triangle edges $u\to v\ne e$ of $D$ satisfy the $k$-TF inequality at $P_0$ with room to spare, and by continuity also at every $P\in\cM_D$ near $P_0$. For the edge $e$ itself, in a triangle that contains it, \cref{lem:M1} gives $|\rho_{P_0}(u,v\mid W)|>0$ while $\lambda^{\mathrm{std}}_e=0$ at $P_0$, so the inequality also holds near $P_0$. Since the eigenvalues of $\Sigma_{P_0}$ are interior, $(D,P)\in\cM_k$ for every $P\in\cM_D$ close enough to $P_0$.
\emph{Membership of the points.} By the previous paragraph, $P_a$, $P_a^{(t)}$ with $|t|\le Ms_{\mathrm{eff}}n^{-1/2}$ for every fixed $M$, and $Q=\Pi_{G'}P_a$ (with $G'$ in place of $G^*$) are within $O(n^{-1/2})$ of $P_0$ in the parameter of the corresponding model, so the pairs $(G^*,P_a)$, $(G^*,P_a^{(t)})$ and $(G',Q)$ belong to $\cM_k$ for $n$ large.
\end{proof}

\paragraph{A bound on a fixed neighbourhood.} Let \DSci\ be the procedure of \cref{def:ds} with $\gamma_3=0$, the parent adjustment rule $A(D)=\pa_D(X)$, and step (1) replaced by the \emph{conditional independence screen}: reject $D$ if the $t$-test of some local Markov constraint $(j,i,\pa_D(j))$ of $D$, $i$ a non-descendant non-parent of $j$, rejects at level $\gamma_1/N_C$, where $N_C$ is the total number of such constraints over $\cG$; \DSciO\ is the same with the optimal adjustment set. \DSci\ is honest at level $1-\alpha$ over $\cM_k$ for every $k>0$ and $n\ge p+4$ (proof of \cref{thm:coverage}). Unlike \cref{thm:K}(ii), the following bound holds on a fixed neighbourhood and for all large $n$, with a constant that does not depend on $P_0$ or $k$, at the price of a worse constant.

\begin{proposition}[Conditional independence screen]\label{prop:CI}
Assume \eqref{eq:K2} and \eqref{eq:K1p} and let $\phi_0\in\Theta_H\setminus\cN$. For every fixed $k>0$ there are $r_0>0$ and $n_1$, depending on $P_0$ and $k$, such that $\Pi^{\mathrm{loc}}_n(C^{\mathrm{CI}}_n;r_0)\le C\,(1+K(P_0))$ for $n\ge n_1$, where $C$ depends only on $(p,\sigu/\sigl,\alpha,\gamma_1,\gamma_2,|\cG|)$, not on $P_0$ or $k$. The radius $r_0$, the threshold $n_1$ and the exponential rate at which cancellation explanations are rejected depend on $k$.
\end{proposition}

\begin{proof}[Proof of \cref{prop:CI}]
We follow the proof of \cref{lem:benign-local}, with the same $r_0$ and $n_1$ up to a further reduction of $r_0$. The lengths of the candidate intervals and the estimation errors contribute at most $C_1(p,\sigu/\sigl,\alpha')|\cG|n^{-1/2}$, and $\orl_n\ge2z_{1-\alpha/2}(\sigl/\sigu)n^{-1/2}$ by \cref{lem:oracle}; cases (a) and (b) contribute exponentially small terms and case (c) nothing. In case (d), on $\cM_{G^*}$ near $P_0$, $\thadj_{G'}-\theta_{G^*}=\lambda_e\tilde q$ and, because $D_{G'}\ge0$ is smooth and vanishes on $\{\lambda_e=0\}$, $D_{G'}=\lambda_e^2\tilde g$ with $\tilde g(\phi_0)=I/2>0$. Hence the ratio $|\thadj_{G'}-\theta_{G^*}|/\sqrt{2D_{G'}}=|\tilde q|/\sqrt{2\tilde g}$ is continuous near $P_0$ with value $|\delta^{\mathrm{adj}}|/\sqrt I$ at $P_0$, and for $r_0$ small it is at most $|\delta^{\mathrm{adj}}|/\sqrt I+1$ on the $r_0$-ball, where also $R_j^2\le1/2$ for a fixed topological order of each candidate. (The additive $1$ covers base points at which $\delta^{\mathrm{adj}}=0$ but $\tilde q\not\equiv0$.) By \cref{lem:Dh} some local Markov constraint of $G'$ has $|\rho|\ge\sqrt{D_{G'}}/(p\sigu/\sigl)$, and \cref{lem:T} with $\sup_{x\ge0}x\Phibar(x-2\bar c)\le2\bar c+1$ gives
\begin{align*}
\Delta_{G'}\Prob(G'\text{ kept})&\le\Bigl(\frac{|\delta^{\mathrm{adj}}|}{\sqrt I}+1\Bigr)\sqrt{2D}\,\Phibar\Bigl(\frac{\sqrt{D}\sqrt{(m+1)/2}}{p\sigu/\sigl}-2\bar c\Bigr)+\text{exp.}\\
&\le\frac{2p(\sigu/\sigl)(2\bar c+1)}{\sqrt{m+1}}\Bigl(\frac{|\delta^{\mathrm{adj}}|}{\sqrt I}+1\Bigr)+\text{exp.}
\end{align*}
Dividing by $\orl_n$ and summing over at most $|\cG|$ candidates bounds the contribution by $C_2|\cG|$ times the maximum of $|\delta^{\mathrm{adj}}|/(s_{\mathrm{eff}}\sqrt I)+1/s_{\mathrm{eff}}$ over the non-equivalent ordered pairs $(G^*,G')$ in $\Xeq$, where $C_2$ depends on $p$, $\sigu/\sigl$ and $\bar c$. The screening tests have at least $n-p\ge p+2$ degrees of freedom once $n_1\ge2p+2$, $t$ quantiles decrease in the degrees of freedom, and $N_C\le|\cG|p^2$; so $\bar c$ is at most the $1-\gamma_1/(2|\cG|p^2)$ quantile of $t_{p+2}$, which depends only on $\gamma_1$, $|\cG|$ and $p$. By \cref{lem:Kadj} and $s_{\mathrm{eff}}\ge\sigl/\sigu$, the maximum is at most $(1+o(1))K(P_0)+(C_\theta+1)\sigu/\sigl$ as $r_0\to0$, where $1+o(1)$ is the ratio of $s_{\mathrm{eff}}$ at $P_0$ and at the points of the ball.

\emph{Invariance.} The same factorizations show that $K_{G^*,G'}(P_0)$ does not depend on the curve: along a smooth curve $\gamma(t)$ in $\Theta_{G^*}$ with $\gamma(0)=\phi_0^{G^*}$ and $\frac{d}{dt}\lambda_e(\gamma(t))|_0=c\ne0$, the discrepancy has derivative $c\,\tilde q_{\mathrm{int}}(\phi_0^{G^*})$ and $D_{G'}$ has second derivative $2c^2\tilde g(\phi_0^{G^*})$, where $\theta_{G'}\circ\Pi_{G'}-\theta_{G^*}=\lambda_e\tilde q_{\mathrm{int}}$ by \eqref{eq:K2}, so the ratio is the same as for $\gamma(t)=\phi_0^{G^*}+te_e$.
\end{proof}

\paragraph{Adjustment rules and screens without estimated directions.}

\begin{remark}[Choice of adjustment set]\label{rem:Oset}
For a rule $A$ of \cref{def:ds}, let $\delta^A_{G^*,G'}$ be the first-order discrepancy of the $A$-adjustment estimand, $\frac{d}{d\varepsilon}[\theta^A_{G'}(P_\varepsilon)-\theta_{G^*}(P_\varepsilon)]_{\varepsilon=0}$, where $\theta^A_D(\Sigma)$ is the coefficient of $X$ in the population regression of $Y$ on $X$ and $A(D)$; \cref{lem:Kadj} applies to it. Let $K^A(P_0)=\max_{(G^*,G')}|\delta^A_{G^*,G'}|/(s_{\mathrm{eff}}(G^*,P_0)\sqrt{I_{G^*,G'}})$, the maximum being over all ordered pairs as for $K$, and $K^{\mathrm{adj}}$, $K^O$ the cases $A(D)=\pa_D(X)$ and $A(D)=O(X,Y,D)$. For a Markov class $S'$ of $\Xeq$ not containing $G^*$ let $w^A_{G^*,S'}$ be the length of the smallest interval that contains $0$ and every $\delta^A_{G^*,G'}$, $G'\in S'$, and
\[
K^A_+(P_0)=\max_{G^*\in\Xeq}\sum_{S'}\frac{w^A_{G^*,S'}}{s_{\mathrm{eff}}(G^*,P_0)\sqrt{I_{G^*,S'}}},\qquad \mathcal R^A(P_0)=\max_{D,G^*\in\Xeq}\frac{s^A_D(P_0)}{s_{\mathrm{eff}}(G^*,P_0)},
\]
where $s^A_D=\sigma_{Y\mid X,A(D)}/\sigma_{X\mid A(D)}$ is the asymptotic standard deviation of the $A(D)$-adjustment estimator. Then $K^A_+\ge K^A$ and $\mathcal R^A\ge\mathcal R$, with $K^A_+=K^A$ when $\Xeq$ has two Markov classes and, for each $G^*$, the numbers $\delta^A_{G^*,G'}$, $G'$ in the other class, have one sign. The proof of \cref{prop:CI} uses only that $\pa_D(X)$ is a valid adjustment set in $D$ that does not contain $Y$, so it applies to every rule $A$, with $K^A$ in place of $K$. The difference $\delta^A-\delta$ is the derivative of the estimand in a direction normal to $\cM_{G'}$, so it vanishes when the $A$-adjustment estimator is efficient in $\cM_{G'}$; then $\mathcal R^A=\mathcal R$ and $K^A_+=K_+$. The optimal adjustment set $O(X,Y,D)$ of \citet{henckel2022graphical} gives the most efficient adjustment estimator (see also \citealp{rotnitzky2020efficient,witte2020efficient}), which is not efficient among all regular estimators in general \citep{guo2022efficient}. In the five-node enumeration (\cref{app:certify}), $K^O=K$ for $186$ of the $240$ pairs $(G,e)$ with $K>0$, and $K^O>0$ only when $K>0$, whereas $K^{\mathrm{adj}}>0$ although $K=0$ for $54$ pairs.
\end{remark}

Let \DSlrO\ and \DSlreff\ be \DSOp\ and \DSeff\ with $\gamma_3=0$, that is, without steps (0) and (1b), and write $C^{\mathrm{LR}}_n$ and $C^{\mathrm{LR,eff}}_n$ for their intervals ($\DS^{\mathrm{LR}}_A$ for a general rule $A$). Let \DSdir\ be \DSeff\ with step (0) replaced by the true $\Xeq$, its classes and the directions $u_{S^*,S'}$ at $P_0$, with $B_{S'}$ an orthonormal basis of their span and $d_{S'}=\dim$, with the directed test of step (1b) computed on all data, and with a level $\gamma_{1,n}\to0$, $n^{-1}\log(1/\gamma_{1,n})\to0$ in step (1a) and $\alpha'_n=\alpha-\gamma_{1,n}-\gamma_2-\gamma_3\to\alpha'=\alpha-\gamma_2-\gamma_3$; it depends on $P_0$ and serves only as a step in the proof of \cref{thm:K}(ii). Write $C^{P_0}_n$ for its interval.

\begin{proposition}[Screens with known or without directions]\label{thm:explicit}
Assume \eqref{eq:K2} and \eqref{eq:K1p}, let $\phi_0\in\Theta_H\setminus\cN$ and $k>0$; if $\Xeq$ has at least two Markov classes then $q\ge1$, and otherwise the terms with $K$ below are absent.
\begin{enumerate}
\item[(a)] $\displaystyle\lim_{r\downarrow0}\limsup_{n\to\infty}\Pi^{\mathrm{loc}}_n(C^{\mathrm{LR}}_n;r)\le\mathcal R^A(P_0)\,\frac{z_{1-\alpha'/2}+\sqrt{2\log|\Xeq|}}{z_{1-\alpha/2}}+\frac{\varpi_q(\gamma_1)}{2z_{1-\alpha/2}}\,K^A_+(P_0)$ for $\DS^{\mathrm{LR}}_A$.
\item[(b)] \DSlreff\ is asymptotically honest at level $1-\alpha$ over $\cM_k$, and its limiting local price is at most $\mathcal R(P_0)(z_{1-\alpha'/2}+e_m)/z_{1-\alpha/2}+\varpi_q(\gamma_1)K_+(P_0)/(2z_{1-\alpha/2})$, with $e_1=0$.
\item[(c)] Along the weak edge, $|u_{S^*,S'}|^2=I_{G^*,S'}$ for the member $G^*$ of $S^*$ that defines $u_{S^*,S'}$. \DSdir\ is asymptotically honest at level $1-\alpha$ over $\cM_k$, and its limiting local price is at most $\mathcal R(P_0)(z_{1-\alpha'/2}+e_m)/z_{1-\alpha/2}+\bar\varpi_{m-1}(\gamma_3)K_+(P_0)/(2z_{1-\alpha/2})$.
\item[(d)] Under the conditions of \cref{thm:K}(ii), \DSA\ has limiting local price at most $\mathcal R^A(P_0)(z_{1-\alpha'/2}+\sqrt{2\log|\Xeq|})/z_{1-\alpha/2}+\bar\varpi_{m-1}(\gamma_3)K^A_+(P_0)/(2z_{1-\alpha/2})$.
\end{enumerate}
\end{proposition}

At $\alpha=0.05$ and $\gamma_1=\gamma_2=0.01$ the slope of (a) and (b) is $0.367$, $0.493$ and $0.618$ for $q=1$, $5$ and $15$, against $c_K=0.195$, and the intercept coefficient of (a) is $1.708$ for $|\Xeq|=2$ and $2.245$ for $|\Xeq|=12$, the largest size of an $\Xeq$ with two or more Markov classes on five nodes; for (b) with $m=2$ the factor between the two sides of \cref{cor:RK} is $3.66$ for $q=1$ and $4.95$ for $q=15$ (\cref{app:certify}).

\begin{proof}[Proof of \cref{thm:explicit}]
(a) The statistic is $T_D=2nD_D(\hat P)$, with $\hat P$ the saturated maximum likelihood estimate, because the minimizer in \cref{lem:Dproj} is the maximum likelihood estimate of $\cM_D$. Suppose the bound fails. Then there are $\varepsilon>0$, $r_m\downarrow0$, $n_m\to\infty$ and $(G^*_m,P_m)\in\cM_k$ with $\|\Sigma_{P_m}-\Sigma_{P_0}\|\le r_m$ and $\E\Leb(C_{n_m})/\orl_{n_m}$ larger than the right-hand side plus $\varepsilon$. Passing to a subsequence, $G^*_m=G^*$, which lies in $\Xeq$ by Step~0 of \cref{app:benign}, $P_m=P(\phi_m)$ with $\phi_m\to\phi_0^{G^*}$ in $\Theta_{G^*}$, and $t_m=\sqrt{n_m}\,\lambda_e(\phi_m)\to c\in[-\infty,\infty]$. Let $\phi_m'$ be $\phi_m$ with $\lambda_e$ set to zero and $\Sigma'_m$ its covariance matrix; $\Sigma'_m\in\cM_H\subseteq\cM_D$ for every $D\in\Xeq$. Write $n$ for $n_m$.

\emph{Candidates outside $\Xeq$.} For $D\notin\cX$, $D_D>0$ near $P_0$, and the critical values converge to $\chi^2_{q,1-\gamma_1}$, so $D$ is kept with exponentially small probability; candidates in $\Xcanc$ are rejected by the $k$-triangle-faithfulness screen with exponentially small probability, as in Step~2 of \cref{app:benign}. These bounds hold uniformly over $P_m$ near $P_0$, by concentration of the sample covariance matrix. On the event that no such candidate is kept, the output is the hull of the kept candidates in $\Xeq$. On its complement, $\sqrt n\Leb(C_n)$ has second moment $O(n)$, because the estimands of candidates outside $\Xeq$ differ from $\theta_{G^*}$ by $O(1)$, while the complement has probability $O(e^{-c n})$ by \cref{lem:T,lem:B}; by Cauchy--Schwarz its contribution to $\sqrt n\,\E\Leb(C_n)$ tends to zero.

\emph{Local limit for finite $c$.} Let $\Delta_m=\sqrt n(\hat\Sigma-\Sigma_m)$. In coordinates that are orthonormal for $g$ at $P_0$, $\Delta_m$ converges in law to a standard normal vector $\Delta$, and $\Sigma_m=\Sigma'_m+n^{-1/2}t_mv_m+O(n^{-1})$ with $v_m\to v$, the tangent vector along the weak edge of \cref{lem:Dproj}. For a Markov class $S'\not\ni G^*$ the proof of \cref{lem:Dproj}, together with smoothness of $\Pi_{S'}$ on a neighbourhood of $P_0$ in the positive definite cone, gives $D_{S'}(\Sigma'+\Delta)=\frac12\|\Delta^{\perp}\|_g^2+O(\|\Delta\|^3)$ uniformly for $\Sigma'\in\cM_{S'}$ near $P_0$, with $\perp$ the normal component at $\Sigma'$, which depends continuously on $\Sigma'$. Hence $T_{S'}\to\|(cv+\Delta)^\perp\|_g^2$, which is $\chi^2_q(c^2I_{G^*,S'})$ because the normal space has dimension $q$ and $\|v^\perp\|_g^2=I_{G^*,S'}$. The null laws of \cref{thm:coverage}(i) converge to $\chi^2_q$, whose distribution function is continuous, so $\Prob(S'\text{ passes step (1)})\to\beta_q(|c|\sqrt{I_{G^*,S'}})$ with $\beta_q(x)=\Prob(\chi^2_q(x^2)\le\chi^2_{q,1-\gamma_1})$. Jointly, $\sqrt n(\hat\theta^A_D-\theta^A_D(\Sigma_m))\to Z_D=d\theta^A_D[\Delta]$ for $D\in\Xeq$, centred Gaussian with standard deviation $s^A_D(P_0)$, where $\theta^A_D(\Sigma)$ is the coefficient of $X$ in the population regression of $Y$ on $X$ and $A(D)$ and $\hat\theta^A_D$ its sample version. For $D\simeq G^*$ this estimand equals $\theta_{G^*}$ on $\cM_{G^*}$: by \eqref{eq:K1p} $\theta_D=\thadj_D=\theta_{G^*}$ there, and $A(D)$ is a valid adjustment set in $D$. For $D\in S'$, on $\{\lambda_e=0\}$ the law lies in $\cM_H\subseteq\cM_D$, where $\theta^A_D=\theta_D=\theta_H=\theta_{G^*}$ because $A(D)$ is valid and by \eqref{eq:K2}; so, as in \cref{lem:K4E}, the estimand differs from $\theta_{G^*}(P_m)$ by $n^{-1/2}t_m(\delta^A_{G^*,D}+o(1))$. The half-widths satisfy $\sqrt nh_D\to t s^A_D$ with $t=z_{1-\alpha'/2}$, and all these quantities have bounded second moments by \cref{lem:M}. Write $Z_{D,m}=\sqrt n(\hat\theta^A_D-\theta^A_D(\Sigma_m))$ and $b_{D,m}=\sqrt n(\theta^A_D(\Sigma_m)-\theta_{G^*}(P_m))$, so that $b_{D,m}=0$ for $D\simeq G^*$ and $b_{D,m}\to\delta^A_{G^*,D}c$ for $D\in S'$. If the kept set $\mathcal K\subseteq\Xeq$ is not empty, then, since step (2) can only remove candidates and the spread of $\{b_{D,m}:D\in\mathcal K\}$ is at most the sum over the classes $S'$ that pass step (1) of the spreads of $\{0\}\cup\{b_{D,m}:D\in S'\}$,
\[
\sqrt n\Leb(C_n)\ \le\ \max_{\Xeq}Z_{\cdot,m}-\min_{\Xeq}Z_{\cdot,m}+\sum_{S'\text{ passes}}\mathrm{spread}_{S'}(b_{\cdot,m})+2\max_{\Xeq}\sqrt nh_D .
\]
The right-hand side has bounded second moments and converges in law, so its expectation converges to that of its limit; only the marginal laws of $Z$ and of each $T_{S'}$ enter, because the spreads are deterministic and the expectation is linear. In the limit $Z_{\cdot,m}$, $b_{\cdot,m}$ and $\sqrt nh_D$ are replaced by $Z$, $\delta^Ac$ and $ts^A$ and the spread of $S'$ becomes $w^A_{G^*,S'}|c|$. The maximal inequality $\E\max_{i\le N}Z_i\le\max_i\mathrm{sd}(Z_i)\sqrt{2\log N}$ for centred Gaussian variables, applied to $Z$ and to $-Z$, gives
\[
\limsup_m\sqrt n\,\E\Leb(C_n)\le2\max_{\Xeq}s^A\bigl(t+\sqrt{2\log|\Xeq|}\bigr)+\sum_{S'}w^A_{G^*,S'}|c|\,\beta_q\bigl(|c|\sqrt{I_{G^*,S'}}\bigr),
\]
and dividing by $\sqrt n\,\orl_n\to2z_{1-\alpha/2}s_{\mathrm{eff}}(G^*,P_0)$ and using $x\beta_q(x)\le\varpi_q(\gamma_1)$ with $x=|c|\sqrt{I_{G^*,S'}}$ bounds the limit by the right-hand side, a contradiction.

\emph{Infinite $c$.} Since $\phi_m\to\phi_0^{G^*}$, $|t_m|\le\varepsilon\sqrt n$ for every $\varepsilon>0$ and $m$ large. On a compact neighbourhood of $P_0$, $\sqrt{2D_{S'}}$ is bounded above and below by constant multiples of the $g$-distance to $\cM_{S'}$, because $D_{S'}$ is smooth and vanishes exactly on $\cM_{S'}$ with non-degenerate normal Hessian. Since $\sqrt{2D_{S'}(\Sigma_m)}\ge c_1|t_m|/\sqrt n$, the bound of \citet[Lemma~1]{laurent2000adaptive}, applied to the quadratic forms $u^\top\hat\Sigma u$, $u\in\{e_i\pm e_j\}$, in whitened coordinates, which are scaled $\chi^2_{n-1}$ variables, gives $\Prob(S'\text{ passes step (1)})\le C\exp(-c_2\min(t_m^2,n))$, so $|t_m|\Prob(S'\text{ passes})\to0$. In the bound on $\sqrt n\Leb(C_n)$ above, the term of $S'$ is its spread, of order $|t_m|$, times the indicator that $S'$ passes, so its contribution to $\sqrt n\,\E\Leb(C_n)$ tends to zero, and the rest of the argument for finite $c$ applies without the classes $S'$. The constants $c_1,c_2,C$ depend on $P_0$ but do not enter the bound.

\emph{(b), honesty.} Let $(G^*,P)\in\cM_k$. As in the proof of \cref{thm:coverage}, $G^*$ survives steps (1) and (2) with probability at least $1-\gamma_1-\gamma_2$, so it suffices that $\Prob_P(\theta_{G^*}(P)\notin W^{\mathrm{eff}}_{G^*})\to\alpha'$ uniformly over $\cM_k$. The eigenvalues of $\Sigma_P$ lie in $[\sigl^2,\sigu^2]$, a compact set $\mathcal C$ of positive definite matrices. Write $(n-1)\hat\Sigma=\Sigma^{1/2}W\Sigma^{1/2}$ with $W$ Wishart with $n-1$ degrees of freedom and identity scale, whose law does not depend on $\Sigma$, and $X_n=\sqrt n(W/(n-1)-I)$, which converges in law to a Gaussian matrix $\mathbb G$. The map $f=\theta_{G^*}\circ\Pi_{G^*}$ is smooth on the positive definite cone and equals $\theta_{G^*}$ on $\cM_{G^*}$, so $\sqrt n(f(\hat\Sigma)-f(\Sigma))=\ell_\Sigma(X_n)+o_P(1)$ uniformly over $\Sigma\in\mathcal C$, with $\ell_\Sigma=df_\Sigma[\Sigma^{1/2}\cdot\Sigma^{1/2}]$; for $\Sigma\in\cM_{G^*}$, $\ell_\Sigma(\mathbb G)$ is centred Gaussian with standard deviation $s_{\mathrm{eff}}(G^*,\Sigma)\ge\sigl/\sigu$ (\cref{lem:oracle}), because the derivative of the maximum likelihood plug-in is the efficient influence function. Also $\hat s_{G^*}/s_{\mathrm{eff}}(G^*,\Sigma)\to1$ in probability uniformly, by continuity of $s_{\mathrm{eff}}\circ\Pi_{G^*}$ on a compact neighbourhood of $\mathcal C$. If the convergence of the coverage were not uniform, there would be $n_m\to\infty$ and $\Sigma_m\in\cM_{G^*}\cap\mathcal C$ violating it, with $\Sigma_m\to\bar\Sigma$ along a subsequence; $\bar\Sigma\in\cM_{G^*}$ because $\cM_{G^*}$ is closed in the positive definite cone (it is defined by vanishing local Markov partial covariances), so $\ell_{\Sigma_m}\to\ell_{\bar\Sigma}$ and $s_{\mathrm{eff}}(G^*,\Sigma_m)\to s_{\mathrm{eff}}(G^*,\bar\Sigma)>0$; but $\ell_{\Sigma_m}(X_{n_m})/s_{\mathrm{eff}}(G^*,\Sigma_m)\to N(0,1)$ in law by Slutsky's lemma, a contradiction.

\emph{(b), the bound.} The argument of (a) applies with three changes. First, members $D,D'$ of one Markov class $S\subseteq\Xeq$ have $\cM_D=\cM_{D'}$, hence $\Pi_D=\Pi_{D'}$, and $\theta_D=\theta_{D'}$ on $\cM_S$ by \eqref{eq:K1p}, so they have the same estimator and the same estimated standard deviation; there are at most $m$ distinct estimators among the members of $\Xeq$. Second, the offsets are $b_{D,m}=\sqrt n(\theta_D(\Pi_D\Sigma_m)-\theta_{G^*}(\Sigma_m))$: for $D\simeq G^*$, $\Pi_D\Sigma_m=\Sigma_m$ and $b_{D,m}=0$; for $D\in S'$, on $\{\lambda_e=0\}$ the law lies in $\cM_H\subseteq\cM_D$, where $\theta_D\circ\Pi_D=\theta_D=\theta_H=\theta_{G^*}$, and the derivative along $e$ is $\delta_{G^*,D}$ by definition, so $b_{D,m}\to c\,\delta_{G^*,S'}$ and the spread of a class is $|b_{S',m}|\to|c\,\delta_{G^*,S'}|$. Third, $Z_D=d(\theta_D\circ\Pi_D)[\Delta]$ has standard deviation $s_{\mathrm{eff}}(D,P_0)$ because $P_0\in\cM_D$, the half-widths satisfy $\sqrt nh_D\to z_{1-\alpha'/2}s_{\mathrm{eff}}(D,P_0)$, and $\E[\max Z-\min Z]\le2e_m\max_Ds_{\mathrm{eff}}(D,P_0)$: for $m=2$, $\E|Z_1-Z_2|=\sqrt{2/\pi}\,\mathrm{sd}(Z_1-Z_2)\le2\sqrt{2/\pi}\max_is_i$, and for $m\ge3$ the maximal inequality is used. The moment bounds needed for uniform integrability hold for $n\ge n_0(p)$, with $n_0(p)$ and the exponent $a$ below determined by the degrees of the rational functions involved: $\theta_D\circ\Pi_D$ and $s_{\mathrm{eff}}\circ\Pi_D$ are rational functions bounded by $C(1+\|\Sigma\|+\|\Sigma^{-1}\|)^a$, so on $\{\|W/(n-1)-I\|\le\varepsilon\}$ the mean value theorem bounds $\sqrt n|f(\hat\Sigma)-f(\Sigma)|$ by $L\|X_n\|$, and the complement has probability $O(e^{-cn})$ while the inverse Wishart moments of the required order are finite for $n\ge n_0(p)$, so H\"older's inequality applies. Dividing the limit bound by $2z_{1-\alpha/2}s_{\mathrm{eff}}(G^*,P_0)$ and maximizing over $G^*$ gives (b); $\max_Ds_{\mathrm{eff}}(D,P_0)/s_{\mathrm{eff}}(G^*,P_0)\le\mathcal R(P_0)$.

\emph{(c), the directions.} Let $\Sigma'\in\cM_D$. Since $\mu_D$ vanishes on $\cM_D$, $\mu_D(\Sigma'+\Delta)=d\mu_D[\Delta]+O(\|\Delta\|^2)$, and $2D_D=\sum_j\log(1+|\mu_{D,j}|^2)$ because $\sigma^2_{j\mid\pa(j)}=\sigma^2_{j\mid\mathrm{pred}(j)}+\gamma_j^\top\Sigma_{E(j)E(j)\cdot\pa(j)}\gamma_j$; so $D_D(\Sigma'+\Delta)=\frac12|d\mu_D[\Delta]|^2+O(\|\Delta\|^3)$. Comparing with $D_D(\Sigma'+\Delta)=\frac12\|\Delta^\perp\|_g^2+O(\|\Delta\|^3)$ from the proof of (a) gives $|d\mu_D[\Delta]|=\|\Delta^\perp\|_g$ for all $\Delta$: $d\mu_D$ vanishes on the tangent space of $\cM_D$ and maps its normal space isometrically onto $\R^q$. In particular $|u_{S^*,S'}|^2=\|v^\perp\|_g^2=I_{G^*,S'}$, and $u_{S^*,S'}$ depends on the member of $S^*$ only through a non-zero factor, because $v$ is determined modulo the tangent space of $\cM_H\subseteq\cM_D$ up to such a factor. Since $\mu_{D,j}$ is invariant under scaling of $\Sigma$, $\hat\tau_{D,j}=\sqrt{n-j}\,\mu_{D,j}(\hat\Sigma)$, and in the local limit of (a), $\hat\tau_{D_{S'}}=\sqrt n\mu_{D_{S'}}(\Sigma_m)+d\mu_{D_{S'}}[\Delta_m]+o_P(1)\to N(cu_{S^*,S'},I_q)$ in law, with $S^*$ the class of $G^*$.

\emph{(c), honesty.} Let $(G^*,P)\in\cM_k$. The LR screen rejects $G^*$ with probability $\gamma_{1,n}$. If $G^*$ is Markov equivalent to the members of a class $S'$ of $\Xeq$, then $P\in\cM_{D_{S'}}$, and the directed test rejects with probability exactly $\gamma_3$ by \cref{thm:coverage}(i). Step (2) rejects $G^*$ with probability at most $\gamma_2$, and the efficient Wald interval of $G^*$ misses $\theta_{G^*}(P)$ with probability tending to $\alpha'$ uniformly over $\cM_k$, by the argument for (b), since $z_{1-\alpha'_n/2}\to z_{1-\alpha'/2}$. Hence the coverage is at least $1-\alpha-o(1)$ uniformly.

\emph{(c), the bound.} The argument of (b) applies with three changes. First, candidates outside $\cX$ are still kept with exponentially small probability: the $1-\gamma_{1,n}$ quantile of the law in \cref{thm:coverage}(i) is $O(q+\log(1/\gamma_{1,n}))=o(n)$, because $\Prob(n\log(1+U_j/V_j)>s)\le\Prob(V_j<n/2)+\Prob(U_j>s/2)$ for $s\le n$, the first term is $e^{-cn}$ and the second at most $Ce^{-s/5}$, while $T_D=2nD_D(\hat P)\ge nD_D(\Sigma_m)$ with probability $1-O(e^{-c n})$ and $D_D$ is bounded away from zero near $P_0$. Second, a class $S'\ne S^*$ passes step (1) with at most the probability that it passes the directed test, which tends to $\Prob(\chi^2_{d}(c^2|B_{S'}^\top u_{S^*,S'}|^2)\le\chi^2_{d,1-\gamma_3})=\beta_d(|c|\sqrt{I_{G^*,S'}})$, $d=d_{S'}$, because $u_{S^*,S'}\in U_{S'}$ and the null laws converge to $\chi^2_d$, whose distribution function is continuous; as in (a), only this marginal law enters, and $x\beta_d(x)\le\varpi_d(\gamma_3)\le\bar\varpi_{m-1}(\gamma_3)$ for $x=|c|\sqrt{I_{G^*,S'}}$. Third, for infinite $c$, $|B_{S'}^\top\sqrt n\mu_{D_{S'}}(\Sigma_m)|\ge|t_m|(\sqrt{I_{G^*,S'}}-C\varepsilon)$ when $|t_m|\le\varepsilon\sqrt n$, and $\hat\tau=\sqrt n\mu(\hat\Sigma)(1+O(n^{-1}))$ with $\mu$ smooth, so the bound of \citet[Lemma~1]{laurent2000adaptive} used in (a) gives $\Prob(S'\text{ passes})\le C\exp(-c_2\min(t_m^2,n))$ and $|t_m|\Prob(S'\text{ passes})\to0$. The rest of the argument for (b) is unchanged.

\emph{(d).} The argument of (a) applies with the changes of the proof of \cref{thm:K}(ii) below for the classes $S'$ other than that of $G^*$. Members of the class of $G^*$ have offset $b_{D,m}=0$, by \eqref{eq:K1p} and validity of $A(D)$. Members of one class $S'\ne S^*$ have different offsets, but step (1b) rejects the whole class at once, so its term is the spread of $\{0\}\cup\{b_{D,m}:D\in S'\}$ times the indicator that $S'$ passes; the spread is at most $C|t_m|$ and tends to $w^A_{G^*,S'}|c|$ for finite $c$. The two cases of the proof of \cref{thm:K}(ii), the event $A_n$ and the directions are the same as for \DSeff, with $|b_{S',m}|$ replaced by this spread, so the probability that $S'$ passes tends to at most $\E\beta_{d}(|c|\sqrt{I_{G^*,S'}})$ with $d\le m-1$ for finite $c$ and $|t_m|\Prob(S'\text{ passes})\to0$ for infinite $c$, and $x\beta_d(x)\le\bar\varpi_{m-1}(\gamma_3)$ replaces $x\beta_q(x)\le\varpi_q(\gamma_1)$.
\end{proof}

\begin{proof}[Proof of \cref{thm:K}(ii)]
We follow the proof of \cref{thm:explicit}(c), whose arguments hold verbatim with the fixed level $\gamma_1$ in place of $\gamma_{1,n}$ (the critical value of step (1a) is then constant), and treat only the differences, which concern step (0). Suppose the bound fails. As in the proof of \cref{thm:explicit}(a) there are $\varepsilon>0$, $r_m\downarrow0$, $n=n_m\to\infty$ and $(G^*,P_m)\in\cM_k$ with $\|\Sigma_{P_m}-\Sigma_{P_0}\|\le r_m$ violating the bound by $\varepsilon$, with $G^*\in\Xeq$, $P_m=P(\phi_m)$, $\phi_m\to\phi_0^{G^*}$, $a_m=\lambda_e(\phi_m)$ and $t_m=\sqrt n\,a_m\to c\in[-\infty,\infty]$. Let $S^*$ be the class of $G^*$. Step (1b) only removes candidates, and the hull of fewer intervals is shorter, so all terms of the bound on $\sqrt n\Leb(C_n)$ in the proof of \cref{thm:explicit}(a) and (b) are handled as there, except the terms $|b_{S',m}|\ind\{S'\text{ passes step (1)}\}$ of the classes $S'\ne S^*$ of $\Xeq$, where $|b_{S',m}|\le C|t_m|$ and $b_{S',m}\to c\,\delta_{G^*,S'}$ for finite $c$. It suffices to show that the expectation of such a term has limit at most $\bar\varpi_{m-1}(\gamma_3)|\delta_{G^*,S'}|/\sqrt{I_{G^*,S'}}$ if $c$ is finite and zero if it is not.

\emph{Geometry.} For $(D,e')\in\cG\times S$, $D_{D-e'}(P_0)=0$ if and only if $e'=e$ and $D-e\simeq H$. Indeed, if $e'\ne e$, the endpoints of $e'$ are non-adjacent in $D-e'$, hence d-separated by some set in it, but adjacent in $H$, so not d-separated in $H$, and $P_0$, which is faithful to $H$, is not in $\cM_{D-e'}$; if $e'=e$, $D-e$ has the skeleton of $H$, and $P_0\in\cM_{D-e}$ if and only if $D-e\simeq H$, as in \cref{lem:Z}. Hence there are a neighbourhood $\mathcal V$ of $\Sigma_{P_0}$ and $d_0>0$ with $D_{D-e'}\ge d_0$ on $\mathcal V$ for all other pairs, while all pairs $(D,e)$ with $D-e\simeq H$ give the function $D_H$. As in the proof of \cref{prop:CI}, on $\cM_{G^*}$ near $P_0$, $D_H=\lambda_e^2\tilde g_H$ and $D_{S'}=\lambda_e^2\tilde g_{S'}$ with continuous $\tilde g_H,\tilde g_{S'}>0$, because the $e$-curve leaves $\cM_H$ transversally and $I_{G^*,S'}>0$; and $\sqrt{2D_H}$ is equivalent to the $g$-distance to $\cM_H$ near $\Sigma_{P_0}$, by the argument of the proof of \cref{thm:explicit}(a) for infinite $c$, which holds for any DAG model, in particular $\cM_H$, so $\sqrt{D_H(\Sigma')}\le C_0(\sqrt{D_H(\Sigma)}+\|\Sigma'-\Sigma\|)$ near $\Sigma_{P_0}$. The sample covariance of $n_1$ observations satisfies $\Prob(\|\hat\Sigma^{(1)}-\Sigma_m\|>s)\le C\exp(-cn_1\min(s,s^2))$, uniformly for $\Sigma_m$ near $\Sigma_{P_0}$, by \citet[Lemma~1]{laurent2000adaptive} applied to whitened quadratic forms as in the proof of \cref{thm:explicit}(a). Passing to a subsequence, one of the following two cases holds for all $m$.

\emph{Case I: $D_H(\Sigma_m)>c_0\tau_n$}, with $c_0=1/(4C_0^2)$. Then $a_m^2\ge c\tau_n$, so $t_m^2\ge cn\tau_n\to\infty$, and $D_{S'}(\Sigma_m)\ge c\,a_m^2$. The critical value of step (1a) is fixed, so, as for infinite $c$ in the proof of \cref{thm:explicit}(a), step (1a) alone keeps $S'$ with probability at most $C\exp(-c_2\min(t_m^2,n))$, and $|t_m|\Prob(S'\text{ passes})\to0$.

\emph{Case II: $D_H(\Sigma_m)\le c_0\tau_n$.} Let $A_n$ be the event that the minimum in step (0) is at most $\tau_n$ and is attained at a pair $(\hat D,\hat e)$ with $\hat e=e$ and $\hat D-e\simeq H$. On $A_n$, $\hat H\simeq H$, so $\hat\cX=\Xeq$ with the same classes and, by the fixed rule, the same members $D_{S'}$ and orders. The complement of $A_n$ requires $D_H(\hat\Sigma^{(1)})>\tau_n$, hence $\|\hat\Sigma^{(1)}-\Sigma_m\|\ge\sqrt{\tau_n}/(2C_0)$, or $D_{D-e'}(\hat\Sigma^{(1)})\le\tau_n<d_0/2$ for some other pair, hence $\|\hat\Sigma^{(1)}-\Sigma_m\|\ge c$ for $m$ large. So $\Prob(A_n^c)\le C\exp(-cn_1\tau_n)+C\exp(-cn_1)$, which is $o(n^{-1/2})$ because $n_1\tau_n\ge n\varepsilon_n\tau_n$ and $n\varepsilon_n\tau_n/\log n\to\infty$; since $|b_{S',m}|\le C\sqrt n$, the contribution of $A_n^c$ vanishes. On $A_n$, $\hat u_{S^{**},S'}=U_{S^{**},S'}(\Pi_H\hat\Sigma^{(1)})$ for every class $S^{**}\ne S'$, with $U_{S^{**},S'}$ smooth near $\Sigma_{P_0}$ and $U_{S^{**},S'}(\Sigma_{P_0})$ a non-zero multiple of $u_{S^{**},S'}$; we write $u$ for $U_{S^*,S'}(\Sigma_{P_0})$, which spans the same line as $u_{S^*,S'}$. Since $\|\Pi_H\hat\Sigma^{(1)}-\Sigma_{P_0}\|\le C(\|\hat\Sigma^{(1)}-\Sigma_m\|+r_m)$, for every $\eta>0$, $\Prob(A_n,\ |\hat u_{S^*,S'}-u|>\eta)\le C\exp(-cn_1\eta^2)$ for $m$ large, which is $o(n^{-1/2})$ because $n_1/\log n\to\infty$.

\emph{Case II, finite $c$.} Given the first part, $\hat B=\hat B_{S'}$ is fixed and $\hat u_{S^*,S'}$ lies in its span. The second part is an independent sample of size $n_2=n(1-\varepsilon_n)$, with $\sqrt{n_2}\,a_m\to c$, so, as in the proof of \cref{thm:explicit}(c), $\hat\tau^{(2)}_{D_{S'}}$ converges in law to $N(c\,u_{S^*,S'},I_q)$. The null law of the directed statistic given $B$ is that of $|B^\top\zeta|^2$ with $\zeta$ as in \cref{thm:coverage}(i), and $\zeta\Rightarrow N(0,I_q)$, so it converges to $\chi^2_{d(B)}$ uniformly over the compact set of $q\times d$ matrices with orthonormal columns, $d\le m-1$. Along a further subsequence $(\hat B,\hat\tau^{(2)}_{D_{S'}})$ converges in law to $(B,Z)$ with $Z\sim N(c\,u_{S^*,S'},I_q)$ independent of $B$, because the two parts are independent, and $BB^\top u=u$ almost surely by the previous paragraph, whatever the dimension $d(B)\le m-1$ of the limit (the rank of the estimated directions is only lower semicontinuous). Hence $\Prob(A_n,\ S'\text{ passes step (1b)})\to\E\beta_{d(B)}(|c|\sqrt{I_{G^*,S'}})$, with $\beta_d(x)=\Prob(\chi^2_d(x^2)\le\chi^2_{d,1-\gamma_3})$ and $|u_{S^*,S'}|^2=I_{G^*,S'}$ (\cref{thm:explicit}(c)). Since $b_{S',m}$ is deterministic with $b_{S',m}\to c\delta_{G^*,S'}$, and $x\beta_d(x)\le\varpi_d(\gamma_3)\le\bar\varpi_{m-1}(\gamma_3)$, the limit of the expected term is at most $\bar\varpi_{m-1}(\gamma_3)|\delta_{G^*,S'}|/\sqrt{I_{G^*,S'}}$.

\emph{Case II, infinite $c$.} Now $a_m\to0$ and $|t_m|\to\infty$. Since $\mu_{D_{S'}}$ vanishes on $\cM_H\subseteq\cM_{D_{S'}}$, $\mu_{D_{S'}}(\Sigma(\phi))=\lambda_e\tilde\mu(\phi)$ on $\Theta_{G^*}$ near $\phi_0^{G^*}$ with $\tilde\mu$ continuous and $\tilde\mu(\phi_0^{G^*})=u_{S^*,S'}$, so $\mu_{D_{S'}}(\Sigma_m)=a_m(u_{S^*,S'}+o(1))$. On $A_n\cap\{|\hat u_{S^*,S'}-u|\le\eta\}$, $|\hat B^\top u_{S^*,S'}|\ge|u_{S^*,S'}|-C\eta$, so $|\hat B^\top\sqrt{n_2}\,\mu_{D_{S'}}(\Sigma_m)|\ge c'|t_m|$ for $\eta$ small and $m$ large, and the argument for infinite $c$ in the proof of \cref{thm:explicit}(c), applied to the second part given the first, gives $\Prob(S'\text{ passes step (1b)}\mid\text{first part})\le C\exp(-c_2\min(t_m^2,n))$ on this event. With the bounds on the other events, $|t_m|\Prob(S'\text{ passes})\to0$.

In all cases the terms of the classes $S'\ne S^*$ are bounded as required, and the rest of the proof of \cref{thm:explicit}(c), with the efficient estimators, the maximal inequality $e_m$ and the division by $\sqrt n\orl_n$, gives a contradiction.
\end{proof}

\begin{proof}[Proof of \cref{cor:RK}]
The lower bound is \cref{thm:K}(i) with $\Psi_\alpha(K)\ge\frac12c_P+c_KK$, together with $\Pi^*_k(P_0)\ge1$ from \cref{thm:free}(i). It holds for every asymptotically honest sequence and every $r>0$: for each $r$, $\limsup_n\Pi^{\mathrm{loc}}_n(C_n;r)\ge\liminf_n\Pi^{\mathrm{loc}}_n(C_n;r)$ is at least the bound, and the limit as $r\downarrow0$ exists because $\Pi^{\mathrm{loc}}_n(C_n;r)$ is non-decreasing in $r$; so it bounds $\Pi^*_k(P_0)$. For the upper bound, let $\gamma_1,\gamma_2>0$ with $\gamma_1+\gamma_2<\alpha-\gamma_3$, and take $\varepsilon_n=\tau_n=n^{-1/3}$. By \cref{thm:coverage,thm:K}(ii), \DSeff\ with these levels is asymptotically honest over $\cM_k$ and has limiting local price at most $\mathcal R(z_{1-(\alpha-\gamma_1-\gamma_2-\gamma_3)/2}+e_m)/z_{1-\alpha/2}+b_mK_+$, so this bounds $\Pi^*_k(P_0)$. The left-hand side does not depend on $\gamma_1,\gamma_2$, and letting $\gamma_1+\gamma_2\downarrow0$ gives the bound with $a_m$; for given $\eta>0$, $\gamma_1+\gamma_2$ small enough gives $(a_m+\eta)\mathcal R$. Also $K_+\le(m-1)K$. For $m=2$, $e_2=\sqrt{2/\pi}$, $K_+=K$ and $d\le1$, so $b_2=\varpi_1(\gamma_3)/(2z_{1-\alpha/2})$. For $m=1$, all members of $\Xeq$ have the same efficient estimator by \eqref{eq:K1p}, so $\mathcal R=1$ and $K=0$, and \cref{thm:explicit}(b) with $\gamma_1+\gamma_2\downarrow0$ gives the upper bound $1$. For the ratio, the maximum on the left is at least $c_P\mathcal R$ and at least $c_KK$, so $(a_m\mathcal R+(m-1)b_mK)/\max\{c_P\mathcal R,\frac12c_P+c_KK\}\le a_m/c_P+(m-1)b_m/c_K$.
\end{proof}

Throughout the rest of this subsection, $f$, $\phi_\lambda$, $\phi_0'$ and $P_\lambda$ are as in \cref{thm:K-lambda}, $(G^*,G')$ is an ordered non-equivalent pair in $\Xeq$ with $e=u\to v$ in $G^*$, and $I(\lambda)=I_{G^*,G'}(P_\lambda)$, $\delta(\lambda)=\delta_{G^*,G'}(P_\lambda)$. We use the condition
\begin{enumerate}
\item[(A)] there are $D^*\simeq G^*$ and $D'\simeq G'$ with $D^*-e\simeq D'-e\simeq H$, $(D^*-e)-f\simeq(D'-e)-f\simeq H-f$ and $D^*-f\simeq D'-f$;
\end{enumerate}
and show that $I(0)=0$ if and only if (A) holds, and that (A) holds if and only if (F) and (C) hold.

\begin{lemma}[Only the Markov classes matter]\label{lem:Kclass}
If $D^*\simeq G^*$ with $D^*-e\simeq H$ and $D'\simeq G'$, then $I_{D^*,D'}(\lambda)=\kappa(\lambda)I_{G^*,G'}(\lambda)$ for a smooth function $\kappa>0$; in particular the two vanish to the same order at $\lambda=0$.
\end{lemma}
\begin{proof}
$D_{G'}$ is the Kullback--Leibler distance to $\cM_{G'}=\cM_{D'}$, so it does not depend on the representative of the class of $G'$. $\cM_{G^*}=\cM_{D^*}$ is a smooth manifold near $P_\lambda$, and $\cM_H=\cM_{G^*-e}=\cM_{D^*-e}$ is a smooth hypersurface in it. The restriction of $D_{G'}$ to $\cM_{G^*}$ is smooth, non-negative and zero on $\cM_H$, so its gradient vanishes on $\cM_H$, and its Hessian at a point of $\cM_H$, which is positive semi-definite because the point is a minimum, has the tangent vectors of $\cM_H$ in its kernel; the second derivative along a curve that leaves $\cM_H$ is therefore the Hessian evaluated at the normal component of its velocity. The $e$-curves of $G^*$ and $D^*$ both lie in $\cM_{G^*}$ and leave $\cM_H$ with non-zero normal velocity, because the structural parametrizations are immersions; the ratio of the squared normal velocities is the smooth positive function $\kappa$.
\end{proof}

\begin{lemma}[One differing v-structure]\label{lem:onev}
Let $G_1,G_2$ be DAGs on the same skeleton with $\mathcal V(G_1)=\mathcal V(G_2)\cup\{\tau\}$, where $\tau=a\to z\leftarrow b\notin\mathcal V(G_2)$ and $\mathcal V$ denotes the set of v-structures. If some $D\simeq G_2$ has $\mathrm{CC}_D(a,z)=\mathrm{CC}_{G_1}(a,z)$, then some $D'\simeq G_2$ has $\mathrm{CC}_{D'}(a,z)=\mathrm{CC}_{G_1}(a,z)$ and $\mathrm{CC}_{D'}(z,b)=\mathrm{CC}_{G_1}(z,b)$. Moreover, $D'$ can be chosen so that reversing one of its edges $a-z$, $z-b$ gives a member of the class of $G_1$ that agrees with $G_1$ on all edges at $z$.
\end{lemma}
\begin{proof}
Write $N(v)$ for the neighbours of $v$ in the skeleton, $P=\pa_{G_1}(z)$ and $U=P\setminus\{a,b\}$; $a\notin N(b)$ because $\tau$ is a v-structure, and $\mathrm{CC}_D(s,t)=\operatorname{ch}_D(s)\cap\operatorname{ch}_D(t)$. Every member of the class of $G_1$ has $a\to z\leftarrow b$. We use three facts. (a) If $s\to t$ and $w\in N(s)\cap N(t)$, then $t\to w$ forces $s\to w$ by acyclicity, so $\mathrm{CC}(s,t)=\operatorname{ch}(t)\cap N(s)$. (b) Reversing a covered edge $x\to y$, one with $\pa(y)=\pa(x)\cup\{x\}$, gives a Markov equivalent DAG (the `if' direction of \citealp[Lemma~1]{chickering1995transformational}). (c) A member of the class of $G_2$ has every v-structure of $G_1$ other than $\tau$, and no others.

\emph{Shape.} Let $D\simeq G_2$ with $\mathrm{CC}_D(a,z)=\mathrm{CC}_{G_1}(a,z)=\operatorname{ch}_{G_1}(z)\cap N(a)$. If $a\to z$ in $D$, then $z\to b$ by (c), $D$ agrees with $G_1$ on the edges $z-w$, $w\in N(a)$, by (a), and for $w\in N(z)\setminus N(a)$, $w\ne b$, both have $w\to z$ if and only if $a\to z\leftarrow w$ is a v-structure, by (c); so $\pa_D(z)=U\cup\{a\}$ (type $\alpha$). If $z\to a$ in $D$, then $\mathrm{CC}_D(a,z)=\operatorname{ch}_D(a)\cap N(z)$ by (a), so for $w\in N(a)\cap N(z)$, $w\to a$ in $D$ if and only if $w\to z$ in $G_1$; for $w\in N(a)\setminus N(z)$, $w\to a$ would make $z\to a\leftarrow w$ a v-structure of $D$, hence of $G_1$, which has $a\to z$; and $U\subseteq N(a)$, since for $u\in U\setminus N(a)$ the v-structure $a\to z\leftarrow u$ of $G_1$ would force $a\to z$ in $D$. So $\pa_D(a)=U\cup\{z\}$ (type $\beta$).

\emph{A covered edge.} For $D\simeq G_2$ let $\Delta(D)$ be the set of edges of $D$ that have the opposite orientation in $G_1$; if $x\to y$ in $D$ is in $\Delta(D)$, call $y$ a head and $x$ a $\Delta$-parent of $y$. Let $y$ be a head none of whose proper ancestors in $D$ is a head, and $x$ a $\Delta$-parent of $y$ none of whose children is a $\Delta$-parent of $y$; the first head in a topological order of $D$ and its last $\Delta$-parent qualify. Then $x\to y$ is covered, unless $x=z$, $y\in\{a,b\}$ and the other of $a,b$ is a parent of $z$ in $D$. Indeed, let $w\in\pa_D(x)$. Then $w\to x$ in $G_1$, as otherwise $x$ would be a head that is a proper ancestor of $y$, so $G_1$ has $w\to x\leftarrow y$. If $w\notin N(y)$, this is a v-structure of $G_1$ but not of $D$, which has $x\to y$, so it is $\tau$, the exception; otherwise $w\to y$ in $D$ by acyclicity. Conversely, let $w\in\pa_D(y)$, $w\ne x$. If $w\notin N(x)$, then $w\to y\leftarrow x$ is a v-structure of $D$, hence of $G_1$, which has $y\to x$; so $w\in N(x)$. If $x\to w$ in $D$, then either $w\to y\in\Delta(D)$, excluded by the choice of $x$, or $G_1$ has $w\to y\to x$ and hence $w\to x$, so that $w$ is a head that is a proper ancestor of $y$; so $w\to x$ in $D$. Hence $\pa_D(y)=\pa_D(x)\cup\{x\}$. The second half does not use the exception. For two equivalent graphs this is the edge selection of \citet[Lemma~2]{chickering1995transformational}, and there is no exception; here the single differing v-structure $\tau$ is the only source of one.

\emph{Descent.} Start from $D$ and repeat the following. If $D$ is of type $\alpha$, it differs from $G_1$ at $z$ only in $z\to b$, so $\Delta(D)\ne\emptyset$ and the only edge of $\Delta(D)$ at $z$ is $z\to b$. If the chosen pair is not $(y,x)=(b,z)$, the exception cannot occur ($y=a$ would require $z\to a$), so $x\to y$ is covered and does not meet $z$, and its reversal gives a member of the class of $G_2$ of type $\alpha$ with one edge fewer in $\Delta$. If it is $(b,z)$, the second half of the argument above gives $\pa_D(b)\setminus\{z\}\subseteq\pa_D(z)\setminus\{a\}=U$; conversely, for $u\in U$, $G_1$ has $u\to z\leftarrow b$, which is not a v-structure because it would force $b\to z$ in $D$, so $u\in N(b)$ and $u\to b$ in $D$. So $\pa_D(b)=U\cup\{z\}$, and reversing $z\to b$ creates no cycle (a second directed path from $z$ to $b$ would end in an edge $u\to b$ with $u\in U\subseteq\pa_D(z)$), removes no v-structure (the other parents of $b$ are adjacent to $z$) and adds exactly $\tau$ (the parents of $z$ other than $a$ are adjacent to $b$); the result is in the class of $G_1$ and has $\pa(z)=P$, and we stop. If $D$ is of type $\beta$, then $z\to a\in\Delta(D)$ and $a$ is a head. If some head $y\ne a$ has no proper ancestor that is a head, the exception cannot occur ($y=b$ would require $a\to z$), so $x\to y$ is covered; and $x\ne a$, since otherwise $\pa_D(y)=U\cup\{z,a\}$, so $G_1$, which has $y\to a\to z$, has $y\to z$, and $y\in U\subseteq\pa_D(a)$ ($y\ne b$ because $y\in N(a)$), contradicting $a\to y$. The reversal keeps $\pa(a)$, hence type $\beta$, and removes one edge from $\Delta$. Otherwise every head is $a$ or a descendant of $a$, because the first head among the ancestors of a head is $a$. An edge $z-w$ in $\Delta(D)$ with $w\notin\{a,b\}$ would have head $z$, or head $w\in U\subseteq\pa_D(a)$ if $z\to w$ in $D$; both are proper ancestors of $a$, which is impossible. So $\pa_D(z)\setminus\{a,b\}=U$. If $b\to z$ in $D$, reversing $z\to a$ creates no cycle, removes no v-structure and adds exactly $\tau$, as in the previous case with $a$ in place of $b$ and $U\subseteq N(a)$; the result is in the class of $G_1$ and has $\pa(z)=P$, and we stop. If $z\to b$, then $\pa_D(a)=\pa_D(z)\cup\{z\}$, so $z\to a$ is covered, and its reversal gives a member of type $\alpha$ with one edge fewer in $\Delta$.

Every step that does not stop removes one edge from $\Delta$, and $\Delta$ never becomes empty, because it contains $z\to b$ for type $\alpha$ and $z\to a$ for type $\beta$; so the descent stops at some $D'\simeq G_2$, and reversing $z\to b$ or $z\to a$ in $D'$ gives $D_1\simeq G_1$ with $\pa_{D_1}(z)=P$. By (a), $\mathrm{CC}_{D_1}(a,z)=\operatorname{ch}_{D_1}(z)\cap N(a)=\mathrm{CC}_{G_1}(a,z)$, and likewise for $(z,b)$. Reversing $a-z$ changes neither set: $\mathrm{CC}(a,z)$ does not involve the edge between $a$ and $z$, and $\mathrm{CC}(z,b)$ cannot gain or lose $a$, because $a\notin N(b)$; the same holds for $z-b$. So $D'$ has the same two sets.
\end{proof}

\paragraph{Remarks on (D).} If $R=\emptyset$, the second requirement of (D) holds if and only if some DAG Markov equivalent to $H$ has $\pa(x_e)=U$ (\cref{lem:Econf}), which, reformulating the local validity condition of \citet[Lemma~3.1]{maathuis2009estimating}, is a local condition on the CPDAG $\mathcal C$ of $H$: $\pa_{\mathcal C}(x_e)\subseteq U$, the vertices of $U\setminus\pa_{\mathcal C}(x_e)$ are joined to $x_e$ by undirected edges of $\mathcal C$, and they are pairwise adjacent (each of them is adjacent to every parent of $x_e$ in $\mathcal C$ by the first orientation rule of \citealp{meek1995causal}). The set $U$ does not depend on the choice of $G_\wedge$ within its Markov class in $\Xeq$, by the argument of the paragraph on common children below. Without the parent-set requirement the condition fails: on six nodes there are configurations that satisfy the other requirements of (D) but in which no member of $[H]_Y$ gives $x_e$ the parent set $U$, for instance because a v-structure of $H$ at $x_e$ fixes its parents, and $K(P_\lambda)$ stays bounded (\cref{app:certify}).

\begin{lemma}[Reversing the weak edge]\label{lem:Econf}
Let $(G^*,G')$ satisfy (F), with $x_e,z,x_f$, $G_\wedge,G_\vee$ and $U=\pa_{G_\wedge}(z)\setminus\{x_e,x_f\}$ as in \cref{sec:gen-constant}. Call $D\in[G_\vee]_Y$ an $e$-configuration if $z\to x_e$, $x_f\to z$, $\pa_D(x_e)=U\cup\{z\}$ and $\pa_D(z)=U\cup\{x_f\}$. Then:
\begin{enumerate}
\item[(i)] if $D$ is an $e$-configuration, then $D\in\Xeq$, $e\notin R$, and the graph $D_\wedge$ obtained from $D$ by reversing $e$ is a member of $\Xeq$ equivalent to $G_\wedge$; if (C) holds, $(D_\wedge-e)-f\simeq H-f$;
\item[(ii)] an $e$-configuration exists if and only if $e\notin R$ and some member of $[H]_Y$ has $\pa(x_e)=U$;
\item[(iii)] if $R=\emptyset$, some member of $[H]_Y$ has $\pa(x_e)=U$ if and only if some DAG Markov equivalent to $H$ has.
\end{enumerate}
\end{lemma}
\begin{proof}
Write $N(v)$ for the neighbours of $v$ and $P=U\cup\{x_e,x_f\}$. The vertex $Y$ is not in $\{x_e,z,x_f\}$: it is a sink of $G_\wedge$, which has $x_e\to z\leftarrow x_f$, and if it were $z$, then $G_\vee$ would also have this v-structure.

(i) Reversing $z\to x_e$ in $D$ creates no cycle, because a second directed path from $z$ to $x_e$ would end in an edge $u\to x_e$ with $u\in U\subseteq\pa_D(z)$. It removes no v-structure at $x_e$, because $z$ is adjacent to every $u\in U$, and at $z$ it adds exactly $x_e\to z\leftarrow x_f$, because $x_e$ is adjacent to every $u\in U$ and not to $x_f$. So $D_\wedge\simeq G_\wedge$ by (F). Since $D$ and $G_\wedge$ both satisfy $\cK$ and orient $e$ differently, $e\notin R$, and $D_\wedge$ satisfies $\cK$. By fact (a) in the proof of \cref{lem:onev}, $\mathrm{CC}_D(x_e,z)=\operatorname{ch}_D(x_e)\cap N(z)=(N(x_e)\cap N(z))\setminus U$, and likewise $\mathrm{CC}_{G_\wedge}(x_e,z)=\operatorname{ch}_{G_\wedge}(z)\cap N(x_e)=(N(x_e)\cap N(z))\setminus U$, using $U\subseteq\pa_D(x_e)\subseteq N(x_e)$ and $x_f\notin N(x_e)$; this set equals $\mathrm{CC}_H(x_e,z)$ because $G_\wedge\in\Xeq$. So $D-e=D_\wedge-e\simeq H$ by the paragraph on common children in the proof of \cref{thm:K-lambda}, and $D,D_\wedge\in\Xeq$. By the same paragraph, $\mathrm{CC}_{D_\wedge}(z,x_f)=\mathrm{CC}_{G_\wedge}(z,x_f)$, which is $\mathrm{CC}_H(z,x_f)$ under (C), so $(D_\wedge-e)-f\simeq H-f$.

(ii) If $D$ is an $e$-configuration, $D-e\in[H]_Y$ has $\pa(x_e)=U$, and $e\notin R$ by (i). Conversely, let $e\notin R$ and let $T\in[H]_Y$ have $\pa_T(x_e)=U$; in particular $U\subseteq N(x_e)$. Start from $D=G_\wedge-e\in[H]_Y$, which has $\pa_D(z)=U\cup\{x_f\}$ and in which $z$ is not an ancestor of $x_e$ or of $U$. As long as some edge of $D$ that $T$ orients the other way has its head in $\mathrm{An}_D(\{x_e\}\cup U)$, the set of ancestors in $D$ of $\{x_e\}\cup U$, including these vertices themselves, choose such a head $y$ none of whose proper ancestors is a head, and $x$ as in the proof of \cref{lem:onev}; there, with $G_1=T$, which is equivalent to $D$, there is no exception, so $x\to y$ is covered, and we reverse it. Both $x$ and $y$ are ancestors of $\{x_e\}\cup U$, and $z$ never becomes one, because a new directed path from $z$ would use the edge $y\to x$, so that $z$ was already an ancestor of $y$, hence of $\{x_e\}\cup U$. Hence $z\notin\{x,y\}$ and $\pa_D(z)=U\cup\{x_f\}$ throughout. Each reversal removes one edge from the set of edges oriented differently in $D$ and $T$, so the descent stops. The edges at $Y$ and those of $R$ have the same orientation in $D$ and $T$ and are never reversed. When no such head is left, every parent of $x_e$ in $D$ is a parent in $T$, hence in $U$, and every $u\in U$ is a parent of $x_e$, since $x_e\to u$ would have head $u$; so $\pa_D(x_e)=U$. Let $G_1$ be $D$ with $x_e\to z$ added. It is acyclic, because a directed path from $z$ to $x_e$ would end in an edge $u\to x_e$ with $u\in U\subseteq\pa_D(z)$. Its v-structures with arm $e$ are $x_e\to z\leftarrow c$ with $c\in\pa_D(z)\setminus N(x_e)=\{x_f\}$, and the same holds for $G_\wedge$, which is $G_\wedge-e$ with $x_e\to z$ added; the other v-structures are, in both graphs, those of $H$ whose non-central vertices are not $x_e,z$, because $D\simeq G_\wedge-e\simeq H$. So $G_1\in[G_\wedge]_Y$, with $\pa_{G_1}(x_e)=U$ and $\pa_{G_1}(z)=P$, and $\mathrm{CC}_{G_1}(x_e,z)=\operatorname{ch}_{G_1}(z)\cap N(x_e)=\mathrm{CC}_{G_\wedge}(x_e,z)=\mathrm{CC}_H(x_e,z)=\mathrm{CC}_{G_\vee}(x_e,z)$, since $\pa_{G_1}(z)=\pa_{G_\wedge}(z)$ and $G_\wedge,G_\vee\in\Xeq$. Now run the descent of the proof of \cref{lem:onev} with this $G_1$, $G_2=G_\vee$, $a=x_e$, $b=x_f$, starting from $D=G_\vee$, and in type $\alpha$ choose a head other than $x_f$ whenever one has no proper ancestor that is a head. Edges at $Y$ and in $R$ are oriented alike in $G_1$ and $G_\vee$ and are never reversed. If the descent stops in type $\beta$, the graph at which it stops is an $e$-configuration. If it stops in type $\alpha$, then $x_f$ is the only head without a head among its proper ancestors, so every head is $x_f$ or one of its descendants; $x_e$ and $U$ are proper ancestors of $x_f$, so they are not heads, and as above $\pa_D(x_e)=U$. Then $x_e\to z$ is covered, since $\pa_D(z)=U\cup\{x_e\}$; after its reversal $z\to x_f$ is covered, since $\pa_D(x_f)=U\cup\{z\}$; and after that reversal $D$ is an $e$-configuration. Neither reversal touches $Y$, $e\notin R$, and $f\notin R$, because $D$ and $G_\wedge$ orient $f$ differently.

(iii) If $R=\emptyset$ and $T$ is only Markov equivalent to $H$, the first descent in (ii) still applies: the sink $Y$ is not an ancestor of $\{x_e\}\cup U$, so the edges at $Y$ are never reversed, and the resulting $D$ is in $[H]_Y$.
\end{proof}

\begin{proof}[Proof of \cref{thm:K-lambda}]
\emph{Base points.} $\cN$ is the zero set of the analytic function $\psi=\prod b_{ijS}\prod a_{G^*,G'}$ on $\Theta_H$ (\cref{app:benign}). The map $\lambda\mapsto\psi(\phi_\lambda)$ is analytic on $\R$ and non-zero at the coefficient of $f$ in $\phi_0$, so its zeros are isolated, and a punctured neighbourhood of $\lambda=0$ avoids them. The eigenvalues of $\Sigma_{P_\lambda}$ are interior for small $\lambda$ by continuity. By \cref{lem:Dproj}, $I(\lambda)$, $\delta(\lambda)$ and $s_{\mathrm{eff}}(G^*,P_\lambda)$ are analytic in $(\lambda,\phi_0')$ near $\lambda=0$, because $M_j$ and $\sigma^2_{j\mid\pa(j)}$ are positive at $P_{\lambda=0}$; in particular $\delta$ is bounded and $s_{\mathrm{eff}}$ is bounded away from zero for small $\lambda$.

\emph{(iii), first implication.} Let $D^*,D'$ be as in (A). By \cref{lem:Kclass} it suffices to show $I_{D^*,D'}(0)=0$. $P_{\lambda=0}\in\cM_{H-f}=\cM_{(D^*-e)-f}$, so the coefficients of $e$ and $f$ in the parameter of $P_{\lambda=0}$ under $D^*$ are zero. Setting the coefficient of $e$ to $\varepsilon$ keeps that of $f$ at zero, so $P_\varepsilon\in\cM_{D^*-f}=\cM_{D'-f}\subseteq\cM_{D'}$ and $D_{G'}(P_\varepsilon)\equiv0$.

\emph{(iii), second implication.} Suppose a witness configuration of the proof of \cref{thm:E} (\cref{app:E}) uses an edge $g\neq f$, and consider its sparse point: unit variances, $\lambda_e=s$, and only $g$ with a non-zero coefficient $\beta$ in $D^*=G^*-e$. Its base point $s=0$ is the law in which only the endpoints of $g$ are dependent. Since $g$ is an edge of $H-f$, this law belongs to $\cM_{H-f}$, so it equals $P_{\lambda=0}$ for some $\phi_0'$, and its parameter under $G^*$ is the sparse one; hence the $e$-curve from it is the sparse curve. By the proof of \cref{thm:E}, some local Markov partial covariance of $G'$ has non-zero $s$-derivative there, so $I(0)>0$ at this $\phi_0'$. By \cref{lem:Dproj}, $I(0)$ is analytic on all of the connected set $\Theta_{H-f}$, because $M_j$ and $\sigma^2_{j\mid\pa(j)}$ are positive at every point of it (the eigenvalue condition of the theorem only restricts to an open subset); so it is non-zero for almost every $\phi_0'$. Hence $I(0)=0$ for almost every $\phi_0'$ forces every witness to use $f$, which is (F).

\emph{(F) in terms of v-structures.} For $D\in\Xeq$ let $\mathcal V(D)$ be its set of v-structures and $\mathcal A_e(D)$ the subset of those with arm $e$. Since $D-e\simeq H$, $\mathcal V(D)\setminus\mathcal A_e(D)$ is the set of v-structures of $H$ other than those whose two non-central vertices are the endpoints of $e$, the same for all members of $\Xeq$. So $\mathcal V(G^*)\triangle\mathcal V(G')=\mathcal A_e(G^*)\triangle\mathcal A_e(G')$. If $e$ has the same orientation $u\to v$ in both graphs, this consists of the v-structures $u\to v\leftarrow c$ with $c\notin\adj(u)$ and $c\to v$ in exactly one of them (case A of \cref{app:E}); otherwise it is the disjoint union of the v-structures $u\to v\leftarrow c$, $c\in V_+(G^*)$, and $v\to u\leftarrow c$, $c\in V_+(G')$ (cases B1 and B2). These are the witness configurations, and each is determined by its witness edge, its arm other than $e$. Hence all witness edges equal $f$ if and only if the symmetric difference is a single v-structure with arms $e$ and $f$, which is (F); its non-central vertices $x_e,x_f$ are non-adjacent.

\emph{Common children.} Let $D$ be equivalent to $G^*$ or to $G'$. By the same decomposition, $D-e\simeq H$ if and only if $\mathrm{CC}_D(x_e,z)=\mathrm{CC}_H(x_e,z)$, the set of centres of the v-structures of $H$ with non-central vertices $x_e,z$; and, deleting $f$ in the same way, if $D-e\simeq H$, then $(D-e)-f\simeq H-f$ if and only if $\mathrm{CC}_D(z,x_f)=\mathrm{CC}_H(z,x_f)$, because $x_e$ is not a common neighbour of $z$ and $x_f$. Now let $D\simeq G_\wedge$ with $D-e\simeq H$, so that $x_e\to z\leftarrow x_f$ in $D$, and let $w$ be a common neighbour of $z$ and $x_f$. If $w\in\adj(x_e)$, then $z\to w$ forces $x_e\to w$ and $x_f\to w$ by acyclicity, so $w\in\mathrm{CC}_D(z,x_f)$ if and only if $w\in\mathrm{CC}_D(x_e,z)=\mathrm{CC}_H(x_e,z)$. If $w\notin\adj(x_e)$, then $w\to z$ would make $x_e\to z\leftarrow w$ a v-structure, so $z-w$ has the same orientation in all members of the class, and $w\in\mathrm{CC}_D(z,x_f)$ if and only if $z\to w$. Hence $\mathrm{CC}_D(z,x_f)=\mathrm{CC}_{G_\wedge}(z,x_f)$ for all such $D$, and the class of $G_\wedge$ has a member $D$ with $D-e\simeq H$ and $(D-e)-f\simeq H-f$ if and only if (C) holds.

\emph{(A) in graphical terms.} Similarly $\mathcal V(D-f)=(\mathcal V(D)\setminus\mathcal A_f(D))\cup\{z\to w\leftarrow x_f:w\in\mathrm{CC}_D(z,x_f)\}$. Suppose that (F) and (C) hold. By (C) and the previous paragraph, $D_\wedge=G_\wedge$ satisfies $D_\wedge-e\simeq H$ and $(D_\wedge-e)-f\simeq H-f$. The pair $G_1=G_\wedge$, $G_2=G_\vee$ satisfies the hypotheses of \cref{lem:onev} with $a=x_e$, $b=x_f$ and $D=G_\vee$, because $\mathrm{CC}_{G_\vee}(x_e,z)=\mathrm{CC}_H(x_e,z)=\mathrm{CC}_{G_\wedge}(x_e,z)$. So some $D_\vee\simeq G_\vee$ has $\mathrm{CC}_{D_\vee}(x_e,z)=\mathrm{CC}_H(x_e,z)$ and $\mathrm{CC}_{D_\vee}(z,x_f)=\mathrm{CC}_{G_\wedge}(z,x_f)=\mathrm{CC}_H(z,x_f)$, that is, $D_\vee-e\simeq H$ and $(D_\vee-e)-f\simeq H-f$. By (F), $\mathcal V(D_\wedge)$ and $\mathcal V(D_\vee)$ differ only in $x_e\to z\leftarrow x_f$, which has arm $f$, and the two graphs have the same common children of $z$ and $x_f$; so $D_\wedge-f\simeq D_\vee-f$, and $D_\wedge,D_\vee$, in the order of $(G^*,G')$, satisfy (A). Conversely, (A) implies $I(P_{\lambda=0})=0$ and hence (F), and its representative of the class of $G_\wedge$ gives (C) by the previous paragraph.

\emph{(iii), (C).} Suppose that (F) holds and (C) fails, and compute $I$ with $D_\wedge=G_\wedge$ and $D_\vee=G_\vee$ (\cref{lem:Kclass}). There is a common neighbour $w$ of $z$ and $x_f$ that belongs to exactly one of $\mathrm{CC}_{D_\wedge}(z,x_f)$ and $\mathrm{CC}_H(z,x_f)$. Let $P^\circ\in\cM_{H-f}$ have unit variances and non-zero coefficients only on the edges $w-z$ and $w-x_f$ of $H-f$. Under $P^\circ$, $x_e$ and every vertex other than $w,z,x_f$ are independent of all other variables; $z\indep x_f$ if $w\in\mathrm{CC}_H(z,x_f)$, and otherwise $z\indep x_f\mid w$ while $z$ and $x_f$ are dependent. By \cref{lem:Dproj}, $I(P^\circ)=0$ if and only if the $e$-velocity lies in $T_{P^\circ}\cM_{G'}$, that is, $T_{P^\circ}\cM_{G^*}=T_{P^\circ}\cM_{G'}$, since $T\cM_{G^*}=T\cM_H\oplus\R v_e$, $T\cM_H\subseteq T\cM_{G'}$ and the two dimensions agree. For a DAG $D$ with $P^\circ\in\cM_D$ all coefficients on edges at $x_e$ vanish, and the generators $A\,E_{ik}\Sigma+\Sigma E_{ki}A^\top$ and $AE_{ii}A^\top$ of $T_{P^\circ}\cM_D$, $A=(I-B)^{-1}$, fall into three groups with disjoint supports. Edges and variances away from $x_e$ span the tangent space of $\cM_{D-x_e}$ in the block of $V\setminus\{x_e\}$; the variance of $x_e$ gives the entry $(x_e,x_e)$; and the edges at $x_e$ give the row vectors $\mathrm{Cov}(\cdot,p)$, $p\in\pa_D(x_e)$, and $\mathrm{Cov}(\cdot,\varepsilon_c)$, $c\in\operatorname{ch}_D(x_e)$, where $\varepsilon_c$ is the residual of $c$ on $\pa_D(c)\setminus\{x_e\}$. Since $D_\wedge-x_e=(D_\wedge-e)-x_e\simeq H-x_e\simeq D_\vee-x_e$, because deleting a vertex from two equivalent graphs leaves equivalent graphs, the first group is the same for $D_\wedge$ and $D_\vee$. Hence the tangent spaces agree if and only if $\mathcal R_H+\R L_\wedge=\mathcal R_H+\R L_\vee$, where $\mathcal R_H$ is the span of the third group for $H$ and $L_\wedge,L_\vee$ are the rows of the generator of $e$: $\mathrm{Cov}(\cdot,z)$ if $z\to x_e$, and $\mathrm{Cov}(\cdot,\varepsilon_z)$ if $x_e\to z$. Identify $\mathrm{Cov}(\cdot,\ell)$ with the linear combination $\ell$ of the variables, which is injective, and write $r(z\mid S)$ for the residual of $z$ on $S$. Under $P^\circ$, $\mathcal R_H$ meets the span of $w,z,x_f$ in $\R\ell_w$, with $\ell_w=0$ if $w\notin\adj(x_e)$, $\ell_w=w$ if $w\to x_e$ in $H$, and $\ell_w=r(w\mid\pa_H(w)\cap\{z,x_f\})$ if $x_e\to w$ in $H$; and $L_\wedge,L_\vee$ are $z$ or $r(z\mid\pa_D(z)\cap\{w,x_f\})$. There are four cases.
(1) $w\notin\adj(x_e)$ and $z\to w$ in $G_\wedge$. Then $w\in\mathrm{CC}_{D_\wedge}(z,x_f)$, so $w\notin\mathrm{CC}_H(z,x_f)$ and $z,x_f$ are dependent; $L_\wedge=r(z\mid x_f)$. Neither $x_e\to z\leftarrow x_f$ nor $x_e\to z\leftarrow w$ is a v-structure of $G_\vee$, so $L_\vee=z$; $\ell_w=0$, and $r(z\mid x_f)\notin\R z$.
(2) $w\notin\adj(x_e)$ and $w\to z$ in $G_\wedge$. Then $x_e\to z\leftarrow w$ is a v-structure of $G_\wedge$, hence of $G_\vee$ by (F), so $D_\vee$ has $x_e\to z\leftarrow w$ and $z\to x_f$. Now $w\in\mathrm{CC}_H(z,x_f)$, so $z,x_f$ are dependent given $w$; $L_\wedge=r(z\mid w,x_f)$ has a non-zero coefficient on $x_f$, $L_\vee=r(z\mid w)$ has none, and $\ell_w=0$.
(3) $w\in\adj(x_e)$ and $w\in\mathrm{CC}_H(x_e,z)$. Then $z\to w$ in $D_\wedge$ and in $D_\vee$, so $w\in\mathrm{CC}_{D_\wedge}(z,x_f)$ and $w\notin\mathrm{CC}_H(z,x_f)$; in $H$, $x_e\to w\leftarrow z$ and hence $w\to x_f$, and $z,x_f$ are dependent. $L_\wedge=r(z\mid x_f)$, $L_\vee=z$ and $\ell_w=r(w\mid z)$. A combination of $L_\wedge$ and $L_\vee$ has coefficient zero on $w$, so it lies in $\R\ell_w$ only if it vanishes, which the coefficient of $r(z\mid x_f)$ on $x_f$ excludes.
(4) $w\in\adj(x_e)$ and $w\notin\mathrm{CC}_H(x_e,z)$. Then $w\to z$ in $D_\wedge$, and in $D_\vee$ if $x_e\to z$, since $z\to w$ would force $x_e\to w$; so $w\notin\mathrm{CC}_{D_\wedge}(z,x_f)$ and $w\in\mathrm{CC}_H(z,x_f)$; in $H$, $z\to w\leftarrow x_f$ and hence $w\to x_e$, so $\ell_w=w$, and $z,x_f$ are dependent given $w$. $L_\wedge=r(z\mid w,x_f)$ has a non-zero coefficient on $x_f$, while $L_\vee\in\{z,r(z\mid w)\}$ and $\ell_w$ have none.
In each case $\mathcal R_H+\R L_\wedge\neq\mathcal R_H+\R L_\vee$, so $I(P^\circ)>0$, and $I(P_{\lambda=0})>0$ for almost every $\phi_0'$ by analyticity.

\emph{(iv).} Let the pair satisfy (F) and (C), so that $I(P_{\lambda=0})=0$ by (iii), let $P_0=P_{\lambda=0}$ and let $T=T_{P_0}\cM_{G^*}=T_{P_0}\cM_{G'}$. For a DAG $D$ let $\theta^{\mathrm{adj}}_D(P)$ be the coefficient of $X$ in the regression of $Y$ on $\{X\}\cup\pa_D(X)$, defined for every $P$, so that $\theta^{\mathrm{adj}}_D=\theta_D$ on $\cM_D$. By \cref{lem:Dproj}, the derivative of $\Pi_{G'}$ at $P_0\in\cM_{G'}$ is the projection onto $T$, which fixes the velocity $v_e\in T$ of the $e$-curve; so $\delta=d(\theta^{\mathrm{adj}}_{G'}-\theta^{\mathrm{adj}}_{G^*})[v_e]$. Let $D^*,D'\in\Xeq$ be equivalent to $G^*,G'$. By \eqref{eq:K1p}, $\theta^{\mathrm{adj}}_{D^*}=\theta_{G^*}=\theta^{\mathrm{adj}}_{G^*}$ on $\cM_{G^*}$, and likewise for $D'$, so the differentials of $\Delta\theta=\theta^{\mathrm{adj}}_{D'}-\theta^{\mathrm{adj}}_{D^*}$ and of $\theta^{\mathrm{adj}}_{G'}-\theta^{\mathrm{adj}}_{G^*}$ agree on $T$; by \eqref{eq:K2} they vanish on $T_{P_0}\cM_H$, and $T=T_{P_0}\cM_H\oplus\R v_e$. Hence $\delta=c\,d\Delta\theta[v]$ with $c\ne0$ for every $v\in T\setminus T_{P_0}\cM_H$. Write $D_\wedge$ for the one of $D^*,D'$ that is equivalent to $G_\wedge$ and $D_\vee$ for the other, and choose them below so that $(D_\wedge-e)-f\simeq H-f$, which gives $P_0\in\cM_K$ for $K=D_\wedge-f$. We take $v$ to be the velocity of the curve $P_\varepsilon\in\cM_K$ that sets the coefficient of $e$ in $K$ to $\varepsilon$ and keeps the others at their values at $P_0$; then $v\in T\setminus T_{P_0}\cM_H$, because $\cM_K\subseteq\cM_{D_\wedge}$ and the structural parametrization of $\cM_{D_\wedge}$ is an immersion. On $\cM_K$, $\theta^{\mathrm{adj}}_{D_\wedge}=\theta_{D_\wedge}=\theta_K$, because the coefficient of $f$ vanishes there.

First, run the descent of the proof of \cref{lem:onev} with $G_1=G_\wedge$, $G_2=G_\vee$, $a=x_e$, $b=x_f$, from $D=G_\vee$, whose hypothesis holds because $\mathrm{CC}_{G_\vee}(x_e,z)=\mathrm{CC}_H(x_e,z)=\mathrm{CC}_{G_\wedge}(x_e,z)$. The edges at $Y$ and in $R$ are oriented alike in $G_\wedge$ and $G_\vee$ and are never reversed, so the final graph $D_\vee$ and the graph $D_\wedge\simeq G_\wedge$ obtained from it by the last reversal satisfy $\cK$ and have the same common children of $x_e$ and $z$ as $G_\wedge$; hence they lie in $\Xeq$ by the paragraph on common children, and $(D_\wedge-e)-f\simeq H-f$ by (C). If the descent stops in type $\alpha$, then $D_\vee=D_\wedge$ with $f$ reversed, so $D_\vee-f=K$, $\theta^{\mathrm{adj}}_{D_\vee}=\theta_K$ on $\cM_K$, $\Delta\theta$ vanishes on $\cM_K$, and $\delta=0$. This happens whenever there is no $e$-configuration. It may also happen when one exists; then both computations below are valid, so $\delta=0$, and hence (D) fails, there as well.

Second, let $D_\vee$ be an $e$-configuration and $D_\wedge$ the graph obtained by reversing $e$, as in \cref{lem:Econf}(i). Then $\pa_K(x_e)=U$ and $\pa_K(z)=U\cup\{x_e\}$, so $e$ is covered in $K$, $K'=D_\vee-f$ is $K$ with $e$ reversed, $\cM_{K'}=\cM_K$, and $\theta^{\mathrm{adj}}_{D_\vee}=\theta_{K'}$ on $\cM_K$. The two parent sets of $X$ in $K$ and $K'$ differ only if $X\in e$. If $X=z$, then on $\cM_K$
\[
\Delta\theta=\beta_{YX\cdot U}-\beta_{YX\cdot U\cup\{x_e\}}=\beta_{Yx_e\cdot U\cup\{X\}}\,\beta_{x_eX\cdot U},
\]
and $\beta_{x_eX\cdot U}=\varepsilon\Var(x_e\mid U)/\Var(X\mid U)$ along the curve; if $X=x_e$, then $\Delta\theta=-\beta_{Yz\cdot U\cup\{X\}}\,\beta_{zX\cdot U}$ with $\beta_{zX\cdot U}=\varepsilon$. So $\delta$ is a non-zero multiple of $\beta_{Yw\cdot U\cup\{X\}}(P_0)$ if $X\in e$ and is zero otherwise. As a function of $\phi_0'$, $\beta_{Yw\cdot U\cup\{X\}}(P_0)$ is rational on $\Theta_{H-f}$; it vanishes identically if $Y$ and $w$ are d-separated given $U\cup\{X\}$ in $H-f$, and otherwise it is not identically zero, because the parameters of a linear Gaussian model on a DAG at which the distribution is not faithful to the DAG form a Lebesgue-null set \citep[Theorem~3.2]{spirtes2000causation}; hence it is non-zero for almost every $\phi_0'$. With \cref{lem:Econf}(ii), this proves (iv).

\emph{Order of vanishing.} Suppose $I(0)$ vanishes for almost every, hence every, $\phi_0'$. By (iii) the pair satisfies (F), and by \cref{thm:E} it has a witness configuration, whose witness edge is $f$. Take $\phi_0'=0$, which is not a generic point but suffices to show that an analytic function is not identically zero: all coefficients of $H$ other than that of $f$ are zero and all variances are one. Then $P_\lambda$ is the law in which only the endpoints of $f$ are dependent, and its parameter under $G^*$ has a single non-zero coefficient $\mu(\lambda)$ on $f$, with $\mu(0)=0$ and $\mu'(0)=1$, and variances that equal one at $\lambda=0$ and are $1+O(\lambda^2)$ near it, so they do not contribute to first derivatives in $\lambda$ (if $f$ points in opposite directions in $H$ and $G^*$, then $\mu=\lambda/(1+\lambda^2)$ and the two variances are $1/(1+\lambda^2)$ and $1+\lambda^2$). Let $c_w$ be the witness partial covariance of $G'$. It vanishes identically when $\mu=0$, because a law in which only the endpoints of $e$ are dependent belongs to $\cM_{G'}$. At unit variances the table of \cref{app:E} gives $\partial_sc_w|_{s=0}=-\beta/(1+\beta^2)$ or $\beta$. Hence $\partial_\lambda\partial_sc_w|_{\lambda=s=0}=\mp1\ne0$, so $\dot c_w(\lambda)=\pm\lambda+O(\lambda^2)$ and, by \cref{lem:Dproj}, $I(\lambda)\ge c\lambda^2$ at $\phi_0'=0$. Writing $I(\lambda)=\lambda^2I_2(\phi_0')+O(\lambda^3)$, the coefficient $I_2=\sum_j(\partial_\lambda\dot c_j)^\top M_j^{-1}(\partial_\lambda\dot c_j)/\sigma^2_{j\mid\pa(j)}$ is analytic in $\phi_0'$ and positive at $\phi_0'=0$, so it is positive for almost every $\phi_0'$. Thus, outside a null set of $\phi_0'$, every pair has either $I(0)>0$ or $I(0)=0$ and $I(\lambda)\ge c\lambda^2$ for small $\lambda$.

\emph{(i) and (ii).} For each pair, $\delta$ is bounded, $s_{\mathrm{eff}}$ is bounded away from zero, and $\sqrt{I(\lambda)}\ge c$ or $\ge c|\lambda|$; so $K_{G^*,G'}(P_\lambda)\le C/|\lambda|$, and the maximum over the finitely many pairs gives (i). If some pair has $I(0)=0$ and $\delta(0)\ne0$, then $I\ge0$ is smooth with $I(0)=0$, so $I'(0)=0$ and $I(\lambda)\le C_I\lambda^2$, and $K(P_\lambda)\ge K_{G^*,G'}(P_\lambda)\ge|\delta(0)|/(2\bar s\sqrt{C_I}\,|\lambda|)$ for small $\lambda$. Otherwise every pair has $I(0)>0$, so that $K_{G^*,G'}$ is bounded, or $I(0)=0$ and $\delta(0)=0$, so that $\delta(\lambda)=O(\lambda)$ and $K_{G^*,G'}(P_\lambda)=O(\lambda)/\sqrt{c\lambda^2}=O(1)$. The statement about the price follows from \cref{thm:K} at each $P_\lambda$.
\end{proof}

\begin{figure}[t]
\centering
\begin{tikzpicture}[>={Stealth[length=2mm]},v/.style={circle,draw,inner sep=1.2pt,minimum size=5.5mm},every node/.style={font=\small}]
\begin{scope}
\node[v] (W) at (0,1.2) {$W$};
\node[v] (X) at (1.4,1.2) {$X$};
\node[v] (Z) at (2.8,1.2) {$Z$};
\node[v] (Y) at (2.1,0) {$Y$};
\draw[->,very thick] (W) -- node[above] {$f$} (X);
\draw[dashed] (X) -- node[above] {$e$} (Z);
\draw[->] (X) -- (Y);
\draw[->] (Z) -- (Y);
\node at (1.4,-0.75) {(a) $C$: $Z\to X$;\ $M$: $X\to Z$};
\end{scope}
\begin{scope}[xshift=5.4cm]
\node[v] (A) at (0,2.0) {$A$};
\node[v] (B) at (0,0) {$B$};
\node[v] (X2) at (1.6,1.0) {$X$};
\node[v] (M) at (3.4,2.0) {$M$};
\node[v] (Y2) at (3.4,0) {$Y$};
\draw[->] (A) -- (X2);
\draw[->,very thick] (B) -- node[below right=-2pt] {$f$} (X2);
\draw[->] (A) -- (M);
\draw[dashed] (X2) -- node[below right=-2pt] {$e$} (M);
\draw[->] (X2) -- (Y2);
\draw[->] (M) -- (Y2);
\node at (1.7,-0.75) {(b) $G_\wedge$: $M\to X$;\ $G_\vee$: $X\to M$};
\end{scope}
\end{tikzpicture}
\caption{Two instances of \cref{thm:K-lambda}. The weak edge $e$ is dashed and the margin edge $f$, whose coefficient $\lambda$ tends to zero, is thick. (a) The four-variable model of \cref{sec:toy}, with $R=\{W\to X\}$: $Z$ is a confounder under $C$ and a mediator under $M$. (b) A five-node instance with only $Y$ a sink as background knowledge: the explanations differ only in the v-structure $M\to X\leftarrow B$ formed by $e$ and $f$, which is (F); $X$ and $B$ have no common neighbour, which gives (C); and $X$ is an endpoint of $e$, $A$ can be the only parent of $M$, and $M\to Y$ connects $M$ to $Y$ given $\{A,X\}$, which is (D). So $K(P_\lambda)\asymp1/\lambda$ in both, the factor being $b/w$ in (a).}
\label{fig:instances}
\end{figure}
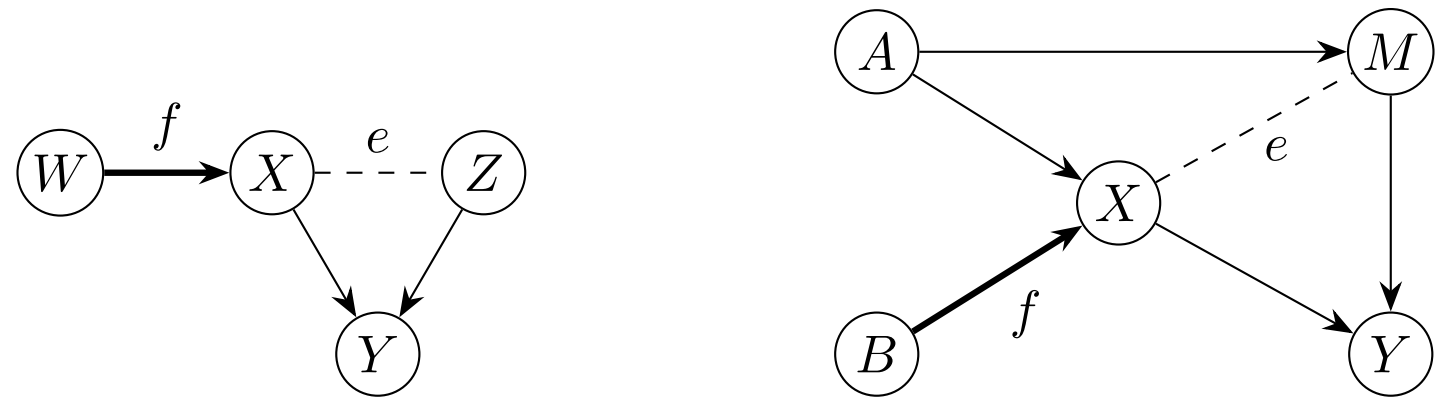

\paragraph{An instance without required edges.} On five nodes (\cref{fig:instances}(b)) let $H$ have the edges $A\to X$, $B\to X$, $A\to M$, $X\to Y$ and $M\to Y$, let $e=X-M$, and let $Y$ be a sink, with no other background knowledge. The explanations are $G_\vee=H+(X\to M)$, in which $M$ mediates, and the members of the class of $G_\wedge=H+(M\to X)$, in which $M$ is a parent of $X$ and $M\to X\leftarrow B$ is the v-structure that separates the two. For the margin edge $f=B\to X$ the pair satisfies (F), (C) and (D) with $x_e=M$, $z=X$, $x_f=B$, $U=\{A\}$ and $w=M$, so $K(P_\lambda)\asymp1/\lambda$; with $f=A\to X$ instead, (F) fails and $K(P_\lambda)$ stays bounded. \Cref{app:sim} gives the values of $\lambda K(P_\lambda)$ at one base point and a simulation along $\lambda$.

\begin{lemma}[Regularization in the target]\label{lem:ureg}
Let $O$ be a random interval, possibly empty, that depends on the data and on external randomization, let $\varepsilon>0$, and let
\[
\Psi(u)=\E\Bigl[\frac1\varepsilon\int_{u-\varepsilon/2}^{u+\varepsilon/2}1\{v\in O+[-\varepsilon,\varepsilon]\}\,dv\ \Big|\ \text{data}\Bigr].
\]
Then (i) $\Psi$ takes values in $[0,1]$ and is $\varepsilon^{-1}$-Lipschitz in $u$; (ii) $\Psi(u)\ge\Prob(u'\in O\mid\text{data})$ whenever $|u-u'|\le\varepsilon/2$; (iii) $\int\Psi(u)\,du\le\E[\Leb(O)\mid\text{data}]+2\varepsilon$. Conversely, every measurable $\Psi$ with values in $[0,1]$ defines the randomized set $\{u:\Psi(u)\ge V\}$, $V\sim U(0,1)$ independent, which contains $u$ with probability $\E\Psi(u)$ and has expected Lebesgue measure $\E\int\Psi$.
\end{lemma}
\begin{proof}
(i) is clear. (ii) If $u'\in O$, then $[u-\varepsilon/2,u+\varepsilon/2]\subseteq[u'-\varepsilon,u'+\varepsilon]\subseteq O+[-\varepsilon,\varepsilon]$, so the average is one; take the conditional expectation. (iii) By Fubini, the integral of the average is $\Leb(O+[-\varepsilon,\varepsilon])$, which is $\Leb(O)+2\varepsilon$ for a non-empty interval and $0$ otherwise. The converse is Fubini's theorem.
\end{proof}

\paragraph{The two directions of $K$.} In the notation of \cref{sec:gen-constant}, the proof of \cref{prop:canonical} below gives
\[
K_{1\to2}=\frac{|\alpha_2\rho-\alpha_1|}{\sqrt{a_0^2+\alpha_1^2}\,\sqrt{1-\rho^2}},\qquad K_{2\to1}=\frac{|\alpha_1\rho-\alpha_2|}{\sqrt{a_0^2+\alpha_2^2}\,\sqrt{1-\rho^2}} .
\]
In the four-variable model with unit variances at $a=0$, $K_{C\to M}=b/w$ and $K_{M\to C}=b\sqrt{1+w^2}/(w\sqrt{1+b^2})$, and the dual lower bound of \cref{tab:price-curve} is not a function of $K$ alone ($1.826$ at $(w,b)=(0.25,1)$ with $K=4$, $1.909$ at $(2,2)$ with $K=1$). If $\Xeq$ has three or more Markov classes, the full limit experiment has more than two lines and is not canonical, and $\Pi^*_{\mathrm{can}}$ of a pair is a lower bound; the values of \cref{tab:price-curve} are lower bounds, or approximations, for the limit price at all benign base points with canonical parameters $(1,b,0,\sqrt v)$.

\paragraph{The identity \eqref{eq:R-identity}.} With $Q=a_0^2(1-\rho^2)-\rho^2(\alpha_1^2+\alpha_2^2)+2\alpha_1\alpha_2\rho$, a direct computation from the formulas for $K_{1\to2}$ and $K_{2\to1}$ above gives $1-K_{1\to2}^2=Q/((a_0^2+\alpha_1^2)(1-\rho^2))$ and the same with the indices exchanged. Dividing gives \eqref{eq:R-identity} when $Q\ne0$, and since $1-K_{1\to2}^2$ and $1-K_{2\to1}^2$ both have the sign of $Q$, $K_{1\to2}>1$ if and only if $K_{2\to1}>1$. At the degenerate point $K_{1\to2}=K_{2\to1}=1$, $Q=0$ and the ratio is free.

\begin{proof}[Proof of \cref{prop:canonical}]
\emph{Tangent spaces.} $\cM_{G^*}$ has dimension $\dim\cM_H+1$, because the coefficient of $e$ is one free parameter, and so does $\cM_{G'}$; hence $L_0$ has codimension one in $L_1$ and in $L_2$, and $L_0\subseteq L_1\cap L_2$. The velocity $v_e$ of the $e$-curve of $G^*$ lies in $L_1\setminus L_0$, and by \cref{thm:E} it is not in $L_2$ (some local Markov partial covariance of $G'$ has a non-zero derivative along it, \cref{lem:Dproj}), so $L_1\cap L_2$, which contains $L_0$ and is a proper subspace of $L_1$, equals $L_0$. By \eqref{eq:K2}, $\theta_{G^*}=\theta_{G'}=\theta_H$ on $\cM_H$, so $a_1$ and $a_2$ define the same functional on $L_0$.

\emph{The constant $K$.} Write $v_e=l_0+cn_1$ with $l_0\in L_0$ and $c\neq0$. The $g$-orthogonal projection of $v_e$ onto $L_2$ is $l_0+c\rho n_2$, so by \cref{lem:Dproj} $v^\perp=c(n_1-\rho n_2)$ and $I=c^2(1-\rho^2)$, and $\delta=\langle a_2,l_0+c\rho n_2\rangle-\langle a_1,l_0+cn_1\rangle=c(\alpha_2\rho-\alpha_1)$, since $a_1,a_2$ agree on $L_0$. Hence $|\delta|/\sqrt I=|\alpha_2\rho-\alpha_1|/\sqrt{1-\rho^2}$, and dividing by $s_{\mathrm{eff},1}=\sqrt{a_0^2+\alpha_1^2}$ gives $K_{1\to2}$; $K_{2\to1}$ is the same computation with the roles exchanged.

\emph{(i) The limit experiment.} For the lower bound we may restrict the class to laws with mean zero: honesty over the class implies honesty over this subclass, and the price over the class is at least that over it. We use a common chart for the two models: points of $\cM_H$ are parametrized through $\Theta_H$, and the normal coordinate through the coefficient of $e$, scaled so that the velocity has normal component $n_g$; at $\tau=0$ both models give the same law. The model of centred Gaussian laws is a regular exponential family, and the charts are smooth, so the experiments $(P_{n,h}^n:h\in L_1\cup L_2)$, $P_{n,h}=P_{h/\sqrt n}$, are locally asymptotically normal with one central sequence: $\log dP^n_{n,h}/dP^n_{n,0}=\langle h,\Delta_n\rangle-\frac12|h|^2+o_P(1)$ uniformly on bounded sets of $h$, where $\Delta_n$ is the orthogonal projection onto $L_1+L_2$ of the normalized score of the model of centred Gaussian laws and $\Delta_n\rightsquigarrow N(0,I)$ under $P_{n,0}$; the limit experiment is $Z\sim N(h,I)$ on $L_1+L_2$, $h\in L_1\cup L_2$, with target $\langle a_g,h\rangle$ \citep[Chapter~7]{vaart1998asymptotic}. Write $h=h_0+\tau n_g$ with $h_0\in L_0$, and let $\eta$ be the coordinate of $h_0$ along $\bar a/a_0$ (along any unit vector of $L_0$ if $a_0=0$), where $\bar a$ is the common projection of $a_1,a_2$ onto $L_0$, and $h_0'$ its component orthogonal to $\bar a$ in $L_0$. The coordinates of $Z$ along $\bar a/a_0$, $n_1$ and the unit vector $n_1^\perp$ of $\mathrm{span}(n_1,n_2)$ orthogonal to $n_1$ are independent with means $\eta$, $\tau\langle n_g,n_1\rangle$, $\tau\langle n_g,n_1^\perp\rangle$; the coordinates in $L_0\ominus\R\bar a$ have mean $h_0'$, which does not enter the target. In the subexperiment $h_0'=0$ the coordinates other than the first three are pure noise, that is, external randomization, so its minimax price over randomized sets is $\Pi^*_{\mathrm{can}}$. The minimax price of the full limit experiment is at least that of the subexperiment, because honesty over the full class implies honesty over the subclass. Conversely, a randomized set that is a function of the first three coordinates and external randomization has coverage and expected length that do not depend on $h_0'$, since neither the law of these coordinates nor the target does; so the two minimax prices are equal. The numerical values of $\Pi^*_{\mathrm{can}}$ in \cref{sec:gen-constant} are computed over sets that are equivariant under the translations $\eta\mapsto\eta+c$, which leave lengths and oracles invariant and move the target by $a_0c$, that is, over sets $a_0Z_0+A(Z)$; this is the reduction of \cref{app:lp} (Hunt--Stein, \citealp{kiefer1957invariance,lecam1986asymptotic}) and is not used in the proof.

\emph{(ii) Transfer to finite samples.} Let $\mathbb F_\varepsilon$ be the set of $\varepsilon^{-1}$-Lipschitz functions $\R\to[0,1]$ with the topology of uniform convergence on compact sets. By the Arzel\`a--Ascoli theorem it is a compact metrizable space, and for every $u$ and $B$ the maps $f\mapsto f(u)$ and $f\mapsto\int_{-B}^Bf$ are continuous. Let $C_n$ be asymptotically honest intervals, so that, by the definition in \cref{sec:setup}, which is uniform over the class, $\inf_{\cM_k}\Prob(\text{cover})\ge1-\alpha_n$ with $\alpha_n'=\max(\alpha_n,\alpha)\to\alpha$. Let $L=\liminf_n\Pi^{\mathrm{loc}}_n(C_n;r)$ and pass to a subsequence along which $\Pi^{\mathrm{loc}}_n\to L$. Fix $M>0$, $\varepsilon>0$, and let $K_M$ be the set of $(g,h)$ with $h\in L_g$, $|h|\le M$.
(a) By the class membership argument of the proof of \cref{thm:K}(i), the pairs $(G_g,P_{n,h})$, $(g,h)\in K_M$, belong to $\cM_k$ and satisfy $\|\Sigma_P-\Sigma_{P_0}\|\le r$ for $n$ large, uniformly on $K_M$. By smoothness, uniformly on $K_M$, $u_n(g,h)=\sqrt n(\theta_{G_g}(P_{n,h})-\theta_0)\to\langle a_g,h\rangle$ and $\sqrt n\,\orl_n(G_g,P_{n,h})\to2z\,s_g$ with $s_g=\sqrt{a_0^2+\alpha_g^2}$.
(b) Let $O_n=\sqrt n(C_n-\theta_0)$ and let $\Psi_n$ be its regularization from \cref{lem:ureg}, a statistic with values in $\mathbb F_\varepsilon$. Since $\mathbb F_\varepsilon$ is compact, $(\Psi_n,\Delta_n)$ is tight under $P_{n,0}$, so along a further subsequence it converges in distribution to some $(\Psi,\Delta)$ with $\Delta\sim N(0,I)$. The limit $\exp(\langle h,\Delta\rangle-\frac12|h|^2)$ of the likelihood ratios has mean one, so $P^n_{n,h}$ and $P^n_{n,0}$ are mutually contiguous by Le Cam's first lemma \citep[Lemma~6.4]{vaart1998asymptotic}. By Le Cam's third lemma \citep[Theorem~6.6]{vaart1998asymptotic}, whose proof applies verbatim to random elements of a Polish space, under $P_{n,h}$ the pair converges to the law with density $\exp(\langle h,\delta\rangle-\frac12|h|^2)$ with respect to the law of $(\Psi,\Delta)$. In that law $\Delta\sim N(h,I)$, and since the density depends on $\delta$ only, the conditional law of $\Psi$ given $\Delta=\delta$ is a kernel $\kappa(\delta,\cdot)$ that does not depend on $h$ (it exists because $\mathbb F_\varepsilon$ is Polish). So $\Psi\sim\kappa(Z,\cdot)$ is a randomized statistic of the limit experiment with values in $\mathbb F_\varepsilon$, and $\Psi_n\rightsquigarrow\Psi$ under every $h\in L_1\cup L_2$.
(c) For $(g,h)\in K_M$ and $n$ large, $|u_n(g,h)-\langle a_g,h\rangle|\le\varepsilon/2$, so \cref{lem:ureg}(ii) gives $\E_{n,h}\Psi_n(\langle a_g,h\rangle)\ge1-\alpha_n'$, and by continuity of evaluation $\E_h\Psi(\langle a_g,h\rangle)\ge1-\alpha$. For every $B$, by \cref{lem:ureg}(iii), $\E_h\int_{-B}^B\Psi=\lim_n\E_{n,h}\int_{-B}^B\Psi_n\le\liminf_n\sqrt n\,\E_{n,h}\Leb(C_n)+2\varepsilon\le2z\,s_g L+2\varepsilon$, using $\sqrt n\,\E_{n,h}\Leb(C_n)\le\Pi^{\mathrm{loc}}_n\sqrt n\,\orl_n(G_g,P_{n,h})$; letting $B\to\infty$, the same bound holds for $\E_h\int\Psi$. The bound is obtained for each $(g,h)$ and then the supremum is taken, so no exchange of limits is needed. By \cref{lem:ureg}, $\Psi$ defines a randomized set of the limit experiment that is honest on $K_M$ and has price at most $L+\varepsilon/(z\min_gs_g)$ there. The canonical parameters $(g,\eta,\tau)$ with $|\eta|,|\tau|\le M/2$, which we call the box of $M$, lie in $K_M$ (with $h_0'=0$), so $L+\varepsilon/(z\min_gs_g)\ge\Pi^*_{\mathrm{can},M}$, the minimax price over randomized sets of the canonical experiment that are honest on this box. Letting $\varepsilon\to0$ gives $L\ge\Pi^*_{\mathrm{can},M}$.
(d) $\Pi^*_{\mathrm{can},M}$ is non-decreasing in $M$ and at most $\Pi^*_{\mathrm{can}}$; it remains to show $\Pi^*_{\mathrm{can}}\le\lim_M\Pi^*_{\mathrm{can},M}$. The sets in this step need not be intervals, so \cref{lem:ureg}(iii) is not available, and we smooth without enlarging. Write $\Pi^*_{\mathrm{can}}[\alpha']$ and $\Pi^*_{\mathrm{can},M}[\alpha']$ for the minimax prices when the required coverage is $1-\alpha'$ and the oracle lengths are kept, and let $X=(Z_0,Z)$, $\theta=(g,\eta,\tau)$, $t(\theta)=a_0\eta+\alpha_g\tau$ and $s=\min_gs_g$. For each $M$ take a randomized set $S_M$ that is honest on the box of $M$ with price at most $\Pi^*_{\mathrm{can},M}[\alpha]+1/M$ there, let $\Psi^0_M(x,u)=\Prob(u\in S_M\mid X=x)$ and $\Psi_M(x,u)=\varepsilon^{-1}\int_{-\varepsilon/2}^{\varepsilon/2}\Psi^0_M(x,u-r)\,dr$. Then $\Psi_M(x,\cdot)\in\mathbb F_\varepsilon$ and $\int\Psi_M(x,u)\,du=\int\Psi^0_M(x,u)\,du$, so expected lengths are unchanged. For coverage, $t(\theta)-r=t(\theta_r)$ with $\theta_r=\theta-r(0,a_0,\alpha_g)/s_g^2$, a parameter of the same model whose mean differs from that of $\theta$ by $|r|/s_g$ in Euclidean norm, so the total variation distance between the laws of $X$ under $\theta$ and $\theta_r$ is at most $|r|/(\sqrt{2\pi}\,s_g)\le\varepsilon/(2s)$ for $|r|\le\varepsilon/2$; each of $\eta$ and $\tau$ moves by at most $\varepsilon/(2s)$, so $\theta_r$ is in the box of $M$ when $\theta$ is in the box of $M-\varepsilon/s$, and hence $\E_\theta\Psi_M(X,t(\theta))\ge1-\alpha-\varepsilon/(2s)$ for $\theta$ in the box of $M-\varepsilon/s$. The canonical experiment does not change with $M$: if $\ell_\theta$ is the density of the law of $X$ under $\theta$ with respect to that under $\theta_0=(1,0,0)$, then $\E_\theta f(\Psi_M,X)=\E_{\theta_0}[f(\Psi_M,X)\ell_\theta(X)]$ for bounded continuous $f$. The pairs $(\Psi_M,X)$ are tight under $\theta_0$, so along a subsequence they converge in distribution to some $(\Psi,X)$, and since $\ell_\theta$ is continuous and the law of $\ell_\theta(X)$ does not depend on $M$, so that the integrands are uniformly integrable, $\E_\theta f(\Psi_M,X)\to\E_{\theta_0}[f(\Psi,X)\ell_\theta(X)]$ for every $\theta$. As in (b), $\Psi$ is a randomized statistic of the canonical experiment with values in $\mathbb F_\varepsilon$ and $\Psi_M\rightsquigarrow\Psi$ under every $\theta$. As in (c), by continuity of evaluation and of truncated integrals, $\Psi$ covers $t(\theta)$ with probability at least $1-\alpha-\varepsilon/(2s)$ at every $\theta$, which lies in the box of $M-\varepsilon/s$ for $M$ large, and has price at most $\lim_M\Pi^*_{\mathrm{can},M}[\alpha]$; by \cref{lem:ureg} it defines a randomized set, so $\Pi^*_{\mathrm{can}}[\alpha+\varepsilon/(2s)]\le\lim_M\Pi^*_{\mathrm{can},M}[\alpha]$. Finally, $\alpha'\mapsto\Pi^*_{\mathrm{can}}[\alpha']$ is convex on $(0,1)$, because a mixture by external randomization of two sets that are honest at levels $\alpha'$ and $\alpha''$ is honest at the mixed level and its expected length is the mixture of theirs, and it is finite, because the union of the two oracle intervals $a_0Z_0+\alpha_g\langle Z,n_g\rangle\pm z_{1-\alpha'/2}\,s_g$, $g=1,2$, is honest at level $\alpha'$ with finite price. So it is continuous, and letting $\varepsilon\to0$ gives $\Pi^*_{\mathrm{can}}[\alpha]\le\lim_M\Pi^*_{\mathrm{can},M}[\alpha]$. With (c), $L\ge\Pi^*_{\mathrm{can}}$.

\emph{Four-variable case.} With $E=D-wP$, the pair $(P,E)$ has independent coordinates with variance $v$; dividing by $\sqrt v$ gives unit variance. Under $M$ the mean is $(t,0)/\sqrt v=\tau_Mn_M$ with $n_M=(1,0)$ and $\tau_M=t/\sqrt v$; under $C$ it is $(t,-wt)/\sqrt v=\tau_Cn_C$ with $n_C=(1,-w)/\sqrt{1+w^2}$ and $\tau_C=t\sqrt{1+w^2}/\sqrt v=t/v$. So $\rho=\langle n_M,n_C\rangle=1/\sqrt{1+w^2}=\sqrt v$. The target offsets are $bt=b\sqrt v\,\tau_M$ under $M$ and $0$ under $C$, and $S_\beta$ has variance $v=a_0^2$, so $(a_0,\alpha_1,\alpha_2,\rho)=(\sqrt v,b\sqrt v,0,\sqrt v)$, and dividing by $a_0$ gives $(1,b,0,\sqrt v)$. For $G^*=M$, $G'=C$ the formulas for $K_{1\to2}$ and $K_{2\to1}$ reduce to $b\sqrt{1+w^2}/(w\sqrt{1+b^2})$ and $b/w$, as in \cref{sec:gen-constant}.
\end{proof}

\begin{proof}[Proof of \cref{prop:Kpos}]
Fix $G\in\Xeq$ with $e=u\to v$, $D=G-e\in[H]_Y$ (\cref{prop:graph}(a)) and $W=\pa_D(v)$. At $P_0$ the parameter of $G$ is that of $D$ with $\lambda_e=0$. Write $\varepsilon_j$ and $\omega_j$ for the noises and noise variances of $D$, $\tau_{a\to b}$ for total effects in $D$, $r_u=u-\sum_{p\in W}b_px_p$ for the residual of $u$ on $W$, and identify tangent vectors with their scores, so that $g(w,w')=\E[\ell_w\ell_{w'}]$. Ancestors of a vertex include the vertex itself, and $u$ and the vertices of $W$ are non-descendants of $v$ in $D$, because $G=D+(u\to v)$ is acyclic.

\emph{Normal direction.} The curve in $\cM_G$ that sets $\lambda_e=t$ and replaces $\lambda_{p\to v}$ by $\lambda_{p\to v}-tb_p$, $p\in W$, adds $tr_u$ to the equation of $v$; its velocity $w_G$ has score $r_u\varepsilon_v/\omega_v$. The tangent space $L_0=T_{P_0}\cM_D=T_{P_0}\cM_H$ is spanned by the scores $x_p\varepsilon_j/\omega_j$, $p\in\pa_D(j)$, and $(\varepsilon_j^2-\omega_j)/(2\omega_j^2)$. Hence $\E[r_u\varepsilon_v]=0$, and by Isserlis' formula $\E[r_u\varepsilon_v\,x_p\varepsilon_j]=\E[r_u\varepsilon_j]\,\E[\varepsilon_vx_p]$ for $j\ne v$. This vanishes: $\E[\varepsilon_vx_p]\ne0$ requires $v$ to be an ancestor of $p$ or $p=v$, hence a proper ancestor of $j$, and $\E[r_u\varepsilon_j]\ne0$ requires $j$ to be an ancestor of $u$ or of $W$, which closes a cycle with $u\to v$ in $G$ or with $W\to v$. For $j=v$ and $p\in W$, $\E[r_u\varepsilon_v^2x_p]=\omega_v\E[r_ux_p]=0$, and the variance scores give $2\E[\varepsilon_v\varepsilon_j]\E[r_u\varepsilon_j]=0$ for $j\ne v$ and a multiple of $\E[r_u\varepsilon_v]=0$ for $j=v$. So $w_G\perp L_0$, $\|w_G\|^2=\Var(r_u)/\omega_v>0$, and $T_{P_0}\cM_G=L_0\oplus\R w_G$ by dimension.

\emph{The normal component.} Since $\theta_G=[(I-\Lambda)^{-1}]_{YX}$ has derivative $\tau_{X\to a}\tau_{j\to Y}$ in $\lambda_{a\to j}$ at $P_0$, and $\Cov(\varepsilon_X,x_a)=\omega_X\tau_{X\to a}$,
\[
d\theta_G[w_G]=A_G:=\tau_{v\to Y}\Bigl(\tau_{X\to u}-\sum_{p\in W}b_p\tau_{X\to p}\Bigr)=\frac{\tau_{v\to Y}}{\omega_X}\Cov(\varepsilon_X,r_u).
\]
By \eqref{eq:K2}, $\theta_G=\theta_H$ on $\cM_H$, so the gradient of $\theta_G$ in $T_{P_0}\cM_G$ is $a_0+\alpha_Gn_G$, where $a_0$ is the gradient of $\theta_H$ in $L_0$, $n_G=w_G/\|w_G\|$ and $\alpha_G=A_G/\|w_G\|$. Along a curve in $\cM_G$ with velocity $n_G$, an estimator with influence function $a_0$ has, by Le Cam's third lemma, local bias $g(a_0,n_G)-d\theta_G[n_G]=-\alpha_G$; this is the statistical reading in \cref{sec:gen-constant}.

\emph{Reduction.} (Here $n_G=w_G/\|w_G\|$ is not normalized to $\rho\ge0$ as in \cref{prop:canonical}; the formula for $K_{1\to2}$ is invariant under $(n_2,\alpha_2,\rho)\mapsto(-n_2,-\alpha_2,-\rho)$.) The score $x_u\varepsilon_v/\omega_v$ of the $e$-curve of $G^*$ is $r_u\varepsilon_v/\omega_v+\sum_{p\in W}b_px_p\varepsilon_v/\omega_v$, so its velocity is $v_e=v_0+w_{G^*}$ with $v_0\in L_0$. By \cref{lem:Dproj} the derivative of $\Pi_{G'}P_\varepsilon$ is the projection of $v_e$ onto $T_{P_0}\cM_{G'}$, so, with $\rho=g(n_{G^*},n_{G'})$,
\[
\delta_{G^*,G'}=g(a_0+\alpha_{G'}n_{G'},v_e)-g(a_0+\alpha_{G^*}n_{G^*},v_e)=\|w_{G^*}\|(\rho\,\alpha_{G'}-\alpha_{G^*}),\qquad I_{G^*,G'}=\|w_{G^*}\|^2(1-\rho^2),
\]
and $|\rho|<1$ outside $\cN$. If all $\alpha_G$ vanish, every $\delta$ does and $K(P_0)=0$. If $\alpha_G\ne0$ and $G'$ is in another Markov class, then $\delta_{G,G'}=\delta_{G',G}=0$ would give $\alpha_G=\rho\alpha_{G'}=\rho^2\alpha_G$; so $K(P_0)>0$.

\emph{Graphical form.} $A_G$ is the product of $\tau_{v\to Y}$ and $\beta_G=\tau_{X\to u}-\sum_pb_p\tau_{X\to p}$, rational functions of the parameter of $D$, which is a rational function of $\phi_0$ because $\cM_D=\cM_H$. As a sum over directed paths, $\tau_{v\to Y}$ is not identically zero if and only if (a) holds. For (b), give $D^+$ the parameters of $D$ except that $X$ receives $\kappa\xi$ with $\Var\xi=\sigma^2$ and an own noise of variance $\omega'$; its marginal on $V$ is the model of $D$ with $\omega_X=\kappa^2\sigma^2+\omega'$, every parameter of $D$ arises in this way, and $\Cov(\xi,x_a)=\kappa\sigma^2\tau_{X\to a}$, so that the partial covariance of $\xi$ and $u$ given $W$ in $D^+$ is $\kappa\sigma^2\beta_G$. It vanishes identically if $\xi$ and $u$ are d-separated given $W$, and otherwise it is not identically zero, by the same faithfulness result applied to the linear Gaussian model on $D^+$ \citep[Theorem~3.2]{spirtes2000causation}; so $\beta_G\not\equiv0$ if and only if (b) holds. As $\Theta_H$ is connected, $A_G$ is then non-zero outside a null set, and the union of these null sets and $\cN$ is the exceptional set. Under \eqref{eq:K1p}, for $G_1\simeq G_2$ in $\Xeq$, $\theta_{G_1}=\thadj_{G_2}=\theta_{G_2}$ on the common model $\cM_c=\cM_{G_1}=\cM_{G_2}$, so the two members have the same gradient on the same tangent space $T_{P_0}\cM_c\supset L_0$. Since $T_{P_0}\cM_c\ominus L_0$ is a line, $|\alpha_{G_1}|=|\alpha_{G_2}|$ is the length of the component of this gradient orthogonal to $L_0$, and $A_G=\alpha_G\|w_G\|$ vanishes identically for both members or for neither; so (a) and (b) hold for all or none of them, although $e$ may have opposite orientations in them.
\end{proof}

\paragraph{The criterion of \cref{prop:Kpos} is not a path condition.} If $X\in W$ and $X$ has no parent in $D$, every path from $\xi$ leaves $X$ along an edge out of $X$ and is blocked there (if $X\in W$ has a parent $a$, the collider $\xi\to X\leftarrow a$ is opened instead), and if $X\notin W$ and $X\to c\leftarrow u$ with $c\in W$, the path $\xi\to X\to c\leftarrow u$ is open although $u$ need not be a descendant of $X$. So condition (b) neither implies nor is implied by $e$ lying on a directed path from $X$ to $Y$; on five nodes the two happen to coincide, on six they differ in both directions (\cref{app:certify}).

\paragraph{When honesty is free.} In the canonical experiment of \cref{prop:canonical} write $s_g=\sqrt{a_0^2+\alpha_g^2}$, $z=z_{1-\alpha/2}$, $T_g=a_0Z_0+\alpha_g\langle Z,n_g\rangle$, $\beta=\alpha_1n_1-\alpha_2n_2$, $c_\Lambda=\Lambda^2/(1+\Lambda^2)$ for $\Lambda>0$ and $c_f=(2\pi)^{-3/2}e^{-1}$, where $s_1,s_2>0$ are efficient standard deviations, and for $d\in\R$ let
\[
G(d)=\int\bigl|\ind\{|v-d|\le zs_1\}-\ind\{|v|\le zs_2\}\bigr|\,\min\{w_1(v-d),w_2(v)\}\,dv,\qquad w_g(v)=\bigl|1-e^{(z^2-v^2/s_g^2)/2}\bigr|,
\]
\[
\delta_{\mathrm{can}}=\sup_{\Lambda>0}\Bigl[\frac{c_f\Lambda\sqrt{c_\Lambda}}{z(s_1+s_2)(1+\Lambda^2)}\int_{|\zeta|\le1}G\bigl(\sqrt{c_\Lambda}\langle\beta,\zeta\rangle\bigr)\,d\zeta-\bigl(1-\sqrt{c_\Lambda}\bigr)\Bigr],
\]
with $\zeta\in\R^2$. Like $\Pi^*_{\mathrm{can}}$, $\delta_{\mathrm{can}}$ is invariant under a common scaling of $(a_0,\alpha_1,\alpha_2)$.

\begin{lemma}[Strict price in the canonical experiment]\label{lem:strict}
$\Pi^*_{\mathrm{can}}(a_0,\alpha_1,\alpha_2,\rho)\ge1+\delta_{\mathrm{can}}$ and $\delta_{\mathrm{can}}\ge0$. If $\alpha_1=\alpha_2=0$, then $\Pi^*_{\mathrm{can}}=1$; otherwise $\delta_{\mathrm{can}}>0$. Moreover $\alpha_1=\alpha_2=0$ if and only if $K_{1\to2}=K_{2\to1}=0$.
\end{lemma}

\begin{proof}
Let $S$ be a randomized confidence set that is honest over the canonical experiment with price at most $1+\delta$, and let $\psi(x,u)=\Prob(u\in S\mid X=x)$ with $X=(Z_0,Z)$, which is jointly measurable because a randomized set is given by a Markov kernel. By Fubini's theorem the coverage at $(g,\eta,\tau)$ is $\E\psi(X,t_g)$ with $t_g=a_0\eta+\alpha_g\tau$, and the expected length is $\E\int\psi(X,u)\,du\le(1+\delta)\,2zs_g$.

\emph{One model.} Fix $g$ and the prior $(\eta,\tau)\sim N(0,\Lambda^2I_2)$. The pair $(Z_0,\langle Z,n_g\rangle)$ is sufficient, and $\langle Z,n_g^\perp\rangle\sim N(0,1)$, with $n_g^\perp$ a unit vector orthogonal to $n_g$, does not depend on the parameter, so the posterior of $t_g$ is $N(c_\Lambda T_g,c_\Lambda s_g^2)$; let $p_g(u\mid x)$ be its density and $m_g$ the marginal law of $X$. Let $\psi^*_g(x,u)=\ind\{|u-c_\Lambda T_g(x)|\le z\sqrt{c_\Lambda}s_g\}=\ind\{p_g\ge k_g\}$ with $k_g=\varphi(z)/(\sqrt{c_\Lambda}s_g)$, and $\bar w_g=|1-p_g/k_g|$. Averaging the coverage over the prior gives $\E_{m_g}\int\psi p_g\,du\ge1-\alpha=\E_{m_g}\int\psi^*_gp_g\,du$, and pointwise $(\psi-\psi^*_g)(1-p_g/k_g)=|\psi-\psi^*_g|\,\bar w_g$. Averaging the length bound over the prior, and writing $\psi-\psi^*_g=(\psi-\psi^*_g)(1-p_g/k_g)+(\psi-\psi^*_g)p_g/k_g$,
\[
(1+\delta)\,2zs_g-2z\sqrt{c_\Lambda}\,s_g\ \ge\ \E_{m_g}\!\int(\psi-\psi^*_g)\,du\ \ge\ \E_{m_g}\!\int|\psi-\psi^*_g|\,\bar w_g\,du .
\tag{B.$\dagger$}
\]
The right-hand side is non-negative, so $\delta\ge\sqrt{c_\Lambda}-1$ for every $\Lambda$, and $\Pi^*_{\mathrm{can}}\ge1$.

\emph{The common region.} Let $\mathcal B=\{(z_0,\zeta)\in\R\times\R^2:|z_0|\le\Lambda,\ |\zeta|\le1\}$. Under $m_g$, $Z_0$ and $\langle Z,n_g\rangle$ are $N(0,1+\Lambda^2)$ and $\langle Z,n_g^\perp\rangle$ is $N(0,1)$, independent, so the density of $m_g$ on $\mathcal B$ is at least $(2\pi)^{-3/2}(1+\Lambda^2)^{-1}\exp\{-(\Lambda^2+1)/(2(1+\Lambda^2))-1/2\}=c_f/(1+\Lambda^2)$, and by (B.$\dagger$)
\[
\int_{\mathcal B}\int|\psi-\psi^*_g|\,\bar w_g\,du\,dx\ \le\ \frac{1+\Lambda^2}{c_f}\,2zs_g\bigl(1+\delta-\sqrt{c_\Lambda}\bigr),\qquad g=1,2.
\]
Pointwise, $|\psi^*_1-\psi^*_2|\min\{\bar w_1,\bar w_2\}\le|\psi-\psi^*_1|\bar w_1+|\psi-\psi^*_2|\bar w_2$. On the left, $c_\Lambda(T_1-T_2)=c_\Lambda\langle\beta,Z\rangle$ does not depend on $z_0$, and the substitution $u=c_\Lambda T_2+\sqrt{c_\Lambda}\,v$ shows that the integral over $u$ is $\sqrt{c_\Lambda}\,G(\sqrt{c_\Lambda}\langle\beta,\zeta\rangle)$ at $x=(z_0,\zeta)$. Integrating over $\mathcal B$,
\[
2\Lambda\sqrt{c_\Lambda}\int_{|\zeta|\le1}G\bigl(\sqrt{c_\Lambda}\langle\beta,\zeta\rangle\bigr)\,d\zeta\ \le\ \frac{1+\Lambda^2}{c_f}\,2z(s_1+s_2)\bigl(1+\delta-\sqrt{c_\Lambda}\bigr),
\]
which is $\delta\ge$ the bracket in the definition of $\delta_{\mathrm{can}}$; taking the supremum over $\Lambda$ gives $\delta\ge\delta_{\mathrm{can}}$. As $\Lambda\to\infty$ the bracket is at least $-(1-\sqrt{c_\Lambda})\to0$, so $\delta_{\mathrm{can}}\ge0$.

\emph{Positivity.} If $\alpha_1=\alpha_2=0$, the two Wald intervals $a_0Z_0\pm za_0$ coincide, are honest and have price one. Otherwise $\langle\beta,\zeta\rangle=(\alpha_1-\rho\alpha_2)\zeta_1-\alpha_2\sqrt{1-\rho^2}\,\zeta_2$ is not identically zero, so it vanishes only on a null set. For $d\ne0$ the two intervals in the definition of $G(d)$ differ, their symmetric difference has positive length, and the weights vanish only at the endpoints, so $G(d)>0$; and $G$ is continuous, by dominated convergence. By Fatou's lemma, $\liminf_{\Lambda\to\infty}\int_{|\zeta|\le1}G(\sqrt{c_\Lambda}\langle\beta,\zeta\rangle)\,d\zeta\ge\int_{|\zeta|\le1}G(\langle\beta,\zeta\rangle)\,d\zeta>0$, while $1-\sqrt{c_\Lambda}\le1/(2\Lambda^2)$, so the bracket is of order $1/\Lambda$ and positive for $\Lambda$ large. Finally, $K_{1\to2}=K_{2\to1}=0$ means $\alpha_1=\rho\alpha_2$ and $\alpha_2=\rho\alpha_1$, hence $\alpha_1=\rho^2\alpha_1$, and both vanish because $|\rho|<1$.
\end{proof}

The bound is far from sharp, because the weights vanish at the endpoints of the intervals and only a bounded region of the data is used: on the $27$ canonical parameters of \cref{app:certify}, the largest value of the bracket over a grid of $\Lambda$, a lower bound for $\delta_{\mathrm{can}}$, is between $1.6\times10^{-8}$ and $4.3\times10^{-4}$, and along the ray $\alpha_1=\alpha_2\to0$ at $\rho=0.5$ it is about $9.3\times10^{-5}K_{1\to2}^4$ (script \texttt{d44c\_strict\_bound.py}).

\emph{The procedure \DSsw.} Let $\gamma_1,\gamma_2,\gamma_3,\alpha'>0$ with $\gamma_1+\gamma_2+\gamma_3+\alpha'=\alpha$, and let $\varepsilon_n,\tau_n$ be as in \cref{thm:K}(ii). \DSsw\ is \DSeff\ with two changes: steps (1)--(3) use only the second part of the sample; and if step (0) returns no estimate, or returns $\hat\cX$ such that $(\hat\cX,\hat e)$ is free, the levels $(\gamma_1,\gamma_2,0,\alpha'+\gamma_3)$ are used, so that step (1b) is omitted, and otherwise $(\gamma_1,\gamma_2,\gamma_3,\alpha')$. With the $t$-intervals of an adjustment rule in step (3) we obtain its version with adjustment intervals.

\begin{proof}[Proof of \cref{thm:free}]
\emph{(i), lower bound.} If $m\ge2$, \cref{prop:canonical}(ii) for any ordered non-equivalent pair and \cref{lem:strict} give $\liminf_n\Pi^{\mathrm{loc}}_n(C_n;r)\ge\Pi^*_{\mathrm{can}}\ge1$ for every $r>0$ and every asymptotically honest sequence, so $\Pi^*_k(P_0)\ge1$. If $m=1$, the proof of \cref{prop:canonical}(ii) applies verbatim, with one model in place of two, to the limit experiment $Z\sim N(h,I)$, $h\in T_{P_0}\cM_D$, of one explanation $D\in\Xeq$, with target $\langle a_D,h\rangle$ and $|a_D|=s_{\mathrm{eff}}(D,P_0)$: in step (d) this Gaussian location experiment, which does not change with the box, takes the place of the canonical experiment, and its Wald interval that of the union of the two oracle intervals; the argument for (B.$\dagger$) with the prior $N(0,\Lambda^2I)$ on $h$ shows that its minimax price is at least $\sqrt{c_\Lambda}$ for every $\Lambda$, hence at least one.

\emph{(ii).} If $K(P_0)>0$, some ordered pair has $K_{1\to2}>0$, so $(\alpha_1,\alpha_2)\ne0$ for it, and \cref{prop:canonical}(ii) and \cref{lem:strict} give $\Pi^*_k(P_0)\ge\Pi^*_{\mathrm{can}}\ge1+\delta_{\mathrm{can}}>1$. This is also the direction ``only if'' of (i).

\emph{(iii).} If $m=1$, then $K=0$ by definition and $\mathcal R=1$ by \eqref{eq:K1p}, and the bound is \cref{thm:explicit}(b), with $e_1=0$. Let $m\ge2$ and $K(P_0)=0$. Every $D\in\Xeq$ has a member $D'$ of another Markov class, and $|\rho|<1$ for the pair $(D,D')$ outside $\cN$ (\cref{thm:E}), so the last paragraph of the proof of \cref{prop:Kpos} applies to it and gives $\alpha_D=0$. Hence the gradient of $\theta_D$ in $T_{P_0}\cM_D$ is the common gradient $a_0\in L_0$ of $\theta_H$, $s_{\mathrm{eff}}(D,P_0)=|a_0|$ and $\mathcal R(P_0)=1$. The derivative of $\Pi_D$ at $P_0\in\cM_D$ is the $g$-orthogonal projection onto $T_{P_0}\cM_D$ (\cref{lem:Dproj}), so the limits $Z_D=d(\theta_D\circ\Pi_D)[\Delta]=g(a_0,\Delta)$ in the proof of \cref{thm:explicit}(b) are equal for all $D$, and the term $\E[\max Z-\min Z]$ vanishes. Also $\delta_{G^*,S'}=0$ for every pair, so the offsets $b_{S',m}$ tend to zero for finite $c$, and the terms of the classes $S'\ne S^*$ vanish for infinite $c$ as there. The limit is therefore at most $2z_{1-\alpha'/2}|a_0|/(2z_{1-\alpha/2}|a_0|)$. The same changes in the proof of \cref{thm:K}(ii), which follows that of \cref{thm:explicit}(c), give the bound for \DSeff. Letting $\gamma_1+\gamma_2\downarrow0$ in the bound for \DSlreff\ gives $\Pi^*_k(P_0)\le1$, which with the lower bound proves ``if'' in (i). The statements on free pairs $(\Xeq,e)$ in (i) are \cref{prop:Kpos} and the definition of $K$ for $m=1$.

\emph{\DSsw, honesty.} Given the first part, the levels and the matrices $\hat B_{S'}$ are fixed, and steps (1)--(3) are those of \DSeff, or of \DSlreff, applied to the second part, an independent sample of size $n_2\to\infty$. The proofs of \cref{thm:coverage}(ii) and \cref{thm:explicit}(b) show that the coverage of these procedures is at least $1-\alpha-o(1)$ uniformly over $\cM_k$ and over the directions, because the directed test of the true graph has exact level given its direction; so the unconditional coverage is at least $1-\alpha-o(1)$ uniformly. With adjustment intervals the conditional coverage is at least $1-\alpha$ whenever $n_2\ge p+4$, by \cref{thm:coverage}(ii).

\emph{\DSsw, the bounds.} Since $n_2/n\to1$, the local limits in the proofs of \cref{thm:K}(ii) and \cref{thm:explicit}(b) are the same for the second part as for the full sample. Let $P_0$ be a base point to which \cref{cor:RK} applies, with $m\ge2$ or $m=1$, and follow the proof of \cref{thm:K}(ii). In Case~II, on the event $A_n$, $\hat\cX=\Xeq$ and $\hat e=e$, so the levels are $(\gamma_1,\gamma_2,0,\alpha'+\gamma_3)$ if $(\Xeq,e)$ is free and $(\gamma_1,\gamma_2,\gamma_3,\alpha')$ otherwise; the bound is that of (iii) for \DSlreff\ in the first case, because $K(P_0)=0$ outside $\cN$ by (i), and that of \cref{thm:K}(ii) in the second, and the complement of $A_n$ contributes nothing, as there. In Case~I, step (1a), whose level is $\gamma_1$ for both choices, keeps the classes $S'\ne S^*$ with exponentially small probability, as in that proof, so only the class of the true graph contributes, and the price tends to at most $z_{1-\alpha''/2}/z_{1-\alpha/2}$ with $\alpha''$ the level of the final interval. At a base point at which $(\Xeq,e)$ is free, step (0) returns either no estimate or a pair with $\hat e=e$ and $\hat D-e\simeq H$, hence $\hat\cX=\Xeq$, except with probability $O(e^{-cn_1})$, by the separation $d_0$ in that proof; so $\alpha''=\alpha'+\gamma_3$ except on an event whose contribution vanishes, because $\sqrt n\Leb(C_n)$ has second moment $O(n)$ and $n_1/\log n\to\infty$. Hence the limiting local price is at most $z_{1-(\alpha-\gamma_1-\gamma_2)/2}/z_{1-\alpha/2}$ at base points at which $(\Xeq,e)$ is free, which is at most $1+\eta$ for $\gamma_1+\gamma_2$ small, and at most the bound of \cref{thm:K}(ii) at the others, which is at most $(a_m+\eta)\mathcal R+b_mK_+$ for $\gamma_1+\gamma_2$ small, as in the proof of \cref{cor:RK}; there $z_{1-\alpha'/2}/z_{1-\alpha/2}\le a_m+\eta$ covers Case~I.
\end{proof}

\subsection{Proof of \texorpdfstring{\cref{thm:uniform}}{the uniformity theorem}}\label{app:uniform}

Throughout, $H_2=H-f$, $\cX_2=\{G^*\in\cG:G^*-e-f\simeq H_2\}$, and for $G^*\in\cX_2$ we use the coordinates $(\phi',\varepsilon,\lambda')$ of $\Theta_{G^*}$, where $\varepsilon=\beta^{G^*}_e$, $\lambda'=\beta^{G^*}_f$ and $\phi'$ collects the other coefficients and the noise variances; $P_{00}$ has coordinates $(\phi'_{00},0,0)$, because $P_{00}\in\cM_{H_2}=\cM_{G^*-e-f}$, and $\phi_{00}\mapsto\phi'_{00}$ is an analytic diffeomorphism $\Theta_{H_2}\to\Theta_{G^*-e-f}$ (\cref{app:benign}), so null sets correspond. For a candidate $G'$ write $\mathcal D(\phi',\varepsilon,\lambda')=D_{G'}(P)$ and $\Delta(\phi',\varepsilon,\lambda')=\theta_{G'}(\Pi_{G'}P)-\theta_{G^*}(P)$; both are analytic near $(\phi'_{00},0,0)$ by \cref{lem:Dproj}, because $M_j$ and $\sigma^2_{j\mid\pa(j)}$ are positive at $P_{00}$.

\paragraph{Consequences of \eqref{eq:M}.} (a) $\cX_2\subseteq\Xeq$: for $G^*\in\cX_2$, $G^*-e-f\in[H_2]_Y$ and $G^*-e$ is obtained from it by adding $f$ in the orientation of $G^*$, which $\cK$ allows, so $G^*-e\in[H]_Y$. (b) For $G^*\in\cX_2$ and $P\in\cM_{H_2}$, $\theta_{G^*}(P)=\theta_H(P)$, by \eqref{eq:K2} and $\cM_{H_2}\subseteq\cM_H$. (c) $f$ is benign at base points on $\cM_{H_2}$, in the sense of \eqref{eq:K2} and \eqref{eq:K1p} with $(f,H_2)$ in place of $(e,H)$, without any admissibility condition on $f$: an explanation $G'$ of such a base point that differs from $H_2$ only in $f$ and in the orientations $f$ leaves undetermined has $G'-f\in[H_2]_Y$, so $G'\in[H]_Y$ by \eqref{eq:M}; on $\cM_{H_2}\subseteq\cM_H=\cM_{G'}$, $\thadj_{G'}=\theta_{G'}=\theta_H=\theta_{H_2}$ by \eqref{eq:K2} for $e$; and all these explanations lie in the single Markov class of $H$. In the configurations of \cref{app:certify} that have an (F) pair with $e$ benign in $H$ and satisfy \eqref{eq:M}, $3157$ on six nodes and $63$ on five, $f$ is benign, as (c) asserts.

\paragraph{A local form of \eqref{eq:M}.} Write $f=a-b$ with $a\to b$ in $H$. For $D\in[H_2]_Y$, the graphs $D+f$ and $H$ have the same skeleton and the same v-structures without arm $f$, because adding $f$ removes from $D$ the v-structures $a\to w\leftarrow b$, whose centres are the common children of $a$ and $b$, the same in all members of $[H_2]_Y$ and in $H$. With $\adj$ the adjacency of $S-e$, the v-structures with arm $f$ are $a\to b\leftarrow x$, $x\in\pa(b)\setminus(\{a\}\cup\adj(a))$, if $f$ is oriented $a\to b$, and $b\to a\leftarrow x$, $x\in\pa(a)\setminus(\{b\}\cup\adj(b))$, if it is oriented $b\to a$; the latter have centre $a$, the former centre $b$. So \eqref{eq:M} holds if and only if, for every $D\in[H_2]_Y$, $\pa_D(b)\setminus\adj(a)=\pa_H(b)\setminus(\{a\}\cup\adj(a))$ whenever $D+(a\to b)$ is acyclic, which $\cK$ always allows because $H$ has $a\to b$, and $\pa_D(a)\setminus\adj(b)=\pa_H(b)\setminus(\{a\}\cup\adj(a))=\emptyset$ whenever $D+(b\to a)$ is acyclic and allowed by $\cK$. It depends only on the parent sets of the endpoints of $f$ in the members of $[H_2]_Y$; the two forms agree in all configurations of \cref{app:certify}.

\paragraph{True graphs near $P_{00}$.} Let $G'\in\cX(P_{00})\setminus\cX_2$, where $\cX(P_{00})=\{G'\in\cG:P_{00}\in\cM_{G'}\}$, and let $G''$ be $G'$ without its edges of zero coefficient at $P_{00}$. The proof of \cref{lem:HC} applies with $H_2$ in place of $H$: $P_{00}$ is Markov and minimal for $G''$, and it is not faithful to $G''$, since otherwise $G''\simeq H_2$, so that the zero edges of $G'$ are exactly $e$ and $f$ and $G'\in\cX_2$. Hence some triangle edge of $G'$ with non-zero coefficient has $\rho_{P_{00}}(u,v\mid W)=0$ for a legal $W$. As in Step~0 of \cref{app:benign}, there is $r_1(k)$ such that every $(G^*,P)\in\cM_k$ with $\|\Sigma_P-\Sigma_{P_{00}}\|\le r_1(k)$ has $G^*\in\cX_2$, and, reducing $r_1$, every such $G'$ has a $k$-TF element with violation margin at least $\mu(k)>0$ at every such $P$. These radii do not depend on $\lambda$ or $A$.

\paragraph{The threshold $k_0$.} The true graphs of the class must satisfy $k$-TF along their triangles. Near $P_{00}$ this holds for all points of $\cM_{G^*}$ when no triangle contains both $e$ and $f$, and for the points at which the coefficient of $e$ vanishes in any case; the proof uses points with non-zero coefficient of $e$ only when $\iota=1$, where no triangle contains both. Let $G^*\in\cX_2$ and let $(u\to v,W)$ be a $k$-TF element of $G^*$ with triangle $\{u,v,c\}$.
(i) If $u-v\notin\{e,f\}$, then $u-v$ is an edge of $G^*-e-f\simeq H_2$, so its coefficient at $P_{00}$ is non-zero and $\rho_{P_{00}}(u,v\mid W)\ne0$, by faithfulness.
(ii) If $u-v\in\{e,f\}$ and the triangle does not contain both $e$ and $f$, then $u-c$ and $v-c$ are edges of $G^*-e-f$, so $u-c-v$ is an unshielded triple of $G^*-e-f$ and of the equivalent $H_2$ with the same collider status as in $G^*$; as $W$ is legal and $u,v\notin W$, the path is open given $W$ in $H_2$, and $\rho_{P_{00}}(u,v\mid W)\ne0$ by faithfulness, while $\lambda^{\mathrm{std}}$ of $u-v$ vanishes at $P_{00}$.
(iii) If the triangle contains both $e$ and $f$, then $e$ and $f$ share an endpoint and their other endpoints are adjacent, so no v-structure has arms $e$ and $f$, no pair satisfies (F), and $\iota=0$. For an element with $u-v=e$, the condition is void at points with $\varepsilon=0$, where $\lambda^{\mathrm{std}}_e=0$. For an element with $u-v=f$, consider the points with $\varepsilon=0$. There $\lambda^{\mathrm{std}}_f=\lambda'\sigma_u/\sigma_v$ exactly, and $\rho(u,v\mid W)$ is analytic in $(\phi',\lambda')$. If $\rho(\cdot,0)$ does not vanish identically in $\phi'$, it is non-zero at almost every $\phi'_{00}$, and then the ratio $|\rho|/|\lambda^{\mathrm{std}}_f|$ tends to infinity as $\lambda'\to0$. If $\rho(\cdot,0)$ vanishes identically, the ratio tends to $\ell(\phi')=|\partial_{\lambda'}\rho(\phi',0)|\,\sigma_v/\sigma_u$, locally uniformly in $\phi'$; the derivative is analytic in $\phi'$ and equals $1$ at the value of $\phi'$ with all coefficients zero and all variances one, where only $u$ and $v$ are dependent and $\rho=\lambda'/\sqrt{1+\lambda'^2}$ for every $W$, so $\ell(\phi'_{00})>0$ for almost every $\phi_{00}$. We set $\ell=\infty$ in the first case.

Let $k_0$ be half the smallest of the ratios $|\rho|/|\lambda^{\mathrm{std}}|$ at $P_{00}$ of the elements in (i), and of the limits $\ell(\phi'_{00})$ in (iii), over $G^*\in\cX_2$. By continuity, for $k<k_0$ there is $r_2$ with the following property. Every $P\in\cM_{G^*}$, $G^*\in\cX_2$, with $\|\Sigma_P-\Sigma_{P_{00}}\|\le r_2$ and $\beta^{G^*}_e(P)=0$ satisfies $k$-TF under $G^*$. If no triangle of $S$ contains both $e$ and $f$, the restriction $\beta^{G^*}_e(P)=0$ can be dropped. Indeed, the elements in (i) are controlled by the choice of $k_0$, those in (ii) because $|\rho|$ stays bounded away from zero while $\lambda^{\mathrm{std}}\to0$, and those in (iii), at points with $\varepsilon=0$, by the limits $\ell$; without a triangle containing $e$ and $f$ there are no elements in (iii). Without the restriction $\varepsilon=0$, the elements of (iii) need not satisfy $k$-TF near $P_{00}$: their ratio depends on the direction $\lambda'/\varepsilon$ and can vanish along a line, a near-cancellation along the triangle that the class excludes. The radius $r_2$ depends on $P_{00}$ and can be much smaller than $k_0$ suggests near non-generic base points (\cref{app:certify}).

\paragraph{The generic set.} Exclude from $\Theta_{H_2}$ the null sets on which $P_{00}$ is not faithful to $H_2$, on which some $\rho(\cdot,0)$ in (iii) that is not identically zero vanishes or some limit $\ell(\phi'_{00})$ vanishes, and on which, for some ordered pair $(G^*,G')$ of non-equivalent members of $\cX_2$, one of the following fails: $I_{G^*,G'}(P_{00})>0$ if the pair does not satisfy (F); $I_2>0$ in the expansion $I_{G^*,G'}(P_\lambda)=\lambda^2I_2+O(\lambda^3)$ of the proof of \cref{thm:K-lambda} if it does; and $\delta_{G^*,G'}(P_{00})\ne0$ if it satisfies (F) and (D). These are the conditions established for almost every $\phi'_0$ in that proof, and by \cref{lem:Kclass} they do not depend on the representatives. Along the curve $P_\lambda$ of that proof the coordinates here are $(\phi'(\lambda),0,\lambda'(\lambda))$ with $\phi'(0)=\phi'_{00}$ and $\lambda'=u(\lambda)\lambda$, $u(0)\ne0$, because the coefficients of $f$ under $H$ and under $G^*-e$ are coordinate functions of analytic charts of $\cM_H$, both vanishing exactly on $\cM_{H_2}$, so each has non-vanishing differential there and their ratio is analytic and non-zero, while $P_\lambda$ crosses $\cM_{H_2}$ transversally. Finally, fix a directed path $\pi$ from $X$ to $Y$ in $H$ and its first edge $X\to b$, and exclude the null set on which $\tau_{b\to Y}$, the derivative of $\theta_H$ in the coefficient of $X\to b$, vanishes identically as a polynomial in the coefficient of $f$; it is a polynomial in the parameters of $H$ that does not involve the coefficient of $X\to b$ and is not identically zero, being $1$ when the edges of $\pi$ have coefficient $1$ and all others $0$.

Under \eqref{eq:M}, $\iota=1$ if and only if $\cX_2$ contains a pair that satisfies (F) and (D): such a pair lies in $\Xeq$ by (a), and (C) follows from \cref{lem:Vf} below; conversely, a pair in $\Xeq$ with (F), (C) and (D) has the representatives $D_\wedge,D_\vee$ of the proof of \cref{thm:K-lambda}(iv) in its classes, with $(D_\wedge-e)-f\simeq H_2$, and also $(D_\vee-e)-f\simeq H_2$, because $D_\vee$ and $D_\wedge$ differ only in the orientation of $e$ or of $f$; so $D_\wedge,D_\vee\in\cX_2$, and (F) and (D) depend only on the classes.

\begin{lemma}[Deleting the margin edge]\label{lem:Vf}
Let $e=x_e-z$ and $f=z-x_f$ share the endpoint $z$, with $x_e$ and $x_f$ not adjacent, and let $G_1,G_2\in\cX_2$ have v-structures that differ at most in $x_e\to z\leftarrow x_f$. Then $G_1-f\simeq G_2-f$, so that $\cM_{G_1-f}\subseteq\cM_{G_2}$.
\end{lemma}
\begin{proof}
Deleting the edge $f$ from a DAG $D$ removes the v-structures with arm $f$ and creates $z\to w\leftarrow x_f$ for the common children $w$ of $z$ and $x_f$; so the v-structures of $D-f$ are those of $D$ without arm $f$, together with these. The v-structures of $G_1$ and $G_2$ without arm $f$ agree, because the only possible difference has arm $f$. By (a), $G_i-e\simeq H$, and by the paragraph on common children in the proof of \cref{thm:K-lambda}, which uses only that $x_e$ is not a common neighbour of $z$ and $x_f$, $(G_i-e)-f\simeq H_2$ if and only if $\mathrm{CC}_{G_i}(z,x_f)=\mathrm{CC}_H(z,x_f)$; as $G_i\in\cX_2$, the two sets of common children agree. The skeletons agree.
\end{proof}
In particular, within $\cX_2$, (F) implies (C).

\begin{lemma}[The discrepancy at zero margin]\label{lem:Af}
In the setting of \cref{lem:Vf}, let $G_1\not\simeq G_2$ and $\Delta\theta=\theta_{G_2}-\theta_{G_1}$ on $\cM_{G_1-f}$. If the pair does not satisfy (D), then $\Delta\theta\equiv0$ on $\cM_{G_1-f}$.
\end{lemma}
\begin{proof}
Let $D_\wedge,D_\vee\in\cX_2$ be the representatives of the classes of $G_\wedge,G_\vee$ constructed in the proof of \cref{thm:K-lambda}(iv), and $K=D_\wedge-f$. By \cref{lem:Vf}, applied also to $G_i$ and the representative of its class, $\cM_{G_1-f}=\cM_{G_2-f}=\cM_K$, and by \eqref{eq:K1p}, $\theta_{G_i}$ equals $\theta_{D_\wedge}$ or $\theta_{D_\vee}$ on $\cM_{G_i}\supseteq\cM_K$, according to the class of $G_i$. So $\pm\Delta\theta$ is the function $\theta_{D_\vee-f}-\theta_{D_\wedge-f}$ on $\cM_K$ of that proof, which vanishes identically if its first descent, that of the proof of \cref{lem:onev}, stops in type $\alpha$, or if $X\notin e$. Otherwise the descent stops in type $\beta$ at an $e$-configuration (proof of \cref{lem:Econf}(ii)), so $e\notin R$ and some member of $[H]_Y$ gives $x_e$ the parent set $U$ (\cref{lem:Econf}); as (D) fails, $Y$ and $w$ are d-separated given $U\cup\{X\}$ in $H_2$, and $\Delta\theta=\pm\beta_{Yw\cdot U\cup\{X\}}\beta_{x_eX\cdot U}$ or $\pm\beta_{Yw\cdot U\cup\{X\}}\beta_{zX\cdot U}$ on $\cM_K$. The d-separation persists in $K$, which is $K-e\simeq H_2$ with the covered edge $x_e\to z$ added, $\pa_K(x_e)=U$ and $\pa_K(z)=U\cup\{x_e\}$. If $X=z$ and $w=x_e$, a path from $x_e$ that uses the new edge continues from $z$, which is in the conditioning set, either along an edge out of $z$, where $z$ is a non-collider, or into $z$ from some $u\in U$, so that $u$ is a non-collider in the conditioning set. If $X=x_e$ and $w=z$, a path from $z$ that uses the new edge passes through $x_e$, a non-collider in the conditioning set. In both cases the new edge changes the sets of descendants only of $x_e$ and of its proper ancestors. The proper ancestors of $x_e$ belong to $U=\pa_K(x_e)$ or are ancestors of $U$, and $U$ lies in the conditioning set $U\cup\{X\}$, so a proper ancestor of $x_e$ either is conditioned on or already had a conditioned descendant; and $x_e$ itself is either the endpoint $w$ of the query ($X=z$) or conditioned on ($X=x_e$). So no collider on a path that avoids the new edge changes its status. Hence $\beta_{Yw\cdot U\cup\{X\}}\equiv0$ on $\cM_K$.
\end{proof}

\paragraph{The local expansion.} Fix $G^*\in\cX_2$ and a candidate $G'$. There is a neighbourhood $V$ of $(\phi'_{00},0,0)$, and constants $c_1>0$, $C_1$, such that on $V$:
\begin{enumerate}
\item[(a)] if $G'\notin\cX(P_{00})$, then $\mathcal D\ge c_1$;
\item[(b)] if $G'\in\cX(P_{00})\setminus\cX_2$, then $G'$ has a $k$-TF element with violation margin at least $\mu(k)$;
\item[(c)] if $G'\in\cX_2$ and $G'\simeq G^*$, then $\mathcal D=\Delta=0$;
\item[(d)] if $G'\in\cX_2$, $G'\not\simeq G^*$ and the pair does not satisfy (F), then $\mathcal D\ge c_1\varepsilon^2$ and $|\Delta|\le C_1|\varepsilon|$;
\item[(e)] if $G'\in\cX_2$ and the pair satisfies (F), then $\mathcal D\ge c_1\varepsilon^2\lambda'^2$, $\mathcal D(\phi',\varepsilon,0)=0$, $|\Delta|\le C_1|\varepsilon|$, and $|\Delta|\le C_1|\varepsilon\lambda'|$ if (D) fails.
\end{enumerate}
(a) is continuity, since $\mathcal D>0$ at $P_{00}$, and (b) was shown above. (c) $\cM_{G'}=\cM_{G^*}$, and by \eqref{eq:K1p} and (a) of the consequences of \eqref{eq:M}, $\theta_{G'}=\theta_{G^*}$ on it. For (d) and (e), at $\varepsilon=0$, $P\in\cM_{G^*-e}=\cM_H=\cM_{G'-e}\subseteq\cM_{G'}$, so $\mathcal D(\cdot,0,\cdot)=0$, and $\Delta(\cdot,0,\cdot)=0$ by \eqref{eq:K2}. As $\mathcal D\ge0$, Hadamard's lemma gives $\mathcal D=\varepsilon^2d$ and $\Delta=\varepsilon h$ with analytic $d$, $h$, and $2d(\phi'_{00},0,\lambda')=I_{G^*,G'}$ at the point with coordinates $(\phi'_{00},0,\lambda')$. In (d), $d(\phi'_{00},0,0)>0$ on the generic set, which gives $\mathcal D\ge c_1\varepsilon^2$ near it. In (e), \cref{lem:Vf} gives $\mathcal D(\cdot,\cdot,0)=0$, so $d(\cdot,\cdot,0)=0$, and $d\ge0$ gives $d=\lambda'^2\psi$ with $\psi$ analytic. Along the curve $P_\lambda$, $I(\lambda)=2d(\phi'(\lambda),0,\lambda'(\lambda))=2u(\lambda)^2\lambda^2\psi(\phi'(\lambda),0,\lambda'(\lambda))$, so $I_2=2u(0)^2\psi(\phi'_{00},0,0)$ and $\psi(\phi'_{00},0,0)>0$ on the generic set, which gives $\mathcal D\ge c_1\varepsilon^2\lambda'^2$. If (D) fails, \cref{lem:Af} gives $h(\cdot,\cdot,0)=0$, so $h=\lambda'h_1$. Consequently, on $V$ and for every $G'\in\cX_2$,
\begin{equation}\label{eq:uniform-star}
|\Delta|\le C_2\Bigl(\sqrt{\mathcal D}+\iota\min\Bigl\{|\varepsilon|,\frac{\sqrt{\mathcal D}}{|\lambda'|}\Bigr\}\Bigr),
\end{equation}
because in (e) with (D), $|\varepsilon|\le\sqrt{\mathcal D}/(\sqrt{c_1}|\lambda'|)$, and in the other cases $|\Delta|\le C_1\sqrt{\mathcal D/c_1}$. By finiteness of $\cG$, one neighbourhood and one set of constants serve all pairs. Take $r\le\min(r_1,r_2)$ so small that, for every $G^*\in\cX_2$, the points of $\cM_{G^*}$ within $r$ of $\Sigma_{P_{00}}$ have coordinates in $V$, which is possible because the coordinates are continuous on $\cM_{G^*}$, and $\lambda_0,A_0$ so small that the points used below stay in $V$.

\begin{proof}[Proof of the upper bound]
Let $(G^*,P)\in\cM_k(\lambda,A;r)$, so that $G^*\in\cX_2$, $|\varepsilon|\le A$ and $|\lambda'|\ge\lambda$. Let $c_{G'}(P)$ be the coefficient of $X$ in the population regression of $Y$ on $X$ and $A(G')$, the centre of the interval of $G'$ in \DSA. If the output is non-empty, its length is at most $\sum_{G'}\Leb(W_{G'})+2\sum_{G'}|\hat c_{G'}-\theta_{G^*}(P)|\ind\{G'\text{ kept}\}$, the sums over all candidates, and as in Step~1 of \cref{app:benign} the widths and the errors $\hat c_{G'}-c_{G'}(P)$ contribute $O(n^{-1/2})$ uniformly. Near $\Sigma_{P_{00}}$, the projection distance is comparable to the squared distance to the closed set $\cM_{G'}$: $c_8\,\mathrm{dist}(\Sigma,\cM_{G'})\le\sqrt{D_{G'}(\Sigma)}\le C_8\,\mathrm{dist}(\Sigma,\cM_{G'})$, because the Kullback--Leibler divergence between Gaussian laws with covariance eigenvalues in a compact subinterval of $(0,\infty)$ is comparable to the squared distance of the covariances, and $\Pi_{G'}\Sigma$ is close to $\Sigma$ when $\Sigma$ is close to $\cM_{G'}$. Since $c_{G'}=\theta_{G'}$ on $\cM_{G'}$ and both are smooth, $|c_{G'}(P)-\theta_{G'}(\Pi_{G'}P)|\le C\sqrt{D_{G'}(P)}$. So it remains to bound $(|\Delta_{G'}|+C\sqrt{D_{G'}})\Prob(G'\text{ kept})$ for each $G'$. For $G'$ in case (b), the $k$-TF screen keeps $G'$ with probability at most $16e^{-c_Bn\mu^2}$ (\cref{lem:B}), and $|\Delta_{G'}|\le2\sigu/\sigl$. For the others, the likelihood ratio screen keeps $G'$ only if $T_{G'}=2nD_{G'}(\hat\Sigma)$, an identity because $T_{G'}$ is $n$ times a sum of logarithms of ratios of residual sums of squares, is at most the $1-\gamma_1$ quantile $q_n$ of its null law, which is bounded in $n$. The sample covariance satisfies $\Prob(\|\hat\Sigma-\Sigma_P\|>t)\le C_3e^{-c_3nt^2}$ for $t\le1$. On the event $\|\hat\Sigma-\Sigma_P\|\le t_0$, $\hat\Sigma$ stays in the neighbourhood, and $\sqrt{D_{G'}(\hat\Sigma)}\ge c_8\bigl(\mathrm{dist}(\Sigma_P,\cM_{G'})-\|\hat\Sigma-\Sigma_P\|\bigr)\ge(c_8/C_8)\sqrt{D_{G'}(P)}-c_8\|\hat\Sigma-\Sigma_P\|$. Hence $\Prob(G'\text{ kept})\le C_4e^{-c_4nD_{G'}(P)}+C_3e^{-c_3nt_0^2}$ whenever $nD_{G'}(P)\ge4(C_8/c_8)^2q_n$, and in case (a) this is exponentially small. With \eqref{eq:uniform-star} and $\sup_{x\ge0}\sqrt xe^{-c_4x}<\infty$,
\[
(|\Delta_{G'}|+C\sqrt{D_{G'}})\Prob(G'\text{ kept})\le C_5\Bigl(n^{-1/2}+\iota\min\Bigl\{A,\frac1{\lambda\sqrt n}\Bigr\}\Bigr),
\]
using $|\varepsilon|\le A$ and $\sqrt{\mathcal D}/|\lambda'|\le\sqrt{\mathcal D}/\lambda$. Summing over the candidates and adding the exponentially small terms, which are $O(n^{-1/2})$ for $n\ge n_0$, gives the upper bound; honesty is \cref{thm:coverage}. The output is empty only if the true graph is rejected, which does not increase the length.
\end{proof}

\begin{proof}[Proof of the lower bound]
\emph{The term $n^{-1/2}$.} Take the true graph $G$, whose coordinates at $\varepsilon=0$ are those of $H=G-e$, and the edge $X\to b$ of the generic set. Let $P_\pm\in\cM_H$ have the parameters of $P_{00}$ under $H$, except that the coefficient of $f$ is $2\lambda_0$ and that of $X\to b$ is shifted by $\pm c\,n^{-1/2}$; if $X\to b$ is $f$, its coefficient is $2\lambda_0\pm c\,n^{-1/2}$. Then $\beta^G_e(P_\pm)=0$ and $|\beta^G_f(P_\pm)|\ge\lambda_0\ge\lambda$ for $n\ge n_0$. For small $\lambda_0$ and large $n_0$, $P_\pm$ lie within $r$ of $\Sigma_{P_{00}}$, and they satisfy $k$-TF under $G$ by the property of $r_2$ with $\beta^G_e=0$, whether or not a triangle contains $e$ and $f$. So $(G,P_\pm)\in\cM_k(\lambda,A;r)$ for every $\lambda\le\lambda_0$ and every $A$. $\KL(P_+\|P_-)\le C_9c^2/n$, and, since a directed path uses $X\to b$ at most once, $\theta_G(P_+)-\theta_G(P_-)=\theta_H(P_+)-\theta_H(P_-)=2c\,n^{-1/2}\tau_{b\to Y}$ exactly, with $\tau_{b\to Y}$ evaluated at coefficient $2\lambda_0$ of $f$. On the generic set it is non-zero for all but finitely many values of $\lambda_0$, which we avoid. With $c$ so small that $\dtv(P_+^n,P_-^n)\le(1-2\alpha)/2$, \cref{lem:two-point} gives a lower bound of order $n^{-1/2}$.

\emph{The term $\iota\min\{A,1/(\lambda\sqrt n)\}$.} Let $\iota=1$, so that no triangle of $S$ contains both $e$ and $f$ ((iii) of the threshold paragraph), and let $(G^*,G')=(D_\vee,D_\wedge)$ or $(D_\wedge,D_\vee)$, the representatives above, which lie in $\cX_2$ and satisfy (F) and (D). Let $P$ have coordinates $(\phi'_{00},\varepsilon,\lambda')$ under $G^*$ and let $Q=\Pi_{G'}P$. By (e), $\KL(P\|Q)=\mathcal D\le C_6\varepsilon^2\lambda'^2$, and $|\theta_{G'}(Q)-\theta_{G^*}(P)|=|\varepsilon h|\ge|\varepsilon\delta_{G^*,G'}(P_{00})|/2$ for small coordinates. Under $G'$, $\beta^{G'}_e(Q)$ vanishes at $\varepsilon=0$, where $Q=P\in\cM_{G'-e}$, and $\beta^{G'}_f(Q)$ vanishes at $\lambda'=0$, where $Q=P\in\cM_{G'-f}$ by \cref{lem:Vf}. So $\beta^{G'}_e(Q)=\varepsilon\kappa_e$ and $\beta^{G'}_f(Q)=\lambda'\kappa_f$ with analytic $\kappa_e,\kappa_f$. These are non-zero at $(\phi'_{00},0,0)$: on $\cM_{G'-f}$ and on $\cM_H$ respectively, the two coefficients are defining functions of the same hypersurface $\cM_{H_2}$, and the structural parametrizations are immersions. Take $\lambda'=\lambda\max\{1,2/|\kappa_f(\phi'_{00},0,0)|\}$ and $\varepsilon=\min\{A\min\{1,1/(2|\kappa_e(\phi'_{00},0,0)|)\},c_7/(\lambda\sqrt n)\}$, with $c_7$ so small that $\dtv(P^n,Q^n)\le\sqrt{n\mathcal D/2}\le(1-2\alpha)/2$ (Pinsker). For small $\lambda_0$, $A_0$ and large $n_0$, both $(G^*,P)$ and $(G',Q)$ belong to $\cM_k(\lambda,A;r)$, by the property of $r_2$ without the restriction $\varepsilon=0$. Then \cref{lem:two-point} gives $\sup\E\Leb(C_n)\ge\frac14(1-2\alpha)|\delta_{G^*,G'}(P_{00})|\,\varepsilon$, which is of order $\min\{A,1/(\lambda\sqrt n)\}$.
\end{proof}

\subsection{The malignant side}\label{app:malignant}

\begin{lemma}\label{lem:M1}
Let $P_0$ be faithful to $H$, $G_i\in\Xeq$, $\{u,v,c\}$ a triangle of $G_i$ with $e=u-v$, and $W$ legal. Then $\rho_{P_0}(u,v\mid W)\ne0$.
\end{lemma}
\begin{proof}
In $G_i-e$ the path $u-c-v$ has the same collider status at $c$ as in the triangle. By legality it is active given $W$, and $u,v\notin W$; so $u$ and $v$ are not d-separated by $W$ in $G_i-e$, which is equivalent to $H$. Faithfulness gives $\rho_{P_0}(u,v\mid W)\ne0$.
\end{proof}

\begin{proof}[Proof of \cref{prop:malignant}]
$k_0(P_0)>0$ because for a triangle edge $u-v\ne e$ of $G_i$, $u$ and $v$ are adjacent in $G_i-e$, hence not d-separated by any $W$ in $H$, and there are finitely many legal $W$. (i) For $k<k_0$, $(G_1,P_0)$ and $(G_2,P_0)$ belong to $\cM_k$: triangle edges other than $e$ satisfy $k$-TF with room to spare, and $e$ has coefficient zero. \Cref{lem:two-point} with $P=Q=P_0$ gives $\E_{P_0}\Leb(C_n)\ge(1-2\alpha)\Delta_0$, and the oracle length is at most $2z_{1-\alpha/2}(\sigu/\sigl)n^{-1/2}$ by \cref{lem:oracle}. (ii) Let $\phi_i$ be the parameter of $P_0$ in $\Theta_{G_i}$ and $\phi_i(t)$ the same parameter with $\lambda_e=t$. Then $\lambda^{\mathrm{std}}_e(\phi_i(t))=t\,r_i(t)$ with $r_i$ smooth and positive; with $c_i=2/r_i(0)$ and $t_i=c_i\lambda$ small enough, $|\lambda^{\mathrm{std}}_e|\ge\lambda$. For small $t$, $k$-TF holds at $\phi_i(t)$: for edges other than $e$ by continuity and $k<k_0$, and for $e$ because by \cref{lem:M1} the partial correlations stay away from zero while $\lambda_e^{\mathrm{std}}=O(t)$; the eigenvalues stay in the interior. Let $P=P_{\phi_1(t_1)}$ and $Q=P_{\phi_2(t_2)}$. Smoothness gives $|\theta_{G_1}(P)-\theta_{G_2}(Q)|\ge\Delta_0-2L\lambda$, and $\KL(P\|P_0),\KL(Q\|P_0)\le C_1\lambda^2$, so by Pinsker's inequality $\dtv(P^n,Q^n)\le2C\lambda\sqrt n$. \Cref{lem:two-point} finishes the proof.
\end{proof}

\begin{proof}[Proof of \cref{prop:K1p-nec}]\label{app:K1p-nec}
\emph{Step 1: the pairs are in the class.} $k_1(P_0)>0$ as in \cref{prop:malignant}. Since $G_i-e\simeq H$, the coefficient of $e$ in the regression parametrization of $P_0$ under $G_i$ is zero. For triangles containing $e$, \cref{lem:M1} and finiteness of the legal sets give $\eta=\min|\rho_{P_0}(u,v\mid W)|>0$ ($\eta=+\infty$ if there are no such triangles). For other triangle edges, $|\rho_{P_0}|>k|\lambda^{\mathrm{std}}_{G_i}(P_0)|$ because $k<k_1(P_0)$. Partial correlations are continuous on the positive definite cone, and for $P\in\cM_{G_1}=\cM_{G_2}$ the regression coefficients of both $G_i$, hence $\lambda^{\mathrm{std}}$, are continuous functions of $\Sigma_P$. As there are finitely many strict inequalities and the eigenvalues of $\Sigma_{P_0}$ are interior, there is $r_1>0$ such that every $P\in\cM_{G_1}$ with $\|\Sigma_P-\Sigma_{P_0}\|\le r_1$ satisfies all $k$-TF inequalities of $G_1$ and $G_2$ and the eigenvalue bounds; for triangles containing $e$ we use $|\rho_P|\ge\eta/2\ge k|\lambda^{\mathrm{std}}_e(P)|$, where the second inequality holds near $P_0$ since $\lambda^{\mathrm{std}}_e(P_0)=0$. Hence $(G_1,P),(G_2,P)\in\cM_k$.
\emph{Step 2: $c(r)>0$.} $q\circ F_{G_1}$ is real analytic on the connected set $\Theta_{G_1}$ and not identically zero, so it does not vanish identically on any non-empty open set. The preimage under the continuous map $F_{G_1}$ of the $\Sigma$-ball of radius $r$ around $\Sigma_{P_0}$ contains an open neighbourhood of $F_{G_1}^{-1}(\Sigma_{P_0})$; hence $c(r)>0$.
\emph{Step 3: two points.} For such $P$ the targets are $\theta_{G_1}(P)$ and $\theta_{G_2}(P)=\thadj_{G_2}(P)$, the latter because $Y$ is a sink and $P\in\cM_{G_2}$. \Cref{lem:two-point} with $P=Q$ gives $\E_P\Leb(C_n)\ge(1-2\alpha)|q(P)|$, and by \cref{lem:oracle} $\orl_n\le2z_{1-\alpha/2}(\sigu/\sigl)n^{-1/2}$. Taking the supremum over $P$ proves the claim. The argument uses finite-sample honesty; for asymptotic honesty replace $1-2\alpha$ by $1-2\alpha-o(1)$.
\end{proof}

\subsection{Proof of the classification}\label{app:dichotomy}

\begin{proof}[Proof of \cref{thm:dichotomy}]
By \cref{prop:graph}(c) and \cref{thm:K1p-graph} the three cases are (K2)$\wedge$(K1$'$), $\neg$(K2), and (K2)$\wedge\neg$(K1$'$). Let $\cN'$ be the union of the exceptional set $\cN$ of \cref{thm:K}, which is defined in \cref{app:benign} without using \eqref{eq:K2} or \eqref{eq:K1p}, and, in case 2, the zero set of $(\theta_{G_1}-\theta_{G_2})^2$ for a pair $G_1,G_2\in\Xeq$ with $\theta_{G_1}-\theta_{G_2}$ not identically zero on $\cM_H$, which exists by \cref{prop:graph}(b); this zero set is null because the function is real analytic on the connected set $\Theta_H$ and not identically zero. Case 1 is \cref{thm:K,cor:RK}, case 2 is \cref{prop:malignant}(i) and case 3 is \cref{prop:K1p-nec}.
\end{proof}

\subsection{Proof of \texorpdfstring{\cref{thm:E}}{the first-order distinguishability theorem}}\label{app:E}

Write $D^*=G^*-e$, $G^*=D^*+(u\to v)$ and $G'=D'+e'$ with $e'\in\{u\to v,v\to u\}$ and $D'=G'-e\simeq D^*$.

\paragraph{Reduction to one point.} $\Theta_{D^*}$ is identified with $\Theta_{G^*}\cap\{\lambda_e=0\}$. Partial covariances are rational functions of $\Sigma$ without poles on the positive definite cone and partial correlations are real analytic there, so $\partial_{\lambda_e}h_{G'}$ is real analytic on $\Theta_{G^*}$; as $\Theta_{D^*}$ is connected, it suffices to find one point of $\Theta_{D^*}$ where it is non-zero. On $\{\lambda_e=0\}$ all partial covariances of local Markov constraints of $G'$ vanish (because $\cM_{D^*}=\cM_{D'}\subseteq\cM_{G'}$), so the derivative of a partial correlation is the derivative of the partial covariance divided by positive residual standard deviations, and it suffices to show that some partial covariance has non-zero derivative.

\paragraph{Unshielded colliders created by $e$.} Let $V_+(G^*)=\{c\in\pa_{D^*}(v):c\notin\adj(u)\}$, the third vertices of unshielded colliders $u\to v\leftarrow c$ with arm $e$; define $V_+(G')$ likewise, with collider $v$ if $e'=u\to v$ and collider $u$ (and the roles of $u,v$ exchanged) if $e'=v\to u$. An unshielded collider of $G^*$ either has no arm equal to $e$, and these are the unshielded colliders of $D^*$ whose parents are not $\{u,v\}$, or it has arm $e$. Since $D^*\simeq D'$, the first kind coincides in $G^*$ and $G'$. As $G^*\not\simeq G'$ and the skeletons agree, the second kind differs: either the colliders differ, or the sets $V_+$ differ.

\paragraph{Witness constraints.} \emph{Case A: $e'=u\to v$.} There is $c\notin\adj(u)$ such that $c\to v$ is in exactly one of $D^*,D'$; the other has $v\to c$. (A1) $c\to v\in D^*$, $v\to c\in D'$: $u\to v\to c$ in $G'$, $u$ is an ancestor and a non-parent of $c$, and $(c,u,\pa_{G'}(c))$ is a local Markov constraint of $G'$ with $v\in\pa_{G'}(c)$. (A2) $v\to c\in D^*$, $c\to v\in D'$: $u\to v\leftarrow c$ in $G'$; choose $(j,i)$ among $\{u,c\}$ with $i$ a non-descendant of $j$; $v$ is a common child, so $v\notin\pa_{G'}(j)$. \emph{Case B: $e'=v\to u$.} Unshielded colliders with arm $e$ have collider $v$ in $G^*$ and $u$ in $G'$, so they differ exactly when $V_+(G^*)$ or $V_+(G')$ is non-empty. (B1) Some $c\in\pa_{D^*}(v)$ with $c\notin\adj(u)$: if $c\to v\in D'$ then $c\to v\to u$ in $G'$ and $(u,c,\pa_{G'}(u))$ is a constraint with $v\in\pa_{G'}(u)$; if $v\to c\in D'$ then $v$ is a common parent of $u,c$ in $G'$; choose $(j,i)$ among $\{u,c\}$ with $i$ a non-descendant of $j$, so that $v\in\pa_{G'}(j)$. (B2) $V_+(G^*)=\emptyset$ and some $c\in\pa_{D'}(u)$ with $c\notin\adj(v)$: choose $(j,i)$ among $\{c,v\}$ with $i$ a non-descendant of $j$; $u$ is not a parent of $j$ in either choice ($c\to u$ and $v\to u$ in $G'$), and $c-u\in D^*$. In every case $i,j$ are non-adjacent, so the constraint is a local Markov constraint of $G'$ and $\{i,j\}\cap\pa_{G'}(j)=\emptyset$.

\paragraph{Sparse points.} Take all noise variances equal to $1$, $\lambda_e=s$, and all coefficients of $D^*$ zero except one, $\beta\ne0$, on the edge $c\to v$ (A1, B1), $v\to c$ (A2) or the edge between $c$ and $u$ (B2). All other vertices are independent standard normals, independent of $(u,v,c)$, and can be dropped from conditioning sets. The partial covariances are $\sigma_{uc\cdot v}=-s\beta/(1+\beta^2+s^2)$ in A1 and B1, with derivative $-\beta/(1+\beta^2)$ at $s=0$; $\sigma_{uc}=s\beta$ in A2 and $\sigma_{cv}=s\beta$ in B2, with derivative $\beta$. In each case the conditioning set intersected with $\{u,v,c\}$ is $\{v\}$ (A1, B1) or empty (A2, B2), as shown above. \qed

\paragraph{Witness edges.} For a non-equivalent pair $(G^*,G')$ with $e=u\to v$ in $G^*$ let $W(G^*,G')$ be the set of edges $g$ for which a configuration above exists: if $e'=u\to v$, the edges $\{c,v\}$ with $c\in\adj(v)\setminus(\{u\}\cup\adj(u))$ such that $c\to v$ is in exactly one of $D^*,D'$ (cases A1 and A2); if $e'=v\to u$, the edges $\{c,v\}$ with $c\in\pa_{D^*}(v)\setminus\adj(u)$ (case B1) and the edges $\{c,u\}$ with $c\in\pa_{D'}(u)\setminus\adj(v)$ (case B2). The sparse computation of B2 does not use $V_+(G^*)=\emptyset$, which is needed only to make the case distinction exhaustive, so every $g\in W(G^*,G')$ carries a sparse point with a non-zero derivative, and $W(G^*,G')\ne\emptyset$ by the proof of \cref{thm:E}.

\section{Computations}\label{app:numerics}

\subsection{Enumeration of graphs and checks of the graphical results}\label{app:certify}

The graphical results of \cref{sec:general} are proved; the computations in this subsection check the implementations and give the counts quoted in the text. Except in the paragraph on required edges and the four-variable instance, where equality of effects is checked by random identity tests, all labels below are decided by exact graph operations or exact rational arithmetic, without random identity tests.

\paragraph{Exhaustive enumeration.} We take $p\in\{4,5\}$ nodes, $X=0$, $Y=p-1$, every skeleton containing $X-Y$ ($32$ and $512$ labelled skeletons), every acyclic orientation with $Y$ a sink, and every edge $e$; this gives $392$ and $26988$ pairs $(G,e)$. For each pair we compute $\Xeq$ and its Markov equivalence classes, $\MP(H)$ and $\MP''(S_c)$ by our implementation of the rules R1--R4 of \cref{app:K1p}, and the bad paths. The closure agrees edge by edge with the orientations common to the corresponding extension sets in all $26988$ cases for $\MP(H)$ and all $19010$ classes for $\MP''(S_c)$, and $\MP''(S_c)$ does not depend on the representative (script \texttt{d22h\_meek\_closure.py}). At $p=5$, \eqref{eq:K2} holds for $12686$ pairs and fails for $14302$. Among the former, \eqref{eq:K1p} fails for $4260$; for each of these, an exact rational computation exhibits two members of $S_c$ with different effects, at a distribution whose only non-zero coefficients lie on one directed path from $X$ through the first vertex of a bad path to $Y$. The stronger condition ``identifiable in $\MP(G')$ for all $G'\in\Xeq$'' fails for $5106$, so $846$ pairs are benign only by \eqref{eq:K1p} (script \texttt{d21c\_certified\_counts.py}). At $p=4$ the corresponding numbers are $237$, $155$, $98$, $104$ and $6$.

\paragraph{Required edges.} For general required edges, $\MP(H)$ need not be a chain graph, about a fifth of the cases on five and six nodes below, and the criterion of \cref{thm:K1p-graph} agreed with exact computation of \eqref{eq:K1p} in every case, also without the chain structure, although a step of its proof can then fail; admissibility failed only when $e\in R$, and we have no polynomial-time test for it. For these checks we drew, for each skeleton, a random member of $\cG$ with $Y$ a sink and either one or two of its edges not at $Y$ as required (``general''), or a random partition of the nodes other than $Y$ into two or three ordered tiers with the cross-tier edges of the skeleton required in the tier order (``tier''); six draws per skeleton at $p=4$ and three at $p=5$. For each resulting class $\cG$, every $(G,e)$ with $G\in\cG$ was classified: admissibility of $e$, equality of the explanation set with $[H]_Y$, the chain property of $\MP(H)$, \eqref{eq:K2} and \eqref{eq:K1p} by exact rational computation at two random rational points per class (equality is only checked; inequality is proved), and the Perkovi\'c criteria on the common orientations of $[H]_Y$ and of each class $S_c$, and on $\MP''(S_c)$ computed by enumeration. For every admissible pair, \eqref{eq:K2} agreed with the Perkovi\'c criterion on $[H]_Y$. At $p=5$ there are $30746$ pairs in the general experiment, $1493$ of which are not admissible, all with $e\in R$, and $19382$ in the tier experiment, all admissible; among admissible pairs, $5951$ general ones have a non-chain $\MP(H)$ and none of the tier ones does (the smallest non-chain example on five nodes is $R=\{3\to0\}$, $e=\{1,2\}$, where $3\to0-2-3$ is a partially directed cycle). The criterion for \eqref{eq:K1p}, in the form with common orientations of $S_c$, had no disagreement with the identity check on $26494$ and $18067$ classes, and in the form $\MP''(S_c)$ none on the $25685$ classes with admissible $e$ in the general experiment; the step $\mathrm{Ext}''=S_c$ in the proof of \cref{thm:K1p-graph} failed only for non-chain $\MP(H)$, in $14$ of $692$ classes with a bad path in $\MP''(S_c)$ and in none of $2315$ with chain $\MP(H)$ (scripts \texttt{d27\_background\_knowledge.py}, \texttt{d27b\_toy\_instance.py}). The four-variable model with $R=\{W\to X\}$ is checked for all eight pairs $(G,e)$ by identity tests at four random rational points, inequalities being proofs. An independent sample on six nodes (scripts in \texttt{d27c\_independent\_sample}, random skeletons with seven to ten edges, $240000$ general and $90000$ tier draws) is not exhaustive; it found no disagreement of either criterion with the identity checks, $\MP(H)$ non-chain in $18.6\%$ of the general draws and in none of the tier draws, and $\mathrm{Ext}''\ne S_c$ with a bad path in $68$ classes, all with non-chain $\MP(H)$. \Cref{lem:tier} is also checked numerically (script \texttt{d37\_tier\_chain.py}) on $50000$ random DAGs $H$ for each $p=4,\dots,10$, with the nodes other than $Y$ assigned at random to two to four ordered tiers and $Y$ alone in the last tier: in all $350000$ cases $\MP(H)$ was a chain graph with chordal chain components, its closure under R1 alone was the same graph, as in Lemma~1 of \citet{bang2023wiser}, and in the $347210$ cases with at most twelve undirected edges in the CPDAG, enumeration showed that $[H]_Y$ is the interventional Markov equivalence class used in the proof and that its common orientations are those of $\MP(H)$. As a control, with a single required edge, an undirected edge of the CPDAG oriented as in $H$, instead of tiers, the closure was not a chain graph in $86031$ of $314474$ draws. \Cref{lem:tier-cut} is checked separately (script \texttt{d39\_cut\_chain.py}), with $50000$ draws for each $p=4,\dots,10$ ($346622$ non-empty graphs) of a random DAG $G$ with $Y$ a sink, one to three random node sets $I_j$, $R$ the edges of $G$ with exactly one endpoint in some $I_j$ (edges at $Y$ included, which does not change $[H]_Y$), and a random edge $e$. In all $343532$ cases with at most twelve undirected edges in the CPDAG, the union of the members of $[H]_Y$ was a chain graph with chordal chain components and equal to the Meek closure. In $157046$ of these cases the weak edge was required and not at $Y$, and in $4075$ of them it was not admissible; the $11745$ such cases with tier knowledge had no failure.

\paragraph{Checks of \cref{thm:E}.} Over the same enumeration, the non-equivalent ordered pairs $(G^*,G')$ in $\Xeq$, counted once per triple $(G^*,G',e)$, number $144$ ($p=4$) and $32256$ ($p=5$). Every pair falls into one of the cases of \cref{app:E}, the witness is a local Markov constraint of $G'$, and the derivative at the sparse point is exactly non-zero. A check that uses neither the case analysis nor sparse points, at a random rational point with all coefficients of $D^*$ non-zero, also succeeds for every pair (script \texttt{d21\_e1\_sparse\_witness.py}). An earlier version of this paper certified the condition on $1204$ triples drawn from the protocol below; \cref{thm:E} supersedes it.

\paragraph{The lower bound function.} $c_P$, $c_K$ and $\Psi_\alpha(K)$ of \cref{thm:K}(i) are computed by numerical quadrature and a grid over $d$ (script \texttt{d26\_lower\_bound\_psi.py}); at $\alpha=0.05$, $c_P=0.8498863$ (the value of \cref{thm:lfp}'s crossing term) and $c_K=0.19526$ (rounded down to $0.1952$ in the text), attained near $d=1.405$. At the fourteen points $(w,b)$ of \cref{tab:price-curve}, $\Psi_{0.05}(K)$ is smaller than the dual lower bound, between $0.849$ and $1.32$ against $1.02$ to $2.73$, which is consistent with $\Psi_\alpha$ being a valid but cruder lower bound; in the five-node configuration of \cref{tab:simK} with $K=5.66$, $\Psi_{0.05}(K)=1.63$, against $1.53$ for the closed form $c_P/2+c_KK$; the four-variable model is an instance of the graph classification (with $R=\{W\to X\}$) but its parameter class, with bounded coefficients and variances in $[1,\kappa]$, differs from $\cM_k$, so this is a plausibility check. An independent recomputation (script \texttt{d26b\_psi\_independent\_check.py}) agrees.

\paragraph{The constants of \cref{thm:explicit}.} $\varpi_q(\gamma)$ is computed from the non-central $\chi^2$ distribution function by a grid search followed by a one-dimensional optimization, and $c_K$ in the same way (script \texttt{d40\_T2\_constants.py}). At $\alpha=0.05$, $\gamma_1=\gamma_2=0.01$, the upper slope $\varpi_q(\gamma_1)/(2z_{1-\alpha/2})$ is $0.3664$, $0.4146$, $0.4467$, $0.4717$, $0.4926$ for $q=1,\dots,5$, $0.5669$ for $q=10$ and $0.6180$ for $q=15$ (rounded up), against $c_K=0.19526$, ratios between $1.87$ and $3.17$ for $q\le15$; the intercept coefficient $(z_{1-\alpha'/2}+\sqrt{2\log N})/z_{1-\alpha/2}$ is $1.1072$, $1.7079$, $1.8635$ and $2.2446$ for $N=1$, $2$, $3$ and $12$. For \cref{cor:RK} (script \texttt{d42\_T5\_constants.py}), $\varpi_d(0.01)$ is non-decreasing in $d=1,\dots,10$, so $\bar\varpi_d=\varpi_d$ there. The factor $a_m/c_P+(m-1)b_m/c_K$ of \cref{cor:RK}, minimized over $\gamma_3$, is $3.3912$ at $\gamma_3=0.0300$ for $m=2$ (the continuous optimum is $0.03003$; at $\gamma_3=0.03$ exactly $a_2=1.59403$ and $b_2=0.29595$), and $5.693$, $7.835$, $10.064$ and $12.379$ for $m=3,\dots,6$ (rounded up), at $\gamma_3$ between $0.039$ and $0.047$; at $\gamma_3=0.01$ it is $3.589$ for $m=2$, and at $\gamma_3=0.03$ it is $5.759$, $8.042$, $10.453$ and $12.975$ for $m=3,\dots,6$. With the constants rounded as in \cref{cor:RK}, $1.595/0.849+0.296/0.195=3.397$. With the levels $\gamma_2=\gamma_3=0.01$ of the simulations, the factor of \cref{thm:explicit}(c) for $m=2$ is $3.658$, and a grid over $(\mathcal R,K)\in[1,60]\times[0,80]$ gives $3.608$, below the closed form as it should.

\paragraph{The canonical experiment.} For \cref{prop:canonical} and the comparison in \cref{sec:gen-constant} we solved, for $27$ values of $(1,\alpha_1,\alpha_2,\rho)$ with $\alpha_1\in\{0.5,1,2\}$, $\alpha_2\in\{-1,0.5,1.5\}$ and $\rho\in\{0.4,0.7,0.9\}$, the discretized minimax linear program, which allows randomization and does not restrict to intervals, so that it approximates $\Pi^*_{\mathrm{can}}$ (grid steps $0.6$ for $Z$, $0.42$ for the offset, $\tau$ on the symmetric grid with step $0.7$ that contains $0$ and covers $[-6,6]$, that is, $|\tau|\le6.3$); it reproduces the four-variable program of \cref{app:lp} for four values of $(w,b)$ (script \texttt{d28\_canonical\_lp.py}, differences at most $0.025$; the two programs differ in grid, the $\tau$-grid of the canonical program being symmetric and containing $\tau=0$, and its steps for $Z$ and the offset being $0.6$ and $0.42$ in canonical units). For each of the $27$ values, $K_{1\to2}$ and $K_{2\to1}$ are matched, in either order because the price is symmetric under exchange of the labels, by a least-squares fit with residual below $10^{-6}$ to a four-variable parameter with $w,b\in[0.05,20]$ (this excludes degenerate matches with $w,b\to0$ and $b/w$ fixed), and the program is solved at both parameters (scripts \texttt{d28b\_depends\_on\_K.py}, \texttt{d28c\_summarize.py}). For $21$ of the $27$ such a match exists; for $19$ of these $21$ the general price is below the four-variable price, the average absolute difference is $0.025$, and the difference exceeds $0.05$ only for the one with the largest $K$, where it is $-0.14$; the Spearman correlation of the price with $\max\{K_{1\to2},K_{2\to1}\}$ over the $27$ cases is $0.88$. None of these parameters is at the degenerate point ($|Q|\ge0.012$), so matched parameters have the same $\mathcal R$. The $43$ parameters of the bound comparison in \cref{sec:gen-constant} are these $27$ and the $16$ of the form $(1,\alpha_1,0,\rho)$ with $\alpha_1\in\{0.5,1,2,4\}$ and $\rho\in\{0.4,0.7,0.9,0.97\}$ (script \texttt{d29b\_dual\_scan.py}); on them the ratio of the program value to $\max\{\Psi_\alpha(K),c_P\mathcal R\}$ has median $1.09$ and maximum $1.56$, the bound never exceeding the program value, and the ratio to $\Psi_\alpha(K)$ alone has median $1.44$ and maximum $3.92$ (script \texttt{d29d\_bound\_ratio.py}), and the dual bound of \cref{thm:lfp} extended to the canonical experiment (script \texttt{d29\_canonical\_dual.py}, numerical, not certified) lies below the program value by $0.01$ to $0.08$, and at $(1,1,0,\sqrt{1/2})$ it is $1.2737$, against the computer-assisted bound $1.2743$ in the four-variable model at $(w,b)=(1,1)$. For the pair with the largest difference ($\alpha_1=2$, $\alpha_2=-1$, $\rho=0.9$) a refinement of the grid in $\tau$, in $Z$ and the offset, and in all three changes each of the two prices by about $0.012$ and the difference only from $-0.134$ to $-0.136$ (scripts in \texttt{d28d\_independent\_checks} of the supplementary material). On the same $27$ parameters, the largest value over a grid of $60$ values of $\Lambda$ in $[0.3,10^6]$ of the bracket in the definition of $\delta_{\mathrm{can}}$ (\cref{lem:strict}), a lower bound for $\delta_{\mathrm{can}}$, computed with a trapezoidal rule in $v$ and Gauss--Legendre quadrature in $\zeta$, is positive in every case and below the program value minus one; it is invariant under scaling of $(a_0,\alpha_1,\alpha_2)$, as it should be (script \texttt{d44c\_strict\_bound.py}). The same directory contains an independent check of the canonical geometry in the full eight-parameter four-variable model, in which $\dim L_0$, $\dim L_1$, $\dim L_2$ are $7$, $8$, $8$.

\paragraph{The identity for $\mathcal R$.} Identity \eqref{eq:R-identity} is checked at $10^4$ random canonical parameters (largest relative error $4\times10^{-12}$, and $K_{1\to2}>1$ exactly when $K_{2\to1}>1$) and on the $43$ parameters of the dual scan; the value $1.2737$ of the canonical dual bound at $(1,1,0,\sqrt{1/2})$ is saved with it (script \texttt{d29c\_R\_identity.py}).

\paragraph{The canonical experiment and $K$.} On the $43$ canonical parameters, the lower bound $\max\{\Psi_\alpha(K),c_P\mathcal R\}$ is below the value of the discretized minimax program by a factor $1.09$ in the median and at most $1.56$, against $1.44$ and $3.92$ for $\Psi_\alpha(K)$ alone; canonical parameters with the same pair $(K_{1\to2},K_{2\to1})$ as a four-variable parameter have prices that differ from it by $0.025$ on average and by up to $0.14$. So $K$ determines the price approximately but not exactly. At the degenerate point $\mathcal R$ is free: $(w,b)=(0.5,0.5)$, $(1,1)$ and $(2,2)$ all have $K=1$, and their dual lower bounds $1.179$, $1.274$ and $1.909$ (\cref{tab:price-curve}) increase with $\mathcal R=\sqrt{1+b^2}$.

\paragraph{The cross-model constant.} For each benign pair $(G,e)$ of the exhaustive enumeration at $p=5$ whose $\Xeq$ has at least two Markov classes ($2472$ pairs), we evaluate $K(P_0)$ of \cref{thm:K} at one random base point, with coefficients $\pm U(0.5,1.5)$ and noise variances $U(0.5,1.5)$, using dual numbers for the first derivatives and the closed forms of \cref{lem:Dproj}. The same code reproduces the four-variable formulas for $K$ to six digits. For each edge $f$ of $H$ we set its coefficient to $\lambda\in\{10^{-1},10^{-2},10^{-3}\}$ and classify the growth of $K(P_\lambda)$ by the two log-log slopes. For \cref{thm:K-lambda}, conditions (A) and (F) are decided by graph operations, (A) over all acyclic orientations of the skeleton in the two Markov classes. For each distinct triple $(\Xeq,e,f)$ and each ordered pair, $I$ and $\delta$ and their first $\lambda$-derivatives at $\lambda=0$ are computed in exact rational arithmetic with numbers of the form $a+b\varepsilon+c\lambda+d\varepsilon\lambda$, $\varepsilon^2=\lambda^2=0$, at up to three random rational points; a non-zero value proves non-vanishing for almost every $\phi_0'$, while three zero values are recorded as vanishing without proof (scripts \texttt{d25\_two\_sided.py}, \texttt{d25b\_classlevel.py}, about one hour on six cores). On five nodes, $K(P_0)>0$ for $240$ of the $2472$ pairs $(G,e)$, that is, $96$ distinct pairs $(\Xeq,e)$, with median $0.71$ and largest value $4.6$ over the $240$ values; over the $2760$ distinct triples $(\Xeq,e,f)$, the $528$ ordered pairs that satisfy (F) but not (A) all violate (C); and the predicted growth of $K(P_\lambda)$, $1/|\lambda|$ for $24$ triples, $|\lambda|$ for $120$ and $1$ for $384$, agrees with the evaluation at $\lambda\in\{10^{-1},10^{-2},10^{-3}\}$ for every triple, while $K\equiv0$ for the remaining $2232$. The computation uses one graph $G$ per triple. For $100$ randomly chosen pairs $(\Xeq,e)$ out of $594$ we repeated it for up to two further graphs $G$ with the same $\Xeq$, using their own $H=G-e$ and evaluating (A) with respect to that $H$: among $8778$ ordered pairs, $I(P_{\lambda=0})>0$ for $7162$, none of which satisfies (A), and $I(P_{\lambda=0})=0$ with order of vanishing one for $1616$, all of which satisfy (A) (script \texttt{d25c\_repcheck.py}). For the graphical form of (iii), (F) is checked against the witness edges and (C) against $(G_\wedge-e)-f\simeq H-f$ on all $24816$ pairs, the tangent-space form of $I(P_{\lambda=0})=0$ (rank of the stacked generators at two random base points) agrees with the exact computation for every pair, and the base point with only the two edges at a common neighbour of $z$ and $x_f$ is a witness for each of the $528$ pairs that satisfy (F) but not (C) (scripts \texttt{d30\_IA\_mechanism.py}, \texttt{d30d\_reformulations.py}, \texttt{d30e\_collider\_side.py}). \Cref{lem:onev} and the intermediate claims of its proof were checked for every pair of DAGs on four and five nodes that differ in one v-structure ($216$ and $21210$ cases, counting every representative of the first class, of which $204$ and $19590$ satisfy the hypothesis) and on random skeletons with six and seven nodes, and the descent was run along every branch that occurs; the last branch of type $\beta$, the switch to type $\alpha$, was never reached, and the proof does not depend on whether it is (scripts \texttt{d30g\_lemma\_Q.py}, \texttt{d30i\_Q\_proof\_checks.py}). Condition (D) of \cref{thm:K-lambda}(iv) was compared with exact values of $\delta(P_{\lambda=0})$ (two random rational points, both orders of each pair) on configurations that satisfy (F), (C), \eqref{eq:K2} and \eqref{eq:K1p}, generated exhaustively on five and six nodes from the two shapes that the descent of \cref{lem:onev} produces ($5443$ on six nodes, counted by pairs of classes and $U$, of which $496$ satisfy (D)) and sampled on seven nodes ($6945$), with no disagreement; the local criterion on the CPDAG agrees with the definition of (D) in all cases. Independently of \cref{lem:onev}, with (F) and (C) decided from their definitions, the exhaustive five-node computation above and $950$ randomly drawn benign pairs $(\Xeq,e)$ with six and seven nodes and $R=\emptyset$ gave no disagreement among $4236$ and $14116$ ordered pairs, $96$ and $370$ of them with $\delta(P_{\lambda=0})\ne0$ (scripts \texttt{d33\_random\_p6\_delta.py}, \texttt{d33b\_Dprime.py}). With one or two random required edges, the weak edge admissible and possibly required, $820$ random benign pairs $(\Xeq,e)$ with six and seven nodes gave no disagreement among $5854$ ordered pairs that satisfy (F) and (C), $376$ of them with $\delta(P_{\lambda=0})\ne0$ and $696$ with $e\in R$, all of which have $\delta(P_{\lambda=0})=0$ (script \texttt{d33c\_Dprime\_R.py}). Without the requirement on the parents of $x_e$ the condition fails: on six nodes $125$ configurations, counted by the graph $H$ ($109$ counted by pairs of classes and $U$), satisfy all other requirements of (D) and have $\delta(P_{\lambda=0})=0$ for every ordered pair of explanations, and $K(P_\lambda)$ is constant in $\lambda$ for the three we evaluated (scripts \texttt{d32e\_criterion.py}, \texttt{d32e\_search.py}, \texttt{d32e\_verify.py}, \texttt{d32e\_Klambda.py}). On random skeletons with six, seven and eight nodes biased towards several common neighbours of $z$ and $x_f$, among $44000$ ordered pairs satisfying (F), $I(P_{\lambda=0})=0$ (tangent-space test) held exactly when (C) held, and the two-edge witness existed whenever it failed (script \texttt{d31\_random\_p67.py}). Counted over all $2472$ pairs $(G,e)$ and edges $f$, without the reduction to distinct triples, the ordered pairs that satisfy $G^*-f\simeq G'-f$ but not $(G^*-e)-f\simeq H-f$ number $3264$ of $27120$; for one of them, at one random base point, $I_{G^*,G'}(P_\lambda)$ equals $0.2146$, $0.2353$, $0.2373$, $0.2375$, $0.2375$ at $\lambda=10^{-1},\dots,10^{-5}$. For \cref{prop:Kpos}, the closed forms of $\delta$ and $I$ in its proof agree exactly with the definitions (dual-number routine above and \cref{lem:Dproj}) on all $1516$ ordered pairs at $p=5$ and on $30725$ sampled pairs at $p=6$. On all skeletons with six nodes there are $114804$ benign pairs $(\Xeq,e)$ with at least two Markov classes; $K(P_0)>0$, in exact arithmetic at three random rational points, for $28000$ of them, exactly those in which some member satisfies (a) and (b), and the criterion takes the same value on all members of each class. Counting single Markov classes, at which $K=0$ by definition, $(\Xeq,e)$ is free in the sense of \cref{thm:free} for $8186$ of the $8426$ benign pairs $(G,e)$ on five nodes ($3722$ of the $3818$ distinct pairs $(\Xeq,e)$) and for $1358485$ of the $1448757$ on six nodes ($546025$ of $574025$; script \texttt{d44b\_free\_fraction.py}, which reads the outputs of \texttt{d34a}). The condition that $e$ lies on a directed path from $X$ to $Y$ in some member of $\Xeq$, which on five nodes is equivalent to $K(P_0)>0$, fails for $3840$ of them: $1296$ have $K>0$ and no such path, and $2544$ have such a path and $K\equiv0$. Random samples with seven and eight nodes ($1900$ pairs $(\Xeq,e)$) and with one or two required edges on six and seven nodes ($700$) gave no disagreement with the criterion (on five nodes the same holds for all $594$ pairs; scripts \texttt{d34a\_Kpos\_exhaustive.py}, \texttt{d34c\_random\_p7.py}, \texttt{d34e\_bk.py}, \texttt{d34b\_inspect.py} for two six-node examples of the two mechanisms in \cref{sec:gen-constant}). A second computation, which evaluates $\delta$ from its definition and the criterion by its own d-separation routine (moralization), and shares only the enumeration of graphs and the benign-case filter with the first, found no disagreement on $2600$ random benign pairs $(\Xeq,e)$ with six and seven nodes, $438$ of them with $K(P_0)>0$, while the directed-path condition failed for $36$ (script \texttt{d35\_Kpos\_crit\_check.py}). For \cref{rem:Oset}, at the same base points, $\delta^O=\delta$ for $40374$ of the $40656$ ordered pairs of explanations; among the $240$ pairs $(G,e)$ with $K>0$, $K^O=K$ for $186$, and $K^O/K$ lies in $[0.93,3.25]$; $K^O>0$ never occurs with $K=0$, whereas $K^{\mathrm{adj}}>0$ with $K=0$ for $54$ pairs (script \texttt{d24f\_oset.py}). Counts are reported both per pair $(G,e)$ and per distinct pair $(\Xeq,e)$, since all $G$ with the same $\Xeq$ give the same problem (scripts \texttt{d24\_K\_local.py}, \texttt{d24c\_merge.py}, \texttt{d24e\_review\_checks.py}). The values of $K$ and the growth exponents are numerical, with one random base point per pair $(G,e)$; the exponents agree across all pairs with the same $(\Xeq,e,f)$.

\paragraph{Uniformity in the margin.} For \cref{thm:uniform} we generate, with the generator of the check of \cref{thm:K-lambda}(iv) (script \texttt{d32e\_criterion.py}), all configurations on five and six nodes, with only $Y$ a sink, in which an ordered pair satisfies (F) and (C) and $e$ is benign in $H$, one per pair of Markov classes and set $U$ (scripts \texttt{d43a\_fbenign\_scan.py}, \texttt{d43f\_condition\_M.py}). On five nodes there are $91$; the $4$ that satisfy (D) all have $f$ malignant in $H_2$, and $63$ satisfy \eqref{eq:M}. On six nodes there are $5443$: of the $496$ that satisfy (D), $64$ satisfy \eqref{eq:M}, $64$ violate it with $f$ benign in $H_2$, and $368$ have $f$ malignant in $H_2$; of the $4947$ that violate (D), $3093$ satisfy \eqref{eq:M}. In every configuration with \eqref{eq:M}, $f$ is benign in $H_2$, as shown in \cref{app:uniform}. For the configurations with $f$ benign in $H_2$ we draw a base point on $\cM_{H_2}$ (coefficients $\pm U(0.4,1.2)$, noise variances $U(0.6,1.6)$) and compute, for every true graph in $\cX_2$ and every candidate $G'$ that explains $P_{00}$, the projection distance and the effect discrepancy at $1460$ points $(\varepsilon,\lambda')$ on two circles of radii $10^{-3}$ and $10^{-2}$ and on the two axes, and classify the pair by the local form of $D$ (script \texttt{d43e\_pair\_types.py}; all $128$ six-node configurations with (D) and $f$ benign in $H_2$, all $3093$ six-node and $63$ five-node configurations with \eqref{eq:M} and without (D), two base points each). Under \eqref{eq:M} every pair in $\cX_2$ has one of the forms (c)--(e) of \cref{app:uniform}, with no degeneration along a line or at higher order; the ratio of $|\Delta|$ to the right side of \eqref{eq:uniform-star} with $C_2=1$ is at most $3.86$, and when (D) fails, $|\Delta|/\sqrt D\le2.99$. Candidates outside $\cX_2$ reach ratios of order $10^{300}$, with $D=0$ to machine precision and $\Delta\neq0$; they are the cancellation explanations, which the $k$-TF screen removes. On the axis $\lambda'=0$, pairs with (F) and (C) have $D\le2\times10^{-29}$, and when (D) fails, $|\Delta|\le4.5\times10^{-15}$ in all $3156$ configurations with \eqref{eq:M}, as \cref{lem:Vf,lem:Af} assert, and also in the $1091$ with $f$ benign in $H_2$ that violate \eqref{eq:M} (script \texttt{d43d\_lambda0\_effect.py}). The modulus $\omega(\eta;\lambda,A)$, the largest $|\Delta|$ over pairs with $D\le\eta^2$, $|\varepsilon|\le A$ and $|\lambda'|\ge\lambda$ whose two members belong to the class with $k=\min\{0.05,k_{00}/2\}$, where $k_{00}$ is the smallest $k$-triangle-faithfulness ratio of the true graphs at $P_{00}$, is computed for $\eta\in[10^{-6},10^{-2}]$, $\lambda\in[10^{-4},10^{-1}]$ and $A\in[10^{-5},10^{-1}]$, two values per decade, from points with $|\varepsilon|\in[10^{-8},0.2]$ and $|\lambda'|\in[10^{-5},0.2]$, five per decade and both signs (script \texttt{d43c\_modulus\_classes.py}, one base point per configuration). For the $64$ configurations with (D) and \eqref{eq:M} and for the four-variable model, $(\eta+\omega)/(\eta+\min\{A,\eta/\lambda\})$ lies in $[0.2956,2.616]$; for $61$ configurations with \eqref{eq:M} and without (D), chosen evenly in the order of the census among those with at most $16$ explanations of $P_{00}$, $(\eta+\omega)/\eta\le1.984$; and for the $64$ configurations with (D) and $f$ benign in $H_2$ that violate \eqref{eq:M}, the largest value of $(\eta+\omega)/(\eta+\min\{A,\eta/\lambda\})$ is between $2.2\times10^3$ and $1.1\times10^4$ according to the configuration.

\paragraph{Other positions of the two edges.} The configurations above all have an (F) pair, so $e$ and $f$ share an endpoint whose other neighbours are not adjacent. To check the other positions we enumerate, on five nodes with only $Y$ a sink, all triples $(H,e,f)$, one per Markov class of $H$, with $e$ benign in $H$ and $f$ satisfying \eqref{eq:M}: $3169$ in which $e$ and $f$ lie in a common triangle, $2936$ in which they have no common endpoint, and $3638$ in which they share an endpoint whose other neighbours are not adjacent (script \texttt{d43g\_geometry\_census.py}, one base point each). No pair in $\cX_2$ degenerates along a line or at higher order, and in the configurations without an (F) pair, $|\Delta|/\sqrt D$ is at most $4.08$. When $e$ and $f$ lie in a common triangle, the ratio of the $k$-TF element of $f$ in that triangle tends, along $\varepsilon=0$, to a finite positive limit $\ell$ in $2916$ of the $3169$ configurations, as in (iii) of \cref{app:uniform}; the threshold built only from the other elements, half their smallest ratio at $P_{00}$, exceeds $\ell$ in $8$ of them, which is why $\ell$ enters $k_0$. With $k=0.99\,k_0$, the points with $\varepsilon=0$ and $\pm\lambda'\in\{10^{-5},10^{-4},10^{-3}\}$ satisfy $k$-TF for every member of $\cX_2$ in all but $13$ of the $9743$ configurations, and in these $13$, whose $k_0$ is below $3\times10^{-3}$, they do so at $\pm\lambda'\in\{10^{-6},10^{-7}\}$; the radius $r_2$ of \cref{app:uniform} depends on the base point.

\paragraph{Sampling protocol for the simulations.} The five-node experiment of \cref{sec:sim} draws its cases from the following protocol. We draw $60$ skeletons on five nodes with $X=1$ and $Y=4$, containing the edge $X-Y$ and each other edge with probability $0.6$ (random seed $7$). For each acyclic orientation $G$ with $Y$ a sink we draw coefficients $\pm U(0.5,1.5)$ and noise variances $U(0.5,1.5)$, keep $G$ if the total effect is identifiable at this point, and consider every edge $e$ of $G$, with base point $P_0$ obtained by setting the coefficient of $e$ to zero. This gives $1202$ pairs $(G,e)$, counted with multiplicity over the $60$ skeletons. The graphical criterion of \cref{prop:graph}, a floating-point test and an exact rational test of \eqref{eq:K2} classify all of them in the same way ($1046$ benign, $156$ malignant; script \texttt{d19\_k2\_graph\_criterion.py}). On the $138$ distinct malignant pairs $(H,e)$, an exact rational computation finds two members of $\Xeq$ with different effects in every case, and none of the $1048$ pairs (triangle containing $e$, legal $W$) has $\rho_{P_0}=0$, in line with \cref{lem:M1} (script \texttt{d16b\_malignant\_m1.py}). Over the true models drawn in this protocol, the smallest $k$-triangle-faithfulness ratio has median $0.039$ (script \texttt{d12d\_g3\_ktf\_check.py}).

\subsection{Dual bound and linear program}\label{app:lp}

The linear program discretizes the limit experiment of \cref{prop:limit} after reduction by translation equivariance in $\tilde\beta$ (Hunt--Stein), with $t$-grid step $0.5$ on $[-5,5]$, a $(P,E)$ grid of step $0.35$ and a $u$ grid of step $0.22$, and is solved with HiGHS (script \texttt{d4\_rho\_star\_lp.py}). Its value is not a bound: off-grid mass is counted as covered and the cells are piecewise constant, which bias it in opposite directions. Its dual multipliers define the witness $(\hat\pi,\hat\mu)$ of \cref{sec:toy-constant}: $\hat\pi$ is the normalized vector of raw weights and $\hat\mu$ is $1.03$ times the raw multipliers, both stored in hexadecimal floating-point form (file \texttt{witness\_w1b1\_fine.json} in the top directory of the released code). The integral in $R(\hat\pi,\hat\mu)$ is evaluated by the midpoint rule on boxes that cover the effective support, with Gaussian tail mass outside subtracted (script \texttt{d8\_lfp\_bound.py}); the three grids $(0.08,0.04)$, $(0.05,0.025)$ and $(0.03,0.015)$ give values that differ by at most $3\cdot10^{-6}$. These midpoint values are not bounds; the computer-assisted bound is given next.

\subsection{Computer-assisted evaluation of the dual lower bound}\label{app:certify-lfp}

Write $D=A-G$ for the integrand of $R(\hat\pi,\hat\mu)$, with $A(p,e)=\sum_ka_k\phi_s(p-t_k)\phi_s(e-m_k)$, $a_k=\hat\pi_k/\orl_{g_k}$, and $G(p,e,u)=\sum_k\hat\mu_k\phi_s(p-t_k)\phi_s(e-m_k)\phi_s(u-\psi_k)$; $31$ of the $42$ witness points carry positive weight. Let $\mathcal B$ be the box that extends $8s$ beyond the extreme centres in each coordinate. The certificate rests on the following exact inequalities.
\begin{enumerate}
\item Outside $\mathcal B$, $\min\{0,D\}\ge-G$, so the outside contribution is at least $-\sum_k\hat\mu_k(1-\Prob_k(\mathcal B))$, a Gaussian box probability; it is $6.3\cdot10^{-15}$.
\item On a box $Q\subseteq\mathcal B$, each factor $\phi_s(x-c)$ attains its maximum at the projection of $c$ onto the interval and its minimum at the farther endpoint, which gives bounds $A_{\mathrm{lo}}\le A\le A_{\mathrm{hi}}$ and $G_{\mathrm{lo}}\le G\le G_{\mathrm{hi}}$, and $\int_QA$, $\int_QG$ are products of Gaussian interval probabilities. With $D_{\mathrm{lo}}=A_{\mathrm{lo}}-G_{\mathrm{hi}}$ and $D_{\mathrm{hi}}=A_{\mathrm{hi}}-G_{\mathrm{lo}}$: if $D_{\mathrm{lo}}\ge0$ the contribution of $Q$ is $0$; if $D_{\mathrm{hi}}\le0$ it is $\int_QD$; otherwise it is at least $\max\{|Q|D_{\mathrm{lo}},\ \int_QD-|Q|D_{\mathrm{hi}},\ -\int_QG\}$, each of which is the integral of a pointwise lower bound of $\min\{0,D\}$. Since also $\int_Q\min\{0,D\}\le\min\{0,\int_QD\}$, each box has a computable gap between lower and upper bound.
\item Starting from boxes of side $0.5$, every box whose gap exceeds $10^{-11}$ is split into eight, until all gaps are below this threshold. This evaluates $6.4\cdot10^7$ boxes and leaves $5.6\cdot10^7$; the total gap is $6.0\cdot10^{-5}$.
\item The constant $(1-\alpha)\sum_k\hat\mu_k$ is computed exactly from the stored binary rationals. All other quantities are computed in double precision, and the bound is reduced by an explicit allowance $10^{-9}S+10^{-14}NW$, where $S$ bounds the sum of absolute values of all summands, $N$ is the number of evaluated boxes and $W=\sum_ka_kL_u+\sum_k\hat\mu_k$ with $L_u$ the width of $\mathcal B$ in $u$. The allowance is $8.9\cdot10^{-6}$; it exceeds by many orders of magnitude the rounding error of double-precision arithmetic with correctly implemented elementary functions. Recomputing all box quantities for $400$ randomly chosen boxes, coarse and refined, with $50$-digit arithmetic (mpmath) gives a largest deviation of $5.6\cdot10^{-17}$.
\end{enumerate}
The result is $1.27438870\le R(\hat\pi,\hat\mu)\le1.27445732$ (script \texttt{d20\_lfp\_certify.py}; $134$ seconds on seven cores), consistent with midpoint evaluations $1.2744552$, $1.2744577$ and $1.2744562$ on three grids. In the optimal solution of the discretized linear program at $(w,b)=(1,1)$, of the $529$ cells of $(P,E)$ with $|p|,|e|\le4$, $522$ carry a single interval and $7$ the empty set, and no cell carries two components. The argument is computer-assisted and relies on the double-precision error model of item~4 rather than on a formally verified interval library.

\subsection{Simulations}\label{app:sim}

\paragraph{The four-variable experiment.} Naive selects $C$ or $M$ according to which conditional independence $t$-statistic is smaller; this is the rule \textup{\textsc{cmp}} of \cref{prop:naive}. The noise variances are the least favourable values $(\sigma_X^2,\sigma_Z^2)=(2,2)$ for $\kappa=4$, with $\sigma_W^2=\sigma_Y^2=1$, $\beta=0$, $b=1$; we vary the true graph, $w\in\{0.05,0.2,1\}$, $a\in\{0,0.1,0.3,1\}$ and $n\in\{500,2000,8000\}$, with $20000$ replications per cell (script \texttt{s1\_sim\_toy.py}). \Cref{tab:simA} in \cref{app:sim} summarizes the $72$ cells. Naive undercovers badly whenever the confounder and mediator models are hard to distinguish, with observed coverage down to $0.489$. The minimum observed coverage of \DSV\ is $0.962$; it is conservative because of the Bonferroni split. Its length ratio grows with $n$ when $w=0.05$, consistent with the malignant regime where $\lambda\sqrt n$ is of order one; when $w=1$ it is bounded and decreases from $2.74$ to $2.40$ as $n$ grows; $w=0.2$ is in transition. The union set has bounded length throughout, as \cref{prop:union} predicts; its largest ratio increases slightly with $n$ because the two intervals overlap less as they shrink. When $a=0$, the benign crossing, the length of \DSV\ relative to the \textsc{Oracle} is $1.95$ under $C$ and $1.13$ under $M$ for every $w$ and $n$; relative to the efficient interval these are multiplied by $\sqrt{1+w^2/2}$, from $1.00$ at $w=0.05$ to $1.22$ at $w=1$, and the observed ratios of the \textsc{Oracle} to the efficient interval, $1.00$, $1.01$ and $1.23$, agree with \cref{lem:oracle-toy}.

\paragraph{One family with $K>0$.} In the five-node cases of \cref{tab:simB}, $\Xeq$ is a single Markov class in $18$ of the $20$ benign cases, and $K$ vanishes in the other two. We therefore also use one benign five-node configuration with $K>0$: $X\to Y$, $Z\to Y$, the collider $W\to X\leftarrow U$, and the weak edge $e=X-Z$, with only the background knowledge that $Y$ is a sink. Here $\Xeq$ consists of the two graphs with $X\to Z$ and with $Z\to X$, which are not Markov equivalent, and $e$ is benign (\cref{sec:toy-model}). The base point has coefficient $\omega$ on $W\to X$ and $U\to X$, $b$ on $Z\to Y$, $1$ on $X\to Y$, and unit noise variances. As $\omega\to0$ the collider that separates the two explanations disappears and $K$, computed from the closed forms of \cref{lem:Dproj}, increases; both collider edges vanish together, which is outside the single-margin setting of \cref{thm:K-lambda}. The true graph is either explanation, the weak edge has coefficient $c/\sqrt n$, and we estimate the local price by the largest ratio of the mean length to the efficient oracle length $2z_{1-\alpha/2}s_{\mathrm{eff}}/\sqrt n$, over the true graph and $c\in\{0,1,2,3,4,6,8,12,16,20,24,32\}$. This maximizes only along the weak edge, so it approximates the local price of the procedure from below. Besides \DSci\ we run \DSciO{}, which uses the optimal adjustment set of each candidate (\cref{rem:Oset}); in this configuration $K^O=K$. There are $2000$ replications per cell (script \texttt{s4\_sim\_K.py}).

\Cref{tab:simK} shows the results. Over the ten configurations $K$ ranges from $0.42$ to $5.66$. At the largest sample size run for each configuration ($n=32000$ for $\omega\in\{0.18,0.25,0.5\}$ at $b=1$, $n=8000$ otherwise), where the maximizing $c$ is inside the grid, the price of \DSciO{} divided by $1+K$ lies between $0.82$ and $1.08$. On the common grid $c\le32$ it lies between $0.86$ and $1.22$ at $n=8000$ and between $0.91$ and $1.68$ at $n=2000$, but for $\omega\le0.25$ these maxima are at the end of the grid, and over the extended grid the ratio is $2.20$ at $\omega=0.18$, $n=8000$, and $4.60$ at $\omega=0.25$, $n=2000$. The price of \DSci\ is larger by the efficiency loss of parent adjustment, a factor of about $3$ at $\omega=2$, and its first-order bias is governed by $K^{\mathrm{adj}}$ instead of $K$. The minimum coverage over all $611$ cells in the class is at least $0.955$ for both procedures. Near the malignant stratum the local behaviour sets in late: for $\omega\le0.25$ and $n\le8000$ the maximizing $c$ is at or near the end of the grid, with prices up to $17.62$ at $\omega=0.25$, $n=2000$, because \DSciO{} keeps candidates that reverse an edge of the collider at $X$; they are separated from the truth only by partial correlations of order $\omega^2$, so heuristically their rejection requires $n$ large compared with $\omega^{-4}$ . At $n=32000$ the ratio is back in the range above. The threshold $n_1$ of \cref{prop:CI} depends on the base point in this way.

\begin{table}[t]
\centering
\caption{Local price along the weak edge for the five-node configuration with $K>0$: largest ratio of mean length to the efficient oracle length over the true graph and $c\le32$, and minimum coverage over all cells. \DSciO{} uses optimal adjustment sets; here $K^O=K$. $2000$ replications per cell. Prices are rounded up and coverages down. $^\dagger$The maximum is at the end of the grid; over $c\le96$ the price of \DSciO{} is $10.84$, at $c=64$.}
\label{tab:simK}
\begin{tabular}{rrrrrcccc}
\toprule
& & & & & \multicolumn{2}{c}{price} & \multicolumn{2}{c}{minimum coverage}\\
\cmidrule(lr){6-7}\cmidrule(lr){8-9}
$\omega$ & $b$ & $K$ & $K^{\mathrm{adj}}$ & $n$ & \DSci & \DSciO{} & \DSci & \DSciO{}\\
\midrule
2 & 1 & 0.75 & 3.18 & 8000 & 8.90 & 1.87 & 0.960 & 0.962\\
1 & 1 & 0.87 & 2.12 & 8000 & 4.32 & 2.02 & 0.964 & 0.963\\
0.5 & 1 & 1.41 & 2.12 & 8000 & 3.24 & 2.42 & 0.964 & 0.964\\
0.5 & 1 & 1.41 & 2.12 & 32000 & 3.24 & 2.42 & 0.964 & 0.961\\
0.35 & 1 & 2.02 & 2.52 & 8000 & 3.34 & 2.84 & 0.963 & 0.965\\
0.25 & 1 & 2.83 & 3.18 & 8000 & 3.63 & 3.32 & 0.961 & 0.965\\
0.25 & 1 & 2.83 & 3.18 & 32000 & 3.60 & 3.29 & 0.967 & 0.965\\
0.18 & 1 & 3.93 & 4.18 & 8000 & 6.05 & 6.00$^\dagger$ & 0.969 & 0.966\\
0.18 & 1 & 3.93 & 4.18 & 32000 & 4.28 & 4.08 & 0.963 & 0.961\\
0.5 & 0.25 & 0.42 & 0.53 & 8000 & 1.78 & 1.40 & 0.956 & 0.957\\
0.5 & 0.5 & 0.77 & 1.06 & 8000 & 2.25 & 1.72 & 0.959 & 0.964\\
0.5 & 2 & 2.83 & 4.24 & 8000 & 5.33 & 3.89 & 0.962 & 0.964\\
0.5 & 4 & 5.66 & 8.49 & 8000 & 9.83 & 7.11 & 0.964 & 0.966\\
\bottomrule
\end{tabular}
\end{table}

\begin{table}[t]
\centering
\caption{Four-variable model, $\kappa=4$: minimum coverage and maximum length ratio, relative to the efficient interval, over the true graph and $a\in\{0,0.1,0.3,1\}$; $20000$ replications per cell (coverage standard error about $0.0015$). Coverages are rounded down and ratios up. Script \texttt{s1b\_table4.py}.}
\label{tab:simA}
\begin{tabular}{rrccccccc}
\toprule
& & \multicolumn{4}{c}{minimum coverage} & \multicolumn{3}{c}{maximum length ratio}\\
\cmidrule(lr){3-6}\cmidrule(lr){7-9}
$w$ & $n$ & Naive & \DSV & Union & Oracle & Naive & \DSV & Union\\
\midrule
0.05 & 500 & 0.489 & 0.977 & 0.947 & 0.947 & 1.39 & 5.59 & 2.64\\
0.05 & 2000 & 0.503 & 0.977 & 0.949 & 0.949 & 1.37 & 9.87 & 2.67\\
0.05 & 8000 & 0.531 & 0.973 & 0.946 & 0.946 & 1.37 & 14.27 & 2.73\\
0.2 & 500 & 0.586 & 0.974 & 0.947 & 0.947 & 1.39 & 4.43 & 2.66\\
0.2 & 2000 & 0.592 & 0.963 & 0.948 & 0.948 & 1.38 & 5.16 & 2.70\\
0.2 & 8000 & 0.692 & 0.963 & 0.946 & 0.946 & 1.39 & 4.52 & 2.75\\
1 & 500 & 0.787 & 0.963 & 0.950 & 0.947 & 1.69 & 2.74 & 3.23\\
1 & 2000 & 0.870 & 0.963 & 0.948 & 0.948 & 1.68 & 2.70 & 3.27\\
1 & 8000 & 0.945 & 0.962 & 0.948 & 0.947 & 1.68 & 2.40 & 3.33\\
\bottomrule
\end{tabular}
\end{table}

All simulations use joint Gaussian data, regressions with intercept, and exact $t$ and $\chi^2$ quantiles. Seeds are fixed in the scripts: the four-variable experiment uses seed $20261031+i$ for cell $i$, and the five-node experiment draws its cases with seed $20261031$ from the enumeration protocol above and uses seed $1000i+17n+c$ for case $i$. The four-variable experiment uses $72$ cells and $20000$ replications per cell (about $13$ minutes on seven cores). The tables round coverages down and length ratios up from the raw simulation output. The length ratios are relative to the efficient interval $2z_{1-\alpha/2}s_{\mathrm{eff}}n^{-1/2}$ of the true graph, computed exactly for the four-variable experiment (script \texttt{s1b\_table4.py}, from the saved cell means of the experiment of 2026-10-31) and with the dual-number routine of \cref{app:certify} for the five-node experiment. The five-node experiment uses $40$ cases, $3$ sample sizes, $3$ values of $c$ and $2000$ replications per cell; cells in which the true model violates $k$-triangle faithfulness for $c>0$ are skipped (between $0$ and $3$ of $20$ cases per row of \cref{tab:simB}). \Cref{tab:simB} comes from the run described below under ``Double screening with estimated directions'' (scripts \texttt{s12\_sim\_p5\_adaptive.py}, \texttt{s12b\_table.py}); its Naive, Hull and Union columns, and \DSlrO, agree with the earlier runs in all $2824$ compared numbers. The largest benign ratio of \DSOp, $7.73$, $10.89$ and $2.94$ at $n=500$, $2000$ and $8000$ over $c\in\{0,1,3\}$, comes from the case in which the weak edge is $X\to Y$, whose candidate class has $18$ graphs.

\paragraph{Earlier runs of the five-node experiment.} The runs of 2026-11-03 and 2026-11-04 on the same data (scripts \texttt{s2\_sim\_p5.py}, \texttt{s3b\_table5.py}, \texttt{s3c\_table5.py}, \texttt{s3d\_table5.py}) include the conditional independence screen. In the $353$ retained cells \DSci\ and \DSci\ without the $k$-triangle-faithfulness screen give the same interval, and \DSci\ returns the empty set with probability at most $0.0025$. A second rerun (tag \texttt{2026-11-04b}) adds \DSlrO\ and \DSlreff, with the null quantiles of $T_D$ from $400000$ draws for each pattern $(q_j)$ and $n$; all earlier columns agree in all $360$ cells. In the retained cells the length of \DSlrO\ is at most $0.998$ times that of \DSciO{}, with median ratio $0.990$, its empty-set rate is at most $0.017$, and the minimum coverage of \DSlreff\ is $0.942$. \DSci\ and \DSciO{} keep the true graph in at least $99.6\%$ of the replications; the efficiency gain of optimal adjustment is small on malignant edges (median ratio of \DSciO{} to \DSci{} $0.98$, against $0.71$ on benign edges); and the largest benign ratio of \DSciO{}, $8.0$, $13.3$ and $7.3$ at $n=500$, $2000$ and $8000$, comes from the same case, of whose $18$ candidates on average $12.0$, $11.2$ and $4.3$ are kept. On benign edges the median length ratio of \DSciO{} and of \DSlrO\ is about $1.11$ for all $n$, the ratio $z_{1-0.03/2}/z_{1-0.05/2}\approx1.107$ of the quantiles used: in the median case only candidates equivalent to the true graph survive. The median of \DSci\ is $1.8$ because parent adjustment is itself $1.57$ times longer than the efficient interval in the median case. In that case at $n=8000$ and $c=0$ the ratio is $7.00$ for \DSciO{} and $2.63$ for \DSlrO, whose single test per candidate rejects the wrong candidates sooner. On malignant edges the ratios of \DSciO{} and \DSlrO\ grow from about $4$ to $10$ as $n$ grows from $500$ to $8000$ at $c=0$. The minimum coverage over the $353$ retained cells is $0.960$ for \DSci, $0.955$ for \DSciO{} and $0.943$ for \DSlrO, whose likelihood ratio test spends the full level $\gamma_1$ on the true graph.

The experiment with $K>0$ (\cref{tab:simK}) uses the eight orientations of its skeleton with $Y$ a sink as candidates, $k=0.05$, the same levels as \DSci, the Fisher $z$ partial correlation interval as in the five-node experiment (as does the simulation of the instance of \cref{thm:K-lambda} below), the efficient oracle $s_{\mathrm{eff}}$ of the true pair, and seed $20261103+i$ for cell $i$. The main grid has $\omega\in\{2,1,0.5,0.35,0.25,0.18\}$ at $b=1$ and $b\in\{0.25,0.5,2,4\}$ at $\omega=0.5$, both true graphs, $c\in\{0,1,2,3,4,6,8,12,16,20,24,32\}$ and $n\in\{2000,8000\}$ ($480$ cells, about $9$ minutes on three cores). Two supplements cover $c\in\{24,32,48,64,96\}$ at $n\in\{2000,8000,32000\}$ and $c\le20$ at $n=32000$, for $\omega\in\{0.18,0.25,0.5\}$ at $b=1$. Cells in which the true model violates $k$-triangle faithfulness are skipped ($611$ of the cells are in the class); a cell run in two grids contributes both runs to the maximum. Script \texttt{s4c\_table.py} merges the grids and produces the table. The five-node cases are drawn from $669$ benign and $135$ malignant pairs, distinct by $(H,e)$, where the benign ones satisfy the stronger condition that the effect is identifiable in $\MP(G')$ for every $G'\in\Xeq$; of the $20$ benign cases drawn, $18$ have a single Markov class in $\Xeq$, and in the other two $K(P_0)$ vanishes at the drawn base point ($0$ and $2\cdot10^{-15}$ in floating point); $353$ of the $360$ cells are in the class (script \texttt{d36\_sim\_cases\_check.py}). In the experiment of \cref{tab:simK}, the largest prices over the extended grid are $17.62$ at $\omega=0.25$, $n=2000$, $c=96$ and $10.84$ at $\omega=0.18$, $n=8000$, $c=64$; in all eight cells in which the price of \DSciO{} exceeds $1.5(1+K)$, the true graph has $X\to Z$, $\omega\le0.25$ and $n\le8000$, and \DSciO{} keeps on average more than four of the eight candidates, although $\Xeq$ has two members; the others reverse an edge of the collider at $X$, and one with $Z\to X$ is kept in at least $82\%$ of the replications. For the instance of \cref{fig:instances}(b), $K(P_\lambda)$ at $\lambda=10^{-1/2},\dots,10^{-3}$, with all other coefficients $0.8$ and unit variances, satisfies $\lambda K(P_\lambda)\to0.943$, against $K\to1$ for the margin edge $A\to X$ and $K\to0.8$ for the control $H'$ (script \texttt{d38\_Klambda\_instance.py}). For $f=A\to X$, condition (F) fails, because the v-structure that separates the two explanations, $M\to X\leftarrow B$, does not contain $f$ ($A$ and $M$ are adjacent, so $M\to X\leftarrow A$ is shielded), and $K(P_\lambda)\to1$. In the graph $H'$ with edges $X\to Y$, $X\to Q$, $C\to Q$, $C\to Z$, $Z\to Y$, weak edge $e'=Q-Z$ and margin edge $f'=X\to Q$, the explanation with $Z\to Q$ has the extra v-structure $X\to Q\leftarrow Z$, so (F) and (C) hold with $z=Q$, but (D) does not, because $X$ is not an endpoint of $e'$, and $K(P_\lambda)\to0.8$, bounded but positive. The simulation uses the twelve orientations of its skeleton with $Y$ a sink as candidates, $k=0.05$, $\lambda\in\{0.5,0.3,0.2,0.1,0.05\}$, both explanations $G_\wedge,G_\vee$ as true graphs, $c\in\{0,1,2,3,4,6,8,12,16,24,32,48,64\}$, $n\in\{8000,32000\}$ and $1000$ replications per cell ($260$ cells, about five minutes on three cores; script \texttt{s5\_sim\_Klambda.py}); $248$ cells are in the class, with minimum coverage $0.954$ for \DSci\ and $0.955$ for \DSciO{}. At $n=32000$ the price of \DSciO{} along the weak edge is $2.39$, $3.07$, $4.16$ and $6.83$ at $\lambda=0.5$, $0.3$, $0.2$ and $0.1$, where $K(P_\lambda)=K^O(P_\lambda)$ is $1.83$, $3.11$, $4.69$ and $9.42$: between $0.65$ and $0.85$ times $1+K$, the ratio decreasing as $\lambda$ decreases. At $n=8000$ the ratio of the price of \DSciO{} to $1+K$ is between $0.69$ and $0.86$ for $\lambda\ge0.1$; at $\lambda=0.05$ the maximizing $c$ is at the end of the grid at both sample sizes. Heuristically, a wrong explanation is kept while its separating statistic, of order $\sqrt{nI}\,|\varepsilon|$, is below a critical value, and its interval is then shifted by about $|\delta\varepsilon|$, of the order of $K$ oracle half-lengths, which suggests why the ratio of the price to $1+K$ is close to one once $n$ is large enough for the local regime; we have no proof. The order $\omega^2$ of the separation of the candidates that reverse an edge of the collider $W\to X\leftarrow U$ is that of the partial correlation of $W$ and $U$ given $X$, which is a product of the two coefficients of order $\omega$.

\paragraph{All five-node classes with $K>0$.} The $96$ pairs $(\Xeq,e)$ of \cref{sec:sim} are the distinct pairs with $K>0$ among the $594$ distinct benign pairs with at least two Markov classes in the enumeration at $p=5$ ($X=0$, $Y=4$, only the requirement that $Y$ is a sink), each represented by the first $(G,e)$ in the enumeration order; every $\Xeq$ has two Markov classes and $2$ or $3$ members, and the skeleton has $8$ or $12$ orientations with $Y$ a sink. The base point has coefficients $\pm U(0.5,1.5)$ and noise variances $U(0.5,1.5)$, redrawn (at most $13$ times) until one representative of each Markov class has $k$-triangle-faithfulness ratio at least $2k$; $K$ and $I$ are recomputed at this point, so they differ from the values of the enumeration, while $K>0$ does not. Each representative is used as the true graph, with its parameters at the base point and the coefficient of $e$ set to $c/\sqrt n$, $c\in\{0,1,2,3,4,6,8,12,16,24,32\}$; members of one Markov class give the same model and, by \eqref{eq:K1p}, the same effect. The candidates are all orientations of the skeleton with $Y$ a sink, $k=0.05$, and the partial correlation interval uses the Fisher $z$ transform ($4224$ cells, $1000$ replications each, about $45$ minutes on five cores; scripts \texttt{s6\_sim\_batchK.py}, \texttt{s6b\_summary.py}). In the class, the true graph is kept with probability at least $0.994$ and the empty set occurs with probability at most $0.005$. For \DSci\ the maximizing $c$ is at the end of the grid for $17$ pairs at $n=8000$ and $18$ at $n=32000$, so only \DSciO{} is summarized. $K$, $K^{\mathrm{adj}}$ and $K^O$ are each maximized over all ordered pairs (script \texttt{s6c\_info\_fix.py}; an earlier version of the table recorded $K^O$ only for the pair that maximizes $K$); $K^O\ne K$ for $22$ pairs, with $0.969\le K^O/K\le2.83$. At $n=32000$ the $70$ pairs with $K^O\le1$ have ratios to $1+K^O$ between $0.88$ and $1.15$, and the $26$ pairs with $K^O>1$ between $0.77$ and $1.18$; as in the instance of \cref{thm:K-lambda}, the ratio tends to decrease as the constant grows. At $c=0$ the price of \DSciO{} ranges from $1.18$ to $2.93$, with median $1.57$.

\paragraph{Earlier procedures on the classes with $K>0$.} Here $K$ ranges from $0.22$ to $5.31$. \DSciO{} adjusts for optimal adjustment sets, so its first-order bias is governed by $K^O$ (\cref{rem:Oset}), between $0.96K$ and $2.83K$; at $n=32000$ its price divided by $1+K^O$ has median $0.97$ and lies between $0.77$ and $1.18$, and the ten pairs with the smallest information have ratios between $0.77$ and $0.98$. On the same data the likelihood ratio procedures were shorter still: \DSlrO\ in every case, with median ratio $0.92$ to \DSciO{}, and \DSlreff\ in $185$ of the $192$ pairs and true graphs, and equal to it up to $0.02\%$ in the ten configurations of \cref{tab:simK}, where $K^O=K$. At $n=32000$ the price of \DSlreff\ was below the bound of \cref{thm:explicit}(b), evaluated for the true graph, in every case, at between $0.57$ and $0.81$ of it, and above the lower bound of \cref{cor:RK}; with the largest price over the true graphs and $\mathcal R$ and $K$ of the base point, it lies between $0.53$ and $0.95$ times $\mathcal R+K$, with median $0.75$. For the $15$ true graphs with $\mathcal R>1.5$, with $\mathcal R$ and $K$ of the true graph, the ratio to $1+K$ ranges from $0.65$ to $1.11$ and the ratio to $\mathcal R+K$ from $0.53$ to $0.81$. The directed procedure \DSdir\ of \cref{thm:explicit}(c), run on the same data, stayed below its bound in every case, at most at $0.79$ of it, and was shorter than \DSlreff\ in $199$ of the $212$ pairs and true graphs, by a median of $1.5\%$ for $q=2$ and $5.6\%$ for $q=5$, the largest $q$ here; its minimum coverage was $0.940$. The likelihood ratio test spends the level $\gamma_1$ on the true graph, so these procedures return the empty set with probability up to $0.02$ and their coverage is closer to the nominal level, with mean $0.967$ against $0.977$ for \DSciO{}; the lowest coverage, $0.936$, is the minimum of $4499$ estimates with standard error about $0.006$, and a rerun of the six lowest cells with $5000$ replications gave between $0.959$ and $0.965$. The minimum coverage of \DSci\ and \DSciO{} over the $4021$ cells is $0.949$ and $0.952$.

\paragraph{Misspecified $k$.} \DSci\ is honest over $\cM_k$ for the $k$ it is given. With base points whose true model has $k$-triangle-faithfulness ratio $r$ near $0.5k$ or $0.25k$, for the benign cases of \cref{tab:simB} in which such points exist and for one pair whose identification rests on the $k$-triangle-faithfulness screen (\cref{rem:ktf-identifies}), the screen was almost never active at $k=0.05$ and $n\le8000$, and misspecification had no visible effect (minimum coverage $0.956$). At $k=0.25$ and $n$ up to $32000$ the screen is active: with $r\ge2k$ coverage stays at least $0.962$, and in the pair whose identification needs it the median length ratio of \DSci\ at $n=32000$ is $1.67$, against $48.10$ without it; with $r$ near $0.5k$ or $0.25k$ the true graph is screened out as $n$ grows and coverage drops to zero, almost only through empty output, the probability of a non-empty interval that misses the target being at most $0.035$ in all $96$ misspecified cells. If $P\in\cM_k$ the empty set has probability at most $\gamma_1+\gamma_2$, so an empty output is evidence against $P\in\cM_k$, though not a diagnosis of $k$, and a non-empty output carries no guarantee when $P\notin\cM_k$.

\paragraph{Design of the experiment with misspecified $k$.} It uses the procedures of the five-node experiment (script \texttt{s7\_sim\_kmis.py}, summary \texttt{s7b\_summary.py}). The pair whose identification needs the screen is the first, in the enumeration order, of the $12$ pairs with five edges at $p=5$ ($X=1$, $Y=4$) that satisfy \eqref{eq:K2} and \eqref{eq:K1p} but whose $\MP(G)$ has a bad path: $G=\{2\to X,3\to X,2\to3,X\to Y,3\to Y\}$ (the fifth node is isolated) with $e=2\to3$, two members of $\Xeq$ in one Markov class and six members of $[G]_Y$. Base points are drawn as in \cref{sec:sim}, by rejection until the ratio $r$ of the true model at $c=0$ lies in $[2k,\infty)$, $[0.4k,0.6k)$ or $[0.2k,0.3k)$; with $k=0.05$, the misspecified bands were reached in $11$ of the twenty benign cases of \cref{tab:simB}; in the other $9$ all $20000$ draws had $r\ge2k$ (in $5$ of them the true graph has no triangle), and with $k=0.25$ we use the four of them in which $r\ge2k$ is also attainable, $n\in\{2000,8000,32000\}$, and four base points per band for the pair above. The median acceptance rate of $r\ge2k$ over the twenty cases is $0.95$ for $k=0.05$; for the four cases at $k=0.25$ it is between $0.02$ and $0.11$. The probability of a non-empty \DSci\ output that misses the target is at most $0.035$ in the misspecified cells and at most $0.037$ in the correctly specified ones; given a non-empty output, the coverage is at least $0.947$ in the $41$ misspecified cells with empty-set rate below $0.95$.

\paragraph{The likelihood ratio screen.} The comparison above runs \DSciO{} and \DSlrO on the same data, with the base points, true graphs, seeds and $c$ grids of the all-classes experiment and of \cref{tab:simK}, $n\in\{8000,32000\}$ and $1000$ replications per cell ($4704$ cells, about one hour on five cores; scripts \texttt{s8\_sim\_LR.py}, \texttt{s8b\_summary.py}); \DSciO{} agrees with the earlier run in all $4224$ cells of the all-classes experiment. The null quantile of $T_D$ is computed from $400000$ draws of the representation in \cref{thm:coverage}(i) for each pattern $(q_j)$ and $n$; two topological orders of the same graph gave $15.10$ and $15.12$ for a $\chi^2_5$ limit of $15.09$, and under a true graph the test rejected in $0.0085$ of $4000$ samples (script \texttt{d40b\_T2\_supplement.py}). The bound is evaluated at the base point for each true graph $G^*$, with $\varpi_q(\gamma_1)$ from \texttt{d40\_T2\_constants.py}; the bound of the theorem, the maximum over the true graphs of the two terms separately, ranges from $1.96$ to $9.83$ over the $106$ base points (script \texttt{d40b\_T2\_supplement.py}, which also gives the largest $|\Xeq|$ with two or more Markov classes at $p=5$, namely $12$); here $q\in\{2,\dots,5\}$, $|\Xeq|\in\{2,3\}$ and $K^O_+$ equals $K^O$ in all cases. The rerun of the six cells with the lowest coverage of \DSlrO uses new seeds and $5000$ replications (script \texttt{s8c\_recheck\_low.py}). The procedure \DSlreff was run in the same way at $n=32000$ ($2352$ cells, about $28$ minutes on five cores; scripts \texttt{s9\_sim\_eff.py}, \texttt{s9b\_summary.py}), with \DSlrO on the same data agreeing with the run above in every cell. Its candidate intervals use the maximum likelihood estimate of each candidate model and the efficient variance formula evaluated there, which at the population covariance reproduces $\theta_D$ and $s_{\mathrm{eff}}(D,\cdot)$ exactly. The bound of \cref{thm:explicit}(b) is evaluated for each true graph, with $m=2$ in all cases; it ranges from $1.60$ to $9.03$, and the bound of the theorem, maximized over the true graphs, from $1.64$ to $9.03$. The price of \DSlreff was at most $1.024$ times that of \DSlrO, the minimum coverage $0.936$ and the mean $0.966$. The directed procedure \DSdir was run on the same data at $n=32000$ ($2352$ cells, about $36$ minutes on four cores; scripts \texttt{s10\_sim\_dir.py}, \texttt{s10b\_summary.py}), with $\gamma_{1,n}=0.001$, $\gamma_2=\gamma_3=0.01$ and $\alpha'_n=0.029$, and with \DSlreff on the same data agreeing with the run above in every cell. The direction $u_{S^*,S'}$ is computed by central differences at the base point, and $|u|^2$ agrees with $I_{G^*,S'}$ to a relative error of $5.5\cdot10^{-10}$ over the $212$ classes; the null quantile of the directed statistic uses $400000$ draws of the representation in \cref{thm:explicit}(c). The price was between $0.55$ and $0.79$ times the bound of \cref{thm:explicit}(c) for the true graph, with $\alpha'_n$ or its limit $\alpha'$, and between $0.46$ and $0.95$ times $\mathcal R+K$ of the base point, with median $0.74$. It was at least $0.843$ and at most $1.011$ times the price of \DSlreff, with medians $0.985$, $0.979$, $0.959$ and $0.944$ for $q=2,3,4,5$; the minimum coverage was $0.940$, the mean $0.967$, and the probability of an empty output at most $0.025$.

\paragraph{Uniformity in the margin.} The configuration of \cref{sec:sim} has $X=0$, $Y=5$ and $H=\{2\to0,\ 2\to3,\ 3\to0,\ 3\to1,\ 4\to3,\ 0\to5,\ 1\to5\}$, with $e=0-1$, $f=2\to0$, $x_e=1$, $z=X$ and $x_f=2$; $G_\wedge=H+(1\to0)$ has the v-structure $1\to0\leftarrow2$, and $G_\vee=H+(0\to1)$. Its candidate class has $36$ graphs, and $\cX_2=\{G_\wedge,G_\vee\}$. The coefficients of $H_2$ are $\pm U(0.5,1.2)$ and the noise variances $U(0.6,1.6)$, drawn with the first seed for which the $k$-triangle-faithfulness ratio of both true graphs at $P_{00}$, without the elements of $e$ and $f$, is at least $3k$ (it is $0.378$). The true graph is $G_\wedge$ or $G_\vee$, with the parameters of $P_{00}$ except that $e$ and $f$ have coefficients $c/\sqrt n$, $c\in\{0,\pm1,\pm2,\pm4,\pm8,\pm16\}$, and $\lambda$. \DSOp\ is run without the directed test ($\gamma_3=0$), with $\gamma_1=\gamma_2=0.01$, $\alpha'=0.03$ and $k=0.05$; the step that \cref{thm:uniform} uses is the likelihood ratio screen. Of the $308$ cells, $17$ are outside the class because the true model violates $k$-triangle faithfulness, all with $G_\wedge$ true and $\lambda\ge0.2$ or $|c|=16$ at $n=2000$; the ratios above are maxima over the cells in the class. The empty set had probability at most $0.015$ (script \texttt{s13\_sim\_uniform.py}, about $8$ minutes on four cores).

\paragraph{Details of the experiment on all classes with $K>0$.} In the experiment of \cref{sec:sim}, the bound for \DSOp\ uses $\mathcal R^O$, $\sqrt{2\log|\Xeq|}$ and $K^O_+$, and the bounds are evaluated for the true graph. The ratio of the price of \DSeff\ to $\mathcal R+K$ uses the largest price over the true graphs and $\mathcal R$ and $K$ of the base point; its median is $0.74$. The price ratio to the procedure with directions from the base point has median $1.002$; the largest ratios occur in the configurations of \cref{tab:simK} with $\omega\le0.35$, in which two edges of a collider are weak together, so that removing either of them costs little, and for $c\le12$ the structure estimate picked a wrong edge or a wrong class of $H$ in up to $0.34$ of the replications there, against at most $0.030$ in the $96$ classes. For larger $c$ the weak edge is no longer below the threshold $\tau_n$ at the size of the first part, and the estimate was correct in $0.83$ and $0.61$ of the replications on average at $c=24$ and $32$ (rounded down), mostly returning no estimate; the price still stayed below the bound. The directed test rejected the true graph in $0.0094$ of the replications on average, against its level $\gamma_3=0.01$. Each cell has standard error about $0.007$ for the coverage; the rerun of the nine lowest cells used new seeds and $5000$ replications. The empty set had probability at most $0.021$.

\paragraph{Double screening with estimated directions.} The runs of \cref{sec:sim} use the base points, true graphs, seeds and $c$ grids of the experiments above, $n=32000$ and $1000$ replications per cell for the classes with $K>0$ and \cref{tab:simK} ($2352$ cells, about $68$ minutes on four cores; scripts \texttt{s11\_sim\_adaptive.py}, \texttt{s11b\_summary.py}), and the cases, seeds and $2000$ replications per cell of the five-node experiment (about $16$ minutes on four cores; scripts \texttt{s12\_sim\_p5\_adaptive.py}, \texttt{s12b\_table.py}); the rerun of the lowest-coverage cells uses \texttt{s11d\_recheck\_low.py}. The first $n_1=\lceil n^{2/3}\rceil$ observations ($1008$ at $n=32000$) form the first part. Step (0) evaluates $D_{D-e'}$ once for each pair of an edge $e'$ and a Markov class of $D-e'$; the directions are central differences at $\Pi_{\hat H}\hat\Sigma^{(1)}$, and the directed statistic is computed on the other $n_2$ observations. For $d=1$ its null quantile given $\hat B$ is computed from the representation of \cref{thm:coverage}(i) in Rao--Blackwell form, $\Prob(|\sum_jw_jt_{n_2-j}|\le y)=\E[2\Phi(y/W^{1/2})-1]$ with $W=\sum_jw_j^2(n_2-j)/V_j$, averaged over $4000$ draws of $V$ and solved by Newton's method, and for $d\ge2$ from $100000$ draws; against $400000$ independent Monte Carlo draws the $d=1$ quantiles differ by at most $0.81\%$, the Monte Carlo error of the latter (script \texttt{s11c\_quantile\_check.py}); the branch $d\ge2$, whose shared draws differ from independent ones by up to $1.8\%$ at the $0.99$ quantile, is not used, because all experiments have $m=2$. In the same runs the procedure with the directions computed from the base point, \DSdir\ with $\gamma_{1,n}=0.001$, agrees with the run of the experiment with the directed procedure in every cell, and in the five-node experiment Naive, Hull, Union and \DSlrO\ agree with the earlier runs in every cell.

\end{document}